\documentclass[10pt]{article}
\usepackage{fullpage}

\usepackage[utf8]{inputenc}
\usepackage[backend=bibtex,sorting=nyt,doi=false,isbn=false,url=false,eprint=false]{biblatex}
\usepackage{latexsym}
\usepackage{amssymb}
\usepackage{amsthm}
\usepackage{amscd}
\usepackage{amsmath}
\usepackage{mathrsfs}
\usepackage{graphicx}
\usepackage{hyperref}
\usepackage{shuffle}
\usepackage{stmaryrd}
\usepackage{enumitem}
\usepackage[noabbrev]{cleveref}
\usepackage{ytableau}

\usepackage[colorinlistoftodos]{todonotes}
\usetikzlibrary{chains,scopes,decorations.markings}

\usepackage{tikz-cd}

\theoremstyle{definition}
\newtheorem* {theorem*}{Theorem}
\newtheorem* {conjecture*}{Conjecture}
\newtheorem{theorem}{Theorem}[section]

\theoremstyle{definition}

\newtheorem* {example*}{Example}

\newtheorem{lemma}[theorem]{Lemma}
\theoremstyle{definition}
\newtheorem{definition}[theorem]{Definition}
\theoremstyle{definition}

\newtheorem{conjecture}[theorem]{Conjecture}
\newtheorem{proposition}[theorem]{Proposition}
\newtheorem{corollary}[theorem]{Corollary}

\newtheorem* {remark*}{Remark}
\newtheorem{remark}[theorem]{Remark}
\theoremstyle{definition}
\newtheorem {example}[theorem]{Example}
\theoremstyle{definition}

\theoremstyle{definition}

\theoremstyle{definition}

\theoremstyle{definition}

\newcommand{\std}{\operatorname{std}}

\newcommand{\Inv}{\mathrm{Inv}}
\newcommand{\Ess}{\operatorname{Ess}}
\newcommand{\rank}{\operatorname{rank}}

\renewcommand{\(}{\left(}
\renewcommand{\)}{\right)}

\newcommand{\CC}{\mathbb{C}}
\newcommand{\QQ}{\mathbb{Q}}

\newcommand{\cK}{\mathcal{K}}
\newcommand{\cR}{\mathcal{R}}

\newcommand{\cX}{\mathcal{X}}
\newcommand{\cW}{\mathcal{W}}

\newcommand{\ZZ}{\mathbb{Z}}

\newcommand{\cyc}{\textsf{cyc}}

\newcommand{\cL}{\mathcal{L}}

\newcommand{\fk}{\mathfrak}

\newcommand{\barr}{\begin{array}}
\newcommand{\earr}{\end{array}}
\newcommand{\ba}{\begin{aligned}}
\newcommand{\ea}{\end{aligned}}
\newcommand{\be}{\begin{equation}}
\newcommand{\ee}{\end{equation}}

\newcommand{\Cyc}{\mathrm{Cyc}}

\newcommand{\qquand}{\qquad\text{and}\qquad}
\newcommand{\quand}{\quad\text{and}\quad}

\newcommand{\I}{\mathsf{Invol}}
\newcommand{\Invol}{I_n}
\newcommand{\Ivex}{I^{\mathsf{vex}}_n}
\newcommand{\cH}{\mathcal H}

\newcommand{\ds}{\displaystyle}

\newcommand{\PP}{\mathbb{Z}_{>0}}
\newcommand{\NN}{\mathbb{Z}_{\geq 0}}

\newcommand{\fkS}{\fk S}

\newcommand{\ben}{\begin{enumerate}}
\newcommand{\een}{\end{enumerate}}

\newcommand{\bei}{\begin{itemize}}
\newcommand{\eei}{\end{itemize}}

\newcommand{\cT}{\mathscr{T}}

\newcommand{\xRightarrow}[2][]{\ext@arrow 0359\Rightarrowfill@{#1}{#2}}

\newcommand{\Sp}{\mathsf{Sp}}

\newcommand{\cA}{\mathcal{A}}

\newcommand{\cG}{\mathcal{G}}

\newcommand{\D}{\mathsf{D}}

\definecolor{darkred}{rgb}{0.7,0,0} 
\newcommand{\defn}[1]{{\color{darkred}\emph{#1}}} 

\usepackage{pbox}
\tikzset{every loop/.style={min distance=10 mm, in=60, out=120, looseness=10}}
\tikzset{
dot/.style = {circle, fill, inner sep=2.2,outer sep=0},
}

\numberwithin{equation}{section}
\allowdisplaybreaks[1]

\newcommand{\bpartial}{\partial^{(\beta)}}

\newcommand{\fkGO}{\mathfrak{G}^{\mathsf{O}}}

\newcommand{\Mat}{\mathsf{Mat}}

\newcommand{\iG}{\widehat{\cG}}
\newcommand{\cU}{\mathcal{U}}
\newcommand{\SD}{\mathsf{SD}}

\newcommand{\word}{\mathsf{word}}

\newcommand{\iD}{\hat D}

\newcommand{\ellinv}{\ell_{\mathsf{inv}}}

\newcommand{\NW}{\mathsf{NW}}
\newcommand{\alignNW}{\mathrm{align}_{\NW}}

\newcommand{\bff}{\mathbf{f}}

\newcommand{\ytabc}[2]{
\ytableausetup{boxsize = #1,aligntableaux=center}
{\begin{ytableau}  #2  \end{ytableau}}
}

\usepackage{nicematrix}
\usepackage{mathtools}
\usepackage{verbatim}

\newcommand{\sdgap}{ \\[-4pt]\\ }

\newcommand{\symDef}{\mathsf{symDef}}
\newcommand{\cWsim}{\sim_\cW}
\newcommand{\fkG}{\mathfrak{G}^{(\beta)}}
\renewcommand{\iG}{\widehat{\mathfrak{G}}^{(\beta)}}

\newcommand{\SPP}{\mathsf{ShPP}}
\newcommand{\SetSPP}{\mathsf{SetShPP}}
\newcommand{\SetPP}{\mathsf{SetPP}}

\newcommand{\MPipes}{\mathsf{MPipes}}

\newcommand{\invRothe}{\mathsf{invRothe}}
\newcommand{\double}{\mathsf{double}}
\newcommand{\halve}{\mathsf{halve}}
\newcommand{\transpose}{\mathsf{transpose}}

\newcommand{\GPD}{\mathsf{GPD}}
\newcommand{\BPD}{\mathsf{BPD}}

\newcommand{\kGPD}{\mathsf{GPD}_K}
\newcommand{\kBPD}{\mathsf{BPD}_K}

\newcommand{\InvPD}{\mathsf{InvGPD}}

\newcommand{\InvBPD}{\mathsf{InvBPD}}
\newcommand{\MInvBPD}{\mathsf{InvMBPD}}

\newcommand{\kInvPD}{\mathsf{InvGPD}_K}
\newcommand{\kInvBPD}{\mathsf{InvBPD}_K}

\newcommand{\kweight}{\weight_K}
\newcommand{\hatkweight}{\widehat\weight_K}
\newcommand{\kpi}{\overline{\pi}}

\newcommand{\NRes}{\mathrm{NRes}}
\newcommand{\NDes}{\mathrm{NDes}}

\newcommand{\circled}[1]{\text{\textcircled{\raisebox{-0.9pt}{#1}}}}

\newcommand{\complexDroop}{\mathsf{complexDroop}}
\newcommand{\droop}{\mathsf{droop}}

\newcommand{\extComplexDroop}{\mathsf{strongComplexDroop}}
\newcommand{\extDroop}{\mathsf{strongDroop}}

\newcommand{\undroop}{\mathsf{undroop}}
\newcommand{\complexUndroop}{\mathsf{complexUndroop}}

\newcommand{\iS}{\widehat{\fkS}}
\newcommand{\half}{\mathsf{half}}

\newcommand{\weight}{\mathsf{wt}}

\newcommand{\blank}{\mathsf{blank}}
\newcommand{\cross}{\mathsf{cross}}
\newcommand{\bump}{\mathsf{bump}}
\newcommand{\upturn}{\mathsf{upturn}}

\renewcommand{\std}[2]{#1\cap #2}

\newcommand{\ASM}{\mathrm{ASM}}
\newcommand{\DASM}{\mathrm{DASM}}

\newcommand{\hatu}{\widehat{u}}
\newcommand{\uPrime}{\ell}
\newcommand{\hatuPrime}{\widehat{\ell}}
\newcommand{\diag}{\mathrm{diag}}

\newcommand{\cToff}{\cT_{\mathrm{off}}}
\newcommand{\cTdiag}{\cT_{\mathrm{diag}}}

\newcommand{\flatten}{\mathsf{flatten}}

\newcommand{\cJdiag}{\mathcal{J}_{\Delta}}
\newcommand{\cJoff}{\mathcal{J}}

\newcommand{\cBoff}{\mathcal{B}}

\newcommand{\cJBdiag}{\mathcal{M}_\Delta}

\newcommand{\init}{\mathrm{init}}
\newcommand{\mon}{\mathbf{m}}
\newcommand{\SymMat}{\Mat^{\mathrm{sym}}}
\newcommand{\MSV}{\mathsf{MSV}}
\newcommand{\SymMSV}{\MSV^{\mathrm{sym}}}
\newcommand{\SymI}{I^{\mathrm{sym}}}
\newcommand{\SymJ}{J^{\mathrm{sym}}}
\newcommand{\cSymG}{\cG^{\mathrm{sym}}}
\newcommand{\mat}[1]{\left[ \begin{array} #1 \end{array}\right]}

\newcommand{\UT}{\nabla}
\newcommand{\LT}{\mathsf{LT}}

\newcommand{\displace}{\mathsf{displace}}
\newcommand{\SymDelta}{\Delta^{\mathrm{sym}}}
\newcommand{\iGconj}{\widetilde{\mathfrak{G}}^{\textsf{O}}}

\newcommand{\DBox}{D^{\mathsf{dbl}}}
\newcommand{\bft}{\mathbf t}
\newcommand{\bftrev}{\bft^{\mathsf{rev}}}
\newcommand{\calT}{\mathcal{T}}
\newcommand{\CaseA}{Case (A)}
\newcommand{\CaseB}{Case (B)}
\newcommand{\CaseC}{Case (C)}

\newcommand{\bpdwidth}{1.2em}
\newcommand{\bpdsmallwidth}{1.0em}

\newcommand{\ETILE}{
\begin{tikzpicture}[x=\bpdwidth,y=\bpdwidth,baseline=(z.base)]
\node at (0.5,1.25) (z) {};
\draw[thick,rounded corners,color=blue] (.5,2)--(.5,1.5)--(0,1.5);
\draw[thick,rounded corners,color=blue] (.5,1)--(.5,1.5)--(1,1.5);
\draw[gray!70, line width=0.7pt](0,1)rectangle(1,2);
\end{tikzpicture}
}

\newcommand{\JTILE}{
\begin{tikzpicture}[x=\bpdwidth,y=\bpdwidth,baseline=(z.base)]
\node at (0.5,1.25) (z) {};
\draw[thick,rounded corners,color=blue] (.5,2)--(.5,1.5)--(0,1.5);
\draw[gray!70, line width=0.7pt](0,1)rectangle(1,2);
\end{tikzpicture}
}

\newcommand{\CTILE}{
\begin{tikzpicture}[x=\bpdwidth,y=\bpdwidth,baseline=(z.base)]
\node at (0.5,1.25) (z) {};
\draw[thick,rounded corners,color=blue] (.5,1)--(.5,1.5)--(1,1.5);
\draw[gray!70, line width=0.7pt](0,1)rectangle(1,2);
\end{tikzpicture}
}

\newcommand{\PTILE}{
\begin{tikzpicture}[x=\bpdwidth,y=\bpdwidth,baseline=(z.base)]
\node at (0.5,1.25) (z) {};
\draw[thick,rounded corners,color=blue] (0,1.5)--(1,1.5);
\draw[thick,rounded corners,color=blue] (.5,1)--(.5,2);
\draw[gray!70, line width=0.7pt](0,1)rectangle(1,2);
\end{tikzpicture}
}

\newcommand{\HTILE}{
\begin{tikzpicture}[x=\bpdwidth,y=\bpdwidth,baseline=(z.base)]
\node at (0.5,1.25) (z) {};
\draw[thick,rounded corners,color=blue] (0,1.5)--(1,1.5);
\draw[gray!70, line width=0.7pt](0,1)rectangle(1,2);
\end{tikzpicture}
}

\newcommand{\VTILE}{
\begin{tikzpicture}[x=\bpdwidth,y=\bpdwidth,baseline=(z.base)]
\node at (0.5,1.25) (z) {};
\draw[thick,rounded corners,color=blue] (.5,1)--(.5,2);
\draw[gray!70, line width=0.7pt](0,1)rectangle(1,2);
\end{tikzpicture}}

\newcommand{\BTILE}{
\begin{tikzpicture}[x=\bpdwidth,y=\bpdwidth,baseline=(z.base)]
\node at (0.5,1.25) (z) {};
\draw[gray!70, line width=0.7pt](0,1)rectangle(1,2);
\end{tikzpicture}
}

\newcommand{\JDIAG}{
\begin{tikzpicture}[x=\bpdwidth,y=\bpdwidth,baseline=(z.base)]
\node at (0.5,1.25) (z) {};
\draw[thick,rounded corners,color=blue] (.5,1.5)--(0,1.5);
\draw[white, line width=0.7pt,fill=white] (1,2)--(1,1)--(0,2)--(1,2);
\draw[gray!70, line width=0.7pt] (0,1)--(1,1)--(0,2)--(0,1);
\end{tikzpicture}
}

\newcommand{\CDIAG}{
\begin{tikzpicture}[x=\bpdwidth,y=\bpdwidth,baseline=(z.base)]
\node at (0.5,1.25) (z) {};
\draw[thick,rounded corners,color=blue] (.5,1)--(.5,1.5);
\draw[white, line width=0.7pt,fill=white] (1,2)--(1,1)--(0,2)--(1,2);
\draw[gray!70, line width=0.7pt] (0,1)--(1,1)--(0,2)--(0,1);
\end{tikzpicture}
}

\newcommand{\PDIAG}{
\begin{tikzpicture}[x=\bpdwidth,y=\bpdwidth,baseline=(z.base)]
\node at (0.5,1.25) (z) {};
\draw[white, line width=0.7pt,fill=white] (1,2)--(1,1)--(0,2)--(1,2);
\draw[gray!70, line width=0.7pt] (0,1)--(1,1)--(0,2)--(0,1);
\draw[thick,rounded corners,color=blue] (0,1.5)--(0.5,1.5);
\draw[thick,rounded corners,color=blue] (.5,1)--(.5,1.5);
\end{tikzpicture}
}

\newcommand{\BDIAG}{
\begin{tikzpicture}[x=\bpdwidth,y=\bpdwidth,baseline=(z.base)]
\node at (0.5,1.25) (z) {};
\draw[white, line width=0.7pt,fill=white] (1,2)--(1,1)--(0,2)--(1,2);
\draw[gray!70, line width=0.7pt] (0,1)--(1,1)--(0,2)--(0,1);
\end{tikzpicture}
}

\newcommand{\EDIAG}{
\begin{tikzpicture}[x=\bpdwidth,y=\bpdwidth,baseline=(z.base)]
\node at (0.5,1.25) (z) {};
\draw[thick,rounded corners,color=blue] (.5,2)--(.5,1.5)--(0,1.5);
\draw[thick,rounded corners,color=blue] (.5,1)--(.5,1.5)--(1,1.5);
\draw[white, line width=0.7pt,fill=white] (1,2)--(1,1)--(0,2)--(1,2);
\draw[gray!70, line width=0.7pt] (0,1)--(1,1)--(0,2)--(0,1);
\end{tikzpicture}
}

\DeclareMathOperator{\MD}{MD}
\DeclareMathOperator{\BZOD}{BZOD}
\DeclareMathOperator{\MOD}{MOD}

\newcommand{\fkSO}{\fkS^{\mathsf{O}}}
 \renewcommand{\O}{\mathsf{O}}

\newtheorem{mainthm}{Theorem}

\begin{document}
\title{Bumpless pipe dreams for involutions and shifted plane partitions}
\author{
    Eric Marberg
    \\
    HKUST \\
    {\tt emarberg@ust.hk}
    \and
    Alejandro H. Morales\\
    UQAM\\
    {\tt morales\_borrero.alejandro@uqam.ca}
    \and
    Tianyi Yu \\
    UQAM\\
    {\tt yu.tianyi@uqam.ca}
}

\date{}

\maketitle

\begin{abstract}
This article studies a triangular variant of bumpless pipe dreams for involutions in the symmetric group. We  
 prove that a natural weight generating function for these objects recovers all involution Schubert polynomials. For vexillary involutions, we derive a $K$-theoretic extension of this result that gives a new formula for orthogonal Grothendieck polynomials.
In the course of showing this formula, we compute a diagonal term order Gr\"obner basis
 for the defining ideals of symmetric matrix Schubert varieties associated to vexillary involutions.
As applications, we prove a conjectured product formula for principal specializations of involution Schubert polynomials. We go on to derive a shifted analogue of the Fomin--Kirillov formula counting bounded plane partitions. This analogue proves the straight-shape case of a  conjecture of Ferroni, Morales, and Panova on the Ehrhart positivity for the order polytopes of shifted partitions.  We also conjecture a  bumpless pipe dream formula for involution Grothendieck polynomials and discuss connections with diagonally symmetric alternating sign matrices and asymptotics of principal specializations.
\end{abstract}


\section{Introduction}

\defn{Bumpless pipe dreams (BPDs)} are tilings of a square grid whose pipes encode permutations and whose weights give \defn{Schubert} and \defn{Grothendieck polynomials}. In this article, we introduce triangular analogues of BPDs indexed by involutions in the symmetric group. These objects give formulas for \defn{involution Schubert polynomials} and, in the vexillary case, \defn{orthogonal Grothendieck polynomials}.
The rest of this introduction reviews our motivations from the permutation
setting and then outlines our main results.

\subsection{Bumpless pipe dreams}

Let $n$ be a 
  positive integer.
Bumpless pipe dreams are certain tilings of an $n\times n$ grid that determine a collection of $n$ distinct ``pipes'' that enter the south side of the grid and exit on the east border.
These objects were introduced in \cite{LLS} and their 
precise definition is reviewed in Section~\ref{bpd-sect}. Here are 
two examples:
\be\label{bpd-ex1}
\begin{tikzpicture}[x=\bpdwidth,y=\bpdwidth,line cap=round,line join=round,baseline=(z.base)]
\node at (0,2.25) (z) {};
\draw[gray!70, line width=0.7pt](0,0)--(5,0);
\draw[gray!70, line width=0.7pt](0,0)--(0,5);
\draw[gray!70, line width=0.7pt](0,1)--(5,1);
\draw[gray!70, line width=0.7pt](1,0)--(1,5);
\draw[gray!70, line width=0.7pt](0,2)--(5,2);
\draw[gray!70, line width=0.7pt](2,0)--(2,5);
\draw[gray!70, line width=0.7pt](0,3)--(5,3);
\draw[gray!70, line width=0.7pt](3,0)--(3,5);
\draw[gray!70, line width=0.7pt](0,4)--(5,4);
\draw[gray!70, line width=0.7pt](4,0)--(4,5);
\draw[gray!70, line width=0.7pt](0,5)--(5,5);
\draw[gray!70, line width=0.7pt](5,0)--(5,5);
\draw[thick,rounded corners,color=blue] (2.5,4)--(2.5,4.5)--(3,4.5);
\draw[thick,rounded corners,color=blue] (1.5,3)--(1.5,3.5)--(2,3.5);
\draw[thick,rounded corners,color=blue] (0.5,2)--(0.5,2.5)--(1,2.5);
\draw[thick,rounded corners,color=blue] (4.5,1)--(4.5,1.5)--(5,1.5);
\draw[thick,rounded corners,color=blue] (3.5,0)--(3.5,0.5)--(4,0.5);
\draw[thick,rounded corners,color=blue] (0.5,0)--(0.5,1);
\draw[thick,rounded corners,color=blue] (0.5,1)--(0.5,2);
\draw[thick,rounded corners,color=blue] (1,2.5)--(2,2.5);
\draw[thick,rounded corners,color=blue] (1.5,2)--(1.5,3);
\draw[thick,rounded corners,color=blue] (2,2.5)--(3,2.5);
\draw[thick,rounded corners,color=blue] (2.5,2)--(2.5,3);
\draw[thick,rounded corners,color=blue] (3,2.5)--(4,2.5);
\draw[thick,rounded corners,color=blue] (4,2.5)--(5,2.5);
\draw[thick,rounded corners,color=blue] (1.5,0)--(1.5,1);
\draw[thick,rounded corners,color=blue] (1.5,1)--(1.5,2);
\draw[thick,rounded corners,color=blue] (2,3.5)--(3,3.5);
\draw[thick,rounded corners,color=blue] (2.5,3)--(2.5,4);
\draw[thick,rounded corners,color=blue] (3,3.5)--(4,3.5);
\draw[thick,rounded corners,color=blue] (4,3.5)--(5,3.5);
\draw[thick,rounded corners,color=blue] (2.5,0)--(2.5,1);
\draw[thick,rounded corners,color=blue] (2.5,1)--(2.5,2);
\draw[thick,rounded corners,color=blue] (3,4.5)--(4,4.5);
\draw[thick,rounded corners,color=blue] (4,4.5)--(5,4.5);
\draw[thick,rounded corners,color=blue] (4,0.5)--(5,0.5);
\draw[thick,rounded corners,color=blue] (4.5,0)--(4.5,1);
\draw[black, line width=1.6pt] (0,5)--(5,5)--(5,0)--(0,0)--(0,5);
\draw[](0.5,-0.4) node {1};
\draw[](1.5,-0.4) node {2};
\draw[](2.5,-0.4) node {3};
\draw[](3.5,-0.4) node {4};
\draw[](4.5,-0.4) node {5};
\draw[](5.4,4.5) node {3};
\draw[](5.4,3.5) node {2};
\draw[](5.4,2.5) node {1};
\draw[](5.4,1.5) node {5};
\draw[](5.4,0.5) node {4};
\end{tikzpicture}
\qquad\qquad
\begin{tikzpicture}[x=\bpdwidth,y=\bpdwidth,line cap=round,line join=round,baseline=(z.base)]
\node at (0,2.25) (z) {};
\draw[gray!70, line width=0.7pt](0,0)--(5,0);
\draw[gray!70, line width=0.7pt](0,0)--(0,5);
\draw[gray!70, line width=0.7pt](0,1)--(5,1);
\draw[gray!70, line width=0.7pt](1,0)--(1,5);
\draw[gray!70, line width=0.7pt](0,2)--(5,2);
\draw[gray!70, line width=0.7pt](2,0)--(2,5);
\draw[gray!70, line width=0.7pt](0,3)--(5,3);
\draw[gray!70, line width=0.7pt](3,0)--(3,5);
\draw[gray!70, line width=0.7pt](0,4)--(5,4);
\draw[gray!70, line width=0.7pt](4,0)--(4,5);
\draw[gray!70, line width=0.7pt](0,5)--(5,5);
\draw[gray!70, line width=0.7pt](5,0)--(5,5);
\draw[thick,rounded corners,color=blue] (2.5,4)--(2.5,4.5)--(3,4.5);
\draw[thick,rounded corners,color=blue] (3.5,0)--(3.5,0.5)--(4,0.5);
\draw[thick,rounded corners,color=blue] (4.5,1)--(4.5,1.5)--(5,1.5);
\draw[thick,rounded corners,color=blue] (3.5,2)--(3.5,1.5)--(3,1.5);
\draw[thick,rounded corners,color=blue] (0.5,1)--(0.5,1.5)--(1,1.5);
\draw[thick,rounded corners,color=blue] (1.5,2)--(1.5,2.5)--(2,2.5);
\draw[thick,rounded corners,color=blue] (3.5,3)--(3.5,3.5)--(4,3.5);
\draw[thick,rounded corners,color=blue] (0.5,0)--(0.5,1);
\draw[thick,rounded corners,color=blue] (1,1.5)--(2,1.5);
\draw[thick,rounded corners,color=blue] (1.5,1)--(1.5,2);
\draw[thick,rounded corners,color=blue] (2,1.5)--(3,1.5);
\draw[thick,rounded corners,color=blue] (2.5,1)--(2.5,2);
\draw[thick,rounded corners,color=blue] (3,2.5)--(4,2.5);
\draw[thick,rounded corners,color=blue] (3.5,2)--(3.5,3);
\draw[thick,rounded corners,color=blue] (4,3.5)--(5,3.5);
\draw[thick,rounded corners,color=blue] (1.5,0)--(1.5,1);
\draw[thick,rounded corners,color=blue] (2,2.5)--(3,2.5);
\draw[thick,rounded corners,color=blue] (2.5,2)--(2.5,3);
\draw[thick,rounded corners,color=blue] (4,2.5)--(5,2.5);
\draw[thick,rounded corners,color=blue] (2.5,0)--(2.5,1);
\draw[thick,rounded corners,color=blue] (2.5,3)--(2.5,4);
\draw[thick,rounded corners,color=blue] (3,4.5)--(4,4.5);
\draw[thick,rounded corners,color=blue] (4,4.5)--(5,4.5);
\draw[thick,rounded corners,color=blue] (4,0.5)--(5,0.5);
\draw[thick,rounded corners,color=blue] (4.5,0)--(4.5,1);
\draw[black, line width=1.6pt] (0,5)--(5,5)--(5,0)--(0,0)--(0,5);
\draw[](0.5,-0.4) node {1};
\draw[](1.5,-0.4) node {2};
\draw[](2.5,-0.4) node {3};
\draw[](3.5,-0.4) node {4};
\draw[](4.5,-0.4) node {5};
\draw[](5.4,4.5) node {3};
\draw[](5.4,3.5) node {2};
\draw[](5.4,2.5) node {1};
\draw[](5.4,1.5) node {5};
\draw[](5.4,0.5) node {4};
\end{tikzpicture}
.
\ee
Each BPD has an associated permutation in the 
symmetric group $S_n$. In both of the preceding examples
the associated permutation is given in one-line notation by $32154 \in S_5$. 

The left BPD in \eqref{bpd-ex1} is \defn{reduced}
in the sense that no two of its pipes cross more than once.
The permutation associated to such a BPD is easy to define:
just label the pipes left to right from $1$ to $n$ on the south border,
then read the pipe labels from top to bottom on the east side to get the one-line notation.
 The right BPD in \eqref{bpd-ex1} is not reduced,
 and for such BPDs the internal pipe labels change at crossings according 
 to certain local rules (see Section~\ref{kbpd-sect}).
 
A primary motivation to consider bumpless pipe dreams 
comes from Schubert calculus.
Recall from \cite{Weigandt} that the \defn{double Grothendieck polynomials} $\fkG_w(x;y)$  are a
family of polynomials in $\ZZ[\beta][x_1,y_1,x_2,y_2,\dots]$
indexed by permutations $w \in S_n$.
At the specializations $\beta=0$ and $\beta=-1$,
these polynomials respectively represent the
equivariant cohomology and $K$-theory classes
of Schubert varieties in the complete flag variety.
The $\beta=0$ specializations $\fkS_w(x;y) := \mathfrak{G}^{(0)}_w(x;-y)$ 
are the better-known family of \defn{double Schubert polynomials}.

Let $\kBPD(w)$ be the set of bumpless pipe dreams 
associated to a permutation $w\in S_n$,
and write $\BPD(w)\subseteq \kBPD(w)$ for the subset of reduced BPDs.
Then one has
\be
\label{double-sch-eq}
\textstyle
\fkS_w(x;y) = \sum_{D \in \BPD(w)} \weight(D)
\qquand
\fkG_w(x;y) = \sum_{D \in \kBPD(w)}  \beta^{|D|-\ell(w)}\kweight(D)
\ee
for certain natural weights $\weight(D)$ and $\kweight(D)$ attached to a BPD. Here $|D|$ 
is the number of blank square tiles $\BTILE$ in $D$ and
$\ell : S_n \to \NN$ is the Coxeter length function,
which assigns to each permutation its number of inversions.
 Lam, Lee, and Shimozono proved the left formula in \cite{LLS},
and this was extended to the right formula by Weigandt \cite{Weigandt}.
Let $x\oplus y = x+y+\beta xy$. Then the weight functions are 
\be\label{wt-eq}
 \weight(D) = \prod_{(i,j)} (x_i-y_j)
 \qquand
\kweight(D) = \prod_{(i,j)} (x_i\oplus y_j) \prod_{(k,l)} (1+\beta x_k)(1+\beta y_l)
\ee
 where the respective products are over $(i,j)$
where $D_{ij} = \BTILE$ is a blank tile
and $(k,l)$ where 
$D_{kl} =\JTILE$.
Results of Weigandt \cite{Weigandt} also identify a
notable formula intermediate between the two identities in \eqref{double-sch-eq}: if $w \in S_n$ is \defn{vexillary} in the sense of avoiding 
the permutation pattern $2143$, then 
\be\label{double-sch-eq2}
\textstyle
\BPD(w) = \kBPD(w)\quand \fkG_w(x;y) = \sum_{D \in \BPD(w)} \kweight(D).
\ee

There are older formulas for Schubert and Grothendieck polynomials
as generating functions for ``classical'' pipe dreams
\cite{BergeronBilley, fomin1994grothendieck,KnutsonMiller},
which are similar tilings that allow the ``bump'' tile $\ETILE$.
The exclusion of this tile accounts for the name ``bumpless'' pipe dream.
Lam, Lee, and Shimozono's initial motivation to consider bumpless pipe dreams in \cite{LLS} came from studying a more general family of  \defn{backstable Schubert polynomials}.
The identities in \eqref{double-sch-eq} transparently extend
to the backstable setting, unlike ``classical'' pipe dream formulas for double Schubert polynomials.

Since their introduction in \cite{LLS},
bumpless pipe dreams have been extensively studied.
They have found applications in combinatorial $K$-theory, bijective Schubert calculus, Gr\"obner degenerations, and the enumeration of alternating sign matrices and plane partitions, for example; see \cite{GaoHuang,HPW,Huang2023,HuangStriker2024,KleinWeigandt,Weigandt}.
Our focus here is on proving analogues of 
\eqref{double-sch-eq} and \eqref{double-sch-eq2}
that pertain to variants of Schubert and Grothendieck polynomials indexed by involutions in the symmetric group.

\subsection{Involution pipe dreams}

We now introduce the triangular analogues
of bumpless pipe dreams that figure in our main results.
Our discussion here is informal; the precise definitions appear in Section~\ref{bpd-sect}.
Throughout, let $\Invol = \left\{z \in S_n : z^{-1}=z\right\}$ be the set of involutions in $S_n$. Recall that these are the permutations
whose cycles each have length one or two.

Our \defn{involution pipe dreams} will be tilings of the lower half triangle of an $n\times n$
matrix. There is a simple halving operation that turns 
the bumpless pipe dreams in the previous section into such triangular tilings:
\be\label{inv-bpd-ex}
\begin{tikzpicture}[x=\bpdwidth,y=\bpdwidth,line cap=round,line join=round,baseline=(z.base)]
\node at (0,2.25) (z) {};
\draw[gray!70, line width=0.7pt](0,0)--(5,0);
\draw[gray!70, line width=0.7pt](0,0)--(0,5);
\draw[gray!70, line width=0.7pt](0,1)--(5,1);
\draw[gray!70, line width=0.7pt](1,0)--(1,5);
\draw[gray!70, line width=0.7pt](0,2)--(5,2);
\draw[gray!70, line width=0.7pt](2,0)--(2,5);
\draw[gray!70, line width=0.7pt](0,3)--(5,3);
\draw[gray!70, line width=0.7pt](3,0)--(3,5);
\draw[gray!70, line width=0.7pt](0,4)--(5,4);
\draw[gray!70, line width=0.7pt](4,0)--(4,5);
\draw[gray!70, line width=0.7pt](0,5)--(5,5);
\draw[gray!70, line width=0.7pt](5,0)--(5,5);
\draw[thick,rounded corners,color=blue] (2.5,4)--(2.5,4.5)--(3,4.5);
\draw[thick,rounded corners,color=blue] (1.5,3)--(1.5,3.5)--(2,3.5);
\draw[thick,rounded corners,color=blue] (0.5,2)--(0.5,2.5)--(1,2.5);
\draw[thick,rounded corners,color=blue] (4.5,1)--(4.5,1.5)--(5,1.5);
\draw[thick,rounded corners,color=blue] (3.5,0)--(3.5,0.5)--(4,0.5);
\draw[thick,rounded corners,color=blue] (0.5,0)--(0.5,1);
\draw[thick,rounded corners,color=blue] (0.5,1)--(0.5,2);
\draw[thick,rounded corners,color=blue] (1,2.5)--(2,2.5);
\draw[thick,rounded corners,color=blue] (1.5,2)--(1.5,3);
\draw[thick,rounded corners,color=blue] (2,2.5)--(3,2.5);
\draw[thick,rounded corners,color=blue] (2.5,2)--(2.5,3);
\draw[thick,rounded corners,color=blue] (3,2.5)--(4,2.5);
\draw[thick,rounded corners,color=blue] (4,2.5)--(5,2.5);
\draw[thick,rounded corners,color=blue] (1.5,0)--(1.5,1);
\draw[thick,rounded corners,color=blue] (1.5,1)--(1.5,2);
\draw[thick,rounded corners,color=blue] (2,3.5)--(3,3.5);
\draw[thick,rounded corners,color=blue] (2.5,3)--(2.5,4);
\draw[thick,rounded corners,color=blue] (3,3.5)--(4,3.5);
\draw[thick,rounded corners,color=blue] (4,3.5)--(5,3.5);
\draw[thick,rounded corners,color=blue] (2.5,0)--(2.5,1);
\draw[thick,rounded corners,color=blue] (2.5,1)--(2.5,2);
\draw[thick,rounded corners,color=blue] (3,4.5)--(4,4.5);
\draw[thick,rounded corners,color=blue] (4,4.5)--(5,4.5);
\draw[thick,rounded corners,color=blue] (4,0.5)--(5,0.5);
\draw[thick,rounded corners,color=blue] (4.5,0)--(4.5,1);
\draw[black, line width=1.6pt] (0,5)--(5,5)--(5,0)--(0,0)--(0,5);
\end{tikzpicture}
\mapsto
\begin{tikzpicture}[x=\bpdwidth,y=\bpdwidth,line cap=round,line join=round,baseline=(z.base)]
\node at (0,2.25) (z) {};
\draw[gray!70, line width=0.7pt](0,0)--(5,0);
\draw[gray!70, line width=0.7pt](0,0)--(0,5);
\draw[gray!70, line width=0.7pt](0,1)--(5,1);
\draw[gray!70, line width=0.7pt](1,0)--(1,5);
\draw[gray!70, line width=0.7pt](0,2)--(5,2);
\draw[gray!70, line width=0.7pt](2,0)--(2,5);
\draw[gray!70, line width=0.7pt](0,3)--(5,3);
\draw[gray!70, line width=0.7pt](3,0)--(3,5);
\draw[gray!70, line width=0.7pt](0,4)--(5,4);
\draw[gray!70, line width=0.7pt](4,0)--(4,5);
\draw[gray!70, line width=0.7pt](0,5)--(5,5);
\draw[gray!70, line width=0.7pt](5,0)--(5,5);
\draw[thick,rounded corners,color=blue] (1.5,3)--(1.5,4);
\draw[thick,rounded corners,color=blue] (0.5,2)--(0.5,2.5)--(1,2.5);
\draw[thick,rounded corners,color=blue] (3.5,0)--(3.5,0.5)--(4,0.5);
\draw[thick,rounded corners,color=blue] (0.5,0)--(0.5,1);
\draw[thick,rounded corners,color=blue] (0.5,1)--(0.5,2);
\draw[thick,rounded corners,color=blue] (1,2.5)--(2,2.5);
\draw[thick,rounded corners,color=blue] (1.5,2)--(1.5,3);
\draw[thick,rounded corners,color=blue] (2,2.5)--(3,2.5);
\draw[thick,rounded corners,color=blue] (2.5,2)--(2.5,3);
\draw[thick,rounded corners,color=blue] (1.5,0)--(1.5,1);
\draw[thick,rounded corners,color=blue] (1.5,1)--(1.5,2);
\draw[thick,rounded corners,color=blue] (2.5,0)--(2.5,1);
\draw[thick,rounded corners,color=blue] (2.5,1)--(2.5,2);
\draw[thick,rounded corners,color=blue] (4,0.5)--(5,0.5);
\draw[thick,rounded corners,color=blue] (4.5,0)--(4.5,1);
\draw[white, line width=0.7pt,fill=white] (0,5)--(5,0)--(5,5)--(0,5);
\draw[black, line width=1.6pt] (0,5)--(5,0)--(0,0)--(0,5);
\end{tikzpicture}
\qquad\qquad
\begin{tikzpicture}[x=\bpdwidth,y=\bpdwidth,line cap=round,line join=round,baseline=(z.base)]
\node at (0,2.25) (z) {};
\draw[gray!70, line width=0.7pt](0,0)--(5,0);
\draw[gray!70, line width=0.7pt](0,0)--(0,5);
\draw[gray!70, line width=0.7pt](0,1)--(5,1);
\draw[gray!70, line width=0.7pt](1,0)--(1,5);
\draw[gray!70, line width=0.7pt](0,2)--(5,2);
\draw[gray!70, line width=0.7pt](2,0)--(2,5);
\draw[gray!70, line width=0.7pt](0,3)--(5,3);
\draw[gray!70, line width=0.7pt](3,0)--(3,5);
\draw[gray!70, line width=0.7pt](0,4)--(5,4);
\draw[gray!70, line width=0.7pt](4,0)--(4,5);
\draw[gray!70, line width=0.7pt](0,5)--(5,5);
\draw[gray!70, line width=0.7pt](5,0)--(5,5);
\draw[thick,rounded corners,color=blue] (2.5,4)--(2.5,4.5)--(3,4.5);
\draw[thick,rounded corners,color=blue] (3.5,0)--(3.5,0.5)--(4,0.5);
\draw[thick,rounded corners,color=blue] (4.5,1)--(4.5,1.5)--(5,1.5);
\draw[thick,rounded corners,color=blue] (3.5,2)--(3.5,1.5)--(3,1.5);
\draw[thick,rounded corners,color=blue] (0.5,1)--(0.5,1.5)--(1,1.5);
\draw[thick,rounded corners,color=blue] (1.5,2)--(1.5,2.5)--(2,2.5);
\draw[thick,rounded corners,color=blue] (3.5,3)--(3.5,3.5)--(4,3.5);
\draw[thick,rounded corners,color=blue] (0.5,0)--(0.5,1);
\draw[thick,rounded corners,color=blue] (1,1.5)--(2,1.5);
\draw[thick,rounded corners,color=blue] (1.5,1)--(1.5,2);
\draw[thick,rounded corners,color=blue] (2,1.5)--(3,1.5);
\draw[thick,rounded corners,color=blue] (2.5,1)--(2.5,2);
\draw[thick,rounded corners,color=blue] (3,2.5)--(4,2.5);
\draw[thick,rounded corners,color=blue] (3.5,2)--(3.5,3);
\draw[thick,rounded corners,color=blue] (4,3.5)--(5,3.5);
\draw[thick,rounded corners,color=blue] (1.5,0)--(1.5,1);
\draw[thick,rounded corners,color=blue] (2,2.5)--(3,2.5);
\draw[thick,rounded corners,color=blue] (2.5,2)--(2.5,3);
\draw[thick,rounded corners,color=blue] (4,2.5)--(5,2.5);
\draw[thick,rounded corners,color=blue] (2.5,0)--(2.5,1);
\draw[thick,rounded corners,color=blue] (2.5,3)--(2.5,4);
\draw[thick,rounded corners,color=blue] (3,4.5)--(4,4.5);
\draw[thick,rounded corners,color=blue] (4,4.5)--(5,4.5);
\draw[thick,rounded corners,color=blue] (4,0.5)--(5,0.5);
\draw[thick,rounded corners,color=blue] (4.5,0)--(4.5,1);
\draw[black, line width=1.6pt] (0,5)--(5,5)--(5,0)--(0,0)--(0,5);
\end{tikzpicture}
\mapsto
\begin{tikzpicture}[x=\bpdwidth,y=\bpdwidth,line cap=round,line join=round,baseline=(z.base)]
\node at (0,2.25) (z) {};
\draw[gray!70, line width=0.7pt](0,0)--(5,0);
\draw[gray!70, line width=0.7pt](0,0)--(0,5);
\draw[gray!70, line width=0.7pt](0,1)--(5,1);
\draw[gray!70, line width=0.7pt](1,0)--(1,5);
\draw[gray!70, line width=0.7pt](0,2)--(5,2);
\draw[gray!70, line width=0.7pt](2,0)--(2,5);
\draw[gray!70, line width=0.7pt](0,3)--(5,3);
\draw[gray!70, line width=0.7pt](3,0)--(3,5);
\draw[gray!70, line width=0.7pt](0,4)--(5,4);
\draw[gray!70, line width=0.7pt](4,0)--(4,5);
\draw[gray!70, line width=0.7pt](0,5)--(5,5);
\draw[gray!70, line width=0.7pt](5,0)--(5,5);
\draw[thick,rounded corners,color=blue] (3.5,0)--(3.5,0.5)--(4,0.5);
\draw[thick,rounded corners,color=blue] (3,1.5)--(4,1.5);
\draw[thick,rounded corners,color=blue] (0.5,1)--(0.5,1.5)--(1,1.5);
\draw[thick,rounded corners,color=blue] (1.5,2)--(1.5,2.5)--(2,2.5);
\draw[thick,rounded corners,color=blue] (0.5,0)--(0.5,1);
\draw[thick,rounded corners,color=blue] (1,1.5)--(2,1.5);
\draw[thick,rounded corners,color=blue] (1.5,1)--(1.5,2);
\draw[thick,rounded corners,color=blue] (2,1.5)--(3,1.5);
\draw[thick,rounded corners,color=blue] (2.5,1)--(2.5,2);
\draw[thick,rounded corners,color=blue] (1.5,0)--(1.5,1);
\draw[thick,rounded corners,color=blue] (2,2.5)--(3,2.5);
\draw[thick,rounded corners,color=blue] (2.5,2)--(2.5,3);
\draw[thick,rounded corners,color=blue] (2.5,0)--(2.5,1);
\draw[thick,rounded corners,color=blue] (4,0.5)--(5,0.5);
\draw[thick,rounded corners,color=blue] (4.5,0)--(4.5,1);
\draw[white, line width=0.7pt,fill=white] (0,5)--(5,0)--(5,5)--(0,5);
\draw[black, line width=1.6pt] (0,5)--(5,0)--(0,0)--(0,5);
\end{tikzpicture}
.
\ee
More precisely, this operation retains the tiles strictly below the diagonal
and applies the following transformations on the main diagonal:
\[
\BTILE \mapsto \BDIAG,\qquad
\JTILE \mapsto\JDIAG, \qquad
\CTILE \mapsto\CDIAG,\qquad
\PTILE \mapsto\PDIAG,\qquad
\HTILE\mapsto\JDIAG,\qquad
\VTILE\mapsto\CDIAG
.\]
The analogue of $\kBPD(w)$ for involutions $z \in \Invol$
is easily described in these terms:
we define $\kInvBPD(z)$ to consist of the lower halves of all elements of $\kBPD(z)$
that are \defn{symmetric} under the usual transpose.
For example, the left pipe dream for $z=32154$ in \eqref{inv-bpd-ex} above is symmetric so its lower half belongs to $\kInvBPD(z)$.
We refer to elements of $\kInvBPD(z)$ as \defn{$K$-theoretic bumpless involution pipe dreams} for $z$.

\begin{remark}
The pipe dream for $z=32154$ on the right in \eqref{inv-bpd-ex} 
is not symmetric but its lower half is coincidentally again in $\kInvBPD(z)$. In general, the lower half 
of an asymmetric element of $\kBPD(w)$ for $w \in S_n$ is always an element of  $\kInvBPD(z)$
for some $z\in \Invol$ but one may have $z\neq w$.
\end{remark}

We are more interested here in an involution version of the set $\BPD(w)$.
Defining this requires one extra ingredient.
Let $\cyc(z)$ be the number of nontrivial cycles in $z\in \Invol$,
and set \be\ellinv(z)  =\tfrac{\ell(z)+\cyc(z)}{2}.\ee
This quantity is always an integer, and 
we define $\InvBPD(z)$ to be the set of pipe dreams $D \in \kInvBPD(z)$
whose number of \defn{blank} tiles (meaning tiles equal to $\BTILE$ or $\BDIAG$) is exactly $\ellinv(z)$.
Every element of $\InvBPD(z)$  is \defn{reduced}, again in the sense that no two pipes cross multiple times, where a $\PDIAG$ tile counts as a crossing. However,
$\InvBPD(z)$ may be strictly smaller than the subset of
reduced elements of $\kInvBPD(z)$.
We refer to the objects in $\InvBPD(z)$ as \defn{bumpless involution pipe dreams} for $z$.

\begin{remark}
This formulation is not the same as the first definition of $\InvBPD(z)$ in Section~\ref{bpd-sect},
but is equivalent by Corollary~\ref{k-sub-cor}.
The left triangle  in 
 \eqref{inv-bpd-ex} has $3 = \ellinv(32154)$ blank tiles, so belongs to
 $\InvBPD(z)$ for $z=32154 \in I_5$, 
 while the right triangle has 4 blank tiles, so belongs to
$\kInvBPD(z)\setminus\InvBPD(z)$.
\end{remark}

Let us now explain  the geometric relevance of bumpless involution pipe dreams.
The orthogonal group $\O_n(\CC)$
of $n\times n$ invertible matrices $M$ over $\CC$ with $M^{-1} = M^\top$
acts on the complete flag variety with finitely many orbits,
which are in bijection with $\Invol$.

Wyser and Yong \cite{WyserYong2017} identified a family of polynomials
$\fkSO_z \in \ZZ[x_1,x_2,x_3,\dots]$ indexed by $z \in \Invol$ that represent 
 the   cohomology classes of the closures of the $\O_n(\CC)$-orbits in 
  the coinvariant algebra.
The same polynomials (without any additional $y$-variables) also represent the corresponding
$S$-equivariant cohomology classes, where $S$ is a maximal torus of $\O_n(\CC)$. 
One can write
$ \fkSO_z = 2^{\cyc(z)} \iS_z$
where $\iS_z$ is a certain multiplicity-free sum of single Schubert polynomials $\fkS_w(x;0)$.
The sum $\iS_z$ is called an \defn{involution Schubert polynomial} in \cite{HMP1}
and is defined precisely in Section~\ref{sch-sect}.
Most of our later results will be stated in terms of $\iS_z$, but we phrase things in this introduction  
using $\fkSO_z$.

Our first main result is the following analogue of the left identity in
\eqref{double-sch-eq}:

\begin{mainthm}[{See Theorem~\ref{main-thm}}]
\label{thmA}
If $z \in \Invol$ is any involution then \[
  \fkSO_z = \sum_{D \in \InvBPD(z)} 2^{\cyc(z) - \mathsf{db}(D)} \prod_{(i,j) \in \blank(D)} (x_i +x_j).
\]
where $\mathsf{db}(D)$ is the number of diagonal blank tiles $\BDIAG$ in an involution pipe dream $D$, and $\blank(D)$ is the set of lower triangular positions $(i,j)$ where
$D_{ij} \in \left\{\BTILE,\BDIAG\right\}$.
\end{mainthm}

The polynomials $\fkSO_z$ have a $K$-theoretic generalization $\fkGO_z \in \ZZ[\beta][x_1,x_2,x_3,\dots]$ which we refer to as \defn{orthogonal Grothendieck polynomials}, following \cite{MarPawK}. 
As the name suggests, these objects are equivariant $K$-theory representatives (when $\beta=-1$) for the closures of the $\O_n(\CC)$-orbits in the complete flag variety. They recover $\fkSO_z$ when $\beta=0$.
The geometric definition for the polynomials $\fkGO_z$ is difficult to compute for arbitrary involutions $z \in \Invol$, 
since the corresponding $\O_n(\CC)$-orbit closures can fail to have rational singularities \cite{WyserYong2017}.
In the case when $z \in \Invol$ is vexillary, however,
 this obstruction does not occur and $\fkGO_z$ can be specified by 
an explicit set of divided difference equations \cite{MarPawK}; see Section~\ref{sch-sect}.

Some formulas have been previously identified for $\fkGO_z$ when $z \in \Invol$ is vexillary \cite{MarPawK,MarWen,WyserYong2017}, but these have been quite complicated outside the case when $z$ is also \defn{dominant} ($132$-avoiding).
Our second main result gives an explicit formula for  $\fkGO_z$
in the vexillary case as a finite sum of   weights over the same   pipe dreams that occur in the cohomological formula. 
This provides an involution analogue of \eqref{double-sch-eq2}:

\begin{mainthm}[See Theorems~\ref{fkGO-conj} and \ref{vex-prop1}]
\label{thmB}
If $z \in \Invol$ is vexillary then
$ \InvBPD(z) = \kInvBPD(z)$ and  \[
\fkGO_z = \sum_{D \in \InvBPD(z)}  \prod_{(i,j)\in \blank(D)} (x_i\oplus x_j) \prod_{(k,l)\in\upturn(D)} (1+\beta x_k)(1+\beta x_l)
\]
where $\upturn(D)$ is the set of  lower triangular positions $(k,l)$
 where $D_{kl}\in \left\{\JTILE,\JDIAG\right\}$.
\end{mainthm}

If $D \in \InvBPD(z)$ then
we always have $\cyc(z) \geq \mathsf{db}(D)$, with equality if $z \in \Invol$ is vexillary (see  
Corollary~\ref{db-cor}
and Theorem~\ref{vex-prop1}). Thus the formula  for $\fkGO_z$ just given
specializes  to the expression for $ \fkSO_z$ in Theorem~\ref{thmA} when $\beta=0$.

Involution Schubert polynomials have a second $K$-theoretic deformation,
which conjecturally
satisfies an analogue of the 
formula for $\fkG_w(x;y)$ in  \eqref{double-sch-eq}.
For any $z \in \Invol$ define 
\[ \iG_z = \sum_{\substack{w \in S_n \\ w^{-1} \circ w = z}} \beta^{\ell(w) -\ellinv(z)} \fkG_w(x;0)\] where $\circ$ denotes the \defn{Demazure product} for $S_n$; see Section~\ref{inv-sect}.
These \defn{involution Grothendieck polynomials} specialize to $\iS_z = \widehat{\mathfrak{G}}^{(0)}_z$ but
 currently do not have any known geometric interpretation. 
 Moreover, unlike with $\fkSO_z= 2^{\cyc(z)}\iS_z$, there is no simple scalar 
 relationship between $\fkGO_z$ and $\iG_z$; for example,
 \[
 \fkGO_{21} = 2x_1 + \beta x_1^2 \quad\text{while}\quad
 \iG_{21} = x_1.\]
 Nevertheless, involution Grothendieck polynomials have nice algebraic properties and are involved in some formulas for $\fkGO_z$ \cite{MarbergPawlowski,MarWen}.
The following has been tested for $n\leq 8$:

\begin{conjecture*}[{See Conjecture~\ref{iG-conj}}]
If $z \in \Invol$ is any involution then \[
  \iG_z = \sum_{D \in \kInvBPD(z)} \beta^{|\blank(D)| -\ellinv(z)}\prod_{(i,j) \in \blank(D)} 
  x_{i\oplus j}
 \prod_{(k,l) \in \upturn(D)} 
 (1+\beta x_{k\oplus l})
\]
where $x_{i\oplus j}$ is defined to be $x_i $ when $i=j$
and to be $x_i \oplus x_j = x_i + x_j + \beta x_ix_j$ when $i\neq j$.
\end{conjecture*}

\begin{remark}
As mentioned earlier, 
double Schubert polynomials (and Grothendieck polynomials)
have two pipe dream formulas:
the bumpless version from \cite{LLS,Weigandt} and the classical formula
from \cite{BergeronBilley, fomin1994grothendieck,KnutsonMiller}.
The same is true of involution Schubert polynomials.
 The results here supply a bumpless formula, while a ``classical'' involution   pipe dream formula for $\iS_z$
appears in \cite{HMP6}.

While the latter formula has an extension that recovers $\iG_z$ \cite[Thm.~4.5]{MarbergPawlowski}, we currently do not know of any ``classical'' pipe dream formula for $\fkGO_z$. 
 This discrepancy may be related to the difficulty of computing 
 Gr\"obner bases for \defn{symmetric matrix Schubert varieties} in anti-diagonal term orders; see \cite{MarPaw2024}.
 On the other hand, the pipe dream formula for $\fkSO_z=2^{\cyc(z)}\iS_z$ in Theorem~\ref{thmA} does extend to $\fkGO_z$ when $z$ is vexillary.
 Another useful feature of our new result is that it can easily be generalized to a formula for
the \defn{backstable} version of $\iS_z$ considered in \cite{MarbergPawlowski,Pawlowski2021Universal}. As this is very similar to the permutation setting, the present article does not discuss this application in detail.
 \end{remark}

Our proof strategies for Theorems~\ref{thmA} and \ref{thmB} are quite different.
For the cohomological formula, we work with the normalization 
$\iS_z = 2^{-\cyc(z)}\fkSO_z$ and  
show that its target pipe dream expression
satisfies the same \defn{transition equations}
as those identified for involution Schubert polynomials in \cite{HMP3}.
This is done by exhibiting a specific bijection
between certain sets of bumpless involution pipe dreams in Section~\ref{S: Proof}. This approach is similar to how Weigandt proves the right side of \eqref{double-sch-eq} in \cite{Weigandt}, but more technically involved.

We prove our formula for $\fkGO_z$ by using a bijection of Weigandt \cite{Weigandt} to show that it is equivalent to a generating function over 
certain \defn{set-valued shifted plane partitions}. The precise statement 
(Theorem~\ref{fkGO-conj2}) is the shifted analogue
of a set-valued tableau formula for $\fkG_w(x;y)$
when $w \in S_n$ is vexillary from \cite{KnutsonMillerYong}.
Knutson, Miller, and Yong prove the latter formula
by identifying $\fkG_w(x;y)$ after a change of variables with the $K$-polynomial
of a \defn{matrix Schubert variety}, and then computing a diagonal term order Gr\"obner basis for the defining ideal of this variety when $w$ is vexillary. This Gr\"obner basis 
leads to an alternating-sum formula for $\fkG_w(x;y)$ in the vexillary case which then can be simplified using the combinatorics of \defn{subword complexes} from
\cite{KnutsonMillerSubword}.

We can follow the same approach to derive 
our expression for $\fkGO_z$ by considering \defn{symmetric matrix Schubert varieties} rather than classical matrix Schubert varieties. 
In particular, we use a result from
\cite{GMN} to identify a Gr\"obner basis in any diagonal term order 
for the ideals of symmetric matrix Schubert varieties of vexillary involutions (Theorem~\ref{gr-thm1}). 
This shows that these ideals are prime and have square-free initial ideals 
(Corollary~\ref{sym-i-msv-cor}), which partially resolves some
conjectures about 
symmetric matrix Schubert varieties from \cite{MarPaw2024}.
We then employ similar techniques involving subword complexes
to attain the formula in Theorem~\ref{thmB}.

 \begin{remark}
 When $n$ is even, the symplectic group $\Sp_n(\CC)$
 also acts on the complete flag variety
 with finitely many orbits, now indexed by the fixed-point-free elements of $\Invol$. This leads to a parallel theory of \defn{fpf-involution Schubert polynomials}, \defn{symplectic Grothendieck polynomials}, and \defn{skew-symmetric matrix Schubert varieties}.
 These objects are often considered alongside the involution case \cite{HMP1,HMP3,HMP6,MarPawK,WyserYong2017},
 and the relevant polynomials also possess bumpless pipe dream formulas.
Such formulas are not discussed here since their derivation has been
independently carried out in concurrent work of Hamaker and Weigandt \cite{HamakerWeigandt}.
 \end{remark}
 
\subsection{Applications and outline}

Besides the two main theorems in the previous section,
we also prove some other results. For example, we explain how to
systematically generate the set $\InvBPD(z)$ in Section~\ref{gen-sect}. 
Specifically, Theorem~\ref{droop-thm} shows that this can be done by applying
a sequence of \defn{simple} and \defn{complex droop moves} to an initial \defn{involution Rothe pipe dream}.
Simple droop moves have straightforward analogues in the permutation case \cite{LLS}, but our complex droop moves do not have any direct counterpart in that setting.

For vexillary involutions $z \in \Invol$, our bumpless pipe dreams 
are closely related to shifted plane partitions. Section~\ref{vex-sect}
explores this connection. There, we describe a bijection
between $\InvBPD(z)$ and flagged shifted plane partitions of a certain shape (Theorem~\ref{bij-thm}). This correspondence extends to a weight- and degree-preserving bijection between \defn{marked bumpless involution pipe dreams} and \defn{set-valued shifted plane partitions} (Theorem~\ref{upsilon-thm}).

Such results let us show that the number of shifted plane partitions
with a given shape and entry set 
is equal to $|\InvBPD(z)| = \fkSO_z(\frac{1}{2},\frac{1}{2},\frac{1}{2},\dots)$
for a particular choice of  vexillary involution $z$
(Theorem~\ref{thm: shifted bijection inv bps and plane partitions}).
Using a result of Proctor \cite{Proctor}, we are then able to confirm some product formulas
for principal specializations of involution Schubert polynomials conjectured in \cite{HMP6} (see Corollary~\ref{cor:plane partition product}).

We also obtain a shifted analogue of the Fomin--Kirillov reduced-word formula for bounded plane partitions. Suppose $\lambda$ is a strict partition and $\SPP_{\lambda}(t)$ is the set of shifted plane partitions of shape $\lambda$ with entries in $\{0,1,\ldots,t\}$. Then Theorem~\ref{fmp-thm} expresses $|\SPP_{\lambda}(t-1)|$ as a weighted sum over reduced words of inverse atoms of a dominant involution, showing that
$
|\SPP_{\lambda}(t-1)|\in \QQ_{\geq0}[t].
$
This proves the straight-shape case of a conjecture of Ferroni, Morales, and Panova \cite[Conj.~7.10]{FMP}, and implies Ehrhart positivity for the order polytopes of shifted straight shapes.

The rest of this paper is organized as follows.
Section~\ref{prelim-sect} collects some necessary preliminaries. Section~\ref{bpd-sect} precisely defines our various kinds of bumpless pipe dreams, and includes our main results. Sections~\ref{gen-sect} and \ref{S: Proof} develop the local moves and transition bijection used to prove Theorem~\ref{thmA}. Section~\ref{vex-sect} considers vexillary involutions and their connections to shifted plane partitions. Section~\ref{final-sect} discusses some future directions related to alternating sign matrices and asymptotics of principal specializations. Appendix~\ref{app-sect2}, finally, proves the orthogonal Grothendieck formula in Theorem~\ref{thmB}.

\subsection*{Acknowledgments}

We thank Dave Anderson, Zachary Hamaker, Ilse Fischer, Brendan Pawlowski, and Anna Weigandt
for helpful conversations.

\subsubsection*{Funding Sources}
The first named author was partially supported by Hong Kong RGC grants 16304122 and 16304625. The second named author was partially supported by the NSERC Discovery grant RGPIN-2024-06246.

\subsubsection*{Statement on AI Use} GPT-6 Astra was used to assist with the preparation of this article, and played a significant role
in discovering how to prove the results  in 
Section~\ref{ice-sect} and Appendix~\ref{app-sect2}.
The arguments and text throughout this work have been reviewed and substantially revised by the authors, who are responsible for the final manuscript.

\section{Preliminaries}\label{prelim-sect}

This section contains some minimal preliminaries
on Coxeter groups, permutations, and (involution) Schubert polynomials that will be needed in later sections.
 
\subsection{Involutions in Coxeter groups}\label{inv-sect}

Let $W$ be a Coxeter group with simple generating set $S$
and length function $\ell : W \to \NN$.
The \defn{Demazure product} is the unique associative operation
$\circ : W \times W \to W$ that has $s\circ s = s$ for all $s \in S$ and $u\circ v = uv$ for all $u,v \in W$ with $\ell(uv) = \ell(u) + \ell(v)$. 
Let $\I(W) = \{ z \in W : z=z^{-1}\}$. If $z \in \I(W)$ then  
\be\label{szs-eq}
s\circ z \circ s = \begin{cases}
sz &\text{if $\ell(sz) > \ell(z)$ and $szs=z$} \\
szs &\text{if $\ell(sz) >\ell(z)$ and $szs\neq z$} \\
z &\text{if $\ell(sz) < \ell(z)$}
\end{cases}
\qquad\text{for all $s \in S$,}
\ee
so it follows by induction on length that $\I(W) = \{ w\circ w^{-1} : w \in W\}$. Define 
\be\label{ellinv-def}
\cH(z) = \{ w \in W : w\circ w^{-1}=z\}\quand
\ellinv(z) = \min \{ \ell(w) : w \in \cH(z)\}
\ee
for each $z \in \I(W)$,
and then set
\be
\label{cW-eq}
\cW(z) = \{ w \in \cH(z) : \ell(w) =\ellinv(z)\}.\ee
This notation follows \cite{CJW}.
We refer to the elements of $\cW(z)$ as \defn{(inverse) atoms} for $z$, since 
they are the inverses of the elements of the sets $\cA(z)$ considered in \cite{HMP1,HMP3,HMP6,Mar2019}.

\subsection{Atoms in permutation groups}\label{Sn-sect}

Fix an integer $n>0$, let $[n] = \{1,2,\dots,n\}$,
and suppose $T = \{a_1<a_2<\dots<a_n\}$ is a totally ordered set of size $n$.
Define $S_T$ to be the group of bijections $T \to T$.
For each $i \in [n-1]$ let $s^{T}_i \in S_T$ be the bijection that interchanges $a_i$ and $a_{i+1}$ while fixing the rest of $T$.
Then $S_T$ is a Coxeter group with simple generators
$$\left\{ s^{T}_1,\ s^{T}_2,\ \dots,\ s^{T}_{n-1}\right\}.$$ Let
$
I_T = \I(S_T) = \{ z \in S_T : z=z^{-1}\}
$
be the set of involutions in this group.
We write $S_n$, $\Invol$, and $s_i$ in place of $S_T$, $I_T$, and $s^{T}_i$ when $T = [n]$.

It is well-known that the Coxeter length function $\ell : S_T \to \NN$ is
\be\label{inv-def}
\ell(w) = |\Inv(w)|\quad\text{where }\Inv(w) = \{(i,j) \in T \times T : i < j \text{ and }w(i) > w(j)\}.
\ee
A \defn{reduced word} for $w \in S_T$ is a sequence $(i_1,i_2,\dots,i_k)$ where $i_j \in [n-1]$ and $k=\ell(w)$
and $w = s^T_{i_1}s^T_{i_2}\cdots s^T_{i_k}$.
The elements of $I_T$ are the permutations whose cycles all have length at most two.
One has
\be\label{cyc-def}
\ellinv(z) = \tfrac{\ell(z) + \cyc(z)}{2}
\ee
for all $z \in I_T$,
where $\Cyc(z) = \{ (b,a) \in T\times T : a < b=z(a)\}$ and $\cyc(z) = |\Cyc(z)|$
 \cite[\S2.1]{HMP1}.

The \defn{one-line representation} of an element $w \in S_T$ is the sequence $w(a_1)w(a_2)\cdots w(a_n)$.
Any word $b_1b_2\cdots b_n$ with $\{a_1,a_2,\dots,a_n\} = \{b_1,b_2,\dots,b_n\}$ is the one-line representation of the permutation $w \in S_{T}$ that sends $a_i\mapsto b_i$ for each $i \in [n]$.
We frequently identify permutations with their one-line representations. This loses no information since 
any finite word with all distinct integer letters is
the one-line representation of a unique element of $S_T$ for a unique choice of $T$.

For distinct $a,b\in T$, let $t_{ab}\in S_T$ be the transposition
swapping $a$ and $b$.  Left multiplication by $t_{ab}$ swaps these
values in the one-line representation.
We write $u\lessdot v$ when $v$ covers $u$ in Bruhat order.
For $a<b$, the relation $u\lessdot t_{ab}u$ holds if and only if in the one-line representation,
$a$ appears before $b$ in $u$ and no $t\in T$ with $a<t<b$ appears between them.

For Coxeter groups of the form $W=S_T$, the sets $\cW(z)$ 
from \eqref{cW-eq} may be described in much more explicit terms. 
First, consider two elements $v,w \in S_T$ and set $v_i = v(a_i)$ and $w_i = w(a_i)$ for $ i \in [n]$.
We write $v\cWsim w$ if there is an index $ j \in [n-2]$ and integers $A<B<C$ such that 
\be
\{\ v_jv_{j+1}v_{j+2},\ w_jw_{j+1}w_{j+2}\ \} = \{\ CAB,\ BCA\ \}\quand v_i = w_i\text{ for all }i \notin \{j,j+1,j+2\}.
\ee
Extend $\cWsim$ to an equivalence relation on $S_T$ by taking the reflexive and transitive closure.

\begin{lemma}[{\cite[Thms.~6.4 and 6.10]{HMP2}}]
\label{CAB-lem}
If $z \in I_T$ then $\cW(z)$ is a single equivalence class under $\cWsim$. Additionally, no $w \in \cW(z)$ has any consecutive $321$-patterns.
Finally, if $v,w \in S_T$ have $v\cWsim w$ then $v\circ v^{-1} = w\circ w^{-1}$
where $\circ $ is the Demazure product for $S_T$.

\end{lemma}

The second claim in the lemma follows from \cite[Thm.~6.4]{HMP2}
since that result implies that if $w\in S_T$ has a consecutive $321$-pattern,
then there exists a shorter permutation $v\in S_T$ with $v\circ v^{-1} = w\circ w^{-1}$.

\begin{example}
For $T=\{2,4,6,8\}$ we have $\cW(8642) = \{8264,\ 6824,\ 6482\}$.
\end{example}

Continue to let $w \in S_T$ be an arbitrary element and define $w_i = w(a_i)$ as above.
An index $ i \in [n-1]$ is a \defn{descent} of $w$ if $w_i > w_{i+1}$. Define the \defn{nested descent set} of $w$ recursively by the formula
\[
\NDes(w) = \begin{cases} 
\varnothing &\text{if $w_1<w_2<\dots<w_n$} \\ 
\{ (w_i,w_{i+1})\} \sqcup \NDes(w_1\cdots w_{i-1}w_{i+2}\cdots w_n) &\text{if $i$ is the smallest descent of $w$.}
\end{cases}
\]
Similarly define the \defn{nested residue set} of $w$ by 
\[
\NRes(w) = \begin{cases} 
\{w_1,w_2,\dots,w_n\} &\text{if $w_1<w_2<\dots<w_n$} \\ 
 \NRes(w_1\cdots w_{i-1}w_{i+2}\cdots w_n) &\text{if $i$ is the smallest descent of $w$.}
\end{cases}
\]
Write $z_T(w) \in I_T$ for the involution with  $\Cyc(z_T(w))=\NDes(w)$. Finally, let 
\be
\textstyle \cW_T = \bigsqcup_{z \in I_T} \cW(z).
\ee
To interpret the following statements,
notice that if $a,b \in T$ and $w\in S_T$ then the inequality $w^{-1}(a) < w^{-1}(b)$
is equivalent to saying that $a$ appears to the left of $b$ in the one-line representation of $w$.

\begin{lemma}
\label{cW-pre-lem}
Let $w \in S_T$ and suppose $(b,a),(y,x) \in \NDes(w)$.
Then it is impossible to have both
\[ b>y>x>a \quand w^{-1}(b) < w^{-1}(y) < w^{-1}(x) < w^{-1}(a).\]
\end{lemma}

\begin{proof}
Suppose $w$ is a counterexample to this claim with $T$ as small as possible. 
By hypothesis $\NDes(w)$ is nonempty so there exists a smallest descent $i$
with $w_i>w_{i+1}$. This descent cannot have $w_i=b$ or $w_i = y$
since  $y$ appears between
$b$ and $a$ in $w$ while we have $b>y$ and  $(b,a) \in \NDes(w)$.
We conclude that 
 $\{w_i,w_{i+1}\} \cap\{ a,b,x,y\} =\varnothing$.
 But then the permutation
 $\hat w = w_1\cdots w_{i-1} w_{i+2}\cdots w_n \in S_{T - \{w_i,w_{i+1}\}}$
 is also a counterexample to the desired claim.
 This contradicts the minimality of $w$, so the lemma must hold.
\end{proof}

\begin{theorem}[{\cite[Thm.~2.5]{CJW}}]
\label{cW-lem}
Let $w \in S_T$. Then $w \in \cW_T$ if and only if the following conditions hold:
\ben
\item[(a)] If $a \in \NRes(w)$ and $(c,b) \in\NDes(w)$ have $a<b<c$ then $w^{-1}(a) < w^{-1}(c) $.

\item[(b)] If $(b,a) \in\NDes(w)$ and $c \in \NRes(w)$ have $a<b<c$ then $ w^{-1}(a) < w^{-1}(c)$.

\item[(c)] If $(b,a),(y,x) \in\NDes(w)$ have $a<x$ and $b<y$ 
then $ w^{-1}(a) < w^{-1}(y) $.

\een
Moreover, if these properties hold then $w \in \cW(z)$ for $z =z_T(w) \in I_T$.
\end{theorem}

\begin{proof}
A brief argument is required to derive this statement from the cited result in \cite{CJW}.
First assume that properties (a), (b), and (c) hold. 
We claim that the following additional properties then also hold:
\ben
\item[(d)] If $a,b \in \NRes(w)$ have $a<b$ then $w^{-1}(a) < w^{-1}(b)$.
\item[(e)] If $(b,a) \in \NDes(w)$ then  $w^{-1}(b) < w^{-1}(a)$ and 
no $x$ with $a<x<b$ has $w^{-1}(b)<w^{-1}(x)<w^{-1}(a)$.
\een
Property (d) and the inequality $w^{-1}(b) < w^{-1}(a)$ in property (e) are both clear by definition.
The rest of (e) 
holds because if $(b,a) \in \NDes(w)$ and some $x$ with $a<x<b$ has $w^{-1}(b)<w^{-1}(x)<w^{-1}(a)$,
then there must exist $(y,x) \in \NDes(w)$ with $w^{-1}(b) < w^{-1}(y) < w^{-1}(x) < w^{-1}(a)$. 
 This is impossible if $a<x<y<b$ by Lemma~\ref{cW-pre-lem}
 and violates condition (c) if $a<x<b<y$.
Thus, our claim follows.

The cited result \cite[Thm.~2.5]{CJW} asserts that if (a), (b), (c), (d), and (e) all hold then 
$w \in \cW(z)$ for $z=z_T(w) \in I_T$.
Conversely, if $w \in \cW(z)$ for some $z \in I_T$,
then   $\NDes(w) = \Cyc(z)$ is the set of non-trivial cycles of $z$
and   $\NRes(w)$ is the set of fixed points of $z$ by \cite[Thm.~7.4]{Mar2019}.
Once we have these identifications, 
 properties (a), (b), and (c) are again consequences of \cite[Thm.~2.5]{CJW}.
\end{proof}

Given $w \in S_T$ and any subset $E\subseteq T$, let 
$
\std{w}{E} \in S_E 
$
be the permutation whose one-line representation is obtained from that of $w$ by omitting all letters not in $E$. Alternatively,
when $E$ and $F$ are totally ordered finite sets with $|E|=|F|$,  write $\phi_{(E\to F)}$ for the order-preserving bijection $E\to F$. Then we have
\be
\std{w}{E} =  w \circ \phi_{(E\to F)} \in S_E\quad\text{for }F=w^{-1}(E),
\ee
where $\circ$ denotes function composition.
When $w(E) = E$, this formula is just the restriction of $w$ to the domain $E$.
In this case the
 cycle representation of $\std{w}{E}$ is obtained 
from that of $w$ by retaining only the cycles contained in $E$.
An equivalent form of the following result was proved in \cite{HMP3}.

\begin{lemma}[{\cite[Cor.~3.19]{HMP3}}]
\label{L: W reduction}
Let $z \in I_T$ and $w \in S_T$. Then the following are equivalent:
\ben
\item[(a)] $w \in \cW(z)$.

\item[(b)] $\std{w}{E} \in \cW(\std{z}{E})$ for all subsets $E\subseteq T$ with $z(E) = E$.

\item[(c)] $\std{w}{E}\in \cW(\std{z}{E})$ for all subsets $E\subseteq T$ of the form $E=\{a,b,z(a),z(b)\}$ for some $a,b \in T$.

\een
\end{lemma}

The \defn{Rothe diagram} of $w \in S_T$ is the set of pairs 
\be
D(w)= \left\{(i,j) \in T\times T:
j < w(i)\text{ and }i < w^{-1}(j)\right\}.
\ee
If $T = [n]$ then $D(w)$ is the complement of the set of positions in an $n\times n$ matrix 
that are weakly below (in the same column) or to the right of (in the same row) some pair $(i,w(i))$ for $i \in [n]$. For example,
\be\label{dom-eq}
D(53124) = \left[\begin{smallmatrix} 
\square & \square & \square & \square & 1 \\
\square & \square & 1 & \cdot & \cdot \\
1 & \cdot & \cdot & \cdot & \cdot \\
\cdot & 1 & \cdot & \cdot & \cdot \\
\cdot & \cdot & \cdot & 1 & \cdot  \end{smallmatrix} \right]= \{ (1,1),(1,2),(1,3),(1,4),(2,1),(2,2)\}.
\ee
If $E\subseteq T$ satisfies $w(E)=E$, then 
\begin{equation}
\label{EQ: Rothe pattern}
D(w\cap E)=D(w)\cap(E\times E).
\end{equation}

For each $z \in I_T$, let 
$\iD(z) =\{ (i,j) \in D(z) : i \geq j\}$ be the subset of positions in the usual Rothe diagram that 
are weakly below the main diagonal.
Then it is known \cite[\S3.2]{HMP1}
that
\be
\ell(w) = |D(w)|\qquand \ellinv(z) = |\iD(z)|.
\ee
Notice also that if $J\subseteq T$ satisfies $z(J)=J$, then
$
\iD(\std{z}{J})
=
\iD(z)\cap(J\times J).
$

\subsection{Schubert and Grothendieck polynomials}\label{sch-sect}

This section reviews the definitions of Schubert polynomials,
Grothendieck polynomials, and their involution versions.
Although we discussed double Schubert and Grothendieck polynomials in the introduction, 
here we will only consider their single versions
$\fkS_w = \fkS_w(x;0)$ and $\fkG_w =\fkG_w(x;0)$.

A \defn{partition} is a weakly decreasing sequence of 
nonnegative integers $\lambda=(\lambda_1\geq \lambda_2\geq \dots)$
with finite sum. A permutation is \defn{dominant} if its Rothe diagram is equal to the \defn{Young diagram}
\be
\D_\lambda := \{ (i,j) \in \PP\times \PP : j \leq \lambda_i\}
\ee
for some partition $\lambda$. In this case, the relevant permutation is said to be dominant of shape $\lambda$. 
Let $\delta = (n-1,\dots,3,2,1)$ and write $\lambda \subseteq \delta$ if $\D_\lambda \subseteq \D_\delta$.
For each partition $\lambda \subseteq \delta$,
there exists a unique dominant permutation of shape $\lambda$.
The example in \eqref{dom-eq} shows this permutation when $n=5$ and $\lambda =(4,2)$.

Given $i \in[n-1] $ and $ f\in \ZZ[\beta][x_1,x_2,\dots,x_n]$ let $s_if$ be the polynomial obtained from $f$ by swapping the variables $x_i$ and $x_{i+1}$, and then consider the \defn{divided difference operators} given by
\be
  \partial_i f = \frac{f -s_i f}{x_i-  x_{i+1}}
 \qquand 
 \bpartial_if = -\beta f + (1+\beta x_i) \partial_i f.
 \ee
Both $\partial_i$ and $\bpartial_i$ are $\ZZ[\beta]$-linear operators on $ \ZZ[\beta][x_1,x_2,\dots,x_n]$
which satisfy the usual braid relations for the symmetric group.
One has $\partial_i\partial_i = 0$ while $\bpartial_i\bpartial_i = -\beta \bpartial_i$.
The following definitions are well-known:

\begin{definition}[See \cite{fomin1994grothendieck,Manivel}]
The \defn{Schubert polynomials} $\{\fkS_w\}_{w \in S_n}$
and \defn{Grothendieck polynomials} $\{\fkG_w\}_{w\in S_n}$
are the unique elements of $\ZZ[\beta][x_1,x_2,\dots,x_n]$ indexed by $S_n$ with the following properties: 
\ben
\item[(a)] If $w$ is dominant then $\fkS_w=\fkG_w = \prod_{(i,j) \in D(w)} x_i$.
\item[(b)] If $i \in [n-1]$ has $w\circ s_i \neq w$ then $\partial_i \fkS_{w\circ s_i} = \fkS_{w}$ and $\bpartial_i \fkG_{w\circ s_i} = \fkG_{w}$.
\een
\end{definition}

Each $\fkS_w$ is homogeneous of degree $\ell(w)$ and is obtained from $\fkG_w$ by setting $\beta=0$.
Less obviously, we always have $\fkS_{w\times 1} = \fkS_{w}$ and $\fkG_{w\times 1} = \fkG_w$
where $w\times 1 \in S_{n+1}$ denotes the image of $w$ under the injective group homomorphism 
$S_n \to S_{n+1}$ sending $s^{[n]}_i \mapsto s^{[n+1]}_i$ for $i \in [n-1]$.

A dominant permutation of shape $\lambda$
is an involution if and only if $\lambda=\lambda^\top$ is \defn{symmetric} in the sense of 
being fixed by the transpose operation.
Thus, for any symmetric partition $\lambda =\lambda^\top\subseteq \delta$ there exists 
a unique $z \in \Invol $ that is dominant of shape $\lambda$.
In the following definition, let $a\oplus b = a + b + \beta ab$. Also define
\be\label{i-oplus-j-eq}
x_{(i,j)} = \begin{cases} x_i &\text{if }i =j \\
x_i + x_j &\text{if }i\neq j
\end{cases}
\quand
x_{i\oplus j} = \begin{cases} x_i &\text{if }i =j \\
x_i \oplus  x_j &\text{if }i\neq j.
\end{cases}
\ee

\begin{definition}[See \cite{MarWen,WyserYong2017}]
The \defn{involution Schubert polynomials} $\{\iS_z\}_{z \in \Invol }$
and \defn{involution Grothendieck polynomials} $\{\iG_z\}_{z\in \Invol }$
are the unique elements of $\ZZ[\beta][x_1,x_2,\dots,x_n]$ indexed by $\Invol $ such that:
\ben
\item[(a)] If $z$ is dominant then $\iS_z=\prod_{(i,j) \in \iD(z)} x_{(i,j)}$
and $\iG_z = \prod_{(i,j) \in \iD(z)} x_{i\oplus j}$.
\item[(b)] If $i \in [n-1]$ has $s_i \circ z \circ s_i \neq z$ then $\partial_i \iS_{s_i \circ z \circ s_i} = \iS_{z}$ and $\bpartial_i \iG_{s_i \circ z \circ s_i} = \iG_{z}$.
\een
\end{definition}

It once again holds that  $\iG_z$ becomes $\iS_z$ on setting $\beta=0$.
These polynomials have an alternate formula that shows their close relationship
to Schubert and Grothendieck polynomials. Namely \cite[Thms.~3.3 and 3.6]{MarWen}, it holds for all $z \in \Invol $ that
\be\label{atom-eq}
\textstyle
\iS_z = \sum_{w \in \cW(z)} \fkS_{w^{-1}}
\qquand
\iG_z = \sum_{w \in \cH(z)} \beta^{\ell(w) - \ellinv(z)} \fkG_{w^{-1}}.
\ee 
There is a notable product formula for the sum of all (involution) Grothendieck polynomials:
\begin{proposition}[\cite{fomin1994grothendieck}]
\label{prop:partition function inv Grothendieck}
If $n \in \PP$ 
then $
\sum_{w \in S_n} \beta^{\ell(w)} \fkG_w = \sum_{z \in \Invol } \beta^{\ellinv(z)} \iG_z =\prod_{i =1}^{n-1} (1+\beta x_i)^{n-i}
$.
\end{proposition}

\begin{proof}
The first two sums are equal since every $w \in S_n$ belongs to $\cH(z)$ for a unique $z \in \Invol$.
Equality is obtained between the  first and third expressions by setting $u_i=\beta$ in \cite[Eq.~(2.2) and Thm.~2.3]{fomin1994grothendieck}.
\end{proof}

Schubert polynomials (respectively, Grothendieck polynomials when $\beta=-1$) are meaningful in geometry as
cohomology (respectively, $K$-theory)
representatives for the \defn{Schubert varieties} in the type A complete flag variety \cite{LasSch}.
Similarly, 
the orthogonal group $\O_n(\CC)$ acts on the type A complete flag variety,
and the orbits of this action are in bijection with $\Invol $.
The involution Schubert polynomials $\iS_z$ (after being rescaled by the factor $2^{\cyc(z)}$) represent the cohomology classes of 
the closures of these orbits  \cite{WyserYong2017}. 

The $K$-theory classes of the closures of the $\O_n(\CC)$-orbits in the complete flag variety also have polynomial representatives \cite{MarPawK},
but these are distinct from  $\iG_z$. They are instead defined in the following way.
This construction only applies 
to the orbit closures of involutions $z \in \Invol $ that 
are
\defn{vexillary} in the sense of being $2143$-avoiding,
meaning we never have $z(b)<z(a) < z(d)<z(c)$ 
for any $1\leq a<b<c<d\leq n$.
Let $\Ivex\subseteq \Invol$ be the set of such involutions. 
Every dominant permutation is also vexillary \cite[\S2.2.1]{Manivel}.

\begin{definition}[See \cite{MarPawK,MarWen}]
The \defn{orthogonal  Grothendieck polynomials} $\{\fkGO_z\}_{z\in \Ivex}$
are the unique elements of $\ZZ[\beta][x_1,x_2,\dots,x_n]$ indexed by $\Ivex $ such that:
\ben
\item[(a)] If $z$ is dominant then $\fkGO_z = \prod_{(i,j) \in \iD(z)} x_i \oplus x_j$.
\item[(b)] If $i \in [n-1]$ has $z\neq  s_i \circ z \circ s_i = s_iz s_i  \in \Ivex$ then $\bpartial_i \fkGO_{s_i  z  s_i} = \fkGO_{z}$.
\een
\end{definition}

The polynomials $\fkS_w$ and $\fkG_w$ have two  ``pipe dream'' generating function formulas \cite{BJS,LLS},
 where a pipe dream is a certain kind of wiring diagram defined in the next section.
 For one version of these formulas, there are known involution analogues that give the polynomials $\iS_z$ and $\iG_z$ \cite{HMP6,MarbergPawlowski}.
In this paper, we derive involution versions of the other ``bumpless'' pipe dream formula for $\fkS_w$.

\subsection{Involution Schubert transitions}
\label{SS: involution transition}

This section recalls the transition formulas for involution Schubert polynomials from \cite{HMP3} that we will use in our proof of Theorem~\ref{main-thm}.
These formulas give a recursive
way to compute all 
involution Schubert polynomials.

Throughout, we fix a finite set $T \subset \PP$.
We first recall certain Bruhat covering operators for involutions.
For this, define the \defn{arc diagram} of an involution in $I_T$ to be the matching on $T$ that has edges connecting the two elements in each of the involution's non-trivial cycles.

\begin{definition}[{See \cite[Table 1]{HMP3}}]
\label{D: tau}
For $i,j \in T$ with $i<j$, 
we define $\tau_{ij}:I_T\to I_T$ as follows. 
When applied to $z \in I_T$, the operator $\tau_{ij}$ only modifies cycles
involving $i$ or $j$, and depends only on
the relative order of the elements of $A=\{i,j,z(i),z(j)\}$. The set $A$
has size 2, 3, or 4, and we denote its elements as either $a<b$ or $a<b<c$ or $a<b<c<d$. 

Now, specifically, $\tau_{ij}(z) \in I_T$ is defined to be the involution that takes the same values as $z$ at all points of $T\setminus A$,
and whose arc diagram restricted to $A$ is formed 
from that of $z$ according to the following cases:
\begin{itemize}
\item[(1)] 
$
\begin{tikzpicture}[baseline=-2pt,scale=0.6]
\node at (0,0) {$a$};
\node at (1,0) {$b$};
\end{tikzpicture}
\rightarrow
\begin{tikzpicture}[baseline=-3pt,scale=0.6]
\node at (0,0)  {$a$};
\node at (1,0)  {$b$};
\draw [transform canvas={yshift=5pt}] (0,0) to[out=60,in=120] (1,0);
\end{tikzpicture}
\qquad\text{for }(i,j)=(a,b)$ when $A=\{a<b\}$.

\item[(2)]
$
\begin{tikzpicture}[baseline=-2pt,scale=0.6]
\node at (0,0) {$a$};
\node at (1,0) {$b$};
\node at (2,0) {$c$};
\draw [transform canvas={yshift=5pt}] (0,0) to[out=60,in=120] (1,0);
\end{tikzpicture}
\rightarrow
\begin{tikzpicture}[baseline=-2pt,scale=0.6]
\node at (0,0) {$a$};
\node at (1,0) {$b$};
\node at (2,0) {$c$};
\draw [transform canvas={yshift=5pt}] (0,0) to[out=60,in=120] (2,0);
\end{tikzpicture}
\qquad\text{for }(i,j)\in\{(b,c),(a,c)\}$ when $A=\{a<b<c\}$.

\item[(3)]
$
\begin{tikzpicture}[baseline=-2pt,scale=0.6]
\node at (0,0) {$a$};
\node at (1,0) {$b$};
\node at (2,0) {$c$};
\draw [transform canvas={yshift=5pt}] (1,0) to[out=60,in=120] (2,0);
\end{tikzpicture}
\rightarrow
\begin{tikzpicture}[baseline=-2pt,scale=0.6]
\node at (0,0) {$a$};
\node at (1,0) {$b$};
\node at (2,0) {$c$};
\draw [transform canvas={yshift=5pt}] (0,0) to[out=60,in=120] (2,0);
\end{tikzpicture}
\qquad\text{for }(i,j)\in\{(a,b),(a,c)\}$ when $A=\{a<b<c\}$.

\item[(4)]
$
\begin{tikzpicture}[baseline=-2pt,scale=0.6]
\node at (0,0) {$a$};
\node at (1,0) {$b$};
\node at (2,0) {$c$};
\node at (3,0) {$d$};
\draw [transform canvas={yshift=5pt}] (0,0) to[out=60,in=120] (1,0);
\draw [transform canvas={yshift=5pt}] (2,0) to[out=60,in=120] (3,0);
\end{tikzpicture}
\rightarrow
\begin{tikzpicture}[baseline=-2pt,scale=0.6]
\node at (0,0) {$a$};
\node at (1,0) {$b$};
\node at (2,0) {$c$};
\node at (3,0) {$d$};
\draw [transform canvas={yshift=5pt}] (0,0) to[out=60,in=120] (2,0);
\draw [transform canvas={yshift=5pt}] (1,0) to[out=60,in=120] (3,0);
\end{tikzpicture}
\qquad\text{for }(i,j)=(b,c)$ when $A=\{a<b<c<d\}$.

\item[(5)] $
\begin{tikzpicture}[baseline=-2pt,scale=0.6]
\node at (0,0) {$a$};
\node at (1,0) {$b$};
\node at (2,0) {$c$};
\node at (3,0) {$d$};
\draw [transform canvas={yshift=5pt}] (0,0) to[out=60,in=120] (1,0);
\draw [transform canvas={yshift=5pt}] (2,0) to[out=60,in=120] (3,0);
\end{tikzpicture}
\rightarrow
\begin{tikzpicture}[baseline=-2pt,scale=0.6]
\node at (0,0) {$a$};
\node at (1,0) {$b$};
\node at (2,0) {$c$};
\node at (3,0) {$d$};
\draw [transform canvas={yshift=5pt}] (0,0) to[out=60,in=120] (3,0);
\end{tikzpicture}
\qquad\text{for }(i,j)\in\{(a,c),(b,d),(a,d)\}$ when $A=\{a<b<c<d\}$.

\item[(6)]
$
\begin{tikzpicture}[baseline=-2pt,scale=0.6]
\node at (0,0) {$a$};
\node at (1,0) {$b$};
\node at (2,0) {$c$};
\node at (3,0) {$d$};
\draw [transform canvas={yshift=5pt}] (0,0) to[out=60,in=120] (2,0);
\draw [transform canvas={yshift=5pt}] (1,0) to[out=60,in=120] (3,0);
\end{tikzpicture}
\rightarrow
\begin{tikzpicture}[baseline=-2pt,scale=0.6]
\node at (0,0) {$a$};
\node at (1,0) {$b$};
\node at (2,0) {$c$};
\node at (3,0) {$d$};
\draw [transform canvas={yshift=5pt}] (0,0) to[out=60,in=120] (3,0);
\draw [transform canvas={yshift=5pt}] (1,0) to[out=60,in=120] (2,0);
\end{tikzpicture}
\qquad\text{for }(i,j)\in\{(a,b),(c,d),(a,d)\}$ when $A=\{a<b<c<d\}$.
\end{itemize}
Outside these cases, we 
set $\tau_{ij}(z)=z$.
\end{definition}

For example,
if $y=(1,3)(2,4)(5,7)$ and $z=(1,7)(2,4)$
 then $\tau_{15}(y)=\tau_{37}(y)=\tau_{17}(y)=z$ and $\tau_{23}(y)=y$.
The $\tau_{ij}$ operators are compatible with restriction, essentially by definition.

\begin{lemma}
\label{L: tau and pattern}
Let $z\in I_T$, let $i,j\in T$ with $i<j$,
and suppose $E\subseteq T$ satisfies $z(E)=E$
and $i,j\in E$.
Then $\tau_{ij}(z)(E)=E$ and
$\tau_{ij}(z)\cap E
=
\tau_{ij}(z\cap E)$.
\end{lemma}

\begin{proof}
Since $z(E)=E$ and $i,j\in E$,
we also have $z(i),z(j)\in E$.
The definition of $\tau_{ij}$ only modifies
the restriction of $z$ to
$\{i,j,z(i),z(j)\}$, so restricting to $E$ before or after
applying $\tau_{ij}$ gives the same result.
\end{proof}

The following lemmas summarize \cite[Cor.~3.5 and Thms.~3.16 and 3.20]{HMP3}, which relate $\tau_{ij}$
to the covering relations for the Bruhat order restricted to involutions.

\begin{lemma}[{\cite{HMP3}}]
\label{L: involution cover atoms}
Suppose $y\in I_T$ and
$z=\tau_{ij}(y)$ for some $i<j$ with
$\ellinv(z)=\ellinv(y)+1$.
Then for each $w\in\cW(z)$,
there exists $u\in\cW(y)$ and $a,b\in \{i,j,y(i),y(j)\}$ with $a<b$ and 
$u\lessdot t_{ab}u=w$ and $z=\tau_{ab}(y)$.
\end{lemma}

\begin{lemma}[{\cite{HMP3}}]
\label{L: atom cover transition---missing}
Let  $y,z\in I_T$ and $u \in \cW(y)$. Suppose some $i<j$ in $T$ exist with $\ell(t_{ij}u)=\ell(u)+1$ and $z=\tau_{ij}(y)\neq y$.
Then $t_{ij}u \in \cW(z)$ and $\ellinv(z)=\ellinv(y)+1$.
\end{lemma}

For $y\in I_T$ and $r\in T$, define
\begin{align*}
\hat\Phi^{+}(y,r)
 & =
\left\{
v\in I_T:
\ellinv(v)=\ellinv(y)+1
\text{ and }
v =\tau_{rj}(y)
\text{ for some }j>r
\right\}, \text{ and } \\
\hat\Phi^{-}(y,r)
& =
\left\{
v\in I_T:
\ellinv(v)=\ellinv(y)+1
\text{ and }
v=\tau_{ir}(y)
\text{ for some }i<r
\right\}.
\end{align*}

The involution analogue of Monk's rule
from~\cite{HMP3} states the following.
Suppose $z\in \Invol$ and $p\leq z(p)=q$.
Recall the definition of $x_{(p,q)}$ from \eqref{i-oplus-j-eq}.
We may embed $\Invol$ into $I_{N}$
for $N > n$. 
Then, when $z$ is viewed as an element 
of $I_N$ for large enough $N$,
 \cite[Thm.~1.5]{HMP3} asserts that
\begin{equation}
\label{EQ: HMP}
x_{(p,q)}\iS_z
=
\sum_{v\in\hat\Phi^{+}(z,q)}
\iS_v
-
\sum_{v\in\hat\Phi^{-}(z,p)}
\iS_v.
\end{equation}

The identity we will use is a special case 
of~\eqref{EQ: HMP} identified in~\cite{HMPstanley}.
We first introduce some more terminology.
We equip $\mathbb Z_{>0}\times\mathbb Z_{>0}$
with the lexicographic order where we first compare
the first coordinate and then the second.
For a non-identity permutation $w\in S_T$,
we call the maximal element of $D(w)$
the \defn{maximal Rothe element} of $w$.

A pair $(i,j) \in T \times T$ is an
\defn{inversion} of $w \in S_T$ if $i<j$ and $w(i)>w(j)$.
Let $\Inv(w)$ denote the set of inversions of $w$.
Notice that $\operatorname{Inv}(w)=\varnothing$
if and only if $w$ is the identity.

\begin{lemma}
\label{L: maximal Rothe and inv}
A pair $(q,p) \in T\times T$ is the maximal Rothe element of a non-identity permutation $w\in S_T$
if and only if $(q,w^{-1}(p))$
is the maximal inversion of $w$
in lexicographic order.
\end{lemma}

\begin{proof}
Suppose $(q,p)$ is the maximal Rothe element of $w$
and $(i,j)$ is the maximal inversion of $w$.
The map $(i',j')\mapsto(i',w(j'))$
is a bijection from $\operatorname{Inv}(w)$ to $D(w)$,
so $q=i$.

Consequently, there are no inversions
whose first coordinate is larger than $q$.
Thus, $w$ preserves the relative order on
$J=\{j' \in T: j'>q\text{ and }w(j')<w(q)\}$.
By definition, $j$ is the largest element of $J$,
while $p$ is the largest of
$\{w(j'):j'\in J\}$.
Therefore $p=w(j)$.
\end{proof}

Following~\cite{HMPstanley},
we say that an inversion $(i,j)$ of $z \in I_T$\
is \defn{visible} if $z(j)\leq i$.

\begin{corollary}
\label{C: maximal Rothe and inv}
A pair $(q,p) \in T\times T$ is the maximal Rothe element of a non-identity involution $z\in I_T$
if and only if $(q,z(p))$
is the maximal visible inversion of $z$
in lexicographic order,
and in this case $p \leq q$.
\end{corollary}
\begin{proof}
If $(i,j)$ is an inversion of $z$ that is not visible,
then $(z(j),z(i))$ is a visible inversion of $z$
that is lexicographically larger than $(i,j)$.
Thus, the maximal inversion of $z$ is visible.
Now apply Lemma~\ref{L: maximal Rothe and inv}.
\end{proof}

Fix a non-identity  $ z\in I_T$
with maximal Rothe element $(q,p)$.
Set $p'=z(p)$ and $q'=z(q)$.
By Corollary~\ref{C: maximal Rothe and inv},
we have
$p\leq q$ and 
$(q,p')$ is the maximal visible inversion of $z$.

\begin{lemma}
\label{L: on w}
We have $z(i)<z(j)$
if $i, j \in T$ and $q < i < j$.
Moreover, $q<p'$ and $p<q'$.
\end{lemma}

\begin{proof}
If $z(i) > z(j)$ and $q<i<j$,
then $(i,z(j))\in D(z)$ would be larger than $(q,p)$,
a contradiction.
This proves the first statement.
The second follows immediately from the fact that
$(q,p)\in D(z)$.
\end{proof}

\begin{lemma}
\label{L: ABC pq}
The four indices $p,p',q,q'$ have the following properties:
\begin{itemize}
\item[(A)] if $p=q$ then $p=q<q'=p'$;
\item[(B)] if $p\ne q=q'$ then $p<q'=q<p'$;
\item[(C)] if $p\neq q\neq q'$ then $p<q'<q<p'$.
\end{itemize}
\end{lemma}

\begin{proof}
The first two statements follow immediately
from Lemma~\ref{L: on w}. 
In the last case, we have $p<q$, $p<q'$, and $q<p'$.
It remains to show that $q'<q$.
This holds as we cannot have $q'>q$, as then $(q,q)\in D(z)$,
which is larger than $(q,p)$,
a contradiction.
\end{proof}

Following~\cite[Section~4]{HMPstanley},
we 
refer to the mutually exclusive situations in 
Lemma~\ref{L: ABC pq} as \CaseA, \CaseB, and \CaseC.

We now define the involution used in the transition formula.

\begin{definition}[{\cite[Prop.-Def.~4.26]{HMPstanley}}]
\label{D: eta(w)}
Define $\eta(z)\in I_T$ to be the involution that takes the same 
values as $z$ outside
$\{p,p',q,q'\}$, and which acts on this set as follows:
\begin{itemize}
\item[] \CaseA:
define $\eta(z)$ to fix both $p=q$ and $p'=q'$;

\item[] \CaseB:
define $\eta(z)$ to interchange $p$ and $q=q'$ while fixing $p'$;

\item[] \CaseC:
define $\eta(z)$ to interchange $p$ and $q$
and to interchange $q'$ and $p'$.
\end{itemize}
\end{definition}

Notice that in every case we have $\eta(z)(q)=p$.
The map $\eta$ is also compatible with restriction.

\begin{lemma}
\label{L: transition and pattern}
Suppose $E\subseteq T$ satisfies $z(E)=E$
and $p,q\in E$.
Then $z\cap E$ has maximal Rothe element $(q,p)$ and 
is in the same Case (A/B/C) as $z$ in 
Lemma~\ref{L: ABC pq}, 
and it holds that
$\eta(z)\cap E=\eta(z\cap E)$.
\end{lemma}

\begin{proof}
Since $p,q\in E$ and $z(E)=E$,
we also have $p',q'\in E$.
By~\eqref{EQ: Rothe pattern},
$(q,p)$ remains the maximal Rothe element of $z\cap E$.
The remaining assertions follow directly from
Lemma~\ref{L: ABC pq} 
 and Definition~\ref{D: eta(w)}.
\end{proof}

A key observation in~\cite{HMPstanley}
is the following.
\begin{lemma}[{\cite[Lem.~4.29]{HMPstanley}}]
\label{L: Phi+ only w}
Suppose $1\neq z\in \Invol$ has maximal Rothe element $(q,p)$.
Then 
\[\hat\Phi^{+}(\eta(z),q)=\{z\},
\quand \ellinv(\eta(z))
=\ellinv(z)-1.\]
\end{lemma}

Combining Lemma~\ref{L: Phi+ only w}
with~\eqref{EQ: HMP} gives
\begin{equation}
\label{EQ: iS Tran}
\textstyle
\iS_z
=
x_{(p,q)}\iS_{\eta(z)}
+
\sum_{v\in\hat\Phi^{-}(\eta(z),p)}
\iS_v.
\end{equation}
We will prove in the next section that this recurrence, combined with the base case $\iS_1=1$, 
uniquely determines all involution Schubert polynomials; see Remark~\ref{rec-rem}.

\subsection{Some technical results}

This final preliminary section contains some new (but somewhat technical) results that will be needed in Section~\ref{S: Proof}.
Continue to let $T$ be a finite set of positive integers.

\begin{lemma}
\label{L: negative transition at qp}
Let $y\in I_T$ and suppose $p\leq q=y(p)$.
If $v\in\hat\Phi^{-}(y,p)$ and $a=v(q)$, then
$a<p$, $y(a)<q$, $v(p)<q$, and $v=\tau_{ap}(y)$.
In particular, $v(r)=y(r)$ for every $r>q$.
\end{lemma}
\begin{proof}
Write $v=\tau_{ip}(y)$ for some $i<p$.
Checking the cases in Definition~\ref{D: tau} gives the first four assertions.
The last assertion follows since
 $y$ and $v$ differ only on
$\{a,p,y(a),q\}$, which are all 
at most $q$.
\end{proof}

The following result shows that the recursion for $\iS_z$ defined by \eqref{EQ: iS Tran}
terminates.

\begin{lemma}
\label{L: inv transition order}
Suppose $ z \in \Invol$ has maximal Rothe element $(q,p)$.
Then every non-identity involution on the right-hand side
of~\eqref{EQ: iS Tran}
has maximal Rothe element strictly smaller than $(q,p)$.
\end{lemma}

\begin{proof}
Set $p'=z(p)$ and $q'=z(q)$ as above. Then let $y=\eta(z)$ and take
$v\in\hat\Phi^{-}(y,p)$.
To prove that the maximal Rothe element of $y$
is smaller than $(q,p)$, it suffices to show that
\begin{equation}
\label{EQ: y increasing after}
y(q+1)<\cdots<y(n)
\qquad\text{and}\qquad
y(q)\leq p.
\end{equation}
The first condition guarantees that
any $(i,j)\in D(y)$ has $i\leq q$,
while the second guarantees that $j<p$
if $i=q$.
Thus, every element of $D(y)$ is strictly smaller
than $(q,p)$.
The same argument applies to $v$.

By Definition~\ref{D: eta(w)}, we have $y(q)=p$.
Recall that
$z(q+1)<\cdots<z(n)$
since $(q,p)$ is the maximal Rothe element of $z$.
For any $q<r\leq n$, we have $y(r)=z(r)$
except when $r=p'$.
Moreover, $z(p')=p$, while
$y(p')=p'$ in Cases (A) and (B)
 and $y(p')=q'$ in \CaseC.
Thus, we only need to check
$y(p')<z(p'+1)$ when $p'<n$.
Set $t=z(p'+1)$.
Since $z(q+1)<\cdots<z(n)$,
we have $p=z(p')<z(p'+1)=t$.
Assume toward contradiction that $t<y(p')$.
Then $t<p'<p'+1$, while $z(t)=p'+1$,
so $(t,t)\in D(z)$.
Since $(q,p)$ is the maximal element of $D(z)$
and $p<t$, we must have $t<q$.

In \CaseA, we have $q=p$,
so $p<t<q$ is already a contradiction.
Now suppose we are in Case (B) or (C).
We claim that $(q,t)\in D(z)$.
To prove this, first note that $z(t)=p'+1>q$.
In \CaseB, $z(q)=q>t$,
while in \CaseC,
$z(q)=q'=y(p')>t$.
Thus, $(q,t)\in D(z)$.
By the maximality of $(q,p)$, we must have $t\leq p$,
contradicting $p<t$.
This proves~\eqref{EQ: y increasing after}.

It remains to prove the analogous statement for $v$.
Set $a=v(q)$.
By Lemma~\ref{L: negative transition at qp},
we have $a<p$, $y(a)<q$, and
$v=\tau_{a,p}(y)$.
Thus, the four indices
$a$, $y(a)$, $p$, $y(p)=q$
are all at most $q$.
Since $\tau_{ap}$ only modifies these indices,
we have $v(r)=y(r)$ for every $r>q$,
so $v(q+1)<\cdots<v(n)$.
Together with $v(q)=a<p$,
this is the analogue of~\eqref{EQ: y increasing after}
for $v$.
\end{proof}

\begin{remark}\label{rec-rem}
Since lexicographic order on $[n]\times[n]$ is well ordered
and the identity is the unique involution with empty Rothe diagram,
combining Lemma~\ref{L: inv transition order} with the base case
$\iS_1=1$ shows that~\eqref{EQ: iS Tran}
recursively determines every involution Schubert polynomial.
\end{remark}

We record a consequence of Lemma~\ref{L: involution cover atoms}.

\begin{lemma}
\label{L: local atom identification}
Take $y\in I_T$ and $i<j$ in $T$, and suppose
$v=\tau_{ij}(y)$ satisfies $\ellinv(v)=\ellinv(y)+1$.
Then for all $u\in\cW(v)$ and all $c<d$ in the set $A=\{i,j,y(i),y(j)\}$, we have
\[
t_{cd}u\in\cW(y)
\quad\Longleftrightarrow\quad
\std{t_{cd}u}{A}\in\cW(\std{y}{A}),
\]
and when these equivalent conditions hold,  $t_{cd}u\lessdot u$.
\end{lemma}

\begin{proof}
The forward implication follows from Lemma~\ref{L: W reduction}.
For the converse, we first note that $\cW(\std{y}{A})$ has exactly one element.
After the order-preserving identification of $A$ with $[|A|]$,
Definition~\ref{D: tau} and Theorem~\ref{cW-lem} give the following possibilities:
\[
\begin{array}{c|ccccc}
\std{y}{A} & 12 & 213 & 132 & 2143 & 3412 \\
\hline
\text{unique inverse atom} & 12 & 213 & 132 & 2143 & 3142
\end{array}
\]
By Lemma~\ref{L: involution cover atoms}, there exist $h\in\cW(y)$
and $e<f$ in $A$ such that $h\lessdot t_{ef}h=u$.
Lemma~\ref{L: W reduction} gives $\std{h}{A}\in\cW(\std{y}{A})$.
Thus, if $\std{t_{cd}u}{A}\in\cW(\std{y}{A})$, the uniqueness just proved
implies $\std{h}{A}=\std{t_{cd}u}{A}$.
Both $h=t_{ef}u$ and $t_{cd}u$ are obtained from $u$ by permuting only
values in $A$, so this equality implies $h=t_{cd}u$.
Hence $t_{cd}u\in\cW(y)$ and $t_{cd}u\lessdot u$.
\end{proof}

As above, 
let $ z\in I_T\setminus\{1\}$
have maximal Rothe element $(q,p)$
and set $p'=z(p)$ and $q'=z(q)$.
\begin{lemma}
\label{L: AB eta(w)}
Suppose   $z$ is in Case (A) or (B) and $u\in\cW(z)$.
Assume that $p$ and $p'$ are adjacent in $u$.
In \CaseB, assume in addition that $q$ appears before $p'$ in $u$.
Then $t_{pp'}u\in\cW(\eta(z))$.
\end{lemma}

\begin{proof}
Set $y=\eta(z)$ and $A=\{p,p',q,q'\}$.
Definitions~\ref{D: tau} and~\ref{D: eta(w)} give $z=\tau_{pp'}(y)$,
and Lemma~\ref{L: Phi+ only w} gives $\ellinv(z)=\ellinv(y)+1$.
By Lemma~\ref{L: local atom identification}, it is enough to show that
$\std{t_{pp'}u}{A}\in\cW(\std{y}{A})$.

In \CaseA, $A=\{p,p'\}$ with $p<p'$.
By Lemma~\ref{L: W reduction}, $\std{u}{A}=p'\,p$.
Swapping these two entries gives the identity permutation on $A$,
which is the unique inverse atom of $\std{y}{A}$.

In \CaseB, $A=\{p,q,p'\}$ with $p<q<p'$.
After order-preserving relabeling, $\std{z}{A}=321$ and
$\std{y}{A}=213$.
Theorem~\ref{cW-lem} gives $\cW(321)=\{231,312\}$ and
$\cW(213)=\{213\}$.
Since $q$ appears before $p'$, Lemma~\ref{L: W reduction} forces
$\std{u}{A}=qp'p$.
Swapping $p$ and $p'$ gives $qpp'$, the unique inverse atom
of $\std{y}{A}$.
This proves the result. 
\end{proof}

\begin{lemma}
\label{L: C eta(w)}
Suppose $z$ is in \CaseC. Let $u\in\cW(z)$ and assume $q'$ is before $p$ in $u$.
Then $t_{pq'}u\in\cW(\eta(z))$.
\end{lemma}

\begin{proof}
Set $y=\eta(z)$ and $A=\{p,p',q,q'\}$.
Definitions~\ref{D: tau} and~\ref{D: eta(w)} give $z=\tau_{pq'}(y)$,
and Lemma~\ref{L: Phi+ only w} gives $\ellinv(z)=\ellinv(y)+1$.
By Lemma~\ref{L: ABC pq}, we have $p<q'<q<p'$.
After order-preserving relabeling, $\std{z}{A}=4321$ and
$\std{y}{A}=3412$.
Theorem~\ref{cW-lem} gives $\cW(4321)=\{3241,3412,4132\}$ and
$\cW(3412)=\{3142\}$.
Since $q'$ appears before $p$, Lemma~\ref{L: W reduction} forces
$\std{u}{A}=qq'p'p$.
Swapping $p$ and $q'$ gives $qpp'q'$, the unique inverse atom
of $\std{y}{A}$.
The result follows from Lemma~\ref{L: local atom identification}.
\end{proof}

We have one last technical statement.

\begin{lemma}
\label{L: reverse Fuse in Phi}
Let $y\in I_T$ and suppose $p\leq q=y(p)$.
Take $v\in\hat\Phi^{-}(y,p)$, set $a=v(q)$, and let $u\in\cW(v)$.
Suppose that $p$ appears before $a$ in $u$.
If $p=q$ and $a<y(a)<p$, assume in addition that $y(a)$ appears
before $a$ in $u$.
Then $t_{ap}u\in\cW(y)$ and $t_{ap}u\lessdot u$.
\end{lemma}

\begin{proof}
By Lemma~\ref{L: negative transition at qp}, we have $a<p$,
$y(a)<q$, and $v=\tau_{ap}(y)$.
Set $S=\{a,p,y(a),q\}$.
Lemma~\ref{L: W reduction} gives $\std{u}{S}\in\cW(\std{v}{S})$.
By Lemma~\ref{L: local atom identification}, it remains to check that
$\std{t_{ap}u}{S}\in\cW(\std{y}{S})$.

Apply the order-preserving identification of $S$ with $[|S|]$.
Definition~\ref{D: tau} gives the six possibilities below.
Theorem~\ref{cW-lem}, together with the order conditions in the statement,
forces the word in the fourth column.
All entries in the following table, including $a,p,q$, are written after this
identification:
\[
\begin{array}{c|c|c|c|c}
\std{y}{S} & \std{v}{S} & (a,p,q) & \std{u}{S}
& \std{t_{ap}u}{S} \\
\hline
12   & 21   & (1,2,2) & 21   & 12 \\
132  & 321  & (1,2,3) & 231  & 132 \\
213  & 321  & (1,3,3) & 231  & 213 \\
2143 & 4231 & (1,3,4) & 2341 & 2143 \\
2143 & 3412 & (2,3,4) & 3142 & 2143 \\
3412 & 4321 & (1,2,4) & 3241 & 3142
\end{array}
\]
For the third row, both $231$ and $312$ have $p$ before $a$.
The additional condition puts $y(a)=2$ before $a=1$ and excludes $312$.
In each of the other rows, the condition that $p$ appears before $a$
already determines the fourth column.
The last column is the unique inverse atom of the first column,
as computed in the proof of Lemma~\ref{L: local atom identification}.
This proves the required membership, and the result follows.
\end{proof}

\section{Bumpless pipe dreams}\label{bpd-sect}

In this section we formally introduce a notion of triangular \defn{bumpless pipe dreams}
for involutions in $S_n$. These are the objects of  primary interest in this paper.
After presenting the relevant definitions, we also  state our  main results (Theorems~\ref{main-thm} and \ref{fkGO-conj}), whose proofs appear in later sections.

\subsection{Triangular pipe dreams}

 Let $\LT_n = \{ (i,j) \in [n] \times [n] : i \geq j\}$ be the set of positions weakly below the diagonal in an $n\times n$ matrix.
A \defn{triangular tiling} is a map $D : \LT_n \to  \cToff\sqcup \cTdiag$, written $(i,j)\mapsto D_{ij}$,
where $\cT$ is some allowable set of tiles.
For our applications, the set of tiles will be partitioned $\cT = \cToff\sqcup \cTdiag$
into a subset of \defn{off-diagonal tiles} and a subset of \defn{diagonal tiles}, and it will be required 
that 
\[D_{ij} \in \cToff\quad\text{when $i>j$}
\qquand
 D_{ij} \in \cTdiag\quad\text{when $i=j$.}
 \]
 The off-diagonal tiles will be unit squares and the diagonal tiles will be unit right triangles.
This makes it possible to draw $D$ by assembling the relevant tiles in
the lower half triangle of an $n\times n$ square.  

Let us be more concrete. 
We will only consider  triangular tilings
with tile sets
\be\label{TPD-eq}
\cToff = \Bigl\{
\BTILE,\
\HTILE,\
\VTILE,\
\JTILE, \
\CTILE,\ 
\PTILE,\
\ETILE
\Bigr\}
\quand
\cTdiag = \Bigl\{
\BDIAG,\
\JDIAG, \
\CDIAG,\ 
\PDIAG,\
\EDIAG
\Bigr\},
\ee
and for the rest of this article we fix $\cToff$ and $\cTdiag$ to have these definitions.
We view each tile as representing a union of up to two wire segments,
which we refer to as \defn{pipes}.
The endpoints of each pipe are called \defn{connections}.
Both $\PTILE$ and $\PDIAG$ contain two pipes, shown as the vertical and horizontal lines,
and these pipes have connections on all sides of each tile.

A tiling $ \LT_n\to \cToff \sqcup \cTdiag$ is \defn{wirable}
if whenever two tiles occupy adjacent positions in the same row or column, the tiles either both have   or both do not have a connection on their shared side.
The connections on the diagonal sides of diagonal tiles are irrelevant to this property. 
Here are three examples when $n=8$:
\be\label{t-wirable-ex}
\begin{tikzpicture}[x=\bpdwidth,y=\bpdwidth,line cap=round,line join=round,baseline=(z.base)]
\node at (0,4.0) (z) {};
\draw[gray!70, line width=0.7pt](0,0)--(8,0);
\draw[gray!70, line width=0.7pt](0,0)--(0,8);
\draw[gray!70, line width=0.7pt](0,1)--(8,1);
\draw[gray!70, line width=0.7pt](1,0)--(1,8);
\draw[gray!70, line width=0.7pt](0,2)--(8,2);
\draw[gray!70, line width=0.7pt](2,0)--(2,8);
\draw[gray!70, line width=0.7pt](0,3)--(8,3);
\draw[gray!70, line width=0.7pt](3,0)--(3,8);
\draw[gray!70, line width=0.7pt](0,4)--(8,4);
\draw[gray!70, line width=0.7pt](4,0)--(4,8);
\draw[gray!70, line width=0.7pt](0,5)--(8,5);
\draw[gray!70, line width=0.7pt](5,0)--(5,8);
\draw[gray!70, line width=0.7pt](0,6)--(8,6);
\draw[gray!70, line width=0.7pt](6,0)--(6,8);
\draw[gray!70, line width=0.7pt](0,7)--(8,7);
\draw[gray!70, line width=0.7pt](7,0)--(7,8);
\draw[gray!70, line width=0.7pt](0,8)--(8,8);
\draw[gray!70, line width=0.7pt](8,0)--(8,8);
\draw[thick,rounded corners,color=blue] (0,0.5)--(1,0.5);
\draw[thick,rounded corners,color=blue] (1,0.5)--(2,0.5);
\draw[thick,rounded corners,color=blue] (1.5,0)--(1.5,1);
\draw[thick,rounded corners,color=blue] (0,2.5)--(1,2.5);
\draw[thick,rounded corners,color=blue] (1.5,1)--(1.5,1.5)--(2,1.5);
\draw[thick,rounded corners,color=blue] (2.5,1)--(2.5,0.5)--(2,0.5);
\draw[thick,rounded corners,color=blue] (2.5,0)--(2.5,0.5)--(3,0.5);
\draw[thick,rounded corners,color=blue] (0,3.5)--(1,3.5);
\draw[thick,rounded corners,color=blue] (1.5,3)--(1.5,2.5)--(1,2.5);
\draw[thick,rounded corners,color=blue] (2,1.5)--(3,1.5);
\draw[thick,rounded corners,color=blue] (2.5,1)--(2.5,2);
\draw[thick,rounded corners,color=blue] (3,0.5)--(4,0.5);
\draw[thick,rounded corners,color=blue] (0.5,5)--(0.5,4.5)--(0,4.5);
\draw[thick,rounded corners,color=blue] (1.5,4)--(1.5,3.5)--(1,3.5);
\draw[thick,rounded corners,color=blue] (1.5,3)--(1.5,3.5)--(2,3.5);
\draw[thick,rounded corners,color=blue] (2.5,2)--(2.5,3);
\draw[thick,rounded corners,color=blue] (3.5,2)--(3.5,1.5)--(3,1.5);
\draw[thick,rounded corners,color=blue] (4.5,1)--(4.5,0.5)--(4,0.5);
\draw[thick,rounded corners,color=blue] (4.5,0)--(4.5,0.5)--(5,0.5);
\draw[thick,rounded corners,color=blue] (0.5,6)--(0.5,5.5)--(0,5.5);
\draw[thick,rounded corners,color=blue] (0.5,5)--(0.5,5.5)--(1,5.5);
\draw[thick,rounded corners,color=blue] (1.5,4)--(1.5,4.5)--(2,4.5);
\draw[thick,rounded corners,color=blue] (2.5,4)--(2.5,3.5)--(2,3.5);
\draw[thick,rounded corners,color=blue] (2.5,3)--(2.5,3.5)--(3,3.5);
\draw[thick,rounded corners,color=blue] (3.5,2)--(3.5,2.5)--(4,2.5);
\draw[thick,rounded corners,color=blue] (4.5,1)--(4.5,2);
\draw[thick,rounded corners,color=blue] (5,0.5)--(6,0.5);
\draw[thick,rounded corners,color=blue] (0.5,7)--(0.5,6.5)--(0,6.5);
\draw[thick,rounded corners,color=blue] (0.5,6)--(0.5,6.5)--(1,6.5);
\draw[thick,rounded corners,color=blue] (1.5,6)--(1.5,5.5)--(1,5.5);
\draw[thick,rounded corners,color=blue] (2,4.5)--(3,4.5);
\draw[thick,rounded corners,color=blue] (2.5,4)--(2.5,5);
\draw[thick,rounded corners,color=blue] (3,3.5)--(4,3.5);
\draw[thick,rounded corners,color=blue] (4.5,3)--(4.5,2.5)--(4,2.5);
\draw[thick,rounded corners,color=blue] (4.5,2)--(4.5,2.5)--(5,2.5);
\draw[thick,rounded corners,color=blue] (6,0.5)--(7,0.5);
\draw[thick,rounded corners,color=blue] (6.5,0)--(6.5,1);
\draw[thick,rounded corners,color=blue] (0.5,7)--(0.5,8);
\draw[thick,rounded corners,color=blue] (1,6.5)--(2,6.5);
\draw[thick,rounded corners,color=blue] (1.5,6)--(1.5,7);
\draw[thick,rounded corners,color=blue] (2.5,5)--(2.5,6);
\draw[thick,rounded corners,color=blue] (3,4.5)--(4,4.5);
\draw[thick,rounded corners,color=blue] (4,3.5)--(5,3.5);
\draw[thick,rounded corners,color=blue] (4.5,3)--(4.5,4);
\draw[thick,rounded corners,color=blue] (5,2.5)--(6,2.5);
\draw[thick,rounded corners,color=blue] (6.5,1)--(6.5,2);
\draw[thick,rounded corners,color=blue] (7,0.5)--(8,0.5);
\draw[white, line width=0.7pt,fill=white] (0,8)--(8,0)--(8,8)--(0,8);
\draw[black, line width=1.6pt] (0,8)--(8,0)--(0,0)--(0,8);
\end{tikzpicture}
\qquad
\begin{tikzpicture}[x=\bpdwidth,y=\bpdwidth,line cap=round,line join=round,baseline=(z.base)]
\node at (0,4.0) (z) {};
\draw[gray!70, line width=0.7pt](0,0)--(8,0);
\draw[gray!70, line width=0.7pt](0,0)--(0,8);
\draw[gray!70, line width=0.7pt](0,1)--(8,1);
\draw[gray!70, line width=0.7pt](1,0)--(1,8);
\draw[gray!70, line width=0.7pt](0,2)--(8,2);
\draw[gray!70, line width=0.7pt](2,0)--(2,8);
\draw[gray!70, line width=0.7pt](0,3)--(8,3);
\draw[gray!70, line width=0.7pt](3,0)--(3,8);
\draw[gray!70, line width=0.7pt](0,4)--(8,4);
\draw[gray!70, line width=0.7pt](4,0)--(4,8);
\draw[gray!70, line width=0.7pt](0,5)--(8,5);
\draw[gray!70, line width=0.7pt](5,0)--(5,8);
\draw[gray!70, line width=0.7pt](0,6)--(8,6);
\draw[gray!70, line width=0.7pt](6,0)--(6,8);
\draw[gray!70, line width=0.7pt](0,7)--(8,7);
\draw[gray!70, line width=0.7pt](7,0)--(7,8);
\draw[gray!70, line width=0.7pt](0,8)--(8,8);
\draw[gray!70, line width=0.7pt](8,0)--(8,8);
\draw[thick,rounded corners,color=blue] (0.5,0)--(0.5,1);
\draw[thick,rounded corners,color=blue] (0.5,1)--(0.5,1.5)--(1,1.5);
\draw[thick,rounded corners,color=blue] (1.5,0)--(1.5,0.5)--(2,0.5);
\draw[thick,rounded corners,color=blue] (1.5,2)--(1.5,1.5)--(1,1.5);
\draw[thick,rounded corners,color=blue] (2,0.5)--(3,0.5);
\draw[thick,rounded corners,color=blue] (2.5,0)--(2.5,1);
\draw[thick,rounded corners,color=blue] (1.5,2)--(1.5,2.5)--(2,2.5);
\draw[thick,rounded corners,color=blue] (2.5,1)--(2.5,2);
\draw[thick,rounded corners,color=blue] (3,0.5)--(4,0.5);
\draw[thick,rounded corners,color=blue] (3.5,0)--(3.5,1);
\draw[thick,rounded corners,color=blue] (2.5,3)--(2.5,2.5)--(2,2.5);
\draw[thick,rounded corners,color=blue] (2.5,2)--(2.5,2.5)--(3,2.5);
\draw[thick,rounded corners,color=blue] (3.5,1)--(3.5,2);
\draw[thick,rounded corners,color=blue] (4.5,1)--(4.5,0.5)--(4,0.5);
\draw[thick,rounded corners,color=blue] (4.5,0)--(4.5,0.5)--(5,0.5);
\draw[thick,rounded corners,color=blue] (2.5,3)--(2.5,4);
\draw[thick,rounded corners,color=blue] (3,2.5)--(4,2.5);
\draw[thick,rounded corners,color=blue] (3.5,2)--(3.5,3);
\draw[thick,rounded corners,color=blue] (4.5,1)--(4.5,1.5)--(5,1.5);
\draw[thick,rounded corners,color=blue] (5.5,1)--(5.5,0.5)--(5,0.5);
\draw[thick,rounded corners,color=blue] (5.5,0)--(5.5,0.5)--(6,0.5);
\draw[thick,rounded corners,color=blue] (2.5,4)--(2.5,5);
\draw[thick,rounded corners,color=blue] (3.5,3)--(3.5,3.5)--(4,3.5);
\draw[thick,rounded corners,color=blue] (4,2.5)--(5,2.5);
\draw[thick,rounded corners,color=blue] (5,1.5)--(6,1.5);
\draw[thick,rounded corners,color=blue] (5.5,1)--(5.5,2);
\draw[thick,rounded corners,color=blue] (6,0.5)--(7,0.5);
\draw[thick,rounded corners,color=blue] (6.5,0)--(6.5,1);
\draw[thick,rounded corners,color=blue] (2.5,5)--(2.5,6);
\draw[thick,rounded corners,color=blue] (4,3.5)--(5,3.5);
\draw[thick,rounded corners,color=blue] (5,2.5)--(6,2.5);
\draw[thick,rounded corners,color=blue] (5.5,2)--(5.5,3);
\draw[thick,rounded corners,color=blue] (6.5,2)--(6.5,1.5)--(6,1.5);
\draw[thick,rounded corners,color=blue] (6.5,1)--(6.5,1.5)--(7,1.5);
\draw[thick,rounded corners,color=blue] (7,0.5)--(8,0.5);
\draw[thick,rounded corners,color=blue] (7.5,0)--(7.5,1);
\draw[white, line width=0.7pt,fill=white] (0,8)--(8,0)--(8,8)--(0,8);
\draw[black, line width=1.6pt] (0,8)--(8,0)--(0,0)--(0,8);
\end{tikzpicture}
\qquad
\begin{tikzpicture}[x=\bpdwidth,y=\bpdwidth,line cap=round,line join=round,baseline=(z.base)]
\node at (0,4.0) (z) {};
\draw[gray!70, line width=0.7pt](0,0)--(8,0);
\draw[gray!70, line width=0.7pt](0,0)--(0,8);
\draw[gray!70, line width=0.7pt](0,1)--(8,1);
\draw[gray!70, line width=0.7pt](1,0)--(1,8);
\draw[gray!70, line width=0.7pt](0,2)--(8,2);
\draw[gray!70, line width=0.7pt](2,0)--(2,8);
\draw[gray!70, line width=0.7pt](0,3)--(8,3);
\draw[gray!70, line width=0.7pt](3,0)--(3,8);
\draw[gray!70, line width=0.7pt](0,4)--(8,4);
\draw[gray!70, line width=0.7pt](4,0)--(4,8);
\draw[gray!70, line width=0.7pt](0,5)--(8,5);
\draw[gray!70, line width=0.7pt](5,0)--(5,8);
\draw[gray!70, line width=0.7pt](0,6)--(8,6);
\draw[gray!70, line width=0.7pt](6,0)--(6,8);
\draw[gray!70, line width=0.7pt](0,7)--(8,7);
\draw[gray!70, line width=0.7pt](7,0)--(7,8);
\draw[gray!70, line width=0.7pt](0,8)--(8,8);
\draw[gray!70, line width=0.7pt](8,0)--(8,8);
\draw[thick,rounded corners,color=blue] (0.5,0)--(0.5,1);
\draw[thick,rounded corners,color=blue] (0.5,1)--(0.5,2);
\draw[thick,rounded corners,color=blue] (1.5,0)--(1.5,0.5)--(2,0.5);
\draw[thick,rounded corners,color=blue] (0.5,2)--(0.5,3);
\draw[thick,rounded corners,color=blue] (2,0.5)--(3,0.5);
\draw[thick,rounded corners,color=blue] (2.5,0)--(2.5,1);
\draw[thick,rounded corners,color=blue] (0.5,3)--(0.5,3.5)--(1,3.5);
\draw[thick,rounded corners,color=blue] (2.5,1)--(2.5,2);
\draw[thick,rounded corners,color=blue] (3,0.5)--(4,0.5);
\draw[thick,rounded corners,color=blue] (3.5,0)--(3.5,1);
\draw[thick,rounded corners,color=blue] (1,3.5)--(2,3.5);
\draw[thick,rounded corners,color=blue] (2.5,2)--(2.5,2.5)--(3,2.5);
\draw[thick,rounded corners,color=blue] (3.5,1)--(3.5,2);
\draw[thick,rounded corners,color=blue] (4,0.5)--(5,0.5);
\draw[thick,rounded corners,color=blue] (4.5,0)--(4.5,1);
\draw[thick,rounded corners,color=blue] (2.5,4)--(2.5,3.5)--(2,3.5);
\draw[thick,rounded corners,color=blue] (3,2.5)--(4,2.5);
\draw[thick,rounded corners,color=blue] (3.5,2)--(3.5,3);
\draw[thick,rounded corners,color=blue] (4.5,1)--(4.5,2);
\draw[thick,rounded corners,color=blue] (5,0.5)--(6,0.5);
\draw[thick,rounded corners,color=blue] (5.5,0)--(5.5,1);
\draw[thick,rounded corners,color=blue] (2.5,4)--(2.5,4.5)--(3,4.5);
\draw[thick,rounded corners,color=blue] (3.5,3)--(3.5,3.5)--(4,3.5);
\draw[thick,rounded corners,color=blue] (4,2.5)--(5,2.5);
\draw[thick,rounded corners,color=blue] (4.5,2)--(4.5,3);
\draw[thick,rounded corners,color=blue] (5.5,1)--(5.5,1.5)--(6,1.5);
\draw[thick,rounded corners,color=blue] (6,0.5)--(7,0.5);
\draw[thick,rounded corners,color=blue] (6.5,0)--(6.5,1);
\draw[thick,rounded corners,color=blue] (3,4.5)--(4,4.5);
\draw[thick,rounded corners,color=blue] (4,3.5)--(5,3.5);
\draw[thick,rounded corners,color=blue] (4.5,3)--(4.5,4);
\draw[thick,rounded corners,color=blue] (5,2.5)--(6,2.5);
\draw[thick,rounded corners,color=blue] (6,1.5)--(7,1.5);
\draw[thick,rounded corners,color=blue] (6.5,1)--(6.5,2);
\draw[thick,rounded corners,color=blue] (7,0.5)--(8,0.5);
\draw[thick,rounded corners,color=blue] (7.5,0)--(7.5,1);
\draw[white, line width=0.7pt,fill=white] (0,8)--(8,0)--(8,8)--(0,8);
\draw[black, line width=1.6pt] (0,8)--(8,0)--(0,0)--(0,8);
\end{tikzpicture}
.
\ee
The positions in such pictures are oriented as in a matrix,
so we have $D_{53}=\ETILE$ if $D$ is the first tiling,
 $D_{53}=\VTILE$ if $D$ is the second tiling,
 and
 $D_{53}=\JTILE$ if $D$ is the third tiling.
 We view a wirable tiling as a collection of pipes, each 
traveling from a connection on the south or west side to
a connection on the diagonal side of the lower half triangle of an $n\times n$ square.

\begin{definition}\label{tpd-def}
A \defn{triangular pipe dream} is a wirable triangular tiling $D : \LT_n \to \cToff\sqcup \cTdiag$ 
such that when $D$ is drawn in 
the lower half of an $n\times n$ square, 
there are no connections on its west side and
  $n$ connections on its south side.
  A triangular pipe dream is \defn{bumpless} if it does not include the \defn{bump tiles} $\ETILE$ or $\EDIAG$.
 \end{definition}

The last two examples in \eqref{t-wirable-ex} 
are both triangular pipe dreams, but only the third example 
is bumpless.
From now on we will only consider triangular pipe dreams
$D$ rather than wirable triangular tilings in general.

\subsection{Involution pipe dreams}

To each triangular pipe dream $D$ we associate
a permutation 
$\pi(D) \in S_n$
in the following way.
Label the $n$ connections on the south side of $D$ by $1,2,3,\dots,n$
from west to east. This assigns a label to each of the pipes in $D$.
We form $\pi(D)$ in one-line notation
by iterating over the connections on the diagonal side of $D$
from northwest to southeast and listing the labels of the incoming pipes.
When two pipes meet at the same connection in a $\PDIAG$ tile,
we list the label of the vertical pipe before the label of the horizontal pipe.
For example, adding these labels to the triangular pipe dreams 
\be\label{t-wirable-ex2}
D=
\begin{tikzpicture}[x=\bpdwidth,y=\bpdwidth,line cap=round,line join=round,baseline=(z.base)]
\node at (0,4.0) (z) {};
\draw[gray!70, line width=0.7pt](0,0)--(8,0);
\draw[gray!70, line width=0.7pt](0,0)--(0,8);
\draw[gray!70, line width=0.7pt](0,1)--(8,1);
\draw[gray!70, line width=0.7pt](1,0)--(1,8);
\draw[gray!70, line width=0.7pt](0,2)--(8,2);
\draw[gray!70, line width=0.7pt](2,0)--(2,8);
\draw[gray!70, line width=0.7pt](0,3)--(8,3);
\draw[gray!70, line width=0.7pt](3,0)--(3,8);
\draw[gray!70, line width=0.7pt](0,4)--(8,4);
\draw[gray!70, line width=0.7pt](4,0)--(4,8);
\draw[gray!70, line width=0.7pt](0,5)--(8,5);
\draw[gray!70, line width=0.7pt](5,0)--(5,8);
\draw[gray!70, line width=0.7pt](0,6)--(8,6);
\draw[gray!70, line width=0.7pt](6,0)--(6,8);
\draw[gray!70, line width=0.7pt](0,7)--(8,7);
\draw[gray!70, line width=0.7pt](7,0)--(7,8);
\draw[gray!70, line width=0.7pt](0,8)--(8,8);
\draw[gray!70, line width=0.7pt](8,0)--(8,8);
\draw[thick,rounded corners,color=blue] (0.5,1)--(0.5,1.5)--(1,1.5);
\draw[thick,rounded corners,color=blue] (1.5,0)--(1.5,0.5)--(2,0.5);
\draw[thick,rounded corners,color=blue] (1.5,2)--(1.5,1.5)--(1,1.5);
\draw[thick,rounded corners,color=blue] (1.5,2)--(1.5,2.5)--(2,2.5);
\draw[thick,rounded corners,color=blue] (2.5,3)--(2.5,2.5)--(2,2.5);
\draw[thick,rounded corners,color=blue] (2.5,2)--(2.5,2.5)--(3,2.5);
\draw[thick,rounded corners,color=blue] (4.5,1)--(4.5,0.5)--(4,0.5);
\draw[thick,rounded corners,color=blue] (4.5,0)--(4.5,0.5)--(5,0.5);
\draw[thick,rounded corners,color=blue] (4.5,1)--(4.5,1.5)--(5,1.5);
\draw[thick,rounded corners,color=blue] (5.5,1)--(5.5,0.5)--(5,0.5);
\draw[thick,rounded corners,color=blue] (5.5,0)--(5.5,0.5)--(6,0.5);
\draw[thick,rounded corners,color=blue] (3.5,3)--(3.5,3.5)--(4,3.5);
\draw[thick,rounded corners,color=blue] (2.5,5)--(2.5,6);
\draw[thick,rounded corners,color=blue] (4,3.5)--(5,3.5);
\draw[thick,rounded corners,color=blue] (6.5,2)--(6.5,1.5)--(6,1.5);
\draw[thick,rounded corners,color=blue] (6.5,1)--(6.5,1.5)--(7,1.5);
\draw[thick,rounded corners,color=blue] (0.5,0)--(0.5,1);
\draw[thick,rounded corners,color=blue] (2.5,3)--(2.5,4);
\draw[thick,rounded corners,color=blue] (2.5,4)--(2.5,5);
\draw[thick,rounded corners,color=blue] (2,0.5)--(3,0.5);
\draw[thick,rounded corners,color=blue] (2.5,0)--(2.5,1);
\draw[thick,rounded corners,color=blue] (3,0.5)--(4,0.5);
\draw[thick,rounded corners,color=blue] (3.5,0)--(3.5,1);
\draw[thick,rounded corners,color=blue] (5,1.5)--(6,1.5);
\draw[thick,rounded corners,color=blue] (5.5,1)--(5.5,2);
\draw[thick,rounded corners,color=blue] (2.5,1)--(2.5,2);
\draw[thick,rounded corners,color=blue] (3,2.5)--(4,2.5);
\draw[thick,rounded corners,color=blue] (3.5,2)--(3.5,3);
\draw[thick,rounded corners,color=blue] (4,2.5)--(5,2.5);
\draw[thick,rounded corners,color=blue] (5,2.5)--(6,2.5);
\draw[thick,rounded corners,color=blue] (5.5,2)--(5.5,3);
\draw[thick,rounded corners,color=blue] (3.5,1)--(3.5,2);
\draw[thick,rounded corners,color=blue] (6,0.5)--(7,0.5);
\draw[thick,rounded corners,color=blue] (6.5,0)--(6.5,1);
\draw[thick,rounded corners,color=blue] (7,0.5)--(8,0.5);
\draw[thick,rounded corners,color=blue] (7.5,0)--(7.5,1);
\draw[white, line width=0.7pt,fill=white] (0,8)--(8,0)--(8,8)--(0,8);
\draw[black, line width=1.6pt] (0,8)--(8,0)--(0,0)--(0,8);
\draw[](0.5,-0.4) node {1};
\draw[](1.5,-0.4) node {2};
\draw[](2.5,-0.4) node {3};
\draw[](3.5,-0.4) node {4};
\draw[](4.5,-0.4) node {5};
\draw[](5.5,-0.4) node {6};
\draw[](6.5,-0.4) node {7};
\draw[](7.5,-0.4) node {8};
\draw[](2.8,5.75) node {1};
\draw[](4.8,3.75) node {4};
\draw[](5.55,3.0) node {5};
\draw[](6.05,2.5) node {3};
\draw[](6.55,2.0) node {2};
\draw[](7.05,1.5) node {7};
\draw[](7.55,1.0) node {8};
\draw[](8.05,0.5) node {6};
\end{tikzpicture}
\qquand
E=\begin{tikzpicture}[x=\bpdwidth,y=\bpdwidth,line cap=round,line join=round,baseline=(z.base)]
\node at (0,4.0) (z) {};
\draw[gray!70, line width=0.7pt](0,0)--(8,0);
\draw[gray!70, line width=0.7pt](0,0)--(0,8);
\draw[gray!70, line width=0.7pt](0,1)--(8,1);
\draw[gray!70, line width=0.7pt](1,0)--(1,8);
\draw[gray!70, line width=0.7pt](0,2)--(8,2);
\draw[gray!70, line width=0.7pt](2,0)--(2,8);
\draw[gray!70, line width=0.7pt](0,3)--(8,3);
\draw[gray!70, line width=0.7pt](3,0)--(3,8);
\draw[gray!70, line width=0.7pt](0,4)--(8,4);
\draw[gray!70, line width=0.7pt](4,0)--(4,8);
\draw[gray!70, line width=0.7pt](0,5)--(8,5);
\draw[gray!70, line width=0.7pt](5,0)--(5,8);
\draw[gray!70, line width=0.7pt](0,6)--(8,6);
\draw[gray!70, line width=0.7pt](6,0)--(6,8);
\draw[gray!70, line width=0.7pt](0,7)--(8,7);
\draw[gray!70, line width=0.7pt](7,0)--(7,8);
\draw[gray!70, line width=0.7pt](0,8)--(8,8);
\draw[gray!70, line width=0.7pt](8,0)--(8,8);
\draw[thick,rounded corners,color=blue] (1.5,0)--(1.5,0.5)--(2,0.5);
\draw[thick,rounded corners,color=blue] (0.5,3)--(0.5,3.5)--(1,3.5);
\draw[thick,rounded corners,color=blue] (2.5,2)--(2.5,2.5)--(3,2.5);
\draw[thick,rounded corners,color=blue] (2.5,4)--(2.5,3.5)--(2,3.5);
\draw[thick,rounded corners,color=blue] (2.5,4)--(2.5,4.5)--(3,4.5);
\draw[thick,rounded corners,color=blue] (3.5,3)--(3.5,3.5)--(4,3.5);
\draw[thick,rounded corners,color=blue] (5.5,1)--(5.5,1.5)--(6,1.5);
\draw[thick,rounded corners,color=blue] (3,4.5)--(4,4.5);
\draw[thick,rounded corners,color=blue] (5,2.5)--(6,2.5);
\draw[thick,rounded corners,color=blue] (0.5,0)--(0.5,1);
\draw[thick,rounded corners,color=blue] (0.5,1)--(0.5,2);
\draw[thick,rounded corners,color=blue] (0.5,2)--(0.5,3);
\draw[thick,rounded corners,color=blue] (1,3.5)--(2,3.5);
\draw[thick,rounded corners,color=blue] (2,0.5)--(3,0.5);
\draw[thick,rounded corners,color=blue] (2.5,0)--(2.5,1);
\draw[thick,rounded corners,color=blue] (3,0.5)--(4,0.5);
\draw[thick,rounded corners,color=blue] (3.5,0)--(3.5,1);
\draw[thick,rounded corners,color=blue] (4,0.5)--(5,0.5);
\draw[thick,rounded corners,color=blue] (4.5,0)--(4.5,1);
\draw[thick,rounded corners,color=blue] (5,0.5)--(6,0.5);
\draw[thick,rounded corners,color=blue] (5.5,0)--(5.5,1);
\draw[thick,rounded corners,color=blue] (6,0.5)--(7,0.5);
\draw[thick,rounded corners,color=blue] (6.5,0)--(6.5,1);
\draw[thick,rounded corners,color=blue] (7,0.5)--(8,0.5);
\draw[thick,rounded corners,color=blue] (7.5,0)--(7.5,1);
\draw[thick,rounded corners,color=blue] (2.5,1)--(2.5,2);
\draw[thick,rounded corners,color=blue] (3,2.5)--(4,2.5);
\draw[thick,rounded corners,color=blue] (3.5,2)--(3.5,3);
\draw[thick,rounded corners,color=blue] (4,2.5)--(5,2.5);
\draw[thick,rounded corners,color=blue] (4.5,2)--(4.5,3);
\draw[thick,rounded corners,color=blue] (3.5,1)--(3.5,2);
\draw[thick,rounded corners,color=blue] (4,3.5)--(5,3.5);
\draw[thick,rounded corners,color=blue] (4.5,3)--(4.5,4);
\draw[thick,rounded corners,color=blue] (4.5,1)--(4.5,2);
\draw[thick,rounded corners,color=blue] (6,1.5)--(7,1.5);
\draw[thick,rounded corners,color=blue] (6.5,1)--(6.5,2);
\draw[white, line width=0.7pt,fill=white] (0,8)--(8,0)--(8,8)--(0,8);
\draw[black, line width=1.6pt] (0,8)--(8,0)--(0,0)--(0,8);
\draw[](0.5,-0.4) node {1};
\draw[](1.5,-0.4) node {2};
\draw[](2.5,-0.4) node {3};
\draw[](3.5,-0.4) node {4};
\draw[](4.5,-0.4) node {5};
\draw[](5.5,-0.4) node {6};
\draw[](6.5,-0.4) node {7};
\draw[](7.5,-0.4) node {8};
\draw[](3.8,4.75) node {1};
\draw[](4.55,4.0) node {5};
\draw[](5.05,3.5) node {4};
\draw[](5.8,2.75) node {3};
\draw[](6.55,2.0) node {7};
\draw[](7.05,1.5) node {6};
\draw[](7.55,1.0) node {8};
\draw[](8.05,0.5) node {2};
\end{tikzpicture}
\ee
from \eqref{t-wirable-ex}, we see that $\pi(D) =14532786 \in S_8$
and $\pi(E) = 15437682$.

The sets of triangular pipe dreams that are of primary interest (along with our main result connecting these objects to involution Schubert polynomials) can be described using only the preceding constructions.

\begin{definition}
A triangular pipe dream is \defn{reduced} if no two of its pipes 
 cross more than once, where a crossing occurs at any $\PTILE$ or $\PDIAG$ tile.
Given $z \in \Invol $,  define 
$ \InvPD(z)$ and $ \InvBPD(z)$
to be the respective sets of   triangular pipe dreams and  bumpless pipe dreams $D$ that are reduced with $\pi(D) \in \cW(z)$.
\end{definition}

We refer to the elements of $ \InvPD(z)$ and $ \InvBPD(z)$ as \defn{(generic) involution pipe dreams}
and \defn{bumpless involution pipe dreams} for $z$.
Both examples in \eqref{t-wirable-ex2} are reduced.
Given any triangular pipe dream $D$, let 
\be
\label{upturn-eq}
\blank(D)=\left\{(i,j) \in \LT_n : D_{ij} =\BTILE\text{ or }\BDIAG\right\}
\quand
\upturn(D)=\left\{(i,j) \in \LT_n : D_{ij} =\JTILE\text{ or }\JDIAG\right\}.
\ee Then define the \defn{weight} of $D$ to be the polynomial
$\weight(D) \in \ZZ[x_1,x_2,\dots,x_n]$ given by
\be 
\textstyle \weight(D) = \prod_{(i,j) \in \blank(D)} 
 x_{(i,j)}
\qquad
\text{where }x_{(i,j)} = \begin{cases} x_i &\text{if $i=j$} 
\\ x_i + x_j&\text{if $i\neq j$.}
\end{cases}
\ee
 The following is our first main theorem. This will be proved in Section~\ref{S: Proof}.

\begin{theorem}\label{main-thm}
For each $z \in \Invol $ it holds that
$
\iS_z = \sum_{D \in \InvBPD(z)}\weight(D).
$
\end{theorem}

This result is the ``bumpless'' analogue of an earlier pipe dream formula for $\iS_z$ proved in \cite{HMP6}.
The ordinary Schubert polynomials $\fkS_w$ also have two distinct ``classical'' and ``bumpless'' pipe dream formulas \cite{BJS,LLS}, which are related by the bijection of Gao and Huang \cite{GaoHuang}.
It is an open problem to find an involution analogue of this bijection
that can relate \cite[Thm.~1.5]{HMP6} to Theorem~\ref{main-thm}.

\begin{example}\label{bpd-is-eq}
For the involution $z=21543=(1,2)(3,5) \in I_5$ the set $\InvBPD(z)$ has four elements
\[
\begin{tikzpicture}[x=\bpdwidth,y=\bpdwidth,line cap=round,line join=round,baseline=(z.base)]
\node at (0,2.5) (z) {};
\draw[gray!70, line width=0.7pt](0,0)--(5,0);
\draw[gray!70, line width=0.7pt](0,0)--(0,5);
\draw[gray!70, line width=0.7pt](0,1)--(5,1);
\draw[gray!70, line width=0.7pt](1,0)--(1,5);
\draw[gray!70, line width=0.7pt](0,2)--(5,2);
\draw[gray!70, line width=0.7pt](2,0)--(2,5);
\draw[gray!70, line width=0.7pt](0,3)--(5,3);
\draw[gray!70, line width=0.7pt](3,0)--(3,5);
\draw[gray!70, line width=0.7pt](0,4)--(5,4);
\draw[gray!70, line width=0.7pt](4,0)--(4,5);
\draw[gray!70, line width=0.7pt](0,5)--(5,5);
\draw[gray!70, line width=0.7pt](5,0)--(5,5);
\draw[thick,rounded corners,color=blue] (0.5,3)--(0.5,3.5)--(1,3.5);
\draw[thick,rounded corners,color=blue] (3.5,1)--(3.5,2);
\draw[thick,rounded corners,color=blue] (2.5,0)--(2.5,0.5)--(3,0.5);
\draw[thick,rounded corners,color=blue] (0.5,0)--(0.5,1);
\draw[thick,rounded corners,color=blue] (0.5,1)--(0.5,2);
\draw[thick,rounded corners,color=blue] (0.5,2)--(0.5,3);
\draw[thick,rounded corners,color=blue] (1,3.5)--(2,3.5);
\draw[thick,rounded corners,color=blue] (1.5,3)--(1.5,4);
\draw[thick,rounded corners,color=blue] (1.5,0)--(1.5,1);
\draw[thick,rounded corners,color=blue] (1.5,1)--(1.5,2);
\draw[thick,rounded corners,color=blue] (1.5,2)--(1.5,3);
\draw[thick,rounded corners,color=blue] (3,0.5)--(4,0.5);
\draw[thick,rounded corners,color=blue] (3.5,0)--(3.5,1);
\draw[thick,rounded corners,color=blue] (4,0.5)--(5,0.5);
\draw[thick,rounded corners,color=blue] (4.5,0)--(4.5,1);
\draw[white, line width=0.7pt,fill=white] (0,5)--(5,0)--(5,5)--(0,5);
\draw[black, line width=1.6pt] (0,5)--(5,0)--(0,0)--(0,5);
\draw[](0.5,-0.4) node {1};
\draw[](1.5,-0.4) node {2};
\draw[](2.5,-0.4) node {3};
\draw[](3.5,-0.4) node {4};
\draw[](4.5,-0.4) node {5};
\draw[](1.55,4.0) node {2};
\draw[](2.05,3.5) node {1};
\draw[](3.8,1.75) node {4};
\draw[](4.55,1.0) node {5};
\draw[](5.05,0.5) node {3};
\end{tikzpicture}
\qquad
\begin{tikzpicture}[x=\bpdwidth,y=\bpdwidth,line cap=round,line join=round,baseline=(z.base)]
\node at (0,2.5) (z) {};
\draw[gray!70, line width=0.7pt](0,0)--(5,0);
\draw[gray!70, line width=0.7pt](0,0)--(0,5);
\draw[gray!70, line width=0.7pt](0,1)--(5,1);
\draw[gray!70, line width=0.7pt](1,0)--(1,5);
\draw[gray!70, line width=0.7pt](0,2)--(5,2);
\draw[gray!70, line width=0.7pt](2,0)--(2,5);
\draw[gray!70, line width=0.7pt](0,3)--(5,3);
\draw[gray!70, line width=0.7pt](3,0)--(3,5);
\draw[gray!70, line width=0.7pt](0,4)--(5,4);
\draw[gray!70, line width=0.7pt](4,0)--(4,5);
\draw[gray!70, line width=0.7pt](0,5)--(5,5);
\draw[gray!70, line width=0.7pt](5,0)--(5,5);
\draw[thick,rounded corners,color=blue] (3.5,1)--(3.5,2);
\draw[thick,rounded corners,color=blue] (2.5,0)--(2.5,0.5)--(3,0.5);
\draw[thick,rounded corners,color=blue] (2,2.5)--(3,2.5);
\draw[thick,rounded corners,color=blue] (0.5,2)--(0.5,2.5)--(1,2.5);
\draw[thick,rounded corners,color=blue] (1.5,3)--(1.5,4);
\draw[thick,rounded corners,color=blue] (0.5,0)--(0.5,1);
\draw[thick,rounded corners,color=blue] (0.5,1)--(0.5,2);
\draw[thick,rounded corners,color=blue] (1,2.5)--(2,2.5);
\draw[thick,rounded corners,color=blue] (1.5,2)--(1.5,3);
\draw[thick,rounded corners,color=blue] (1.5,0)--(1.5,1);
\draw[thick,rounded corners,color=blue] (1.5,1)--(1.5,2);
\draw[thick,rounded corners,color=blue] (3,0.5)--(4,0.5);
\draw[thick,rounded corners,color=blue] (3.5,0)--(3.5,1);
\draw[thick,rounded corners,color=blue] (4,0.5)--(5,0.5);
\draw[thick,rounded corners,color=blue] (4.5,0)--(4.5,1);
\draw[white, line width=0.7pt,fill=white] (0,5)--(5,0)--(5,5)--(0,5);
\draw[black, line width=1.6pt] (0,5)--(5,0)--(0,0)--(0,5);
\draw[](0.5,-0.4) node {1};
\draw[](1.5,-0.4) node {2};
\draw[](2.5,-0.4) node {3};
\draw[](3.5,-0.4) node {4};
\draw[](4.5,-0.4) node {5};
\draw[](1.8,3.75) node {2};
\draw[](2.8,2.75) node {1};
\draw[](3.8,1.75) node {4};
\draw[](4.55,1.0) node {5};
\draw[](5.05,0.5) node {3};
\end{tikzpicture}
\qquad
\begin{tikzpicture}[x=\bpdwidth,y=\bpdwidth,line cap=round,line join=round,baseline=(z.base)]
\node at (0,2.5) (z) {};
\draw[gray!70, line width=0.7pt](0,0)--(5,0);
\draw[gray!70, line width=0.7pt](0,0)--(0,5);
\draw[gray!70, line width=0.7pt](0,1)--(5,1);
\draw[gray!70, line width=0.7pt](1,0)--(1,5);
\draw[gray!70, line width=0.7pt](0,2)--(5,2);
\draw[gray!70, line width=0.7pt](2,0)--(2,5);
\draw[gray!70, line width=0.7pt](0,3)--(5,3);
\draw[gray!70, line width=0.7pt](3,0)--(3,5);
\draw[gray!70, line width=0.7pt](0,4)--(5,4);
\draw[gray!70, line width=0.7pt](4,0)--(4,5);
\draw[gray!70, line width=0.7pt](0,5)--(5,5);
\draw[gray!70, line width=0.7pt](5,0)--(5,5);
\draw[thick,rounded corners,color=blue] (3.5,1)--(3.5,2);
\draw[thick,rounded corners,color=blue] (2.5,0)--(2.5,0.5)--(3,0.5);
\draw[thick,rounded corners,color=blue] (0.5,2)--(0.5,2.5)--(1,2.5);
\draw[thick,rounded corners,color=blue] (2.5,2)--(2.5,1.5)--(2,1.5);
\draw[thick,rounded corners,color=blue] (1.5,1)--(1.5,1.5)--(2,1.5);
\draw[thick,rounded corners,color=blue] (0.5,0)--(0.5,1);
\draw[thick,rounded corners,color=blue] (0.5,1)--(0.5,2);
\draw[thick,rounded corners,color=blue] (1,2.5)--(2,2.5);
\draw[thick,rounded corners,color=blue] (2,2.5)--(3,2.5);
\draw[thick,rounded corners,color=blue] (2.5,2)--(2.5,3);
\draw[thick,rounded corners,color=blue] (1.5,0)--(1.5,1);
\draw[thick,rounded corners,color=blue] (3,0.5)--(4,0.5);
\draw[thick,rounded corners,color=blue] (3.5,0)--(3.5,1);
\draw[thick,rounded corners,color=blue] (4,0.5)--(5,0.5);
\draw[thick,rounded corners,color=blue] (4.5,0)--(4.5,1);
\draw[white, line width=0.7pt,fill=white] (0,5)--(5,0)--(5,5)--(0,5);
\draw[black, line width=1.6pt] (0,5)--(5,0)--(0,0)--(0,5);
\draw[](0.5,-0.4) node {1};
\draw[](1.5,-0.4) node {2};
\draw[](2.5,-0.4) node {3};
\draw[](3.5,-0.4) node {4};
\draw[](4.5,-0.4) node {5};
\draw[](2.55,3.0) node {2};
\draw[](3.05,2.5) node {1};
\draw[](3.8,1.75) node {4};
\draw[](4.55,1.0) node {5};
\draw[](5.05,0.5) node {3};
\end{tikzpicture}
\qquad
\begin{tikzpicture}[x=\bpdwidth,y=\bpdwidth,line cap=round,line join=round,baseline=(z.base)]
\node at (0,2.5) (z) {};
\draw[gray!70, line width=0.7pt](0,0)--(5,0);
\draw[gray!70, line width=0.7pt](0,0)--(0,5);
\draw[gray!70, line width=0.7pt](0,1)--(5,1);
\draw[gray!70, line width=0.7pt](1,0)--(1,5);
\draw[gray!70, line width=0.7pt](0,2)--(5,2);
\draw[gray!70, line width=0.7pt](2,0)--(2,5);
\draw[gray!70, line width=0.7pt](0,3)--(5,3);
\draw[gray!70, line width=0.7pt](3,0)--(3,5);
\draw[gray!70, line width=0.7pt](0,4)--(5,4);
\draw[gray!70, line width=0.7pt](4,0)--(4,5);
\draw[gray!70, line width=0.7pt](0,5)--(5,5);
\draw[gray!70, line width=0.7pt](5,0)--(5,5);
\draw[thick,rounded corners,color=blue] (3.5,1)--(3.5,2);
\draw[thick,rounded corners,color=blue] (2.5,0)--(2.5,0.5)--(3,0.5);
\draw[thick,rounded corners,color=blue] (1.5,3)--(1.5,4);
\draw[thick,rounded corners,color=blue] (2.5,2)--(2.5,1.5)--(2,1.5);
\draw[thick,rounded corners,color=blue] (0.5,1)--(0.5,1.5)--(1,1.5);
\draw[thick,rounded corners,color=blue] (2.5,2)--(2.5,3);
\draw[thick,rounded corners,color=blue] (0.5,0)--(0.5,1);
\draw[thick,rounded corners,color=blue] (1,1.5)--(2,1.5);
\draw[thick,rounded corners,color=blue] (1.5,1)--(1.5,2);
\draw[thick,rounded corners,color=blue] (1.5,0)--(1.5,1);
\draw[thick,rounded corners,color=blue] (1.5,2)--(1.5,3);
\draw[thick,rounded corners,color=blue] (3,0.5)--(4,0.5);
\draw[thick,rounded corners,color=blue] (3.5,0)--(3.5,1);
\draw[thick,rounded corners,color=blue] (4,0.5)--(5,0.5);
\draw[thick,rounded corners,color=blue] (4.5,0)--(4.5,1);
\draw[white, line width=0.7pt,fill=white] (0,5)--(5,0)--(5,5)--(0,5);
\draw[black, line width=1.6pt] (0,5)--(5,0)--(0,0)--(0,5);
\draw[](0.5,-0.4) node {1};
\draw[](1.5,-0.4) node {2};
\draw[](2.5,-0.4) node {3};
\draw[](3.5,-0.4) node {4};
\draw[](4.5,-0.4) node {5};
\draw[](1.8,3.75) node {2};
\draw[](2.8,2.75) node {1};
\draw[](3.8,1.75) node {4};
\draw[](4.55,1.0) node {5};
\draw[](5.05,0.5) node {3};
\end{tikzpicture}
\]
whose weights are $x_1x_3(x_4+x_3)$,
$x_1(x_2+x_1)(x_4+x_3)$, $x_1x_2(x_2+x_1)$, and
$x_1(x_2+x_1)(x_3+x_1)$. These weights sum to
$\iS_{21543} = x_1 x_3 x_4 + x_1 x_3^2 + x_1 x_2 x_4 + 2x_1 x_2 x_3 + x_1 x_2^2 + x_1^2 x_4 + 2x_1^2 x_3 + 2x_1^2 x_2 + x_1^3.
$
\end{example}

Our second main theorem uses the same bumpless pipe dreams $D$ with the $K$-theoretic weight function
\be \label{kweight-eq}
 \kweight(D) = \prod_{(i,j) \in \blank(D)} 
  x_i \oplus x_j
 \prod_{(i,j) \in \upturn(D)} 
 (1+\beta x_{i}) (1+\beta x_{j})
\quad
\text{where }a\oplus b = a+b+ \beta ab
\ee

\begin{theorem}\label{fkGO-conj}
For each $z \in \Ivex $ it holds that
$
\fkGO_z = \sum_{D \in \InvBPD(z)} \kweight(D).
$
\end{theorem}

To prove this theorem, we show that it
has an equivalent formulation
in terms of \defn{set-valued shifted plane partitions} in  Section~\ref{svpp-sect}.
We then derive this second statement (Theorem~\ref{fkGO-conj2}) in Appendix~\ref{app-sect2}.

\begin{example}\label{fkGO-ex}
If $z = 526413 = (1,5)(3,6) \in I^{\mathsf{vex}}_6$ then $\InvBPD(z)$ has three elements
\[
\begin{tikzpicture}[x=\bpdwidth,y=\bpdwidth,line cap=round,line join=round,baseline=(z.base)]
\node at (0,2.75) (z) {};
\draw[gray!70, line width=0.7pt](0,0)--(6,0);
\draw[gray!70, line width=0.7pt](0,0)--(0,6);
\draw[gray!70, line width=0.7pt](0,1)--(6,1);
\draw[gray!70, line width=0.7pt](1,0)--(1,6);
\draw[gray!70, line width=0.7pt](0,2)--(6,2);
\draw[gray!70, line width=0.7pt](2,0)--(2,6);
\draw[gray!70, line width=0.7pt](0,3)--(6,3);
\draw[gray!70, line width=0.7pt](3,0)--(3,6);
\draw[gray!70, line width=0.7pt](0,4)--(6,4);
\draw[gray!70, line width=0.7pt](4,0)--(4,6);
\draw[gray!70, line width=0.7pt](0,5)--(6,5);
\draw[gray!70, line width=0.7pt](5,0)--(5,6);
\draw[gray!70, line width=0.7pt](0,6)--(6,6);
\draw[gray!70, line width=0.7pt](6,0)--(6,6);
\draw[thick,rounded corners,color=blue] (3.5,2)--(3.5,3);
\draw[thick,rounded corners,color=blue] (0.5,1)--(0.5,1.5)--(1,1.5);
\draw[thick,rounded corners,color=blue] (2.5,0)--(2.5,0.5)--(3,0.5);
\draw[thick,rounded corners,color=blue] (2.5,3)--(2.5,2.5)--(2,2.5);
\draw[thick,rounded corners,color=blue] (2.5,3)--(2.5,4);
\draw[thick,rounded corners,color=blue] (1.5,2)--(1.5,2.5)--(2,2.5);
\draw[thick,rounded corners,color=blue] (0.5,0)--(0.5,1);
\draw[thick,rounded corners,color=blue] (1,1.5)--(2,1.5);
\draw[thick,rounded corners,color=blue] (1.5,1)--(1.5,2);
\draw[thick,rounded corners,color=blue] (2,1.5)--(3,1.5);
\draw[thick,rounded corners,color=blue] (3,1.5)--(4,1.5);
\draw[thick,rounded corners,color=blue] (3.5,1)--(3.5,2);
\draw[thick,rounded corners,color=blue] (4,1.5)--(5,1.5);
\draw[thick,rounded corners,color=blue] (4.5,1)--(4.5,2);
\draw[thick,rounded corners,color=blue] (1.5,0)--(1.5,1);
\draw[thick,rounded corners,color=blue] (3,0.5)--(4,0.5);
\draw[thick,rounded corners,color=blue] (3.5,0)--(3.5,1);
\draw[thick,rounded corners,color=blue] (4,0.5)--(5,0.5);
\draw[thick,rounded corners,color=blue] (4.5,0)--(4.5,1);
\draw[thick,rounded corners,color=blue] (5,0.5)--(6,0.5);
\draw[thick,rounded corners,color=blue] (5.5,0)--(5.5,1);
\draw[white, line width=0.7pt,fill=white] (0,6)--(6,0)--(6,6)--(0,6);
\draw[black, line width=1.6pt] (0,6)--(6,0)--(0,0)--(0,6);
\end{tikzpicture}
\qquad
\begin{tikzpicture}[x=\bpdwidth,y=\bpdwidth,line cap=round,line join=round,baseline=(z.base)]
\node at (0,2.75) (z) {};
\draw[gray!70, line width=0.7pt](0,0)--(6,0);
\draw[gray!70, line width=0.7pt](0,0)--(0,6);
\draw[gray!70, line width=0.7pt](0,1)--(6,1);
\draw[gray!70, line width=0.7pt](1,0)--(1,6);
\draw[gray!70, line width=0.7pt](0,2)--(6,2);
\draw[gray!70, line width=0.7pt](2,0)--(2,6);
\draw[gray!70, line width=0.7pt](0,3)--(6,3);
\draw[gray!70, line width=0.7pt](3,0)--(3,6);
\draw[gray!70, line width=0.7pt](0,4)--(6,4);
\draw[gray!70, line width=0.7pt](4,0)--(4,6);
\draw[gray!70, line width=0.7pt](0,5)--(6,5);
\draw[gray!70, line width=0.7pt](5,0)--(5,6);
\draw[gray!70, line width=0.7pt](0,6)--(6,6);
\draw[gray!70, line width=0.7pt](6,0)--(6,6);
\draw[thick,rounded corners,color=blue] (3.5,2)--(3.5,3);
\draw[thick,rounded corners,color=blue] (0.5,1)--(0.5,1.5)--(1,1.5);
\draw[thick,rounded corners,color=blue] (2.5,0)--(2.5,0.5)--(3,0.5);
\draw[thick,rounded corners,color=blue] (2,3.5)--(3,3.5);
\draw[thick,rounded corners,color=blue] (1.5,3)--(1.5,3.5)--(2,3.5);
\draw[thick,rounded corners,color=blue] (0.5,0)--(0.5,1);
\draw[thick,rounded corners,color=blue] (1,1.5)--(2,1.5);
\draw[thick,rounded corners,color=blue] (1.5,1)--(1.5,2);
\draw[thick,rounded corners,color=blue] (2,1.5)--(3,1.5);
\draw[thick,rounded corners,color=blue] (3,1.5)--(4,1.5);
\draw[thick,rounded corners,color=blue] (3.5,1)--(3.5,2);
\draw[thick,rounded corners,color=blue] (4,1.5)--(5,1.5);
\draw[thick,rounded corners,color=blue] (4.5,1)--(4.5,2);
\draw[thick,rounded corners,color=blue] (1.5,0)--(1.5,1);
\draw[thick,rounded corners,color=blue] (1.5,2)--(1.5,3);
\draw[thick,rounded corners,color=blue] (3,0.5)--(4,0.5);
\draw[thick,rounded corners,color=blue] (3.5,0)--(3.5,1);
\draw[thick,rounded corners,color=blue] (4,0.5)--(5,0.5);
\draw[thick,rounded corners,color=blue] (4.5,0)--(4.5,1);
\draw[thick,rounded corners,color=blue] (5,0.5)--(6,0.5);
\draw[thick,rounded corners,color=blue] (5.5,0)--(5.5,1);
\draw[white, line width=0.7pt,fill=white] (0,6)--(6,0)--(6,6)--(0,6);
\draw[black, line width=1.6pt] (0,6)--(6,0)--(0,0)--(0,6);
\end{tikzpicture}
\qquad
\begin{tikzpicture}[x=\bpdwidth,y=\bpdwidth,line cap=round,line join=round,baseline=(z.base)]
\node at (0,2.75) (z) {};
\draw[gray!70, line width=0.7pt](0,0)--(6,0);
\draw[gray!70, line width=0.7pt](0,0)--(0,6);
\draw[gray!70, line width=0.7pt](0,1)--(6,1);
\draw[gray!70, line width=0.7pt](1,0)--(1,6);
\draw[gray!70, line width=0.7pt](0,2)--(6,2);
\draw[gray!70, line width=0.7pt](2,0)--(2,6);
\draw[gray!70, line width=0.7pt](0,3)--(6,3);
\draw[gray!70, line width=0.7pt](3,0)--(3,6);
\draw[gray!70, line width=0.7pt](0,4)--(6,4);
\draw[gray!70, line width=0.7pt](4,0)--(4,6);
\draw[gray!70, line width=0.7pt](0,5)--(6,5);
\draw[gray!70, line width=0.7pt](5,0)--(5,6);
\draw[gray!70, line width=0.7pt](0,6)--(6,6);
\draw[gray!70, line width=0.7pt](6,0)--(6,6);
\draw[thick,rounded corners,color=blue] (1.5,4)--(1.5,5);
\draw[thick,rounded corners,color=blue] (3.5,2)--(3.5,3);
\draw[thick,rounded corners,color=blue] (0.5,1)--(0.5,1.5)--(1,1.5);
\draw[thick,rounded corners,color=blue] (2.5,0)--(2.5,0.5)--(3,0.5);
\draw[thick,rounded corners,color=blue] (0.5,0)--(0.5,1);
\draw[thick,rounded corners,color=blue] (1,1.5)--(2,1.5);
\draw[thick,rounded corners,color=blue] (1.5,1)--(1.5,2);
\draw[thick,rounded corners,color=blue] (2,1.5)--(3,1.5);
\draw[thick,rounded corners,color=blue] (3,1.5)--(4,1.5);
\draw[thick,rounded corners,color=blue] (3.5,1)--(3.5,2);
\draw[thick,rounded corners,color=blue] (4,1.5)--(5,1.5);
\draw[thick,rounded corners,color=blue] (4.5,1)--(4.5,2);
\draw[thick,rounded corners,color=blue] (1.5,0)--(1.5,1);
\draw[thick,rounded corners,color=blue] (1.5,2)--(1.5,3);
\draw[thick,rounded corners,color=blue] (1.5,3)--(1.5,4);
\draw[thick,rounded corners,color=blue] (3,0.5)--(4,0.5);
\draw[thick,rounded corners,color=blue] (3.5,0)--(3.5,1);
\draw[thick,rounded corners,color=blue] (4,0.5)--(5,0.5);
\draw[thick,rounded corners,color=blue] (4.5,0)--(4.5,1);
\draw[thick,rounded corners,color=blue] (5,0.5)--(6,0.5);
\draw[thick,rounded corners,color=blue] (5.5,0)--(5.5,1);
\draw[white, line width=0.7pt,fill=white] (0,6)--(6,0)--(6,6)--(0,6);
\draw[black, line width=1.6pt] (0,6)--(6,0)--(0,0)--(0,6);
\end{tikzpicture}
\]
and one can check that the right side of Theorem~\ref{fkGO-conj} gives the same result as 
\[\textstyle
\bpartial_2 \(\prod_{(i,j)\in\iD(y)} x_i \oplus x_j\) = \bpartial_2 \fkGO_{y} = \fkGO_z
\]
for the dominant involution $y=s_2 z s_2 = 563412 = (1,5)(2,6) \in I^{\mathsf{vex}}_6$ which has shape $(4,4,2,2)$.
\end{example}

\subsection{K-theoretic pipe dreams}\label{kbpd-sect}

We now discuss a second way of associating a permutation  to any triangular pipe dream, which will be relevant to a second (conjectural) $K$-theoretic analogue of Theorem~\ref{main-thm}.
Continue to suppose $D$ is a triangular pipe dream, which is not necessarily
bumpless or reduced.
We form a permutation $\kpi(D) \in S_n$ from $D$ in 
the same way as $\pi(D)$,
except we change how the pipes are labeled whenever they cross
at a $\PTILE$ tile or meet at a $\PDIAG$ tile according to the rules
\be\label{internal-eq}
\begin{tikzpicture}[x=\bpdwidth,y=\bpdwidth,baseline=(z.base)]
\node at (0.5,1.25) (z) {};
\node at (-0.3,1.5) {$a$};
\node at (1.3,1.5) {$a$};
\node at (0.5,0.65) {$b$};
\node at (0.5,2.4) {$b$};
\draw[thick,rounded corners,color=blue] (0,1.5)--(1,1.5);
\draw[thick,rounded corners,color=blue] (.5,1)--(.5,2);
\end{tikzpicture}
\qquad
\begin{tikzpicture}[x=\bpdwidth,y=\bpdwidth,baseline=(z.base)]
\node at (0.5,1.25) (z) {};
\node at (-0.3,1.5) {$b$};
\node at (1.3,1.5) {$a$};
\node at (0.5,0.65) {$a$};
\node at (0.5,2.4) {$b$};
\draw[thick,rounded corners,color=blue] (0,1.5)--(1,1.5);
\draw[thick,rounded corners,color=blue] (.5,1)--(.5,2);
\end{tikzpicture}
\qquad
\begin{tikzpicture}[x=\bpdwidth,y=\bpdwidth,baseline=(z.base)]
\node at (0.5,1.25) (z) {};
\node at (-0.3,1.5) {$a$};
\node at (1,1.5) {$a$};
\node at (0.5,0.65) {$b$};
\node at (0.5,2) {$b$};
\draw[thick,rounded corners,color=blue] (0,1.5)--(0.5,1.5);
\draw[thick,rounded corners,color=blue] (.5,1)--(.5,1.5);
\end{tikzpicture}
\qquad
\begin{tikzpicture}[x=\bpdwidth,y=\bpdwidth,baseline=(z.base)]
\node at (0.5,1.25) (z) {};
\node at (-0.3,1.5) {$b$};
\node at (1,1.5) {$a$};
\node at (0.5,0.65) {$a$};
\node at (0.5,2) {$b$};
\draw[thick,rounded corners,color=blue] (0,1.5)--(0.5,1.5);
\draw[thick,rounded corners,color=blue] (.5,1)--(.5,1.5);
\end{tikzpicture}
\qquad\text{when $a<b$.}\ee
Then
 $\kpi(D)$ is given in one-line notation
by reading the labels of the diagonal connections of $D$ from northwest to southeast as before. The last two rules in \eqref{internal-eq}
mean that when two pipes meet at the same connection in a $\PDIAG$ tile,
we list the larger label first (which may not necessarily correspond to the vertical pipe).

We have $\pi(D) = \kpi(D)$
if  $D$ is reduced.
For non-trivial examples of $\kpi(D)$, consider 
\be
F=
\begin{tikzpicture}[x=\bpdwidth,y=\bpdwidth,line cap=round,line join=round,baseline=(z.base)]
\node at (0,4.0) (z) {};
\draw[gray!70, line width=0.7pt](0,0)--(8,0);
\draw[gray!70, line width=0.7pt](0,0)--(0,8);
\draw[gray!70, line width=0.7pt](0,1)--(8,1);
\draw[gray!70, line width=0.7pt](1,0)--(1,8);
\draw[gray!70, line width=0.7pt](0,2)--(8,2);
\draw[gray!70, line width=0.7pt](2,0)--(2,8);
\draw[gray!70, line width=0.7pt](0,3)--(8,3);
\draw[gray!70, line width=0.7pt](3,0)--(3,8);
\draw[gray!70, line width=0.7pt](0,4)--(8,4);
\draw[gray!70, line width=0.7pt](4,0)--(4,8);
\draw[gray!70, line width=0.7pt](0,5)--(8,5);
\draw[gray!70, line width=0.7pt](5,0)--(5,8);
\draw[gray!70, line width=0.7pt](0,6)--(8,6);
\draw[gray!70, line width=0.7pt](6,0)--(6,8);
\draw[gray!70, line width=0.7pt](0,7)--(8,7);
\draw[gray!70, line width=0.7pt](7,0)--(7,8);
\draw[gray!70, line width=0.7pt](0,8)--(8,8);
\draw[gray!70, line width=0.7pt](8,0)--(8,8);
\draw[thick,rounded corners,color=blue] (0.5,0)--(0.5,0.5)--(1,0.5);
\draw[thick,rounded corners,color=blue] (2.5,1)--(2.5,0.5)--(2,0.5);
\draw[thick,rounded corners,color=blue] (2.5,0)--(2.5,0.5)--(3,0.5);
\draw[thick,rounded corners,color=blue] (2.5,2)--(2.5,2.5)--(3,2.5);
\draw[thick,rounded corners,color=blue] (3.5,1)--(3.5,1.5)--(4,1.5);
\draw[thick,rounded corners,color=blue] (5.5,1)--(5.5,0.5)--(5,0.5);
\draw[thick,rounded corners,color=blue] (5.5,0)--(5.5,0.5)--(6,0.5);
\draw[thick,rounded corners,color=blue] (1.5,5)--(1.5,5.5)--(2,5.5);
\draw[thick,rounded corners,color=blue] (4.5,3)--(4.5,2.5)--(4,2.5);
\draw[thick,rounded corners,color=blue] (4.5,2)--(4.5,2.5)--(5,2.5);
\draw[thick,rounded corners,color=blue] (6.5,1)--(6.5,0.5)--(6,0.5);
\draw[thick,rounded corners,color=blue] (6.5,0)--(6.5,0.5)--(7,0.5);
\draw[thick,rounded corners,color=blue] (2,5.5)--(3,5.5);
\draw[thick,rounded corners,color=blue] (4.5,3)--(4.5,4);
\draw[thick,rounded corners,color=blue] (5.5,3)--(5.5,2.5)--(5,2.5);
\draw[thick,rounded corners,color=blue] (5.5,2)--(5.5,2.5)--(6,2.5);
\draw[thick,rounded corners,color=blue] (7.5,1)--(7.5,0.5)--(7,0.5);
\draw[thick,rounded corners,color=blue] (7.5,0)--(7.5,0.5)--(8,0.5);
\draw[thick,rounded corners,color=blue] (1,0.5)--(2,0.5);
\draw[thick,rounded corners,color=blue] (1.5,0)--(1.5,1);
\draw[thick,rounded corners,color=blue] (2.5,1)--(2.5,2);
\draw[thick,rounded corners,color=blue] (3,2.5)--(4,2.5);
\draw[thick,rounded corners,color=blue] (1.5,1)--(1.5,2);
\draw[thick,rounded corners,color=blue] (1.5,2)--(1.5,3);
\draw[thick,rounded corners,color=blue] (1.5,3)--(1.5,4);
\draw[thick,rounded corners,color=blue] (1.5,4)--(1.5,5);
\draw[thick,rounded corners,color=blue] (3,0.5)--(4,0.5);
\draw[thick,rounded corners,color=blue] (3.5,0)--(3.5,1);
\draw[thick,rounded corners,color=blue] (4,0.5)--(5,0.5);
\draw[thick,rounded corners,color=blue] (4.5,0)--(4.5,1);
\draw[thick,rounded corners,color=blue] (5,1.5)--(6,1.5);
\draw[thick,rounded corners,color=blue] (5.5,1)--(5.5,2);
\draw[thick,rounded corners,color=blue] (4,1.5)--(5,1.5);
\draw[thick,rounded corners,color=blue] (4.5,1)--(4.5,2);
\draw[thick,rounded corners,color=blue] (6,1.5)--(7,1.5);
\draw[thick,rounded corners,color=blue] (6.5,1)--(6.5,2);
\draw[white, line width=0.7pt,fill=white] (0,8)--(8,0)--(8,8)--(0,8);
\draw[black, line width=1.6pt] (0,8)--(8,0)--(0,0)--(0,8);
\draw[](0.5,-0.4) node {1};
\draw[](1.5,-0.4) node {2};
\draw[](2.5,-0.4) node {3};
\draw[](3.5,-0.4) node {4};
\draw[](4.5,-0.4) node {5};
\draw[](5.5,-0.4) node {6};
\draw[](6.5,-0.4) node {7};
\draw[](7.5,-0.4) node {8};
\draw[](2.8,5.75) node {2};
\draw[](4.8,3.75) node {1};
\draw[](5.55,3.0) node {5};
\draw[](6.05,2.5) node {4};
\draw[](6.55,2.0) node {6};
\draw[](7.05,1.5) node {3};
\draw[](7.55,1.0) node {7};
\draw[](8.05,0.5) node {8};
\end{tikzpicture}
\qquand
G = \begin{tikzpicture}[x=\bpdwidth,y=\bpdwidth,line cap=round,line join=round,baseline=(z.base)]
\node at (0,4.0) (z) {};
\draw[gray!70, line width=0.7pt](0,0)--(8,0);
\draw[gray!70, line width=0.7pt](0,0)--(0,8);
\draw[gray!70, line width=0.7pt](0,1)--(8,1);
\draw[gray!70, line width=0.7pt](1,0)--(1,8);
\draw[gray!70, line width=0.7pt](0,2)--(8,2);
\draw[gray!70, line width=0.7pt](2,0)--(2,8);
\draw[gray!70, line width=0.7pt](0,3)--(8,3);
\draw[gray!70, line width=0.7pt](3,0)--(3,8);
\draw[gray!70, line width=0.7pt](0,4)--(8,4);
\draw[gray!70, line width=0.7pt](4,0)--(4,8);
\draw[gray!70, line width=0.7pt](0,5)--(8,5);
\draw[gray!70, line width=0.7pt](5,0)--(5,8);
\draw[gray!70, line width=0.7pt](0,6)--(8,6);
\draw[gray!70, line width=0.7pt](6,0)--(6,8);
\draw[gray!70, line width=0.7pt](0,7)--(8,7);
\draw[gray!70, line width=0.7pt](7,0)--(7,8);
\draw[gray!70, line width=0.7pt](0,8)--(8,8);
\draw[gray!70, line width=0.7pt](8,0)--(8,8);
\draw[thick,rounded corners,color=blue] (0.5,1)--(0.5,1.5)--(1,1.5);
\draw[thick,rounded corners,color=blue] (2.5,0)--(2.5,0.5)--(3,0.5);
\draw[thick,rounded corners,color=blue] (1.5,2)--(1.5,2.5)--(2,2.5);
\draw[thick,rounded corners,color=blue] (2.5,2)--(2.5,1.5)--(2,1.5);
\draw[thick,rounded corners,color=blue] (4.5,1)--(4.5,1.5)--(5,1.5);
\draw[thick,rounded corners,color=blue] (4.5,3)--(4.5,2.5)--(4,2.5);
\draw[thick,rounded corners,color=blue] (2.5,5)--(2.5,6);
\draw[thick,rounded corners,color=blue] (3.5,4)--(3.5,5);
\draw[thick,rounded corners,color=blue] (4.5,3)--(4.5,4);
\draw[thick,rounded corners,color=blue] (5.5,2)--(5.5,3);
\draw[thick,rounded corners,color=blue] (0.5,0)--(0.5,1);
\draw[thick,rounded corners,color=blue] (1,1.5)--(2,1.5);
\draw[thick,rounded corners,color=blue] (1.5,1)--(1.5,2);
\draw[thick,rounded corners,color=blue] (2,2.5)--(3,2.5);
\draw[thick,rounded corners,color=blue] (2.5,2)--(2.5,3);
\draw[thick,rounded corners,color=blue] (2.5,3)--(2.5,4);
\draw[thick,rounded corners,color=blue] (2.5,4)--(2.5,5);
\draw[thick,rounded corners,color=blue] (1.5,0)--(1.5,1);
\draw[thick,rounded corners,color=blue] (3,2.5)--(4,2.5);
\draw[thick,rounded corners,color=blue] (3.5,2)--(3.5,3);
\draw[thick,rounded corners,color=blue] (3,0.5)--(4,0.5);
\draw[thick,rounded corners,color=blue] (3.5,0)--(3.5,1);
\draw[thick,rounded corners,color=blue] (4,0.5)--(5,0.5);
\draw[thick,rounded corners,color=blue] (4.5,0)--(4.5,1);
\draw[thick,rounded corners,color=blue] (5,0.5)--(6,0.5);
\draw[thick,rounded corners,color=blue] (5.5,0)--(5.5,1);
\draw[thick,rounded corners,color=blue] (6,0.5)--(7,0.5);
\draw[thick,rounded corners,color=blue] (6.5,0)--(6.5,1);
\draw[thick,rounded corners,color=blue] (7,0.5)--(8,0.5);
\draw[thick,rounded corners,color=blue] (7.5,0)--(7.5,1);
\draw[thick,rounded corners,color=blue] (3.5,1)--(3.5,2);
\draw[thick,rounded corners,color=blue] (3.5,3)--(3.5,4);
\draw[thick,rounded corners,color=blue] (5,1.5)--(6,1.5);
\draw[thick,rounded corners,color=blue] (5.5,1)--(5.5,2);
\draw[thick,rounded corners,color=blue] (6,1.5)--(7,1.5);
\draw[thick,rounded corners,color=blue] (6.5,1)--(6.5,2);
\draw[white, line width=0.7pt,fill=white] (0,8)--(8,0)--(8,8)--(0,8);
\draw[black, line width=1.6pt] (0,8)--(8,0)--(0,0)--(0,8);
\draw[](0.5,-0.4) node {1};
\draw[](1.5,-0.4) node {2};
\draw[](2.5,-0.4) node {3};
\draw[](3.5,-0.4) node {4};
\draw[](4.5,-0.4) node {5};
\draw[](5.5,-0.4) node {6};
\draw[](6.5,-0.4) node {7};
\draw[](7.5,-0.4) node {8};
\draw[](2.8,5.75) node {2};
\draw[](3.8,4.75) node {4};
\draw[](4.8,3.75) node {1};
\draw[](5.8,2.75) node {6};
\draw[](6.55,2.0) node {7};
\draw[](7.05,1.5) node {5};
\draw[](7.55,1.0) node {8};
\draw[](8.05,0.5) node {3};
\end{tikzpicture}
.
\ee
The diagonal labels in the above pictures are 
given according to the convention defining $\kpi$,
and so we have 
\[\kpi(F) = 21546378 \neq 21536478 = \pi(F)
\qquand \kpi(G) = 24167583 \neq 14267583 = \pi(G).\]

\begin{definition}
For each $z \in \Invol $  define 
$ \kInvPD(z)$ and $ \kInvBPD(z) $
to be the respective sets of   triangular pipe dreams and  bumpless pipe dreams $D$ with $\kpi(D) \circ \kpi(D)^{-1}=z$.
\end{definition}

We refer to the elements of these sets as \defn{$K$-theoretic (generic/bumpless) involution pipe dreams} for $z$.
Consider the alternate $K$-theoretic weight of $D$
defined by
\be 
 \hatkweight(D) = \prod_{(i,j) \in \blank(D)} 
  x_{i\oplus j}
 \prod_{(i,j) \in \upturn(D)} 
 (1+\beta x_{i\oplus j})
\quad
\text{where }x_{i\oplus j} = \begin{cases} x_i &\text{if $i=j$} 
\\ x_i \oplus x_j &\text{if $i\neq j$.}
\end{cases}
\ee
Notice that if $i\neq j$ then $1+\beta x_{i\oplus j} = (1+\beta x_i)(1+\beta x_j)$.
We pose the following conjecture.

\begin{conjecture}\label{iG-conj}
For each $z \in \Invol $ it holds that
$
\iG_z = \sum_{D \in \kInvBPD(z)} \beta^{ |\blank(D)|-\ellinv(z)}\cdot \hatkweight(D).
$
\end{conjecture}

This identity would give a bumpless variant of \cite[Thm.~4.5]{MarbergPawlowski}.
One difficulty in finding a proof is that {transition equations} in the style of \cite{Lenart2003,MarPawSp} have not yet been formulated
 for $\iG_z$.

\begin{example}
When  $z=21543=(1,2)(3,5) \in I_5$
is the involution from Example~\ref{bpd-is-eq}, we have 
\[
\kInvBPD(z)
=
\InvBPD(z)
\sqcup\left\{\
\begin{tikzpicture}[x=\bpdwidth,y=\bpdwidth,line cap=round,line join=round,baseline=(z.base)]
\node at (0,2.5) (z) {};
\draw[gray!70, line width=0.7pt](0,0)--(5,0);
\draw[gray!70, line width=0.7pt](0,0)--(0,5);
\draw[gray!70, line width=0.7pt](0,1)--(5,1);
\draw[gray!70, line width=0.7pt](1,0)--(1,5);
\draw[gray!70, line width=0.7pt](0,2)--(5,2);
\draw[gray!70, line width=0.7pt](2,0)--(2,5);
\draw[gray!70, line width=0.7pt](0,3)--(5,3);
\draw[gray!70, line width=0.7pt](3,0)--(3,5);
\draw[gray!70, line width=0.7pt](0,4)--(5,4);
\draw[gray!70, line width=0.7pt](4,0)--(4,5);
\draw[gray!70, line width=0.7pt](0,5)--(5,5);
\draw[gray!70, line width=0.7pt](5,0)--(5,5);
\draw[thick,rounded corners,color=blue] (3.5,1)--(3.5,2);
\draw[thick,rounded corners,color=blue] (2.5,0)--(2.5,0.5)--(3,0.5);
\draw[thick,rounded corners,color=blue] (2.5,2)--(2.5,1.5)--(2,1.5);
\draw[thick,rounded corners,color=blue] (1.5,2)--(1.5,2.5)--(2,2.5);
\draw[thick,rounded corners,color=blue] (0.5,1)--(0.5,1.5)--(1,1.5);
\draw[thick,rounded corners,color=blue] (0.5,0)--(0.5,1);
\draw[thick,rounded corners,color=blue] (1,1.5)--(2,1.5);
\draw[thick,rounded corners,color=blue] (1.5,1)--(1.5,2);
\draw[thick,rounded corners,color=blue] (2,2.5)--(3,2.5);
\draw[thick,rounded corners,color=blue] (2.5,2)--(2.5,3);
\draw[thick,rounded corners,color=blue] (1.5,0)--(1.5,1);
\draw[thick,rounded corners,color=blue] (3,0.5)--(4,0.5);
\draw[thick,rounded corners,color=blue] (3.5,0)--(3.5,1);
\draw[thick,rounded corners,color=blue] (4,0.5)--(5,0.5);
\draw[thick,rounded corners,color=blue] (4.5,0)--(4.5,1);
\draw[white, line width=0.7pt,fill=white] (0,5)--(5,0)--(5,5)--(0,5);
\draw[black, line width=1.6pt] (0,5)--(5,0)--(0,0)--(0,5);
\end{tikzpicture}
\ \right\}
.
\]
The $K$-theoretic weight of the unreduced element of $\kInvBPD(z)$ 
is \[x_1x_2(x_2+x_1+\beta x_2x_1)(x_3+x_1+\beta x_3x_1)(1+\beta x_4)(1+\beta x_3).\]
We have checked that Conjecture~\ref{iG-conj} holds for this example
(and more generally for all elements of $ \Invol $ up to $n=8$) 
but the value of $\iG_{21543}$ is a sum of 49 monomial terms which is too large to display here.
\end{example}

As each $w \in \cW(z)$ has $w\circ w^{-1}=z$,
we have
$\InvPD(z)\subseteq \kInvPD(z)$
and
$\InvBPD(z)\subseteq \kInvBPD(z)$.
However, $\InvPD(z)$ and $\InvBPD(z)$
may be strictly smaller than the sets of reduced 
elements in 
$ \kInvPD(z)$
and
$\kInvBPD(z)$.
For instance, the triangular bumpless pipe dream
\[
H=
\begin{tikzpicture}[x=\bpdwidth,y=\bpdwidth,line cap=round,line join=round,baseline=(z.base)]
\node at (0,2.25) (z) {};
\draw[gray!70, line width=0.7pt](0,0)--(5,0);
\draw[gray!70, line width=0.7pt](0,0)--(0,5);
\draw[gray!70, line width=0.7pt](0,1)--(5,1);
\draw[gray!70, line width=0.7pt](1,0)--(1,5);
\draw[gray!70, line width=0.7pt](0,2)--(5,2);
\draw[gray!70, line width=0.7pt](2,0)--(2,5);
\draw[gray!70, line width=0.7pt](0,3)--(5,3);
\draw[gray!70, line width=0.7pt](3,0)--(3,5);
\draw[gray!70, line width=0.7pt](0,4)--(5,4);
\draw[gray!70, line width=0.7pt](4,0)--(4,5);
\draw[gray!70, line width=0.7pt](0,5)--(5,5);
\draw[gray!70, line width=0.7pt](5,0)--(5,5);
\draw[thick,rounded corners,color=blue] (3.5,0)--(3.5,0.5)--(4,0.5);
\draw[thick,rounded corners,color=blue] (3,1.5)--(4,1.5);
\draw[thick,rounded corners,color=blue] (1.5,2)--(1.5,2.5)--(2,2.5);
\draw[thick,rounded corners,color=blue] (0.5,1)--(0.5,1.5)--(1,1.5);
\draw[thick,rounded corners,color=blue] (0.5,0)--(0.5,1);
\draw[thick,rounded corners,color=blue] (1,1.5)--(2,1.5);
\draw[thick,rounded corners,color=blue] (1.5,1)--(1.5,2);
\draw[thick,rounded corners,color=blue] (2,1.5)--(3,1.5);
\draw[thick,rounded corners,color=blue] (2.5,1)--(2.5,2);
\draw[thick,rounded corners,color=blue] (1.5,0)--(1.5,1);
\draw[thick,rounded corners,color=blue] (2,2.5)--(3,2.5);
\draw[thick,rounded corners,color=blue] (2.5,2)--(2.5,3);
\draw[thick,rounded corners,color=blue] (2.5,0)--(2.5,1);
\draw[thick,rounded corners,color=blue] (4,0.5)--(5,0.5);
\draw[thick,rounded corners,color=blue] (4.5,0)--(4.5,1);
\draw[white, line width=0.7pt,fill=white] (0,5)--(5,0)--(5,5)--(0,5);
\draw[black, line width=1.6pt] (0,5)--(5,0)--(0,0)--(0,5);
\draw[](0.5,-0.4) node {1};
\draw[](1.5,-0.4) node {2};
\draw[](2.5,-0.4) node {3};
\draw[](3.5,-0.4) node {4};
\draw[](4.5,-0.4) node {5};
\draw[](2.55,3.0) node {3};
\draw[](3.05,2.5) node {2};
\draw[](3.8,1.75) node {1};
\draw[](4.55,1.0) node {5};
\draw[](5.05,0.5) node {4};
\end{tikzpicture}
\]
is reduced, but $\pi(H) = \kpi(H) = 32154 = \kpi(H)\circ \kpi(H)^{-1}$ is not an atom of any involution.
In this example we have $H \in \kInvBPD(32154) $ but $H \notin \InvBPD(z)$ for all $z \in\I(S_5)$.

\subsection{Pipe dream words}

Continue to let $D$ be any triangular pipe dream.
We mention an alternate definition of $\pi(D)$ and $\kpi(D)$,
which involves
 a word with letters in $[n-1]$ associated to  $D$ in the following way.

Start with the empty word and iterate over the northwest-to-southeast diagonals of $D$ from west to east. Read each diagonal in the same west-to-east order. Then,
each time we encounter a $\PTILE$ or $\PDIAG$ tile, add up the number 
 of wires that appear strictly northwest of the given position in the same diagonal, call this number $\sigma$, and append  $\sigma+1$ to our output word. Call the result of this algorithm $\word(D)$.
 
 \begin{example}\label{word-ex}
When $n=8$,   the triangular pipe dream
$D$ from \eqref{t-wirable-ex2} has $D_{ij} \in \left\{ \PTILE, \PDIAG\right\}$
when $(i,j) = (8,3),(8,4), (6,4), (7,6), (8,7), (6,6),(8,8)$
and
$\word(D)=(2,3,2,4,6,3,7)$.
\end{example}

 For a more formal definition of $\word(D)$,
first construct an $n\times n$ matrix 
 $M=M(D)$ whose entry $M_{ij}$ is the number of pipes of $D$ passing through position $(i,j)$. Notice that we have 
\be\label{Mdiag-eq}
\sum_{\substack{(i,j) \in \LT_n \\ i-j= n-k}} M_{ij} =k\quad\text{for each $k \in [n]$.}
\ee
Next, define $\Sigma = \Sigma(D)$ to be the $n\times n$ matrix with 
\be\textstyle
 \Sigma_{ij} = 1 + \sum_{r=1}^{\min\{i,j\}-1} M_{i-r,j-r} \quad\text{if }D_{ij} \in \left\{ \PTILE, \PDIAG\right\}
 \ee
 and zeros everywhere else.
 When $D$ is a triangular pipe dream, both $M$ and $\Sigma$ are lower triangular, and 
the entry
$\Sigma_{ij}$ is either zero or one plus the total number of pipes appearing in $D$ strictly northwest of $(i,j)$ in the same diagonal. 
 
Finally, write $(i,j)\prec (k,l)$
if (a) the two positions are on the same diagonal with $i<k$ or 
(b) $j-i < l-k$.
If $\alpha_1\prec \alpha_2\prec \dots\prec \alpha_k$ are the nonzero positions of $\Sigma$
then $\word(D)$  has the alternate  formula
\be\label{word-eq}
\word(D) = (\Sigma_{\alpha_1},\Sigma_{\alpha_2},\cdots, \Sigma_{\alpha_k}).
\ee

\begin{example}
If $D$ is as in \eqref{t-wirable-ex2},
then (circling the entries where $D_{ij}=\PTILE$ or $\PDIAG$) we have 
\[
M = \begin{bNiceMatrix}[columns-width=auto, margin]
0 \\ 
0 & 0 \\
0 & 0 & 1 \\ 
0 & 0 & 1 & 0 \\
0 & 0 & 1 & 1 & 1 \\
0 & 1 & 2 & \circled{2} & 1 & \circled{2} \\
1 & 1 & 1 & 1 & 1 & \circled{2} & 2 \\
1 & 1 & \circled{2} & \circled{2} & 2 & 2 & \circled{2} &\circled{2} 
\end{bNiceMatrix}
\quand
\Sigma = \begin{bNiceMatrix}[columns-width=auto, margin]
0 \\ 
0 & 0 \\
0 & 0 & 0 \\ 
0 & 0 & 0 & 0 \\
0 & 0 & 0 & 0 & 0 \\
0 & 0 & 0 & \circled{2} & 0 & \circled{3} \\
0 & 0 & 0 & 0 & 0 & \circled{4} & 0 \\
0 & 0 & \circled{2} & \circled{3} & 0 & 0 & \circled{6} & \circled{7} 
\end{bNiceMatrix},
\]
and the resulting word is again $\word(D) = (2,3,2,4,6,3,7)$.
\end{example}

The following identities
relate  $\word(D)$ to $\pi(D)$ and $\kpi(D)$.
Here, we fix $z \in \Invol $.

\begin{proposition}\label{word-pd-prop}
Let $D$ be a triangular pipe dream with $\word(D) = (i_1,i_2,\dots,i_k)$. Then:
\ben
\item[(a)] It always holds that $\pi(D) = s_{i_1} s_{i_2}\cdots s_{i_k}$
and
$\kpi(D) = s_{i_1} \circ s_{i_2}\circ \cdots \circ s_{i_k}$.

\item[(b)] The sequence $\word(D)$ is a reduced word for $\pi(D)$
 if and only if $D$ is reduced.

 \item[(c)] One has $D \in \InvPD(z)$ if and only if $\word(D) $ is a reduced word for some element of $\cW(z)$.
  
 \item[(d)] One has $D \in \kInvPD(z)$ if and only if $z = s_{i_1} \circ s_{i_2}\circ \cdots \circ s_{i_k} \circ \cdots \circ s_{i_2} \circ s_{i_1}$.

\een
\end{proposition}

\begin{proof}
One can deduce this result from the 
discussion of
 \defn{planar histories} in \cite[\S2.3]{Weigandt}.
In more detail, set $\pi(D : -n)=1\in S_n$ and for each $-n<k\leq 0$ form $\pi(D:k) \in S_n$ in one-line notation
 by reading (from northwest to southeast)  the pipe labels on the north and east sides of tiles $D_{1-k,1}$, $D_{2-k,2}$, \dots, $D_{n,n+k}$ and then adding $(n+k+1)(n+k+2)\cdots n$.

Then we have
 $\pi(D:k) = \pi(D:k-1) s_{j_1}s_{j_2}\cdots s_{j_r}$
 where $(j_1,j_2,\dots,j_r)$ is the  subword of $\word(D)$ contributed by the $k$th diagonal.
 Moreover, one has $\ell(\pi(D:k)) = \ell(\pi(D:k-1)) + r$ 
 if and only if no two pipes of $D$ that cross in diagonal $k$
 also cross an odd number of times strictly below this diagonal.
Given these observations, part (b) and
the claim that $\pi(D) = s_{i_1} s_{i_2}\cdots s_{i_k}$
follow by induction on $k$ as $\pi(D) = \pi(D:0)$.
 Combining these facts implies part (c).
 
It remains only to prove that $\kpi(D) = s_{i_1} \circ s_{i_2}\circ \cdots \circ s_{i_k}$ as this implies part (d).
Set $\kpi(D : -n)=1\in S_n$ and for each $-n<k\leq 0$ form $\kpi(D:k) \in S_n$ in the same way as $\pi(D:k)$
but 
using the crossing rules in \eqref{internal-eq}.
Then   $\kpi(D:k) = \kpi(D:k-1) \circ s_{j_1}\circ s_{j_2}\circ \cdots \circ s_{j_r}$
 where $(j_1,j_2,\dots,j_r)$ is the  subword of $\word(D)$ contributed by the $k$th diagonal, so the desired identity follows again by induction.
 \end{proof}
 
 Recall that $\mathsf{db}(D)$ is the number of rows $i \in [n]$
 with $D_{ii} =\BDIAG$.
 
 \begin{corollary}\label{db-cor}
 Suppose $D \in \kInvBPD(z)$. Then $\mathsf{db}(D) \leq \cyc(z)$.
 \end{corollary}
 
 \begin{proof}
 By \eqref{Mdiag-eq} with $k=n$,
 the value of $\mathsf{db}(D)$ is also the number of diagonal tiles $D_{ii} = \PDIAG$. The letters contributed by these tiles give
  the final part of $\word(D)$ and all differ by at least two. We can therefore write $\kpi(D) = u\circ v$ where $u,v\in S_n$ and $v$ is a product of $\mathsf{db}(D)$ simple transpositions that pairwise commute.
  Then $\cyc(v) = \mathsf{db}(D)$ and it follows from \eqref{szs-eq}
  that $z=\kpi(D)\circ \kpi(D)^{-1} = u\circ v\circ u^{-1}$ has $\cyc(z) \geq \cyc(v)$.
 \end{proof}
 
 Here is a good place to note a technical lemma that will be used later.
 
 \begin{lemma}
\label{L: local reconnection}
Let $D$ be a triangular pipe dream.
Suppose $E$ is obtained by replacing a $\PTILE$ tile by a $\ETILE$ tile,
or a $\PDIAG$ tile by a $\EDIAG$ tile. Suppose the pipes involved in these tiles are labeled $a$ and $b$.
\ben\item[(a)] Then $\pi(E)=t_{ab}\pi(D)$.
\item[(b)] If $D$ is reduced, then $E$ is reduced 
if and only if $\ell(\pi(E))=\ell(\pi(D))-1$.
\item[(c)] If $E$ is reduced, then $D$ is reduced if and only if
$\ell(\pi(D))=\ell(\pi(E))+1$.
\een
\end{lemma}
\begin{proof}
The first statement is immediate. 
The remaining statements hold since Proposition~\ref{word-pd-prop} implies that a triangular pipe dream $D$ is reduced
if and only if its number of $\PTILE$
and $\PDIAG$ equals $\ell(\pi(D))$.
\end{proof}

 Given a triangular pipe dream $D$, 
 let $\cross(D)$ be the set of $(i,j) \in \LT_n$ with $D_{ij} = \PTILE$ or $\PDIAG$,
 and let $\bump(D)$ be the set of $(i,j) \in \LT_n$ with $D_{ij} = \ETILE$ or $\EDIAG$.
Then we have
 \be\label{cross-bump-blank-eq}
 |\cross(D)| + |\bump(D)| = |\blank(D)|
 \ee
as a consequence of \eqref{Mdiag-eq}.
We can now be more precise about which subsets of
$K$-theoretic involution pipe dreams 
give $ \InvPD(z)$ and $ \InvBPD(z)$.
As noted at the end of the previous subsection,
these subsets are not always obtained by just restricting to reduced elements.

 \begin{corollary}\label{k-sub-cor}
 For each $z \in \Invol $ the following properties hold:
 \ben
 \item[(a)]
$ \InvPD(z) = \{ D \in \kInvPD(z) : |\cross(D)| = \ellinv(z)\}.$

 \item[(b)]
$ \InvBPD(z) = \{ D \in \kInvBPD(z) : |\cross(D)| = \ellinv(z)\}= \{ D \in \kInvBPD(z) :  |\blank(D)| = \ellinv(z)\}.$
\een
 \end{corollary}
 
 \begin{proof}
 If $D \in \InvPD(z)$ then $|\cross(D)| = |\word(D)| = \ell(\pi(D)) = \ellinv(z)$.
If 
 $D \in \kInvPD(z)$ then 
 \[ \ellinv(z) \leq \ell(\kpi(D))
 \quand 
 \ell(\kpi(D)) \leq |\word(D)| = |\cross(D)|
 \]
  by the definitions of the Demazure product and $\ellinv$ in Section~\ref{inv-sect}.
Equality holds in the second expression if and only if $D$ is reduced, in which case $\kpi(D) = \pi(D)$.
 Thus if  $D \in \kInvPD(z)$ has $|\cross(D)| = \ellinv(z)$
 then $D$ is reduced and $\pi(D) = \kpi(D) \in \cH(z)$ has length $\ellinv(z)$, so $\pi(D) \in \cW(z)$ and $D \in  \InvPD(z)$.

This proves part (a). Part (b) follows from part (a) and 
  \eqref{cross-bump-blank-eq}
  as $|\bump(D)|=0$ if $D$ is bumpless.
  \end{proof}

\subsection{Symmetric pipe dreams}\label{sym-pd-sect}

Triangular pipe dreams are closely related to the \defn{(generic) pipe dreams} previously studied \cite{KZJ,LLS,Weigandt},
 in essentially the same way that strict partitions are related to symmetric partitions,
 as we now explain.
 
Consider a \defn{square tiling} given by a map
$
D : [n]\times [n] \to \cToff=
\Bigl\{
\BTILE,\
\HTILE,\
\VTILE,\
\JTILE, \
\CTILE,\ 
\PTILE,\
\ETILE
\Bigr\}.
$
Such a tiling is defined to be \defn{wirable} in the same way as a triangular tiling. When this holds, we can view $D$ as a collection of distinct pipes that go between connections on the boundary of an $n\times n$ square.

\begin{definition}
A \defn{generic pipe dream} is a wirable tiling $D : [n]\times [n] \to \cToff$ 
such that when $D$ is drawn in 
an $n\times n$ square, 
there are no connections on the west or north sides and
  $n$ connections on both the south and east sides.
  A \defn{bumpless pipe dream} is a generic pipe dream that does not include the \defn{bump tile} $\ETILE$.
 \end{definition} 
 
 When $D$ is a generic pipe dream, the sets $\blank(D)$, $\cross(D)$, $\bump(D)$, and $\upturn(D)$ are defined in the same way as for triangular pipe dreams, but using positions in the entire square $[n]\times[n]$. 
 
  \begin{remark}
 The definition of generic pipe dreams in \cite{KZJ}
 does not include our restrictions on the boundary connections,
 and so gives a simultaneous generalization of both
  the bumpless pipe dreams from \cite{LLS,Weigandt} and the classical \defn{pipe dreams} (or \defn{rc-graphs}) from \cite{BergeronBilley,KnutsonMiller}.
While the latter objects also have involution analogues \cite{HMP6},  these are tilings of a different region of $[n]\times [n]$ compared to our triangular pipe dreams.
  \end{remark}

Given a generic pipe dream $D$, we define a permutation $\kpi(D) \in S_n$
by extending the definition for triangular pipe dreams in the following way.
Again label the $n$ connections on the south side of $D$ by $1,2,3,\dots,n$ from west to east, 
and then update how pipes are labeled whenever they cross
at a $\PTILE$ tile according to the rules in 
\eqref{internal-eq}. Finally, form $\kpi(D)$ in one-line notation by reading the labels on the east side of $D$ from north to south.
For example,
\be\label{pd-ex}
D=\begin{tikzpicture}[x=\bpdwidth,y=\bpdwidth,line cap=round,line join=round,baseline=(z.base)]
\node at (0,2.5) (z) {};
\draw[gray!70, line width=0.7pt](0,0)--(5,0);
\draw[gray!70, line width=0.7pt](0,0)--(0,5);
\draw[gray!70, line width=0.7pt](0,1)--(5,1);
\draw[gray!70, line width=0.7pt](1,0)--(1,5);
\draw[gray!70, line width=0.7pt](0,2)--(5,2);
\draw[gray!70, line width=0.7pt](2,0)--(2,5);
\draw[gray!70, line width=0.7pt](0,3)--(5,3);
\draw[gray!70, line width=0.7pt](3,0)--(3,5);
\draw[gray!70, line width=0.7pt](0,4)--(5,4);
\draw[gray!70, line width=0.7pt](4,0)--(4,5);
\draw[gray!70, line width=0.7pt](0,5)--(5,5);
\draw[gray!70, line width=0.7pt](5,0)--(5,5);
\draw[thick,rounded corners,color=blue] (0.5,0)--(0.5,0.5)--(1,0.5);
\draw[thick,rounded corners,color=blue] (1,0.5)--(2,0.5);
\draw[thick,rounded corners,color=blue] (1.5,0)--(1.5,1);
\draw[thick,rounded corners,color=blue] (1.5,1)--(1.5,1.5)--(2,1.5);
\draw[thick,rounded corners,color=blue] (2,0.5)--(3,0.5);
\draw[thick,rounded corners,color=blue] (2.5,0)--(2.5,1);
\draw[thick,rounded corners,color=blue] (2,1.5)--(3,1.5);
\draw[thick,rounded corners,color=blue] (2.5,1)--(2.5,2);
\draw[thick,rounded corners,color=blue] (3,0.5)--(4,0.5);
\draw[thick,rounded corners,color=blue] (3.5,0)--(3.5,1);
\draw[thick,rounded corners,color=blue] (2.5,2)--(2.5,3);
\draw[thick,rounded corners,color=blue] (3.5,2)--(3.5,1.5)--(3,1.5);
\draw[thick,rounded corners,color=blue] (3.5,1)--(3.5,1.5)--(4,1.5);
\draw[thick,rounded corners,color=blue] (4,0.5)--(5,0.5);
\draw[thick,rounded corners,color=blue] (4.5,0)--(4.5,1);
\draw[thick,rounded corners,color=blue] (2.5,3)--(2.5,4);
\draw[thick,rounded corners,color=blue] (3.5,2)--(3.5,2.5)--(4,2.5);
\draw[thick,rounded corners,color=blue] (4.5,2)--(4.5,1.5)--(4,1.5);
\draw[thick,rounded corners,color=blue] (4.5,1)--(4.5,1.5)--(5,1.5);
\draw[thick,rounded corners,color=blue] (2.5,4)--(2.5,4.5)--(3,4.5);
\draw[thick,rounded corners,color=blue] (4,2.5)--(5,2.5);
\draw[thick,rounded corners,color=blue] (4.5,2)--(4.5,3);
\draw[thick,rounded corners,color=blue] (3,4.5)--(4,4.5);
\draw[thick,rounded corners,color=blue] (4.5,3)--(4.5,3.5)--(5,3.5);
\draw[thick,rounded corners,color=blue] (4,4.5)--(5,4.5);
\draw[black, line width=1.6pt] (0,5)--(5,5)--(5,0)--(0,0)--(0,5);
\draw[](0.5,-0.4) node {1};
\draw[](1.5,-0.4) node {2};
\draw[](2.5,-0.4) node {3};
\draw[](3.5,-0.4) node {4};
\draw[](4.5,-0.4) node {5};
\draw[](5.4,4.5) node {3};
\draw[](5.4,3.5) node {4};
\draw[](5.4,2.5) node {2};
\draw[](5.4,1.5) node {5};
\draw[](5.4,0.5) node {1};
\end{tikzpicture}
\qquand
E=\begin{tikzpicture}[x=\bpdwidth,y=\bpdwidth,line cap=round,line join=round,baseline=(z.base)]
\node at (0,2.5) (z) {};
\draw[gray!70, line width=0.7pt](0,0)--(5,0);
\draw[gray!70, line width=0.7pt](0,0)--(0,5);
\draw[gray!70, line width=0.7pt](0,1)--(5,1);
\draw[gray!70, line width=0.7pt](1,0)--(1,5);
\draw[gray!70, line width=0.7pt](0,2)--(5,2);
\draw[gray!70, line width=0.7pt](2,0)--(2,5);
\draw[gray!70, line width=0.7pt](0,3)--(5,3);
\draw[gray!70, line width=0.7pt](3,0)--(3,5);
\draw[gray!70, line width=0.7pt](0,4)--(5,4);
\draw[gray!70, line width=0.7pt](4,0)--(4,5);
\draw[gray!70, line width=0.7pt](0,5)--(5,5);
\draw[gray!70, line width=0.7pt](5,0)--(5,5);
\draw[thick,rounded corners,color=blue] (0.5,0)--(0.5,1);
\draw[thick,rounded corners,color=blue] (0.5,1)--(0.5,2);
\draw[thick,rounded corners,color=blue] (1.5,0)--(1.5,1);
\draw[thick,rounded corners,color=blue] (0.5,2)--(0.5,2.5)--(1,2.5);
\draw[thick,rounded corners,color=blue] (1.5,1)--(1.5,1.5)--(2,1.5);
\draw[thick,rounded corners,color=blue] (2.5,0)--(2.5,1);
\draw[thick,rounded corners,color=blue] (1,2.5)--(2,2.5);
\draw[thick,rounded corners,color=blue] (2,1.5)--(3,1.5);
\draw[thick,rounded corners,color=blue] (2.5,1)--(2.5,2);
\draw[thick,rounded corners,color=blue] (3.5,0)--(3.5,0.5)--(4,0.5);
\draw[thick,rounded corners,color=blue] (2,2.5)--(3,2.5);
\draw[thick,rounded corners,color=blue] (2.5,2)--(2.5,3);
\draw[thick,rounded corners,color=blue] (3.5,2)--(3.5,1.5)--(3,1.5);
\draw[thick,rounded corners,color=blue] (4,0.5)--(5,0.5);
\draw[thick,rounded corners,color=blue] (4.5,0)--(4.5,1);
\draw[thick,rounded corners,color=blue] (2.5,3)--(2.5,3.5)--(3,3.5);
\draw[thick,rounded corners,color=blue] (3,2.5)--(4,2.5);
\draw[thick,rounded corners,color=blue] (3.5,2)--(3.5,3);
\draw[thick,rounded corners,color=blue] (4.5,1)--(4.5,1.5)--(5,1.5);
\draw[thick,rounded corners,color=blue] (3,3.5)--(4,3.5);
\draw[thick,rounded corners,color=blue] (3.5,3)--(3.5,4);
\draw[thick,rounded corners,color=blue] (4,2.5)--(5,2.5);
\draw[thick,rounded corners,color=blue] (3.5,4)--(3.5,4.5)--(4,4.5);
\draw[thick,rounded corners,color=blue] (4,3.5)--(5,3.5);
\draw[thick,rounded corners,color=blue] (4,4.5)--(5,4.5);
\draw[black, line width=1.6pt] (0,5)--(5,5)--(5,0)--(0,0)--(0,5);
\draw[](0.5,-0.4) node {1};
\draw[](1.5,-0.4) node {2};
\draw[](2.5,-0.4) node {3};
\draw[](3.5,-0.4) node {4};
\draw[](4.5,-0.4) node {5};
\draw[](5.4,4.5) node {3};
\draw[](5.4,3.5) node {2};
\draw[](5.4,2.5) node {1};
\draw[](5.4,1.5) node {5};
\draw[](5.4,0.5) node {4};
\end{tikzpicture}
\ee
are generic pipe dreams when $n=5$, where $E$ is bumpless,
and we have $\kpi(D)=34251$ and $\kpi(E)=32154$.
Alternatively, 
if we define $\word(D) = (i_1,i_2,\dots,i_k)$ by repeating the definition \eqref{word-eq}, so that  
\[\word(D) = (1, 2, 1, 3, 4, 2)
\quand
\word(E)=(2, 1, 4, 2, 1)\]
for the examples above, then 
 $\kpi(D)  = s_{i_1} \circ s_{i_2}\circ \cdots \circ s_{i_k}$.
 
 A generic pipe dream $D$ is \defn{reduced} if no two of its pipes  cross more than once.
 In this case $\word(D)$ is a reduced word for $\kpi(D)$, and 
$\ell(\kpi(D))$ is equal to the number of $\PTILE$ tiles in $D$.

 \begin{definition}
Given $w \in S_n$,  define 
$ \kGPD(w)$ and $ \kBPD(w)$
to be the respective sets of   generic pipe dreams and  bumpless pipe dreams $D$ with $\kpi(D)=w$.
Let $\GPD(w)\subseteq  \kGPD(w)$ and $\BPD(w)\subseteq  \kBPD(w)$ be the subsets of reduced elements in these families.
\end{definition}

Define the \defn{transpose} of each tile in our square tile set by 
\[
\transpose(\HTILE) = \VTILE,\quad \transpose(\VTILE) = \HTILE,
\quand
\transpose(\alpha) = \alpha\quad\text{for }\alpha \notin \Bigl\{
\HTILE,\
\VTILE
\Bigr\}.
\]
The \defn{transpose} of a generic pipe dream $D$ is then the generic pipe dream $\transpose(D)$ with 
\[ \transpose(D)_{ij} =\transpose(D_{ji})
\qquad\text{for all $(i,j) \in [n]\times [n]$.}\] 
For the  examples seen in \eqref{pd-ex}, we have
\[
\transpose(D)=
\begin{tikzpicture}[x=\bpdwidth,y=\bpdwidth,line cap=round,line join=round,baseline=(z.base)]
\node at (0,2.5) (z) {};
\draw[gray!70, line width=0.7pt](0,0)--(5,0);
\draw[gray!70, line width=0.7pt](0,0)--(0,5);
\draw[gray!70, line width=0.7pt](0,1)--(5,1);
\draw[gray!70, line width=0.7pt](1,0)--(1,5);
\draw[gray!70, line width=0.7pt](0,2)--(5,2);
\draw[gray!70, line width=0.7pt](2,0)--(2,5);
\draw[gray!70, line width=0.7pt](0,3)--(5,3);
\draw[gray!70, line width=0.7pt](3,0)--(3,5);
\draw[gray!70, line width=0.7pt](0,4)--(5,4);
\draw[gray!70, line width=0.7pt](4,0)--(4,5);
\draw[gray!70, line width=0.7pt](0,5)--(5,5);
\draw[gray!70, line width=0.7pt](5,0)--(5,5);
\draw[thick,rounded corners,color=blue] (4.5,4)--(4.5,4.5)--(5,4.5);
\draw[thick,rounded corners,color=blue] (3.5,3)--(3.5,3.5)--(4,3.5);
\draw[thick,rounded corners,color=blue] (3.5,2)--(3.5,1.5)--(3,1.5);
\draw[thick,rounded corners,color=blue] (3.5,1)--(3.5,1.5)--(4,1.5);
\draw[thick,rounded corners,color=blue] (2.5,1)--(2.5,1.5)--(3,1.5);
\draw[thick,rounded corners,color=blue] (3.5,1)--(3.5,0.5)--(3,0.5);
\draw[thick,rounded corners,color=blue] (3.5,0)--(3.5,0.5)--(4,0.5);
\draw[thick,rounded corners,color=blue] (0.5,2)--(0.5,2.5)--(1,2.5);
\draw[thick,rounded corners,color=blue] (1.5,0)--(1.5,0.5)--(2,0.5);
\draw[thick,rounded corners,color=blue] (0.5,0)--(0.5,1);
\draw[thick,rounded corners,color=blue] (0.5,1)--(0.5,2);
\draw[thick,rounded corners,color=blue] (1,2.5)--(2,2.5);
\draw[thick,rounded corners,color=blue] (2,2.5)--(3,2.5);
\draw[thick,rounded corners,color=blue] (3,2.5)--(4,2.5);
\draw[thick,rounded corners,color=blue] (3.5,2)--(3.5,3);
\draw[thick,rounded corners,color=blue] (4,2.5)--(5,2.5);
\draw[thick,rounded corners,color=blue] (4.5,2)--(4.5,3);
\draw[thick,rounded corners,color=blue] (2,0.5)--(3,0.5);
\draw[thick,rounded corners,color=blue] (2.5,0)--(2.5,1);
\draw[thick,rounded corners,color=blue] (4,1.5)--(5,1.5);
\draw[thick,rounded corners,color=blue] (4.5,1)--(4.5,2);
\draw[thick,rounded corners,color=blue] (4,3.5)--(5,3.5);
\draw[thick,rounded corners,color=blue] (4.5,3)--(4.5,4);
\draw[thick,rounded corners,color=blue] (4,0.5)--(5,0.5);
\draw[thick,rounded corners,color=blue] (4.5,0)--(4.5,1);
\draw[black, line width=1.6pt] (0,5)--(5,5)--(5,0)--(0,0)--(0,5);
\draw[](0.5,-0.4) node {1};
\draw[](1.5,-0.4) node {2};
\draw[](2.5,-0.4) node {3};
\draw[](3.5,-0.4) node {4};
\draw[](4.5,-0.4) node {5};
\draw[](5.4,4.5) node {5};
\draw[](5.4,3.5) node {3};
\draw[](5.4,2.5) node {1};
\draw[](5.4,1.5) node {2};
\draw[](5.4,0.5) node {4};
\end{tikzpicture}
\qquand
\transpose(E) = 
\begin{tikzpicture}[x=\bpdwidth,y=\bpdwidth,line cap=round,line join=round,baseline=(z.base)]
\node at (0,2.5) (z) {};
\draw[gray!70, line width=0.7pt](0,0)--(5,0);
\draw[gray!70, line width=0.7pt](0,0)--(0,5);
\draw[gray!70, line width=0.7pt](0,1)--(5,1);
\draw[gray!70, line width=0.7pt](1,0)--(1,5);
\draw[gray!70, line width=0.7pt](0,2)--(5,2);
\draw[gray!70, line width=0.7pt](2,0)--(2,5);
\draw[gray!70, line width=0.7pt](0,3)--(5,3);
\draw[gray!70, line width=0.7pt](3,0)--(3,5);
\draw[gray!70, line width=0.7pt](0,4)--(5,4);
\draw[gray!70, line width=0.7pt](4,0)--(4,5);
\draw[gray!70, line width=0.7pt](0,5)--(5,5);
\draw[gray!70, line width=0.7pt](5,0)--(5,5);
\draw[thick,rounded corners,color=blue] (2.5,4)--(2.5,4.5)--(3,4.5);
\draw[thick,rounded corners,color=blue] (3.5,3)--(3.5,3.5)--(4,3.5);
\draw[thick,rounded corners,color=blue] (4.5,1)--(4.5,1.5)--(5,1.5);
\draw[thick,rounded corners,color=blue] (3.5,2)--(3.5,1.5)--(3,1.5);
\draw[thick,rounded corners,color=blue] (1.5,2)--(1.5,2.5)--(2,2.5);
\draw[thick,rounded corners,color=blue] (3.5,0)--(3.5,0.5)--(4,0.5);
\draw[thick,rounded corners,color=blue] (0.5,1)--(0.5,1.5)--(1,1.5);
\draw[thick,rounded corners,color=blue] (0.5,0)--(0.5,1);
\draw[thick,rounded corners,color=blue] (1,1.5)--(2,1.5);
\draw[thick,rounded corners,color=blue] (1.5,1)--(1.5,2);
\draw[thick,rounded corners,color=blue] (2,1.5)--(3,1.5);
\draw[thick,rounded corners,color=blue] (2.5,1)--(2.5,2);
\draw[thick,rounded corners,color=blue] (3,2.5)--(4,2.5);
\draw[thick,rounded corners,color=blue] (3.5,2)--(3.5,3);
\draw[thick,rounded corners,color=blue] (4,3.5)--(5,3.5);
\draw[thick,rounded corners,color=blue] (1.5,0)--(1.5,1);
\draw[thick,rounded corners,color=blue] (2,2.5)--(3,2.5);
\draw[thick,rounded corners,color=blue] (2.5,2)--(2.5,3);
\draw[thick,rounded corners,color=blue] (4,2.5)--(5,2.5);
\draw[thick,rounded corners,color=blue] (2.5,0)--(2.5,1);
\draw[thick,rounded corners,color=blue] (2.5,3)--(2.5,4);
\draw[thick,rounded corners,color=blue] (3,4.5)--(4,4.5);
\draw[thick,rounded corners,color=blue] (4,4.5)--(5,4.5);
\draw[thick,rounded corners,color=blue] (4,0.5)--(5,0.5);
\draw[thick,rounded corners,color=blue] (4.5,0)--(4.5,1);
\draw[black, line width=1.6pt] (0,5)--(5,5)--(5,0)--(0,0)--(0,5);
\draw[](0.5,-0.4) node {1};
\draw[](1.5,-0.4) node {2};
\draw[](2.5,-0.4) node {3};
\draw[](3.5,-0.4) node {4};
\draw[](4.5,-0.4) node {5};
\draw[](5.4,4.5) node {3};
\draw[](5.4,3.5) node {2};
\draw[](5.4,2.5) node {1};
\draw[](5.4,1.5) node {5};
\draw[](5.4,0.5) node {4};
\end{tikzpicture}
.
\]
It always holds that $\kpi(\transpose(D)) = \kpi(D)^{-1}$.
We say that $D$ is \defn{symmetric} if $D=\transpose(D)$.

Finally, define the \defn{double} of a triangular pipe dream $E$ to be the generic pipe dream $\double(E)$ with 
\be\label{double-def}
\double(E)_{ij} = E_{ij}
\qquand 
\double(E)_{ji} = \transpose(E_{ij})\qquad\text{for }1\leq j < i \leq n,
\ee
and whose diagonal tiles $\double(E)_{ii} $ are formed from the diagonal tiles $E_{ii}$ by applying the transformation
\be\label{double-eq}
\BDIAG \mapsto \BTILE,
\qquad
\JDIAG \mapsto \JTILE,
\qquad
\CDIAG \mapsto \CTILE,
\qquad
\PDIAG \mapsto \PTILE,
\qquad
\EDIAG \mapsto \ETILE.
\ee
Applying this operation to the triangular pipe dreams 
in Example~\ref{bpd-is-eq} gives \[
\begin{tikzpicture}[x=\bpdwidth,y=\bpdwidth,line cap=round,line join=round,baseline=(z.base)]
\node at (0,2.5) (z) {};
\draw[gray!70, line width=0.7pt](0,0)--(5,0);
\draw[gray!70, line width=0.7pt](0,0)--(0,5);
\draw[gray!70, line width=0.7pt](0,1)--(5,1);
\draw[gray!70, line width=0.7pt](1,0)--(1,5);
\draw[gray!70, line width=0.7pt](0,2)--(5,2);
\draw[gray!70, line width=0.7pt](2,0)--(2,5);
\draw[gray!70, line width=0.7pt](0,3)--(5,3);
\draw[gray!70, line width=0.7pt](3,0)--(3,5);
\draw[gray!70, line width=0.7pt](0,4)--(5,4);
\draw[gray!70, line width=0.7pt](4,0)--(4,5);
\draw[gray!70, line width=0.7pt](0,5)--(5,5);
\draw[gray!70, line width=0.7pt](5,0)--(5,5);
\draw[thick,rounded corners,color=blue] (4.5,2)--(4.5,2.5)--(5,2.5);
\draw[thick,rounded corners,color=blue] (3.5,1)--(3.5,1.5)--(4,1.5);
\draw[thick,rounded corners,color=blue] (2.5,0)--(2.5,0.5)--(3,0.5);
\draw[thick,rounded corners,color=blue] (2.5,3)--(2.5,2.5)--(2,2.5);
\draw[thick,rounded corners,color=blue] (0.5,2)--(0.5,2.5)--(1,2.5);
\draw[thick,rounded corners,color=blue] (2.5,4)--(2.5,4.5)--(3,4.5);
\draw[thick,rounded corners,color=blue] (1.5,3)--(1.5,3.5)--(2,3.5);
\draw[thick,rounded corners,color=blue] (0.5,0)--(0.5,1);
\draw[thick,rounded corners,color=blue] (0.5,1)--(0.5,2);
\draw[thick,rounded corners,color=blue] (1,2.5)--(2,2.5);
\draw[thick,rounded corners,color=blue] (1.5,2)--(1.5,3);
\draw[thick,rounded corners,color=blue] (2,3.5)--(3,3.5);
\draw[thick,rounded corners,color=blue] (2.5,3)--(2.5,4);
\draw[thick,rounded corners,color=blue] (3,4.5)--(4,4.5);
\draw[thick,rounded corners,color=blue] (4,4.5)--(5,4.5);
\draw[thick,rounded corners,color=blue] (1.5,0)--(1.5,1);
\draw[thick,rounded corners,color=blue] (1.5,1)--(1.5,2);
\draw[thick,rounded corners,color=blue] (3,3.5)--(4,3.5);
\draw[thick,rounded corners,color=blue] (4,3.5)--(5,3.5);
\draw[thick,rounded corners,color=blue] (3,0.5)--(4,0.5);
\draw[thick,rounded corners,color=blue] (3.5,0)--(3.5,1);
\draw[thick,rounded corners,color=blue] (4,0.5)--(5,0.5);
\draw[thick,rounded corners,color=blue] (4.5,0)--(4.5,1);
\draw[thick,rounded corners,color=blue] (4,1.5)--(5,1.5);
\draw[thick,rounded corners,color=blue] (4.5,1)--(4.5,2);
\draw[black, line width=1.6pt] (0,5)--(5,5)--(5,0)--(0,0)--(0,5);
\end{tikzpicture}
\qquad\begin{tikzpicture}[x=\bpdwidth,y=\bpdwidth,line cap=round,line join=round,baseline=(z.base)]
\node at (0,2.5) (z) {};
\draw[gray!70, line width=0.7pt](0,0)--(5,0);
\draw[gray!70, line width=0.7pt](0,0)--(0,5);
\draw[gray!70, line width=0.7pt](0,1)--(5,1);
\draw[gray!70, line width=0.7pt](1,0)--(1,5);
\draw[gray!70, line width=0.7pt](0,2)--(5,2);
\draw[gray!70, line width=0.7pt](2,0)--(2,5);
\draw[gray!70, line width=0.7pt](0,3)--(5,3);
\draw[gray!70, line width=0.7pt](3,0)--(3,5);
\draw[gray!70, line width=0.7pt](0,4)--(5,4);
\draw[gray!70, line width=0.7pt](4,0)--(4,5);
\draw[gray!70, line width=0.7pt](0,5)--(5,5);
\draw[gray!70, line width=0.7pt](5,0)--(5,5);
\draw[thick,rounded corners,color=blue] (4.5,2)--(4.5,2.5)--(5,2.5);
\draw[thick,rounded corners,color=blue] (3.5,1)--(3.5,1.5)--(4,1.5);
\draw[thick,rounded corners,color=blue] (2.5,0)--(2.5,0.5)--(3,0.5);
\draw[thick,rounded corners,color=blue] (3.5,3)--(3.5,2.5)--(3,2.5);
\draw[thick,rounded corners,color=blue] (3.5,3)--(3.5,3.5)--(4,3.5);
\draw[thick,rounded corners,color=blue] (0.5,2)--(0.5,2.5)--(1,2.5);
\draw[thick,rounded corners,color=blue] (2.5,2)--(2.5,1.5)--(2,1.5);
\draw[thick,rounded corners,color=blue] (2.5,4)--(2.5,4.5)--(3,4.5);
\draw[thick,rounded corners,color=blue] (1.5,1)--(1.5,1.5)--(2,1.5);
\draw[thick,rounded corners,color=blue] (0.5,0)--(0.5,1);
\draw[thick,rounded corners,color=blue] (0.5,1)--(0.5,2);
\draw[thick,rounded corners,color=blue] (1,2.5)--(2,2.5);
\draw[thick,rounded corners,color=blue] (2,2.5)--(3,2.5);
\draw[thick,rounded corners,color=blue] (2.5,2)--(2.5,3);
\draw[thick,rounded corners,color=blue] (4,3.5)--(5,3.5);
\draw[thick,rounded corners,color=blue] (1.5,0)--(1.5,1);
\draw[thick,rounded corners,color=blue] (2.5,3)--(2.5,4);
\draw[thick,rounded corners,color=blue] (3,4.5)--(4,4.5);
\draw[thick,rounded corners,color=blue] (4,4.5)--(5,4.5);
\draw[thick,rounded corners,color=blue] (3,0.5)--(4,0.5);
\draw[thick,rounded corners,color=blue] (3.5,0)--(3.5,1);
\draw[thick,rounded corners,color=blue] (4,0.5)--(5,0.5);
\draw[thick,rounded corners,color=blue] (4.5,0)--(4.5,1);
\draw[thick,rounded corners,color=blue] (4,1.5)--(5,1.5);
\draw[thick,rounded corners,color=blue] (4.5,1)--(4.5,2);
\draw[black, line width=1.6pt] (0,5)--(5,5)--(5,0)--(0,0)--(0,5);
\end{tikzpicture}
\qquad\begin{tikzpicture}[x=\bpdwidth,y=\bpdwidth,line cap=round,line join=round,baseline=(z.base)]
\node at (0,2.5) (z) {};
\draw[gray!70, line width=0.7pt](0,0)--(5,0);
\draw[gray!70, line width=0.7pt](0,0)--(0,5);
\draw[gray!70, line width=0.7pt](0,1)--(5,1);
\draw[gray!70, line width=0.7pt](1,0)--(1,5);
\draw[gray!70, line width=0.7pt](0,2)--(5,2);
\draw[gray!70, line width=0.7pt](2,0)--(2,5);
\draw[gray!70, line width=0.7pt](0,3)--(5,3);
\draw[gray!70, line width=0.7pt](3,0)--(3,5);
\draw[gray!70, line width=0.7pt](0,4)--(5,4);
\draw[gray!70, line width=0.7pt](4,0)--(4,5);
\draw[gray!70, line width=0.7pt](0,5)--(5,5);
\draw[gray!70, line width=0.7pt](5,0)--(5,5);
\draw[thick,rounded corners,color=blue] (4.5,2)--(4.5,2.5)--(5,2.5);
\draw[thick,rounded corners,color=blue] (3.5,1)--(3.5,1.5)--(4,1.5);
\draw[thick,rounded corners,color=blue] (2.5,0)--(2.5,0.5)--(3,0.5);
\draw[thick,rounded corners,color=blue] (1.5,3)--(1.5,3.5)--(2,3.5);
\draw[thick,rounded corners,color=blue] (2.5,2)--(2.5,1.5)--(2,1.5);
\draw[thick,rounded corners,color=blue] (0.5,1)--(0.5,1.5)--(1,1.5);
\draw[thick,rounded corners,color=blue] (3.5,3)--(3.5,2.5)--(3,2.5);
\draw[thick,rounded corners,color=blue] (3.5,4)--(3.5,4.5)--(4,4.5);
\draw[thick,rounded corners,color=blue] (2.5,2)--(2.5,2.5)--(3,2.5);
\draw[thick,rounded corners,color=blue] (0.5,0)--(0.5,1);
\draw[thick,rounded corners,color=blue] (1,1.5)--(2,1.5);
\draw[thick,rounded corners,color=blue] (1.5,1)--(1.5,2);
\draw[thick,rounded corners,color=blue] (3,3.5)--(4,3.5);
\draw[thick,rounded corners,color=blue] (3.5,3)--(3.5,4);
\draw[thick,rounded corners,color=blue] (4,4.5)--(5,4.5);
\draw[thick,rounded corners,color=blue] (1.5,0)--(1.5,1);
\draw[thick,rounded corners,color=blue] (1.5,2)--(1.5,3);
\draw[thick,rounded corners,color=blue] (2,3.5)--(3,3.5);
\draw[thick,rounded corners,color=blue] (4,3.5)--(5,3.5);
\draw[thick,rounded corners,color=blue] (3,0.5)--(4,0.5);
\draw[thick,rounded corners,color=blue] (3.5,0)--(3.5,1);
\draw[thick,rounded corners,color=blue] (4,0.5)--(5,0.5);
\draw[thick,rounded corners,color=blue] (4.5,0)--(4.5,1);
\draw[thick,rounded corners,color=blue] (4,1.5)--(5,1.5);
\draw[thick,rounded corners,color=blue] (4.5,1)--(4.5,2);
\draw[black, line width=1.6pt] (0,5)--(5,5)--(5,0)--(0,0)--(0,5);
\end{tikzpicture}
\qquad\begin{tikzpicture}[x=\bpdwidth,y=\bpdwidth,line cap=round,line join=round,baseline=(z.base)]
\node at (0,2.5) (z) {};
\draw[gray!70, line width=0.7pt](0,0)--(5,0);
\draw[gray!70, line width=0.7pt](0,0)--(0,5);
\draw[gray!70, line width=0.7pt](0,1)--(5,1);
\draw[gray!70, line width=0.7pt](1,0)--(1,5);
\draw[gray!70, line width=0.7pt](0,2)--(5,2);
\draw[gray!70, line width=0.7pt](2,0)--(2,5);
\draw[gray!70, line width=0.7pt](0,3)--(5,3);
\draw[gray!70, line width=0.7pt](3,0)--(3,5);
\draw[gray!70, line width=0.7pt](0,4)--(5,4);
\draw[gray!70, line width=0.7pt](4,0)--(4,5);
\draw[gray!70, line width=0.7pt](0,5)--(5,5);
\draw[gray!70, line width=0.7pt](5,0)--(5,5);
\draw[thick,rounded corners,color=blue] (1.5,4)--(1.5,4.5)--(2,4.5);
\draw[thick,rounded corners,color=blue] (0.5,3)--(0.5,3.5)--(1,3.5);
\draw[thick,rounded corners,color=blue] (4.5,2)--(4.5,2.5)--(5,2.5);
\draw[thick,rounded corners,color=blue] (3.5,1)--(3.5,1.5)--(4,1.5);
\draw[thick,rounded corners,color=blue] (2.5,0)--(2.5,0.5)--(3,0.5);
\draw[thick,rounded corners,color=blue] (0.5,0)--(0.5,1);
\draw[thick,rounded corners,color=blue] (0.5,1)--(0.5,2);
\draw[thick,rounded corners,color=blue] (0.5,2)--(0.5,3);
\draw[thick,rounded corners,color=blue] (1,3.5)--(2,3.5);
\draw[thick,rounded corners,color=blue] (1.5,3)--(1.5,4);
\draw[thick,rounded corners,color=blue] (2,3.5)--(3,3.5);
\draw[thick,rounded corners,color=blue] (3,3.5)--(4,3.5);
\draw[thick,rounded corners,color=blue] (4,3.5)--(5,3.5);
\draw[thick,rounded corners,color=blue] (1.5,0)--(1.5,1);
\draw[thick,rounded corners,color=blue] (1.5,1)--(1.5,2);
\draw[thick,rounded corners,color=blue] (1.5,2)--(1.5,3);
\draw[thick,rounded corners,color=blue] (2,4.5)--(3,4.5);
\draw[thick,rounded corners,color=blue] (3,4.5)--(4,4.5);
\draw[thick,rounded corners,color=blue] (4,4.5)--(5,4.5);
\draw[thick,rounded corners,color=blue] (3,0.5)--(4,0.5);
\draw[thick,rounded corners,color=blue] (3.5,0)--(3.5,1);
\draw[thick,rounded corners,color=blue] (4,0.5)--(5,0.5);
\draw[thick,rounded corners,color=blue] (4.5,0)--(4.5,1);
\draw[thick,rounded corners,color=blue] (4,1.5)--(5,1.5);
\draw[thick,rounded corners,color=blue] (4.5,1)--(4.5,2);
\draw[black, line width=1.6pt] (0,5)--(5,5)--(5,0)--(0,0)--(0,5);
\end{tikzpicture}.
\]
Notice that $E$ is bumpless if and only if $\double(E)$ is bumpless.

\begin{proposition}\label{sym-bij-prop}
For each $z \in \Invol $,  the operation $\double$ is
a bijection 
from 
 $\kInvPD(z) $ and $\kInvBPD(z) $
to the subsets of symmetric elements in 
$\kGPD(z)$ and $\kBPD(z)$,  respectively.
\end{proposition}

\begin{proof}
Given a symmetric generic pipe dream $D$, let $\halve(D)$ be the triangular pipe
dream formed by forgetting $D_{ij}$ for $i<j$ and applying the inverse of the transformations in  \eqref{double-eq} to all tiles $D_{ii}$ on the main diagonal.
It is clear that $\double$ and $\halve$ are inverse operations.

To establish the claimed bijections, it is enough to show that if $D \in \kInvPD(z)$ and $E=\transpose(E) \in \kGPD(z)$
then $\double(D) \in \kGPD(z)$ and $\halve(E) \in \kInvPD(z)$.
To prove this, let $\word(D: k)$ be the (increasing) consecutive subsequence of $\word(D)$ contributed by the tiles $D_{ij}$ with $j-i=k$.
Then 
\[
\ba
\word(D) &= \word(D:-n+1) \cdots \word(D:-2) \word(D:-1) \word(D:0) \quand
\\
\word(\double(D)) &= \word(D) \word(D:-1) \word(D:-2)\cdots \word(D:-n+1).
\ea
\]
All letters in $\word(D:k)$ differ by at least two. Since we have $s_i \circ s_j = s_j\circ s_i$ whenever $|i-j|\geq 2$ along with $s_i\circ s_i = s_i$,
 when we compute the Demazure product of the simple transpositions indexed by $\word(\double(D))$,
we can insert an extra copy of $\word(D:0) $ in the middle and then reverse the order of the factors indexed by the second occurrence of each $\word(D:k)$ without changing the result. It follows that 
\[\kpi(\double(D)) = \kpi(D) \circ \kpi(D)^{-1} = z\] so $\double(D) \in \kGPD(z)$.
By similar reasoning, if $E=\transpose(E) \in \kGPD(z)$ then we have \[  \kpi(\halve(E)) \circ  \kpi(\halve(E))^{-1} = \kpi(E) = z\] so $\halve(E) \in \kInvPD(z)$
 as needed.
\end{proof}

\begin{remark}
If $D $ is any symmetric generic pipe dream then  the 
proof above shows that 
\be\kpi(D) =  \kpi(\halve(D)) \circ \kpi(\halve(D))^{-1} \in \Invol .\ee Hence, 
Proposition~\ref{sym-bij-prop}  holds for all $z \in S_n$ if we 
define $\kInvPD(z)=\kInvBPD(z) = \varnothing$ when $z \neq z^{-1}$.
\end{remark}

Fix $z \in \Invol $ and a symmetric generic pipe dream $D=\transpose(D) \in \kGPD(z)$. Then let
\be
\symDef(D) = |\cross(D)|  + | \{ (i,j) \in \cross(D) : i=j\}| - |\Inv(z)| - |\Cyc(z)|,
\ee
where $\Inv(z)$ and  $\Cyc(z)$ are as in \eqref{inv-def} and \eqref{cyc-def}. One has
\be
\symDef(D) = |\blank(D)|  + | \{ (i,j) \in \blank(D) : i=j\}| - |\Inv(z)| - |\Cyc(z)|
\ee
when $D \in \kBPD(z)$ in view of \eqref{Mdiag-eq}.
We refer to $\symDef(D)$ as the \defn{symmetric defect} of $D$. Note that
\[
\tfrac{1}{2} \symDef(D) = |\word(\halve(D))| - \ellinv(z) \geq \ell(\kpi(\halve(D))) - \ellinv(z) \geq 0
\]
since $\kpi(\halve(D)) \in \cH(z)$.
Comparing this equation with Corollary~\ref{k-sub-cor}
shows that $\halve(D) \in \InvPD(z)$ if and only if $ \symDef(D) = 0$.
Thus, we get the following corollary from Proposition~\ref{sym-bij-prop}:

\begin{corollary}\label{sym-bij-cor}
For each $z \in \Invol $,  the operation $\double$ is
a bijection 
from 
 $\InvPD(z) $ and $\InvBPD(z) $, respectively,
to the subsets of symmetric elements in 
$\kGPD(z)$ and $\kBPD(z)$ with symmetric defect zero.
\end{corollary}

\section{Generating involution pipe dreams}\label{gen-sect}

Fix an element $z \in \Invol $.
This section explains how to generate the set
$\InvBPD(z)$. This can be done by 
 applying certain \defn{droop moves} to a specific pipe dream derived from the Rothe diagram.
A similar process constructs  $\BPD(w)$ for $w \in S_n$ \cite[\S5.2]{LLS}.

\subsection{Rothe pipe dreams and droop moves}

We begin with a useful observation.
The \defn{bend tiles} of a (square or triangular) pipe dream
consist of all instances of $\JTILE$, $\CTILE$, $\ETILE$, $\JDIAG$, $\CDIAG$, and $\EDIAG$.

\begin{lemma}
A (square or triangular) pipe dream  is uniquely determined
by its bend tiles.
\end{lemma}

Using this property, we can characterize one reduced involution pipe dream for $z$.
Namely,
define the \defn{involution Rothe pipe dream} of $z$
to be the triangular pipe dream $\invRothe(z)$ 
whose 
bend tiles all have the form $\CTILE$ and $\CDIAG$ and occur in the positions $(i,j) \in \LT_n$ with $i\geq j=z(i)$.
More explicitly, the tiles of $\invRothe(z)$ are given for $i \in [n]$ and $j \in [i-1]$ by
\be\label{inv-rothe}
\invRothe(z)_{ij} = \begin{cases}
\CTILE &\text{if }z(i) =j \\[-10pt]\\
\PTILE &\text{if }j>z(i)\text{ and }i>z(j)\\[-10pt]\\
\HTILE &\text{if }j>z(i)\text{ and }i<z(j)\\[-10pt]\\
\VTILE &\text{if }j<z(i)\text{ and }i>z(j)\\[-10pt]\\
\BTILE &\text{if }j<z(i)\text{ and }i<z(j)
\end{cases}
\quand
\invRothe(z)_{ii} = \begin{cases}
\CDIAG &\text{if }z(i) =i \\[-10pt]\\
\PDIAG &\text{if }i>z(i)\\[-10pt]\\
\BDIAG &\text{if }i<z(i).
\end{cases}
\ee
We see from this description that $\blank(\invRothe(z)) = \iD(z)$ is the involution Rothe diagram of $z$, which explains the name.
The first triangular pipe dream in Example~\ref{bpd-is-eq} is  $\invRothe(21543)$.

\begin{proposition}
For each $z \in \Invol $ it holds that $\invRothe(z)\in\InvBPD(z)$.
\end{proposition}

\begin{proof}
Let $b_1<b_2<\dots<b_k$ be the elements $b \in [n]$ with $z(b) \leq b$.
 Then the one-line representation of $\pi(\invRothe(z))$
 is formed from the word $b_1b_2\cdots b_k$ by inserting $z(b_i)$ directly after each $b_i$ with $z(b_i)<b_i$. 
 For example, when $z = 21543$ we have $\{b_1<b_2<\dots<b_k\} = \{ 2 < 4 < 5\}$ and 
 $\pi(\invRothe(z)) = 21453$.
 Using this description and Theorem~\ref{cW-lem}, it is straightforward to check that $\pi(\invRothe(z)) \in \cW(z)$;
 alternatively, this holds by \cite[Thm.~6.10]{HMP2}.
 Thus $\invRothe(z)$ is  an element of $\InvBPD(z)$ as it is clearly reduced.
\end{proof}

We now fix a triangular pipe dream $D$, which does not need to be reduced. Let $(a,b),(i,j) \in \LT_n$ be  lower triangular positions with $a< i$ and $b < j$.
A \defn{simple droop move} from $(a,b)$ to $(i,j)$ 
refers to a way of transforming $D$ to another triangular pipe dream by modifying a single pipe as in the pictures
\be\label{simple-eq}
\barr{c}
\begin{tikzpicture}[x=\bpdsmallwidth,y=\bpdsmallwidth,line cap=round,line join=round,baseline=(z.base)]
\node at (0,1.75) (z) {};
\draw[color=red!50,fill=red!50] (0,4)--(1,4)--(1,3)--(0,3)--(0,4);
\draw[thick,rounded corners,color=blue] (0.5,3)--(0.5,3.5)--(1,3.5);
\draw[thick,rounded corners,color=blue] (0.5,0)--(0.5,1);
\draw[thick,rounded corners,color=blue] (0.5,1)--(0.5,2);
\draw[thick,rounded corners,color=blue] (0.5,2)--(0.5,3);
\draw[thick,rounded corners,color=blue] (1,3.5)--(2,3.5);
\draw[thick,rounded corners,color=blue] (2,3.5)--(3,3.5);
\draw[thick,rounded corners,color=blue] (3,3.5)--(4,3.5);
\draw[gray, dotted, line width=1pt] (0,4)--(4,4)--(4,0)--(0,0)--(0,4);
\end{tikzpicture}
\ 
\mapsto
\ 
\begin{tikzpicture}[x=\bpdsmallwidth,y=\bpdsmallwidth,line cap=round,line join=round,baseline=(z.base)]
\node at (0,1.75) (z) {};
\draw[color=red!50,fill=red!50] (3,0)--(4,0)--(4,1)--(3,1)--(3,0);
\draw[thick,rounded corners,color=blue] (0.5,0)--(0.5,0.5)--(1,0.5);
\draw[thick,rounded corners,color=blue] (1,0.5)--(2,0.5);
\draw[thick,rounded corners,color=blue] (2,0.5)--(3,0.5);
\draw[thick,rounded corners,color=blue] (3.5,1)--(3.5,0.5)--(3,0.5);
\draw[thick,rounded corners,color=blue] (3.5,1)--(3.5,2);
\draw[thick,rounded corners,color=blue] (3.5,2)--(3.5,3);
\draw[thick,rounded corners,color=blue] (3.5,3)--(3.5,3.5)--(4,3.5);
\draw[gray, dotted, line width=1pt] (0,4)--(4,4)--(4,0)--(0,0)--(0,4);
\end{tikzpicture}
\sdgap
b < j < a< i
\earr
\qquad
\barr{c}
\begin{tikzpicture}[x=\bpdsmallwidth,y=\bpdsmallwidth,line cap=round,line join=round,baseline=(z.base)]
\node at (0,1.75) (z) {};
\draw[color=red!50,fill=red!50] (0,4)--(1,4)--(1,3)--(0,3)--(0,4);
\draw[thick,rounded corners,color=blue] (0.5,3)--(0.5,4);
\draw[thick,rounded corners,color=blue] (0.5,0)--(0.5,1);
\draw[thick,rounded corners,color=blue] (0.5,1)--(0.5,2);
\draw[thick,rounded corners,color=blue] (0.5,2)--(0.5,3);
\draw[thick,rounded corners,color=blue] (1,3.5)--(2,3.5);
\draw[thick,rounded corners,color=blue] (2,3.5)--(3,3.5);
\draw[thick,rounded corners,color=blue] (3,3.5)--(4,3.5);
\draw[gray, dotted, line width=1pt] (0,4)--(4,4)--(4,0)--(0,0)--(0,4);
\draw[white, line width=0.7pt,fill=white] (0,4)--(4,0)--(4,4)--(0,4);
\draw[black, line width=1.6pt] (0,4)--(4,0);
\end{tikzpicture}
\
\mapsto
\
\begin{tikzpicture}[x=\bpdsmallwidth,y=\bpdsmallwidth,line cap=round,line join=round,baseline=(z.base)]
\node at (0,1.75) (z) {};
\draw[color=red!50,fill=red!50] (3,0)--(4,0)--(4,1)--(3,1)--(3,0);
\draw[thick,rounded corners,color=blue] (0.5,0)--(0.5,0.5)--(1,0.5);
\draw[thick,rounded corners,color=blue] (1,0.5)--(2,0.5);
\draw[thick,rounded corners,color=blue] (2,0.5)--(3,0.5);
\draw[thick,rounded corners,color=blue] (4,0.5)--(3,0.5);
\draw[thick,rounded corners,color=blue] (3.5,1)--(3.5,2);
\draw[thick,rounded corners,color=blue] (3.5,2)--(3.5,3);
\draw[thick,rounded corners,color=blue] (3.5,3)--(3.5,3.5)--(4,3.5);
\draw[gray, dotted, line width=1pt] (0,4)--(4,4)--(4,0)--(0,0)--(0,4);
\draw[white, line width=0.7pt,fill=white] (0,4)--(4,0)--(4,4)--(0,4);
\draw[black, line width=1.6pt] (0,4)--(4,0);
\end{tikzpicture}
\sdgap
 b = a < j = i
\earr
\qquad
\barr{c}
\begin{tikzpicture}[x=\bpdsmallwidth,y=\bpdsmallwidth,line cap=round,line join=round,baseline=(z.base)]
\node at (0,2.25) (z) {};
\draw[color=red!50,fill=red!50] (0,4)--(0,5)--(1,5)--(1,4)--(0,4);
\draw[thick,rounded corners,color=blue] (0.5,3)--(0.5,5);
\draw[thick,rounded corners,color=blue] (0.5,0)--(0.5,1);
\draw[thick,rounded corners,color=blue] (0.5,1)--(0.5,2);
\draw[thick,rounded corners,color=blue] (0.5,2)--(0.5,3);
\draw[gray, dotted, line width=1pt] (0,5)--(3,5)--(3,0)--(0,0)--(0,5);
\draw[white, line width=0.7pt,fill=white] (0,5)--(3,2)--(3,5)--(0,5);
\draw[black, line width=1.6pt] (0,5)--(3,2);
\end{tikzpicture}
\
\mapsto
\
\begin{tikzpicture}[x=\bpdsmallwidth,y=\bpdsmallwidth,line cap=round,line join=round,baseline=(z.base)]
\node at (0,2.25) (z) {};
\draw[color=red!50,fill=red!50] (2,0)--(2,1)--(3,1)--(3,0)--(2,0);
\draw[thick,rounded corners,color=blue] (0.5,0)--(0.5,0.5)--(1,0.5);
\draw[thick,rounded corners,color=blue] (1,0.5)--(2,0.5);
\draw[thick,rounded corners,color=blue] (2.5,1)--(2.5,0.5)--(2,0.5);
\draw[thick,rounded corners,color=blue] (2.5,1)--(2.5,2);
\draw[thick,rounded corners,color=blue] (2.5,2)--(2.5,3);
\draw[thick,rounded corners,color=blue] (2.5,3)--(2.5,3.5)--(3,3.5);
\draw[gray, dotted, line width=1pt] (0,5)--(3,5)--(3,0)--(0,0)--(0,5);
\draw[white, line width=0.7pt,fill=white] (0,5)--(3,2)--(3,5)--(0,5);
\draw[black, line width=1.6pt] (0,5)--(3,2);
\end{tikzpicture}
\sdgap
 b = a < j < i
\earr,
\ee
or as in the following examples (where a second pipe must meet the modified pipe at the diagonal):
\be\label{simple-eq2}
\barr{c}
\begin{tikzpicture}[x=\bpdsmallwidth,y=\bpdsmallwidth,line cap=round,line join=round,baseline=(z.base)]
\node at (0,2.0) (z) {};
\draw[color=red!50,fill=red!50] (0,4)--(1,4)--(1,3)--(0,3)--(0,4);
\draw[thick,rounded corners,color=blue] (2.5,0)--(2.5,3.5);
\draw[thick,rounded corners,color=blue] (0.5,3)--(0.5,3.5)--(1,3.5);
\draw[thick,rounded corners,color=blue] (0.5,0)--(0.5,1);
\draw[thick,rounded corners,color=blue] (0.5,1)--(0.5,2);
\draw[thick,rounded corners,color=blue] (0.5,2)--(0.5,3);
\draw[thick,rounded corners,color=blue] (1,3.5)--(2,3.5);
\draw[thick,rounded corners,color=blue] (2,3.5)--(3,3.5);
\draw[thick,rounded corners,color=blue] (3,3.5)--(4,3.5);
\draw[gray, dotted, line width=1pt] (0,4)--(4,4)--(4,0)--(0,0)--(0,4);
\draw[white, line width=0.7pt,fill=white] (2,4)--(4,2)--(4,4)--(2,4);
\draw[black, line width=1.6pt] (2,4)--(4,2);
\end{tikzpicture}
\
\mapsto
\
\begin{tikzpicture}[x=\bpdsmallwidth,y=\bpdsmallwidth,line cap=round,line join=round,baseline=(z.base)]
\node at (0,2.0) (z) {};
\draw[color=red!50,fill=red!50] (3,0)--(4,0)--(4,1)--(3,1)--(3,0);
\draw[thick,rounded corners,color=blue] (2.5,0)--(2.5,3.5);
\draw[thick,rounded corners,color=blue] (0.5,0)--(0.5,0.5)--(1,0.5);
\draw[thick,rounded corners,color=blue] (1,0.5)--(2,0.5);
\draw[thick,rounded corners,color=blue] (2,0.5)--(3,0.5);
\draw[thick,rounded corners,color=blue] (3.5,1)--(3.5,0.5)--(3,0.5);
\draw[thick,rounded corners,color=blue] (3.5,1)--(3.5,2);
\draw[thick,rounded corners,color=blue] (3.5,2)--(3.5,3);
\draw[thick,rounded corners,color=blue] (3.5,3)--(3.5,3.5)--(4,3.5);
\draw[gray, dotted, line width=1pt] (0,4)--(4,4)--(4,0)--(0,0)--(0,4);
\draw[white, line width=0.7pt,fill=white] (2,4)--(4,2)--(4,4)--(2,4);
\draw[black, line width=1.6pt] (2,4)--(4,2);
\end{tikzpicture}
\sdgap
 b < a < j < i
\earr
\qquad
\barr{c}
\begin{tikzpicture}[x=\bpdsmallwidth,y=\bpdsmallwidth,line cap=round,line join=round,baseline=(z.base)]
\node at (0,1.5) (z) {};
\draw[color=red!50,fill=red!50] (0,2)--(0,3)--(1,3)--(1,2)--(0,2);
\draw[thick,rounded corners,color=blue] (2.5,0)--(2.5,2.5);
\draw[thick,rounded corners,color=blue] (0.5,2)--(0.5,2.5)--(1,2.5);
\draw[thick,rounded corners,color=blue] (0.5,0)--(0.5,1);
\draw[thick,rounded corners,color=blue] (0.5,1)--(0.5,2);
\draw[thick,rounded corners,color=blue] (1,2.5)--(2,2.5);
\draw[thick,rounded corners,color=blue] (2,2.5)--(3,2.5);
\draw[thick,rounded corners,color=blue] (3,2.5)--(4,2.5);
\draw[gray, dotted, line width=1pt] (0,3)--(5,3)--(5,0)--(0,0)--(0,3);
\draw[white, line width=0.7pt,fill=white] (2,3)--(5,0)--(5,3)--(2,3);
\draw[black, line width=1.6pt] (2,3)--(5,0);
\end{tikzpicture}
\
\mapsto
\
\begin{tikzpicture}[x=\bpdsmallwidth,y=\bpdsmallwidth,line cap=round,line join=round,baseline=(z.base)]
\node at (0,1.5) (z) {};
\draw[color=red!50,fill=red!50] (4,0)--(4,1)--(5,1)--(5,0)--(4,0);
\draw[thick,rounded corners,color=blue] (2.5,0)--(2.5,2.5);
\draw[thick,rounded corners,color=blue] (0.5,0)--(0.5,0.5)--(1,0.5);
\draw[thick,rounded corners,color=blue] (1,0.5)--(2,0.5);
\draw[thick,rounded corners,color=blue] (2,0.5)--(3,0.5);
\draw[thick,rounded corners,color=blue] (5,0.5)--(3,0.5);
\draw[gray, dotted, line width=1pt] (0,3)--(5,3)--(5,0)--(0,0)--(0,3);
\draw[white, line width=0.7pt,fill=white] (2,3)--(5,0)--(5,3)--(2,3);
\draw[black, line width=1.6pt] (2,3)--(5,0);
\end{tikzpicture}
\sdgap
 b < a < j = i
\earr.
\ee
The highlighted squares show $(a,b)$ on the left and $(i,j)$ on the right. 
The blank space may contain other pipes in $D$ 
as long as these do not involve any other bend tiles. 
Such pipes must therefore consist of vertical or horizontal lines that meet at the diagonal in the dotted region.
To be precise:

\begin{definition}[\defn{Simple droop moves}]
\label{droop-def}
Suppose the following properties hold:
\ben
\item[(1)] $D_{ab} = \CTILE$ (when $a>b$) or $D_{ab} =\CDIAG$ (when $a=b$);

\item[(2)] $D_{ij} = \BTILE$ (when $i>j$) or $D_{ij}=\BDIAG$ (when $i=j$); and

\item[(3)] besides $(a,b)$, 
no other position of $D$ in the rectangle $[a,i]\times[b,j]$
contains a bend tile. 

\een
Then we form $\droop^{ab}_{ij}(D)$ as the triangular pipe dream
whose bend tiles are obtained from    $D$ by   
\bei
\item[(a)] removing the bend tile in position $(a,b)$;

\item[(b)] adding $\JTILE$ (when $j<i$) or $\JDIAG$ (when $i=j$) to position $(i,j)$;

\item[(c)] adding $\CTILE$ to position $(i,b)$ and to position $(a,j)$ when $j<a$; and

\item[(d)] adding $\CDIAG$ to position $(a,a)$ when $b<a<j$ and  to position $(j,j)$ when $a<j<i$.

\eei
\end{definition}

When $j>a>b$ we are in the situation
of \eqref{simple-eq2}, where $D$ must have a vertical pipe
to prevent a $\JDIAG$ bend tile from occurring at the diagonal.
Then we must have $D_{aa} = \PDIAG$ and this changes to $\CDIAG$ in $\droop^{ab}_{ij}(D)$.

\begin{remark}
The pictures in \eqref{simple-eq} and \eqref{simple-eq2} account for 5 of the 6 possible
 configurations of $a\geq b$ and $i\geq j$ satisfying $a<i$ and $b<j$.
 The remaining possibility corresponds to the picture
\be\label{simple-eq3}
\barr{c}
\begin{tikzpicture}[x=\bpdsmallwidth,y=\bpdsmallwidth,line cap=round,line join=round,baseline=(z.base)]
\node at (0,2.0) (z) {};
\draw[color=red!50,fill=red!50] (0,4)--(1,4)--(1,3)--(0,3)--(0,4);
\draw[thick,rounded corners,color=blue] (0.5,3)--(0.5,3.5)--(1,3.5);
\draw[thick,rounded corners,color=blue] (0.5,0)--(0.5,1);
\draw[thick,rounded corners,color=blue] (0.5,1)--(0.5,2);
\draw[thick,rounded corners,color=blue] (0.5,2)--(0.5,3);
\draw[thick,rounded corners,color=blue] (1,3.5)--(2,3.5);
\draw[thick,rounded corners,color=blue] (2,3.5)--(3,3.5);
\draw[thick,rounded corners,color=blue] (3,3.5)--(4,3.5);
\draw[gray, dotted, line width=1pt] (0,4)--(4,4)--(4,0)--(0,0)--(0,4);
\draw[white, line width=0.7pt,fill=white] (3,4)--(4,3)--(4,4)--(3,4);
\draw[black, line width=1.6pt] (3,4)--(4,3);
\end{tikzpicture}
\
\mapsto
\
\begin{tikzpicture}[x=\bpdsmallwidth,y=\bpdsmallwidth,line cap=round,line join=round,baseline=(z.base)]
\node at (0,2.0) (z) {};
\draw[color=red!50,fill=red!50] (3,0)--(4,0)--(4,1)--(3,1)--(3,0);
\draw[thick,rounded corners,color=blue] (0.5,0)--(0.5,0.5)--(1,0.5);
\draw[thick,rounded corners,color=blue] (1,0.5)--(2,0.5);
\draw[thick,rounded corners,color=blue] (2,0.5)--(3,0.5);
\draw[thick,rounded corners,color=blue] (3.5,1)--(3.5,0.5)--(3,0.5);
\draw[thick,rounded corners,color=blue] (3.5,1)--(3.5,2);
\draw[thick,rounded corners,color=blue] (3.5,2)--(3.5,3);
\draw[thick,rounded corners,color=blue] (3.5,3)--(3.5,4);
\draw[gray, dotted, line width=1pt] (0,4)--(4,4)--(4,0)--(0,0)--(0,4);
\draw[white, line width=0.7pt,fill=white] (3,4)--(4,3)--(4,4)--(3,4);
\draw[black, line width=1.6pt] (3,4)--(4,3);
\end{tikzpicture}
\sdgap
 b <a= j < i
\earr.
\ee
This is \textbf{not} a simple droop move, but is an instance of the \defn{strong droop moves} that are defined below.
\end{remark}

Simple droop moves are straightforward analogues
of droop moves on bumpless pipe dreams \cite[\S4.3]{Weigandt}.
To generate all bumpless involution pipe dreams,
we require one new type of move without an obvious parallel in the classical setting.
This \defn{complex droop move}  $\complexDroop^{ab}_{ij}$ 
is indexed by positions $(a,b),(i,j) \in \LT_n$
with $b<a<j<i$.
Like simple droop moves, it transforms a triangular pipe dream by modifying a ``single'' pipe (technically, two pipes), but now as in the following picture:
\be\label{complex-eq}
\begin{tikzpicture}[x=\bpdwidth,y=\bpdwidth,line cap=round,line join=round,baseline=(z.base)]
\node at (0,3.0) (z) {};
\draw[color=red!50,fill=red!50] (0,5)--(0,6)--(1,6)--(1,5)--(0,5);
\draw[dotted,gray, line width=1pt](0,0)--(6,0)--(6,2);
\draw[dotted,gray, line width=1pt](0,0)--(0,6)--(2,6);
\draw[thick,rounded corners,color=blue] (0.5,0)--(0.5,1);
\draw[thick,rounded corners,color=blue] (0.5,1)--(0.5,2);
\draw[thick,rounded corners,color=blue] (0.5,2)--(0.5,3);
\draw[thick,rounded corners,color=blue] (0.5,3)--(0.5,4);
\draw[thick,rounded corners,color=blue] (0.5,4)--(0.5,5);
\draw[thick,rounded corners,color=blue] (0.5,5)--(0.5,5.5)--(1,5.5);
\draw[thick,rounded corners,color=blue] (1,5.5)--(2,5.5);
\draw[thick,rounded corners,color=blue] (2,5.5)--(3,5.5);
\draw[thick,rounded corners,color=blue] (2.5,5)--(2.5,6);
\draw[thick,rounded corners,color=blue] (2.5,4)--(2.5,5);
\draw[thick,rounded corners,color=blue] (2.5,3)--(2.5,4);
\draw[thick,rounded corners,color=blue] (2.5,2)--(2.5,3);
\draw[thick,rounded corners,color=blue] (2.5,1)--(2.5,2);
\draw[thick,rounded corners,color=blue] (2.5,0)--(2.5,1);
\draw[gray!50, line width=0.7pt,fill=gray!50] (0,0)--(2,0)--(2,2)--(0,2);
\draw[gray!50, line width=0.7pt,fill=gray!50] (0,0)--(2,0)--(2,2)--(0,2);
\draw[white, line width=0.7pt,fill=white] (2,6)--(6,2)--(6,6)--(2,6);
\draw[black, line width=1.6pt] (2,6)--(6,2);
\end{tikzpicture}
\
\mapsto
\
\begin{tikzpicture}[x=\bpdwidth,y=\bpdwidth,line cap=round,line join=round,baseline=(z.base)]
\node at (0,3.0) (z) {};
\draw[color=red!50,fill=red!50] (6,0)--(6,1)--(5,1)--(5,0)--(6,0);
\draw[dotted,gray, line width=1pt](0,0)--(6,0)--(6,2);
\draw[dotted,gray, line width=1pt](0,0)--(0,6)--(2,6);
\draw[thick,rounded corners,color=blue] (0.5,0)--(0.5,1);
\draw[thick,rounded corners,color=blue] (0.5,1)--(0.5,2);
\draw[thick,rounded corners,color=blue] (0.5,2)--(0.5,2.5)--(1,2.5);
\draw[thick,rounded corners,color=blue] (2,2.5)--(3,2.5);
\draw[thick,rounded corners,color=blue] (2.5,0)--(2.5,0.5)--(3,0.5);
\draw[thick,rounded corners,color=blue] (3,0.5)--(4,0.5);
\draw[thick,rounded corners,color=blue] (4,0.5)--(5,0.5);
\draw[thick,rounded corners,color=blue] (5.5,1)--(5.5,0.5)--(5,0.5);
\draw[thick,rounded corners,color=blue] (5.5,1)--(5.5,2);
\draw[thick,rounded corners,color=blue] (5,2.5)--(6,2.5);
\draw[thick,rounded corners,color=blue] (5.5,2)--(5.5,3);
\draw[thick,rounded corners,color=blue] (4,2.5)--(5,2.5);
\draw[thick,rounded corners,color=blue] (3,2.5)--(4,2.5);
\draw[thick,rounded corners,color=blue] (1,2.5)--(2,2.5);
\draw[gray!50, line width=0.7pt,fill=gray!50] (0,0)--(2,0)--(2,2)--(0,2);
\draw[white, line width=0.7pt,fill=white] (2,6)--(6,2)--(6,6)--(2,6);
\draw[black, line width=1.6pt] (2,6)--(6,2);
\end{tikzpicture}.
\ee
Here again, the highlighted square on the left is  $(a,b)$
while the highlighted square on the right is  $(i,j)$.
The grey space is unconstrained, but 
the white blank space  below the diagonal must not contain any other bend tiles. Any pipes in this region must consist of vertical or horizontal lines,
which may meet at the diagonal, and these pipes are unchanged by the droop operation.
Here is a detailed definition:

\begin{definition}[\defn{Complex droop moves}]
\label{complex-droop-def}
Suppose $1\leq b<a<j<i\leq n$. 
Assume that 
\[D_{ab} = \CTILE,\qquad D_{aa} =\PDIAG,
\qquad D_{ij}=\BTILE,
\]
and  
no other lower triangular position of $D$ in $[a,j]\times[b,i]$
or $[b,i]\times[a,j]$
contains a bend tile.
Then we form
 $\complexDroop^{ab}_{ij}(D)$ as the triangular pipe dream
whose bend tiles are obtained from  $D$ by  
removing $\CTILE$ from position $(a,b)$,
adding $\JTILE$  to position $(i,j)$, and
adding $\CTILE$ to positions $(j,b)$ and $(i,a)$.

\end{definition}

\begin{example}\label{tri-droop-ex}
It may be possible to apply both $\droop^{ab}_{ij}$ and $\complexDroop^{ab}_{ij}$ 
 to the 
same triangular pipe dream (with different results), as in the following examples
for $(a,b) = (2,1)$ and $(i,j) = (4,3)$:
\[
\begin{tikzpicture}[x=\bpdwidth,y=\bpdwidth,line cap=round,line join=round,baseline=(z.base)]
\node at (0,2.25) (z) {};
\draw[color=red!50,fill=red!50] (0,3)--(0,4)--(1,4)--(1,3);
\draw[gray!70, line width=0.7pt](0,0)--(5,0);
\draw[gray!70, line width=0.7pt](0,0)--(0,5);
\draw[gray!70, line width=0.7pt](0,1)--(5,1);
\draw[gray!70, line width=0.7pt](1,0)--(1,5);
\draw[gray!70, line width=0.7pt](0,2)--(5,2);
\draw[gray!70, line width=0.7pt](2,0)--(2,5);
\draw[gray!70, line width=0.7pt](0,3)--(5,3);
\draw[gray!70, line width=0.7pt](3,0)--(3,5);
\draw[gray!70, line width=0.7pt](0,4)--(5,4);
\draw[gray!70, line width=0.7pt](4,0)--(4,5);
\draw[gray!70, line width=0.7pt](0,5)--(5,5);
\draw[gray!70, line width=0.7pt](5,0)--(5,5);
\draw[thick,rounded corners,color=blue] (0.5,3)--(0.5,3.5)--(1,3.5);
\draw[thick,rounded corners,color=blue] (3.5,1)--(3.5,2);
\draw[thick,rounded corners,color=blue] (2.5,0)--(2.5,0.5)--(3,0.5);
\draw[thick,rounded corners,color=blue] (0.5,0)--(0.5,1);
\draw[thick,rounded corners,color=blue] (0.5,1)--(0.5,2);
\draw[thick,rounded corners,color=blue] (0.5,2)--(0.5,3);
\draw[thick,rounded corners,color=blue] (1,3.5)--(2,3.5);
\draw[thick,rounded corners,color=blue] (1.5,3)--(1.5,4);
\draw[thick,rounded corners,color=blue] (1.5,0)--(1.5,1);
\draw[thick,rounded corners,color=blue] (1.5,1)--(1.5,2);
\draw[thick,rounded corners,color=blue] (1.5,2)--(1.5,3);
\draw[thick,rounded corners,color=blue] (3,0.5)--(4,0.5);
\draw[thick,rounded corners,color=blue] (3.5,0)--(3.5,1);
\draw[thick,rounded corners,color=blue] (4,0.5)--(5,0.5);
\draw[thick,rounded corners,color=blue] (4.5,0)--(4.5,1);
\draw[white, line width=0.7pt,fill=white] (0,5)--(5,0)--(5,5)--(0,5);
\draw[black, line width=1.6pt] (0,5)--(5,0)--(0,0)--(0,5);
\end{tikzpicture}
\ \mapsto\
\begin{tikzpicture}[x=\bpdwidth,y=\bpdwidth,line cap=round,line join=round,baseline=(z.base)]
\node at (0,2.25) (z) {};
\draw[color=red!50,fill=red!50] (2,1)--(2,2)--(3,2)--(3,1);
\draw[gray!70, line width=0.7pt](0,0)--(5,0);
\draw[gray!70, line width=0.7pt](0,0)--(0,5);
\draw[gray!70, line width=0.7pt](0,1)--(5,1);
\draw[gray!70, line width=0.7pt](1,0)--(1,5);
\draw[gray!70, line width=0.7pt](0,2)--(5,2);
\draw[gray!70, line width=0.7pt](2,0)--(2,5);
\draw[gray!70, line width=0.7pt](0,3)--(5,3);
\draw[gray!70, line width=0.7pt](3,0)--(3,5);
\draw[gray!70, line width=0.7pt](0,4)--(5,4);
\draw[gray!70, line width=0.7pt](4,0)--(4,5);
\draw[gray!70, line width=0.7pt](0,5)--(5,5);
\draw[gray!70, line width=0.7pt](5,0)--(5,5);
\draw[thick,rounded corners,color=blue] (3.5,1)--(3.5,2);
\draw[thick,rounded corners,color=blue] (2.5,0)--(2.5,0.5)--(3,0.5);
\draw[thick,rounded corners,color=blue] (1.5,3)--(1.5,4);
\draw[thick,rounded corners,color=blue] (2.5,2)--(2.5,1.5)--(2,1.5);
\draw[thick,rounded corners,color=blue] (0.5,1)--(0.5,1.5)--(1,1.5);
\draw[thick,rounded corners,color=blue] (2.5,2)--(2.5,3);
\draw[thick,rounded corners,color=blue] (0.5,0)--(0.5,1);
\draw[thick,rounded corners,color=blue] (1,1.5)--(2,1.5);
\draw[thick,rounded corners,color=blue] (1.5,1)--(1.5,2);
\draw[thick,rounded corners,color=blue] (1.5,0)--(1.5,1);
\draw[thick,rounded corners,color=blue] (1.5,2)--(1.5,3);
\draw[thick,rounded corners,color=blue] (3,0.5)--(4,0.5);
\draw[thick,rounded corners,color=blue] (3.5,0)--(3.5,1);
\draw[thick,rounded corners,color=blue] (4,0.5)--(5,0.5);
\draw[thick,rounded corners,color=blue] (4.5,0)--(4.5,1);
\draw[white, line width=0.7pt,fill=white] (0,5)--(5,0)--(5,5)--(0,5);
\draw[black, line width=1.6pt] (0,5)--(5,0)--(0,0)--(0,5);
\end{tikzpicture}
\qquand
\begin{tikzpicture}[x=\bpdwidth,y=\bpdwidth,line cap=round,line join=round,baseline=(z.base)]
\node at (0,2.25) (z) {};
\draw[color=red!50,fill=red!50] (0,3)--(0,4)--(1,4)--(1,3);
\draw[gray!70, line width=0.7pt](0,0)--(5,0);
\draw[gray!70, line width=0.7pt](0,0)--(0,5);
\draw[gray!70, line width=0.7pt](0,1)--(5,1);
\draw[gray!70, line width=0.7pt](1,0)--(1,5);
\draw[gray!70, line width=0.7pt](0,2)--(5,2);
\draw[gray!70, line width=0.7pt](2,0)--(2,5);
\draw[gray!70, line width=0.7pt](0,3)--(5,3);
\draw[gray!70, line width=0.7pt](3,0)--(3,5);
\draw[gray!70, line width=0.7pt](0,4)--(5,4);
\draw[gray!70, line width=0.7pt](4,0)--(4,5);
\draw[gray!70, line width=0.7pt](0,5)--(5,5);
\draw[gray!70, line width=0.7pt](5,0)--(5,5);
\draw[thick,rounded corners,color=blue] (0.5,3)--(0.5,3.5)--(1,3.5);
\draw[thick,rounded corners,color=blue] (3.5,1)--(3.5,2);
\draw[thick,rounded corners,color=blue] (2.5,0)--(2.5,0.5)--(3,0.5);
\draw[thick,rounded corners,color=blue] (0.5,0)--(0.5,1);
\draw[thick,rounded corners,color=blue] (0.5,1)--(0.5,2);
\draw[thick,rounded corners,color=blue] (0.5,2)--(0.5,3);
\draw[thick,rounded corners,color=blue] (1,3.5)--(2,3.5);
\draw[thick,rounded corners,color=blue] (1.5,3)--(1.5,4);
\draw[thick,rounded corners,color=blue] (1.5,0)--(1.5,1);
\draw[thick,rounded corners,color=blue] (1.5,1)--(1.5,2);
\draw[thick,rounded corners,color=blue] (1.5,2)--(1.5,3);
\draw[thick,rounded corners,color=blue] (3,0.5)--(4,0.5);
\draw[thick,rounded corners,color=blue] (3.5,0)--(3.5,1);
\draw[thick,rounded corners,color=blue] (4,0.5)--(5,0.5);
\draw[thick,rounded corners,color=blue] (4.5,0)--(4.5,1);
\draw[white, line width=0.7pt,fill=white] (0,5)--(5,0)--(5,5)--(0,5);
\draw[black, line width=1.6pt] (0,5)--(5,0)--(0,0)--(0,5);
\end{tikzpicture}
\ \mapsto\ 
\begin{tikzpicture}[x=\bpdwidth,y=\bpdwidth,line cap=round,line join=round,baseline=(z.base)]
\node at (0,2.25) (z) {};
\draw[color=red!50,fill=red!50] (2,1)--(2,2)--(3,2)--(3,1);
\draw[gray!70, line width=0.7pt](0,0)--(5,0);
\draw[gray!70, line width=0.7pt](0,0)--(0,5);
\draw[gray!70, line width=0.7pt](0,1)--(5,1);
\draw[gray!70, line width=0.7pt](1,0)--(1,5);
\draw[gray!70, line width=0.7pt](0,2)--(5,2);
\draw[gray!70, line width=0.7pt](2,0)--(2,5);
\draw[gray!70, line width=0.7pt](0,3)--(5,3);
\draw[gray!70, line width=0.7pt](3,0)--(3,5);
\draw[gray!70, line width=0.7pt](0,4)--(5,4);
\draw[gray!70, line width=0.7pt](4,0)--(4,5);
\draw[gray!70, line width=0.7pt](0,5)--(5,5);
\draw[gray!70, line width=0.7pt](5,0)--(5,5);
\draw[thick,rounded corners,color=blue] (3.5,1)--(3.5,2);
\draw[thick,rounded corners,color=blue] (2.5,0)--(2.5,0.5)--(3,0.5);
\draw[thick,rounded corners,color=blue] (0.5,2)--(0.5,2.5)--(1,2.5);
\draw[thick,rounded corners,color=blue] (2.5,2)--(2.5,1.5)--(2,1.5);
\draw[thick,rounded corners,color=blue] (1.5,1)--(1.5,1.5)--(2,1.5);
\draw[thick,rounded corners,color=blue] (0.5,0)--(0.5,1);
\draw[thick,rounded corners,color=blue] (0.5,1)--(0.5,2);
\draw[thick,rounded corners,color=blue] (1,2.5)--(2,2.5);
\draw[thick,rounded corners,color=blue] (2,2.5)--(3,2.5);
\draw[thick,rounded corners,color=blue] (2.5,2)--(2.5,3);
\draw[thick,rounded corners,color=blue] (1.5,0)--(1.5,1);
\draw[thick,rounded corners,color=blue] (3,0.5)--(4,0.5);
\draw[thick,rounded corners,color=blue] (3.5,0)--(3.5,1);
\draw[thick,rounded corners,color=blue] (4,0.5)--(5,0.5);
\draw[thick,rounded corners,color=blue] (4.5,0)--(4.5,1);
\draw[white, line width=0.7pt,fill=white] (0,5)--(5,0)--(5,5)--(0,5);
\draw[black, line width=1.6pt] (0,5)--(5,0)--(0,0)--(0,5);
\end{tikzpicture}.
\]
\end{example}

We set $ \droop^{ab}_{ij}(D)=D$ and  $\complexDroop^{ab}_{ij}(D)=D$
in any cases when the preceding definitions do not apply (due to the presence of 
extra bend tiles or when $a,b,i,j \in [n]$ are not in the required order).
When $D$ and $E$ are triangular pipe dreams and
$E=\droop^{ab}_{ij}(D)$ (respectively, $E=\complexDroop^{ab}_{ij}(D)$),
we say that $E$ is obtained from $D$ by applying a \defn{simple droop move}
(respectively, a \defn{complex droop move}).

As droop moves can change where pipes meet the diagonal, it is not obvious that 
$D$, $ \droop^{ab}_{ij}(D)$, and  $\complexDroop^{ab}_{ij}(D)$ are triangular pipe dreams for the same involution.
However,
 this turns out to be so.
 
\begin{lemma}\label{droop-lem1}
Suppose $D$ is a triangular pipe dream and $E$ is obtained 
from $D$ by applying any (simple or complex)  droop move.
Then $D\in\InvBPD(z)$
for a given $z \in \Invol $ 
  if and only if $E\in\InvBPD(z)$.
\end{lemma}

We omit the proof of this lemma, since it will be subsumed
by a stronger statement (Proposition~\ref{droops-preserve-prop}),
which is given a self-contained proof in the next section.
The following is our  main result here.

\begin{theorem}\label{droop-thm}
Let $z \in \Invol $. Then $\InvBPD(z)$ is exactly 
the set of
triangular pipe dreams  obtained by applying sequences of 
(simple or complex) droop moves to $\invRothe(z)$.
\end{theorem}

\begin{proof}
Our argument is similar to, but slightly more complicated than the proof of \cite[Prop.~5.3]{LLS}.
Given $E \in \InvBPD(z)$, consider the set of positions $\upturn(E)$ from \eqref{upturn-eq}.
If this set is empty then $E =\invRothe(z)$.
Otherwise, let $(i,j) $ be a position in $\upturn(E)$ that is maximally northwest.

First assume that $E_{jj} \neq \PDIAG$. 
Then the pipe of $E$ in position $(i,j)$ must pass through a $\CTILE$ tile in some position $(i,b)$
with $b<j$. 
Let $a$ be the largest integer with $b<a<j$ and $E_{aa} = \CDIAG$,
or else set $a=b$ when no such integer exists.
The northwest maximality of $(i,j)$ implies that $E_{ab} $ is either $ \BTILE$ or $\BDIAG$, and by inspecting \eqref{simple-eq} and \eqref{simple-eq2} one sees that there exists a triangular pipe dream $D$ with $E = \droop^{ab}_{ij}(D)$.

If instead $E_{jj} = \PDIAG$, then $i>j$ and the two pipes that meet in $E_{jj}$
must pass through $\CTILE$ tiles in some positions $(j,b)$ and $(i,a)$
 with $b<a<j<i$.
 The northwest maximality of $(i,j)$ again implies that $E_{ab} = \BTILE$
 and by inspecting \eqref{complex-eq}  one checks that there exists a triangular pipe dream $D$ with $E = \complexDroop^{ab}_{ij}(D)$.
 
 In either case, we have $D \in \InvBPD(z)$ by Lemma~\ref{droop-lem1},
 and this pipe dream has $|\upturn(D)| < |\upturn(E)|$.
 The theorem thus follows by induction on the number of upturn tiles in
$E$.
\end{proof}

\subsection{Strong droops and undroops}

 It will be convenient to have some more general droop
moves that preserve the set of triangular pipe dreams 
associated to an involution.
These   moves  will include pictures like \eqref{simple-eq3} and
\[
\begin{tikzpicture}[x=\bpdwidth,y=\bpdwidth,line cap=round,line join=round,baseline=(z.base)]
\node at (0,2.0) (z) {};
\draw[gray!70, line width=0.7pt](0,0)--(4,0);
\draw[gray!70, line width=0.7pt](0,0)--(0,4);
\draw[color=red!50,fill=red!50] (0,4)--(1,4)--(1,3)--(0,3)--(0,4);
\draw[thick,rounded corners,color=blue] (0,3.5)--(0.5,3.5)--(0.5,4);
\draw[thick,rounded corners,color=blue] (0.5,3)--(0.5,3.5)--(1,3.5);
\draw[thick,rounded corners,color=blue] (0,0.5)--(0.5,0.5)--(0.5,1);
\draw[thick,rounded corners,color=blue] (0.5,1)--(0.5,2);
\draw[thick,rounded corners,color=blue] (0.5,2)--(0.5,3);
\draw[thick,rounded corners,color=blue] (1,3.5)--(2,3.5);
\draw[thick,rounded corners,color=blue] (2,3.5)--(3,3.5);
\draw[thick,rounded corners,color=blue] (3,3.5)--(3.5,3.5)--(3.5,4);
\draw[thick,rounded corners,color=blue] (4,0.5)--(3.5,0.5)--(3.5,0);
\draw[gray, dotted, line width=1pt] (0,4)--(4,4)--(4,0)--(0,0)--(0,4);
\end{tikzpicture}
\
\mapsto
\
\begin{tikzpicture}[x=\bpdwidth,y=\bpdwidth,line cap=round,line join=round,baseline=(z.base)]
\node at (0,2.0) (z) {};
\draw[gray!70, line width=0.7pt](0,0)--(4,0);
\draw[gray!70, line width=0.7pt](0,0)--(0,4);
\draw[color=red!50,fill=red!50] (3,0)--(4,0)--(4,1)--(3,1)--(3,0);
\draw[thick,rounded corners,color=blue] (0,3.5)--(0.5,3.5)--(0.5,4);
\draw[thick,rounded corners,color=blue] (0,0.5)--(1,0.5);
\draw[thick,rounded corners,color=blue] (1,0.5)--(2,0.5);
\draw[thick,rounded corners,color=blue] (2,0.5)--(3,0.5);
\draw[thick,rounded corners,color=blue] (3.5,1)--(3.5,0.5)--(3,0.5);
\draw[thick,rounded corners,color=blue] (4,0.5)--(3.5,0.5)--(3.5,0);
\draw[thick,rounded corners,color=blue] (3.5,1)--(3.5,2);
\draw[thick,rounded corners,color=blue] (3.5,2)--(3.5,3);
\draw[thick,rounded corners,color=blue] (3.5,3)--(3.5,4);
\draw[gray, dotted, line width=1pt] (0,4)--(4,4)--(4,0)--(0,0)--(0,4);
\end{tikzpicture}
\qquand
\begin{tikzpicture}[x=\bpdwidth,y=\bpdwidth,line cap=round,line join=round,baseline=(z.base)]
\node at (0,3.0) (z) {};
\draw[color=red!50,fill=red!50] (0,5)--(0,6)--(1,6)--(1,5)--(0,5);
\draw[dotted,gray, line width=1pt](0,0)--(6,0)--(6,2);
\draw[dotted,gray, line width=1pt](0,0)--(0,6)--(2,6);
\draw[thick,rounded corners,color=blue] (0.5,0)--(0.5,1);
\draw[thick,rounded corners,color=blue] (0.5,1)--(0.5,2);
\draw[thick,rounded corners,color=blue] (0,2.5)--(0.5,2.5)--(0.5,3);
\draw[thick,rounded corners,color=blue] (0.5,3)--(0.5,4);
\draw[thick,rounded corners,color=blue] (0.5,4)--(0.5,5);
\draw[thick,rounded corners,color=blue] (0,5.5)--(0.5,5.5)--(0.5,6);
\draw[thick,rounded corners,color=blue] (0.5,5)--(0.5,5.5)--(1,5.5);
\draw[thick,rounded corners,color=blue] (1,5.5)--(2,5.5);
\draw[thick,rounded corners,color=blue] (2,5.5)--(3,5.5);
\draw[thick,rounded corners,color=blue] (2.5,5)--(2.5,6);
\draw[thick,rounded corners,color=blue] (2.5,4)--(2.5,5);
\draw[thick,rounded corners,color=blue] (2.5,3)--(2.5,4);
\draw[thick,rounded corners,color=blue] (2.5,2)--(2.5,3);
\draw[thick,rounded corners,color=blue] (2.5,1)--(2.5,2);
\draw[thick,rounded corners,color=blue] (2,0.5)--(2.5,0.5)--(2.5,1);
\draw[thick,rounded corners,color=blue] (5.5,0)--(5.5,0.5)--(6,0.5);
\draw[gray!50, line width=0.7pt,fill=gray!50] (0,0)--(2,0)--(2,2)--(0,2);
\draw[gray!50, line width=0.7pt,fill=gray!50] (0,0)--(2,0)--(2,2)--(0,2);
\draw[white, line width=0.7pt,fill=white] (2,6)--(6,2)--(6,6)--(2,6);
\draw[black, line width=1.6pt] (2,6)--(6,2);
\end{tikzpicture}
\
\mapsto
\
\begin{tikzpicture}[x=\bpdwidth,y=\bpdwidth,line cap=round,line join=round,baseline=(z.base)]
\node at (0,3.0) (z) {};
\draw[color=red!50,fill=red!50] (6,0)--(6,1)--(5,1)--(5,0)--(6,0);
\draw[dotted,gray, line width=1pt](0,0)--(6,0)--(6,2);
\draw[dotted,gray, line width=1pt](0,0)--(0,6)--(2,6);
\draw[thick,rounded corners,color=blue] (0,5.5)--(0.5,5.5)--(0.5,6);
\draw[thick,rounded corners,color=blue] (0.5,0)--(0.5,1);
\draw[thick,rounded corners,color=blue] (0.5,1)--(0.5,2);
\draw[thick,rounded corners,color=blue] (0,2.5)--(1,2.5);
\draw[thick,rounded corners,color=blue] (2,2.5)--(3,2.5);
\draw[thick,rounded corners,color=blue] (2,0.5)--(3,0.5);
\draw[thick,rounded corners,color=blue] (3,0.5)--(4,0.5);
\draw[thick,rounded corners,color=blue] (4,0.5)--(5,0.5);
\draw[thick,rounded corners,color=blue] (5.5,1)--(5.5,0.5)--(5,0.5);
\draw[thick,rounded corners,color=blue] (5.5,0)--(5.5,0.5)--(6,0.5);
\draw[thick,rounded corners,color=blue] (5.5,1)--(5.5,2);
\draw[thick,rounded corners,color=blue] (5,2.5)--(6,2.5);
\draw[thick,rounded corners,color=blue] (5.5,2)--(5.5,3);
\draw[thick,rounded corners,color=blue] (4,2.5)--(5,2.5);
\draw[thick,rounded corners,color=blue] (3,2.5)--(4,2.5);
\draw[thick,rounded corners,color=blue] (1,2.5)--(2,2.5);
\draw[gray!50, line width=0.7pt,fill=gray!50] (0,0)--(2,0)--(2,2)--(0,2);
\draw[white, line width=0.7pt,fill=white] (2,6)--(6,2)--(6,6)--(2,6);
\draw[black, line width=1.6pt] (2,6)--(6,2);
\end{tikzpicture}.
\]
To give a precise definition, 
continue to
fix a triangular pipe dream $D$
and lower triangular positions
$(a,b),(i,j) \in \LT_n$ with $a< i$ and $b < j$.
The \defn{interior}
of the rectangle $[a,i]\times[b,j]$
refers to all positions in the rectangle 
distinct from the corners $(a,b)$, $(i,j)$, $(i,b)$, and $(a,j)$.

\begin{definition}[\defn{Strong droop moves}]
\label{extended-droop-def}
Let $D$ be a triangular pipe dream.
\ben
\item[(a)] Assume
that $D$ has no bend tiles in the interior of the
rectangle $[a,i]\times[b,j]$.
Suppose the
 tiles of $D$
 also match the inputs to the transformations
\[
\ba
D_{ab} &= \begin{cases}
\CTILE \mapsto \BTILE \\[-10pt]\\
\ETILE \mapsto \JTILE  \\[-10pt]\\
\CDIAG  \mapsto \BDIAG  \\[-10pt]\\
\EDIAG \mapsto \JDIAG
\end{cases}
&
D_{ib} &=\begin{cases}
\VTILE \mapsto \CTILE \\[-10pt]\\
\JTILE  \mapsto \HTILE
\end{cases}
&
D_{ij} &= \begin{cases}
\BTILE \mapsto \JTILE \\[-10pt]\\
\CTILE  \mapsto \ETILE  \\[-10pt]\\
\BDIAG \mapsto \JDIAG \\[-10pt]\\
\CDIAG \mapsto \EDIAG 
\end{cases}
\ea
\]
along with
\[
D_{aj} = \begin{cases}
\HTILE \mapsto \CTILE \\[-10pt]\\
\JTILE  \mapsto \VTILE\\[-10pt]\\
\JDIAG \mapsto \CDIAG
\end{cases}
\quad\text{when }a\geq  j
\qquand
\ba 
D_{aa} &= 
\PDIAG \mapsto \CDIAG
&\text{when }b<a< j,
\\[-10pt]\\
D_{jj} &=
\BDIAG \mapsto  \CDIAG 
&\text{when }a<j<i.
\ea
\]
Then we form $\extDroop^{ab}_{ij}(D)$
as the triangular pipe dream whose bend tiles are the same as the ones 
that occur in $D$ after we apply the above transformations
to $D_{ab}$, $D_{ij}$, $D_{ib}$ and to
 $D_{aj}$, $D_{aa}$, $D_{jj}$ under the listed conditions.

\item[(b)] 
Assume instead
that $D$ has no bend tiles in the interior of the
rectangle $[a,j]\times[b,i]$
or in the interior of the transposed rectangle  $[b,i]\times [a,j]$.
If $b<a<j<i$ and $D_{aa}=\PDIAG$ and the
 tiles of $D$ 
 match the inputs to the transformations
\[
\ba
D_{ab} &= \begin{cases}
\CTILE \mapsto \BTILE \\[-10pt]\\
\ETILE \mapsto \JTILE  
\end{cases}
&
D_{jb} &=\begin{cases}
\VTILE \mapsto \CTILE \\[-10pt]\\
\JTILE  \mapsto \HTILE
\end{cases}
&
D_{ij} &= \begin{cases}
\BTILE \mapsto \JTILE \\[-10pt]\\
\CTILE  \mapsto \ETILE  
\end{cases}
&
D_{ia} &= \begin{cases}
\VTILE \mapsto \CTILE \\[-10pt]\\
\JTILE  \mapsto \HTILE,
\end{cases}
\ea
\]
then we form $\extComplexDroop^{ab}_{ij}(D)$
as the triangular pipe dream whose bend tiles are the same as the ones 
that occur in $D$ after we apply these transformations to $D_{ab}$, $D_{jb}$, $D_{ij}$, 
and $D_{ia}$.

\een
\end{definition}

We set $\extDroop^{ab}_{ij}(D)=D$
and  $\extComplexDroop^{ab}_{ij}(D)=D$
in any cases when the preceding definition does not apply.
If $D$ and $E$ are triangular pipe dreams
with $E =\extDroop^{ab}_{ij}(D)$ or
$E =  \extComplexDroop^{ab}_{ij}(D)$, then
we say that $E$ is obtained from $D$ by
a \defn{strong (simple or complex) droop move}.
Whenever $\droop^{ab}_{ij}(D) \neq D$ or  
 $\complexDroop^{ab}_{ij}(D) \neq D$
 we respectively have 
 \be\droop^{ab}_{ij}(D) = \extDroop^{ab}_{ij}(D)
\quand
\complexDroop^{ab}_{ij}(D) = \extComplexDroop^{ab}_{ij}(D).
 \ee

\begin{proposition}\label{droops-preserve-prop}
Suppose $D$ is a triangular pipe dream and $E$ is obtained 
from $D$ by applying a strong (simple or complex) droop move.
Then $D$ is in $\InvPD(z)$
for $z \in \Invol$ 
if and only if $E$ is in the same set.
\end{proposition}

\begin{proof}
 We may assume that $D\neq E$. 
 By Lemma~\ref{CAB-lem} and Proposition~\ref{word-pd-prop},
 it is enough to show that 
 $
\pi(D)\cWsim\pi(E)
$ and 
$|\cross(D)|=|\cross(E)|$ when $D$ or $E$ belongs to  $\InvPD(z)$.

First suppose $E = \extDroop^{ab}_{ij}(D)$. If $j\leq a$
then tracing the affected pipes gives
$
\pi(D)=\pi(E)
$ and 
$|\cross(D)|=|\cross(E)|$, and the conclusion of  the result follows.
When $a<j$, let $x$ be the label of the pipe of $D$ entering $(a,b)$ from the south. 
The diagonal connection of this pipe is transferred in $E$ from position 
$(a,a)$ to position $(j,j)$.
List the diagonal crossings between these positions as
$
(c_1,c_1),\ldots,(c_k,c_k)$ where $ a<c_1<\cdots<c_k<j$,
and let $q_t>p_t$ be their vertical and horizontal pipe labels, respectively. These crossings and their labels are the same in $D$ and $E$, and we can write
\be\label{strong-eq1}
\pi(D)=\cdots xq_1p_1q_2p_2\cdots q_kp_k\cdots 
\qquand
\pi(E)=\cdots q_1p_1q_2p_2\cdots q_kp_kx\cdots 
\ee
where the initial (respectively, final) ellipses mask identical subwords in each permutation.

When $E = \extComplexDroop^{ab}_{ij}(D)$, define $x$ as before and let $y>x$ be the vertical pipe label at $(a,a)$ in $D$. The complex droop move transfers the diagonal crossing of these two pipes from $(a,a)$ to $(j,j)$, so 
\be\label{strong-eq2}
\pi(D)=\cdots yxq_1p_1q_2p_2\cdots q_kp_k\cdots 
\qquand
\pi(E)=\cdots q_1p_1q_2p_2\cdots q_kp_kyx\cdots
\ee
with the same conventions for the initial and final ellipses.

In both cases, it is easy to see that 
$|\cross(D)|=|\cross(E)|$.
It remains to show that if $D\in\InvPD(z)$ or $E \in \InvPD(z)$,
so that one of $\pi(D)$ or $\pi(E)$ is in $\cW(z)$,
then the two sides of \eqref{strong-eq1} (in the simple case) or \eqref{strong-eq2} (in the complex case)
are related by $\cWsim$. 
When $D \in \InvPD(z)$ reducedness implies that $x>p_t$ for all $t$ and when $E \in \InvPD(z)$ we likewise have $q_t>x$ (in the simple case) or $q_t>y>x$ (in the complex case).
From these properties, proving the desired equivalence under $\cWsim$ is 
straightforward
using the fact that $\cW(z)$ is closed under  $\cWsim$
and no permutation in this set has any consecutive $321$-patterns by Lemma~\ref{CAB-lem}.
\end{proof}

It is convenient to introduce the following additional notation.

\begin{definition}\label{undroop-def}
For a triangular pipe dream $E$ and distinct pairs $(a,b),(i,j) \in \LT_n$ with $a\leq i$ and $b\leq j$,
define 
$\undroop^{ab}_{ij}(E) = D
$
if there is a triangular pipe dream $D$ with 
\be\label{undroop-eq}
D\neq \extDroop^{ab}_{ij}(D) = E.\ee 
Equivalently, if 
the tiles of $E$  match the outputs to the transformations in  Definition~\ref{extended-droop-def}(a),
and applying the inverses of those transformations results in a triangular pipe dream $D$ that has no bend tiles in the region specified by the same definition,
then we set
$\undroop^{ab}_{ij}(E)=D
$.
Otherwise, we have $\undroop^{ab}_{ij}(E) = E
$. One defines the operation $\complexUndroop^{ab}_{ij}$ analogously.
\end{definition}

The \defn{droop moves} for square pipe dreams introduced 
 in \cite[\S5.2]{LLS} are defined in exactly the same way as our simple droop moves for triangular pipe dreams.
Weigandt \cite[\S4.3]{Weigandt} considers a more general class
of \defn{$K$-theoretic droop moves}
that span and preserve the larger set of (not necessarily reduced) pipe dreams $\kBPD(w)$ for any $w \in S_n$.

Such $K$-theoretic droop moves have the following 
relationship to our pipe dream operations:
when two triangular pipe dreams $D$ and $E$ are connected by a 
simple or complex droop move,
their symmetrized versions $\double(D)$ and $\double(E)$ 
are connected by a sequence of $K$-theoretic droop moves.
One can combine this property with Proposition~\ref{sym-bij-prop} to show that $\kInvBPD(z)$ is  preserved by all strong droop moves that do not create any bump tiles $\ETILE$ or $\EDIAG$.

More generally, one can formulate  \defn{$K$-theoretic droop moves}
for triangular pipe dreams with the same property that generate all elements of $\kInvBPD(z)$ from $\invRothe(z)$.
However, we will not attempt a rigorous treatment of these ideas here.

\begin{example}
The symmetrized versions of the moves $\droop^{21}_{43}$ and $\complexDroop^{21}_{43}$ shown
in Example~\ref{tri-droop-ex} are realized via $K$-theoretic droops
(as defined in \cite[\S4.3]{Weigandt}) by
\[
\begin{tikzpicture}[x=\bpdwidth,y=\bpdwidth,line cap=round,line join=round,baseline=(z.base)]
\node at (0,2.25) (z) {};
\draw[gray!70, line width=0.7pt](0,0)--(5,0);
\draw[gray!70, line width=0.7pt](0,0)--(0,5);
\draw[gray!70, line width=0.7pt](0,1)--(5,1);
\draw[gray!70, line width=0.7pt](1,0)--(1,5);
\draw[gray!70, line width=0.7pt](0,2)--(5,2);
\draw[gray!70, line width=0.7pt](2,0)--(2,5);
\draw[gray!70, line width=0.7pt](0,3)--(5,3);
\draw[gray!70, line width=0.7pt](3,0)--(3,5);
\draw[gray!70, line width=0.7pt](0,4)--(5,4);
\draw[gray!70, line width=0.7pt](4,0)--(4,5);
\draw[gray!70, line width=0.7pt](0,5)--(5,5);
\draw[gray!70, line width=0.7pt](5,0)--(5,5);
\draw[thick,rounded corners,color=blue] (1.5,4)--(1.5,4.5)--(2,4.5);
\draw[thick,rounded corners,color=blue] (0.5,3)--(0.5,3.5)--(1,3.5);
\draw[thick,rounded corners,color=blue] (4.5,2)--(4.5,2.5)--(5,2.5);
\draw[thick,rounded corners,color=blue] (3.5,1)--(3.5,1.5)--(4,1.5);
\draw[thick,rounded corners,color=blue] (2.5,0)--(2.5,0.5)--(3,0.5);
\draw[thick,rounded corners,color=blue] (0.5,0)--(0.5,1);
\draw[thick,rounded corners,color=blue] (0.5,1)--(0.5,2);
\draw[thick,rounded corners,color=blue] (0.5,2)--(0.5,3);
\draw[thick,rounded corners,color=blue] (1,3.5)--(2,3.5);
\draw[thick,rounded corners,color=blue] (1.5,3)--(1.5,4);
\draw[thick,rounded corners,color=blue] (2,3.5)--(3,3.5);
\draw[thick,rounded corners,color=blue] (3,3.5)--(4,3.5);
\draw[thick,rounded corners,color=blue] (4,3.5)--(5,3.5);
\draw[thick,rounded corners,color=blue] (1.5,0)--(1.5,1);
\draw[thick,rounded corners,color=blue] (1.5,1)--(1.5,2);
\draw[thick,rounded corners,color=blue] (1.5,2)--(1.5,3);
\draw[thick,rounded corners,color=blue] (2,4.5)--(3,4.5);
\draw[thick,rounded corners,color=blue] (3,4.5)--(4,4.5);
\draw[thick,rounded corners,color=blue] (4,4.5)--(5,4.5);
\draw[thick,rounded corners,color=blue] (3,0.5)--(4,0.5);
\draw[thick,rounded corners,color=blue] (3.5,0)--(3.5,1);
\draw[thick,rounded corners,color=blue] (4,0.5)--(5,0.5);
\draw[thick,rounded corners,color=blue] (4.5,0)--(4.5,1);
\draw[thick,rounded corners,color=blue] (4,1.5)--(5,1.5);
\draw[thick,rounded corners,color=blue] (4.5,1)--(4.5,2);
\draw[black, line width=1.6pt] (0,5)--(5,5)--(5,0)--(0,0)--(0,5);
\end{tikzpicture}
\ \xmapsto[43]{21}
\begin{tikzpicture}[x=\bpdwidth,y=\bpdwidth,line cap=round,line join=round,baseline=(z.base)]
\node at (0,2.25) (z) {};
\draw[gray!70, line width=0.7pt](0,0)--(5,0);
\draw[gray!70, line width=0.7pt](0,0)--(0,5);
\draw[gray!70, line width=0.7pt](0,1)--(5,1);
\draw[gray!70, line width=0.7pt](1,0)--(1,5);
\draw[gray!70, line width=0.7pt](0,2)--(5,2);
\draw[gray!70, line width=0.7pt](2,0)--(2,5);
\draw[gray!70, line width=0.7pt](0,3)--(5,3);
\draw[gray!70, line width=0.7pt](3,0)--(3,5);
\draw[gray!70, line width=0.7pt](0,4)--(5,4);
\draw[gray!70, line width=0.7pt](4,0)--(4,5);
\draw[gray!70, line width=0.7pt](0,5)--(5,5);
\draw[gray!70, line width=0.7pt](5,0)--(5,5);
\draw[thick,rounded corners,color=blue] (1.5,4)--(1.5,4.5)--(2,4.5);
\draw[thick,rounded corners,color=blue] (4.5,2)--(4.5,2.5)--(5,2.5);
\draw[thick,rounded corners,color=blue] (3.5,1)--(3.5,1.5)--(4,1.5);
\draw[thick,rounded corners,color=blue] (2.5,0)--(2.5,0.5)--(3,0.5);
\draw[thick,rounded corners,color=blue] (2.5,2)--(2.5,1.5)--(2,1.5);
\draw[thick,rounded corners,color=blue] (2.5,3)--(2.5,3.5)--(3,3.5);
\draw[thick,rounded corners,color=blue] (0.5,1)--(0.5,1.5)--(1,1.5);
\draw[thick,rounded corners,color=blue] (0.5,0)--(0.5,1);
\draw[thick,rounded corners,color=blue] (1,1.5)--(2,1.5);
\draw[thick,rounded corners,color=blue] (1.5,1)--(1.5,2);
\draw[thick,rounded corners,color=blue] (2.5,2)--(2.5,3);
\draw[thick,rounded corners,color=blue] (3,3.5)--(4,3.5);
\draw[thick,rounded corners,color=blue] (4,3.5)--(5,3.5);
\draw[thick,rounded corners,color=blue] (1.5,0)--(1.5,1);
\draw[thick,rounded corners,color=blue] (1.5,2)--(1.5,3);
\draw[thick,rounded corners,color=blue] (1.5,3)--(1.5,4);
\draw[thick,rounded corners,color=blue] (2,4.5)--(3,4.5);
\draw[thick,rounded corners,color=blue] (3,4.5)--(4,4.5);
\draw[thick,rounded corners,color=blue] (4,4.5)--(5,4.5);
\draw[thick,rounded corners,color=blue] (3,0.5)--(4,0.5);
\draw[thick,rounded corners,color=blue] (3.5,0)--(3.5,1);
\draw[thick,rounded corners,color=blue] (4,0.5)--(5,0.5);
\draw[thick,rounded corners,color=blue] (4.5,0)--(4.5,1);
\draw[thick,rounded corners,color=blue] (4,1.5)--(5,1.5);
\draw[thick,rounded corners,color=blue] (4.5,1)--(4.5,2);
\draw[black, line width=1.6pt] (0,5)--(5,5)--(5,0)--(0,0)--(0,5);
\end{tikzpicture}
\ \xmapsto[23]{12}
\begin{tikzpicture}[x=\bpdwidth,y=\bpdwidth,line cap=round,line join=round,baseline=(z.base)]
\node at (0,2.25) (z) {};
\draw[gray!70, line width=0.7pt](0,0)--(5,0);
\draw[gray!70, line width=0.7pt](0,0)--(0,5);
\draw[gray!70, line width=0.7pt](0,1)--(5,1);
\draw[gray!70, line width=0.7pt](1,0)--(1,5);
\draw[gray!70, line width=0.7pt](0,2)--(5,2);
\draw[gray!70, line width=0.7pt](2,0)--(2,5);
\draw[gray!70, line width=0.7pt](0,3)--(5,3);
\draw[gray!70, line width=0.7pt](3,0)--(3,5);
\draw[gray!70, line width=0.7pt](0,4)--(5,4);
\draw[gray!70, line width=0.7pt](4,0)--(4,5);
\draw[gray!70, line width=0.7pt](0,5)--(5,5);
\draw[gray!70, line width=0.7pt](5,0)--(5,5);
\draw[thick,rounded corners,color=blue] (4.5,2)--(4.5,2.5)--(5,2.5);
\draw[thick,rounded corners,color=blue] (3.5,1)--(3.5,1.5)--(4,1.5);
\draw[thick,rounded corners,color=blue] (2.5,0)--(2.5,0.5)--(3,0.5);
\draw[thick,rounded corners,color=blue] (2.5,2)--(2.5,1.5)--(2,1.5);
\draw[thick,rounded corners,color=blue] (0.5,1)--(0.5,1.5)--(1,1.5);
\draw[thick,rounded corners,color=blue] (2.5,4)--(2.5,4.5)--(3,4.5);
\draw[thick,rounded corners,color=blue] (1.5,3)--(1.5,3.5)--(2,3.5);
\draw[thick,rounded corners,color=blue] (0.5,0)--(0.5,1);
\draw[thick,rounded corners,color=blue] (1,1.5)--(2,1.5);
\draw[thick,rounded corners,color=blue] (1.5,1)--(1.5,2);
\draw[thick,rounded corners,color=blue] (2.5,2)--(2.5,3);
\draw[thick,rounded corners,color=blue] (2,3.5)--(3,3.5);
\draw[thick,rounded corners,color=blue] (2.5,3)--(2.5,4);
\draw[thick,rounded corners,color=blue] (3,4.5)--(4,4.5);
\draw[thick,rounded corners,color=blue] (4,4.5)--(5,4.5);
\draw[thick,rounded corners,color=blue] (1.5,0)--(1.5,1);
\draw[thick,rounded corners,color=blue] (1.5,2)--(1.5,3);
\draw[thick,rounded corners,color=blue] (3,3.5)--(4,3.5);
\draw[thick,rounded corners,color=blue] (4,3.5)--(5,3.5);
\draw[thick,rounded corners,color=blue] (3,0.5)--(4,0.5);
\draw[thick,rounded corners,color=blue] (3.5,0)--(3.5,1);
\draw[thick,rounded corners,color=blue] (4,0.5)--(5,0.5);
\draw[thick,rounded corners,color=blue] (4.5,0)--(4.5,1);
\draw[thick,rounded corners,color=blue] (4,1.5)--(5,1.5);
\draw[thick,rounded corners,color=blue] (4.5,1)--(4.5,2);
\draw[black, line width=1.6pt] (0,5)--(5,5)--(5,0)--(0,0)--(0,5);
\end{tikzpicture}
\ \xmapsto[34]{13}
\begin{tikzpicture}[x=\bpdwidth,y=\bpdwidth,line cap=round,line join=round,baseline=(z.base)]
\node at (0,2.25) (z) {};
\draw[gray!70, line width=0.7pt](0,0)--(5,0);
\draw[gray!70, line width=0.7pt](0,0)--(0,5);
\draw[gray!70, line width=0.7pt](0,1)--(5,1);
\draw[gray!70, line width=0.7pt](1,0)--(1,5);
\draw[gray!70, line width=0.7pt](0,2)--(5,2);
\draw[gray!70, line width=0.7pt](2,0)--(2,5);
\draw[gray!70, line width=0.7pt](0,3)--(5,3);
\draw[gray!70, line width=0.7pt](3,0)--(3,5);
\draw[gray!70, line width=0.7pt](0,4)--(5,4);
\draw[gray!70, line width=0.7pt](4,0)--(4,5);
\draw[gray!70, line width=0.7pt](0,5)--(5,5);
\draw[gray!70, line width=0.7pt](5,0)--(5,5);
\draw[thick,rounded corners,color=blue] (4.5,2)--(4.5,2.5)--(5,2.5);
\draw[thick,rounded corners,color=blue] (3.5,1)--(3.5,1.5)--(4,1.5);
\draw[thick,rounded corners,color=blue] (2.5,0)--(2.5,0.5)--(3,0.5);
\draw[thick,rounded corners,color=blue] (2.5,2)--(2.5,1.5)--(2,1.5);
\draw[thick,rounded corners,color=blue] (0.5,1)--(0.5,1.5)--(1,1.5);
\draw[thick,rounded corners,color=blue] (1.5,3)--(1.5,3.5)--(2,3.5);
\draw[thick,rounded corners,color=blue] (3.5,3)--(3.5,2.5)--(3,2.5);
\draw[thick,rounded corners,color=blue] (3.5,4)--(3.5,4.5)--(4,4.5);
\draw[thick,rounded corners,color=blue] (2.5,2)--(2.5,2.5)--(3,2.5);
\draw[thick,rounded corners,color=blue] (0.5,0)--(0.5,1);
\draw[thick,rounded corners,color=blue] (1,1.5)--(2,1.5);
\draw[thick,rounded corners,color=blue] (1.5,1)--(1.5,2);
\draw[thick,rounded corners,color=blue] (3,3.5)--(4,3.5);
\draw[thick,rounded corners,color=blue] (3.5,3)--(3.5,4);
\draw[thick,rounded corners,color=blue] (4,4.5)--(5,4.5);
\draw[thick,rounded corners,color=blue] (1.5,0)--(1.5,1);
\draw[thick,rounded corners,color=blue] (1.5,2)--(1.5,3);
\draw[thick,rounded corners,color=blue] (2,3.5)--(3,3.5);
\draw[thick,rounded corners,color=blue] (4,3.5)--(5,3.5);
\draw[thick,rounded corners,color=blue] (3,0.5)--(4,0.5);
\draw[thick,rounded corners,color=blue] (3.5,0)--(3.5,1);
\draw[thick,rounded corners,color=blue] (4,0.5)--(5,0.5);
\draw[thick,rounded corners,color=blue] (4.5,0)--(4.5,1);
\draw[thick,rounded corners,color=blue] (4,1.5)--(5,1.5);
\draw[thick,rounded corners,color=blue] (4.5,1)--(4.5,2);
\draw[black, line width=1.6pt] (0,5)--(5,5)--(5,0)--(0,0)--(0,5);
\end{tikzpicture}
\]
and
\[
\begin{tikzpicture}[x=\bpdwidth,y=\bpdwidth,line cap=round,line join=round,baseline=(z.base)]
\node at (0,2.25) (z) {};
\draw[gray!70, line width=0.7pt](0,0)--(5,0);
\draw[gray!70, line width=0.7pt](0,0)--(0,5);
\draw[gray!70, line width=0.7pt](0,1)--(5,1);
\draw[gray!70, line width=0.7pt](1,0)--(1,5);
\draw[gray!70, line width=0.7pt](0,2)--(5,2);
\draw[gray!70, line width=0.7pt](2,0)--(2,5);
\draw[gray!70, line width=0.7pt](0,3)--(5,3);
\draw[gray!70, line width=0.7pt](3,0)--(3,5);
\draw[gray!70, line width=0.7pt](0,4)--(5,4);
\draw[gray!70, line width=0.7pt](4,0)--(4,5);
\draw[gray!70, line width=0.7pt](0,5)--(5,5);
\draw[gray!70, line width=0.7pt](5,0)--(5,5);
\draw[thick,rounded corners,color=blue] (1.5,4)--(1.5,4.5)--(2,4.5);
\draw[thick,rounded corners,color=blue] (0.5,3)--(0.5,3.5)--(1,3.5);
\draw[thick,rounded corners,color=blue] (4.5,2)--(4.5,2.5)--(5,2.5);
\draw[thick,rounded corners,color=blue] (3.5,1)--(3.5,1.5)--(4,1.5);
\draw[thick,rounded corners,color=blue] (2.5,0)--(2.5,0.5)--(3,0.5);
\draw[thick,rounded corners,color=blue] (0.5,0)--(0.5,1);
\draw[thick,rounded corners,color=blue] (0.5,1)--(0.5,2);
\draw[thick,rounded corners,color=blue] (0.5,2)--(0.5,3);
\draw[thick,rounded corners,color=blue] (1,3.5)--(2,3.5);
\draw[thick,rounded corners,color=blue] (1.5,3)--(1.5,4);
\draw[thick,rounded corners,color=blue] (2,3.5)--(3,3.5);
\draw[thick,rounded corners,color=blue] (3,3.5)--(4,3.5);
\draw[thick,rounded corners,color=blue] (4,3.5)--(5,3.5);
\draw[thick,rounded corners,color=blue] (1.5,0)--(1.5,1);
\draw[thick,rounded corners,color=blue] (1.5,1)--(1.5,2);
\draw[thick,rounded corners,color=blue] (1.5,2)--(1.5,3);
\draw[thick,rounded corners,color=blue] (2,4.5)--(3,4.5);
\draw[thick,rounded corners,color=blue] (3,4.5)--(4,4.5);
\draw[thick,rounded corners,color=blue] (4,4.5)--(5,4.5);
\draw[thick,rounded corners,color=blue] (3,0.5)--(4,0.5);
\draw[thick,rounded corners,color=blue] (3.5,0)--(3.5,1);
\draw[thick,rounded corners,color=blue] (4,0.5)--(5,0.5);
\draw[thick,rounded corners,color=blue] (4.5,0)--(4.5,1);
\draw[thick,rounded corners,color=blue] (4,1.5)--(5,1.5);
\draw[thick,rounded corners,color=blue] (4.5,1)--(4.5,2);
\draw[black, line width=1.6pt] (0,5)--(5,5)--(5,0)--(0,0)--(0,5);
\end{tikzpicture}
\xmapsto[34]{21}
\begin{tikzpicture}[x=\bpdwidth,y=\bpdwidth,line cap=round,line join=round,baseline=(z.base)]
\node at (0,2.25) (z) {};
\draw[gray!70, line width=0.7pt](0,0)--(5,0);
\draw[gray!70, line width=0.7pt](0,0)--(0,5);
\draw[gray!70, line width=0.7pt](0,1)--(5,1);
\draw[gray!70, line width=0.7pt](1,0)--(1,5);
\draw[gray!70, line width=0.7pt](0,2)--(5,2);
\draw[gray!70, line width=0.7pt](2,0)--(2,5);
\draw[gray!70, line width=0.7pt](0,3)--(5,3);
\draw[gray!70, line width=0.7pt](3,0)--(3,5);
\draw[gray!70, line width=0.7pt](0,4)--(5,4);
\draw[gray!70, line width=0.7pt](4,0)--(4,5);
\draw[gray!70, line width=0.7pt](0,5)--(5,5);
\draw[gray!70, line width=0.7pt](5,0)--(5,5);
\draw[thick,rounded corners,color=blue] (1.5,4)--(1.5,4.5)--(2,4.5);
\draw[thick,rounded corners,color=blue] (4.5,2)--(4.5,2.5)--(5,2.5);
\draw[thick,rounded corners,color=blue] (3.5,1)--(3.5,1.5)--(4,1.5);
\draw[thick,rounded corners,color=blue] (2.5,0)--(2.5,0.5)--(3,0.5);
\draw[thick,rounded corners,color=blue] (3.5,3)--(3.5,2.5)--(3,2.5);
\draw[thick,rounded corners,color=blue] (3.5,3)--(3.5,3.5)--(4,3.5);
\draw[thick,rounded corners,color=blue] (0.5,2)--(0.5,2.5)--(1,2.5);
\draw[thick,rounded corners,color=blue] (0.5,0)--(0.5,1);
\draw[thick,rounded corners,color=blue] (0.5,1)--(0.5,2);
\draw[thick,rounded corners,color=blue] (1,2.5)--(2,2.5);
\draw[thick,rounded corners,color=blue] (1.5,2)--(1.5,3);
\draw[thick,rounded corners,color=blue] (2,2.5)--(3,2.5);
\draw[thick,rounded corners,color=blue] (4,3.5)--(5,3.5);
\draw[thick,rounded corners,color=blue] (1.5,0)--(1.5,1);
\draw[thick,rounded corners,color=blue] (1.5,1)--(1.5,2);
\draw[thick,rounded corners,color=blue] (1.5,3)--(1.5,4);
\draw[thick,rounded corners,color=blue] (2,4.5)--(3,4.5);
\draw[thick,rounded corners,color=blue] (3,4.5)--(4,4.5);
\draw[thick,rounded corners,color=blue] (4,4.5)--(5,4.5);
\draw[thick,rounded corners,color=blue] (3,0.5)--(4,0.5);
\draw[thick,rounded corners,color=blue] (3.5,0)--(3.5,1);
\draw[thick,rounded corners,color=blue] (4,0.5)--(5,0.5);
\draw[thick,rounded corners,color=blue] (4.5,0)--(4.5,1);
\draw[thick,rounded corners,color=blue] (4,1.5)--(5,1.5);
\draw[thick,rounded corners,color=blue] (4.5,1)--(4.5,2);
\draw[black, line width=1.6pt] (0,5)--(5,5)--(5,0)--(0,0)--(0,5);
\end{tikzpicture}
\xmapsto[43]{12}
\begin{tikzpicture}[x=\bpdwidth,y=\bpdwidth,line cap=round,line join=round,baseline=(z.base)]
\node at (0,2.25) (z) {};
\draw[gray!70, line width=0.7pt](0,0)--(5,0);
\draw[gray!70, line width=0.7pt](0,0)--(0,5);
\draw[gray!70, line width=0.7pt](0,1)--(5,1);
\draw[gray!70, line width=0.7pt](1,0)--(1,5);
\draw[gray!70, line width=0.7pt](0,2)--(5,2);
\draw[gray!70, line width=0.7pt](2,0)--(2,5);
\draw[gray!70, line width=0.7pt](0,3)--(5,3);
\draw[gray!70, line width=0.7pt](3,0)--(3,5);
\draw[gray!70, line width=0.7pt](0,4)--(5,4);
\draw[gray!70, line width=0.7pt](4,0)--(4,5);
\draw[gray!70, line width=0.7pt](0,5)--(5,5);
\draw[gray!70, line width=0.7pt](5,0)--(5,5);
\draw[thick,rounded corners,color=blue] (4.5,2)--(4.5,2.5)--(5,2.5);
\draw[thick,rounded corners,color=blue] (3.5,1)--(3.5,1.5)--(4,1.5);
\draw[thick,rounded corners,color=blue] (2.5,0)--(2.5,0.5)--(3,0.5);
\draw[thick,rounded corners,color=blue] (3.5,3)--(3.5,2.5)--(3,2.5);
\draw[thick,rounded corners,color=blue] (3.5,3)--(3.5,3.5)--(4,3.5);
\draw[thick,rounded corners,color=blue] (0.5,2)--(0.5,2.5)--(1,2.5);
\draw[thick,rounded corners,color=blue] (2.5,2)--(2.5,1.5)--(2,1.5);
\draw[thick,rounded corners,color=blue] (2.5,4)--(2.5,4.5)--(3,4.5);
\draw[thick,rounded corners,color=blue] (1.5,1)--(1.5,1.5)--(2,1.5);
\draw[thick,rounded corners,color=blue] (0.5,0)--(0.5,1);
\draw[thick,rounded corners,color=blue] (0.5,1)--(0.5,2);
\draw[thick,rounded corners,color=blue] (1,2.5)--(2,2.5);
\draw[thick,rounded corners,color=blue] (2,2.5)--(3,2.5);
\draw[thick,rounded corners,color=blue] (2.5,2)--(2.5,3);
\draw[thick,rounded corners,color=blue] (4,3.5)--(5,3.5);
\draw[thick,rounded corners,color=blue] (1.5,0)--(1.5,1);
\draw[thick,rounded corners,color=blue] (2.5,3)--(2.5,4);
\draw[thick,rounded corners,color=blue] (3,4.5)--(4,4.5);
\draw[thick,rounded corners,color=blue] (4,4.5)--(5,4.5);
\draw[thick,rounded corners,color=blue] (3,0.5)--(4,0.5);
\draw[thick,rounded corners,color=blue] (3.5,0)--(3.5,1);
\draw[thick,rounded corners,color=blue] (4,0.5)--(5,0.5);
\draw[thick,rounded corners,color=blue] (4.5,0)--(4.5,1);
\draw[thick,rounded corners,color=blue] (4,1.5)--(5,1.5);
\draw[thick,rounded corners,color=blue] (4.5,1)--(4.5,2);
\draw[black, line width=1.6pt] (0,5)--(5,5)--(5,0)--(0,0)--(0,5);
\end{tikzpicture}
.
\]
\end{example}

\section{Transition arguments}
\label{S: Proof}

This section proves Theorem~\ref{main-thm}.
Recall that the transition recurrence~\eqref{EQ: iS Tran}, together with the base case $\iS_z = 1$
when $z$ is the identity, uniquely determines
the polynomial $\iS_z$.
Thus, our strategy is to show that 
the weight generating functions of bumpless involution pipe dreams
satisfy the same recurrence. 

\subsection{Fixed regions of bumpless pipe dreams}
\label{SS: fixed regions}

We first record some elementary facts about (square) bumpless pipe dreams for permutations.
Recall that such pipe dreams are not assumed to be reduced.

Let $w\in S_n$ be a non-identity permutation with maximal Rothe element $(k,p)$.
Define the \defn{fixed region} of $w$ to consist of
all positions in rows $k+1,\ldots,n$,
together with the positions $(k,j)$ with $j>p$.
For the rest of this subsection, we fix a bumpless pipe dream $D$ for $w$. 

\begin{lemma}
\label{L: BPD fixed bottom}
In row $r$ of the fixed region for $w$,
all pipes in $D$ travel straight vertically except for pipe $w(r)$.
Pipe $w(r)$ of $D$ enters row $r$ from below at position $(r,w(r))$,
bends to the right,
and then travels horizontally to the east boundary.
\end{lemma}

\begin{proof}
Since $(k,p)$ is the maximal Rothe element of $w$,
there are no inversions whose first coordinate is larger than $k$.
Thus,
$w(k+1)<w(k+2)<\cdots<w(n)$.

We argue by induction on $r=n,n-1,\ldots,k$.
The case $r=n$ is immediate.
Now take $k\leq r<n$.
By our inductive hypothesis, 
for every $i$ with $w^{-1}(i)\leq r$,
pipe $i$ travels vertically below row $r$.
In particular, pipe $w(r)$ enters row $r$
from below at $(r,w(r))$.
Since this pipe exits through the east boundary of row $r$,
it bends there and then travels horizontally
through the rest of the row.

It remains to show that no other pipes bend in the fixed region of row $r$.
Suppose this is not the case. 
Then we can find a $\JTILE$ tile at some position $(r,c)$ in the fixed region.
Since pipe $w(r)$ occupies the tiles weakly to the
right of $(r,w(r))$, we have $c<w(r)$.
Moreover, no pipe enters $(r,c)$ from below.
Thus, $c\neq w(i)$ for every $i\leq r$,
so $c\in\{w(r+1),\ldots,w(n)\}$.
If $r>k$, then
$w(r)<w(r+1)<\cdots<w(n)$,
contradicting $c<w(r)$.
If $r=k$, then $c>p$ since $(k,c)$ lies in the fixed region.
Then we have $c<w(k)$ and $w^{-1}(c)>k$,
implying that $(k,c)\in D(w)$.
But $(k,c)$ is larger than $(k,p)$, contradicting 
the definition of the maximal Rothe element. 
\end{proof}

The following is immediate from
Lemma~\ref{L: BPD fixed bottom}.

\begin{corollary}
\label{C: fixed regions description}
Consider a position $(r,j)$ in the fixed region of $w$.
If $D$ has a pipe traveling horizontally in this position, then 
it is pipe $w(r)$.
If $D$ has a pipe traveling vertically in this position,
then it is pipe $j$.
\end{corollary}


We can now characterize all crossing tiles of $D$ 
in the fixed region.

\begin{corollary}
\label{C: BPD crossing fixed region}
Consider a position $(r,j)$ in the fixed region of $w$.
The pipe dream $D$ contains a crossing tile in this position if and only if
$w(r)<j$ and $w^{-1}(j)<r$.
In this case, the crossing involves
pipe $w(r)$ and pipe $j$.
\end{corollary}

\begin{proof}
By Corollary~\ref{C: fixed regions description},
the only pipe that can travel horizontally (respectively, vertically) through $(r,j)$
is pipe $w(r)$ (respectively, $j$).
By Lemma~\ref{L: BPD fixed bottom},
this pipe travels horizontally through $(r,j)$
if and only if $j>w(r)$.
Pipe $j$ travels vertically through $(r,j)$
if and only if it exits in a row strictly above $r$,
that is, if and only if $w^{-1}(j)<r$.
\end{proof}

We also record the possibilities for the position
of the maximal Rothe element itself.

\begin{corollary}
\label{C: j or blank}
The pipe dream contains in position $(k,p)$ either $\BTILE$ or $\JTILE$.
\end{corollary}

\begin{proof}
The position immediately to the right of $(k,p)$
lies in the fixed region.
By Lemma~\ref{L: BPD fixed bottom},
there is no connection from this position
entering $(k,p)$ from the right.
Thus, no pipe exits $(k,p)$ to the right.

There is also no pipe entering $(k,p)$ from below.
If such a pipe exists, by Lemma~\ref{L: BPD fixed bottom},
it travels vertically under row $k$,
so it is pipe $p$.
However, $(k,p)\in D(w)$ implies $w^{-1}(p)>k$,
so pipe $p$ cannot appear in row $k$.
Hence $(k,p)$ has no connection to the right or the bottom, which implies our claim. 
\end{proof}

We finish with a fact about cutting a particular crossing.
By \defn{cutting}, we mean replacing a $\PTILE$ 
by $\ETILE$.

\begin{lemma}
\label{L: No double crossing after cut}
Set $r=w^{-1}(p)$ and $c=w(k)$.
Then position $(r,c)$ of $D$ is a crossing of pipe $p$ and pipe $c$.
After cutting this crossing,
neither of these two pipes crosses any other pipe twice.
\end{lemma}

\begin{figure}
    \centering
    \begin{tikzpicture}[x=0.8cm,y=0.8cm,line cap=round,line join=round]
  \def\m{7}
  \def\n{6}

  \fill[blue!20]   (0,2) rectangle (3,3);
  \fill[blue!20]   (0,0) rectangle (1,2);

  \fill[green!20]  (4,2) rectangle (7,3);

  \fill[yellow!30] (3,0) rectangle (4,2);

  \fill[red!20]    (3,3) rectangle (4,5);
  \fill[red!20]    (3,5) rectangle (7,6);

  \foreach \k in {0,...,\m} {
    \draw[gray!70, line width=0.7pt] (\k,0) -- (\k,\n);
  }
  \foreach \k in {0,...,\n} {
    \draw[gray!70, line width=0.7pt] (0,\k) -- (\m,\k);
  }

  \draw[blue] (0.5,0) -- (0.5,2.5) -- (7,2.5);
  \draw[blue] (3.5,0) -- (3.5,5.5) -- (7,5.5);

  \draw[line width=1.6pt] (0,0) rectangle (\m,\n);
\end{tikzpicture}
    \caption{Regions corresponding to the four sets $A_1$, $A_2$, $B_1$, and $B_2$ defined in the proof of Lemma~\ref{L: No double crossing after cut}}
    \label{F: (proof) two pipes}
\end{figure}
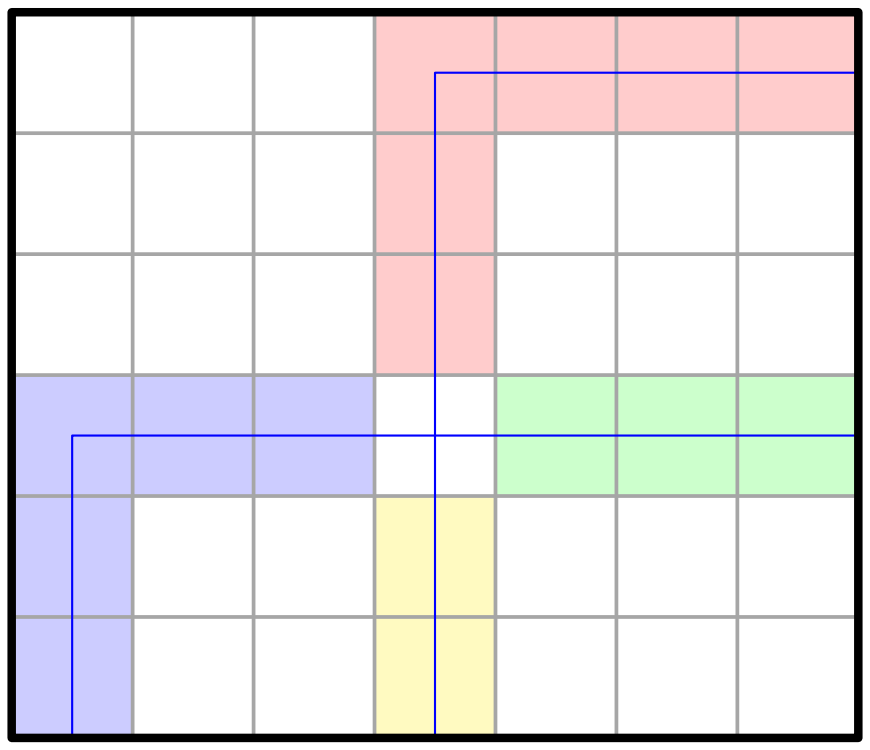

\begin{proof}
Since $(k,p)\in D(w)$, we have $k<r$ and $p<c$.
By Corollary~\ref{C: BPD crossing fixed region},
the position $(r,c)$ is a crossing of pipe $p$
and pipe $c=w(k)$.

Both of these pipes lie entirely in the fixed region
by Lemma~\ref{L: BPD fixed bottom}.
Let $A_1$ and $A_2$ be the sets of labels of pipes in $D$ that 
cross pipe $p$ respectively before and after $(r,c)$,
and define $B_1$ and $B_2$ similarly for pipe $c$.
These four sets are the sets of labels
of the pipes encountered by pipe $p$ and pipe $c$
in the four colored regions in Figure~\ref{F: (proof) two pipes}.
After cutting the crossing at $(r,c)$,
a new double crossing can occur only if
$A_1\cap B_2\neq\varnothing$ or
$A_2\cap B_1\neq\varnothing$.
We show that both intersections are empty.

Using Corollary~\ref{C: BPD crossing fixed region}, we have
\[
\ba
A_1
&=
\{j:j<p,\ r<w^{-1}(j)\}
\sqcup
\{j:p<j<c,\ w^{-1}(j)<r\},
\\
A_2&=\{j:c<j,\ w^{-1}(j)<r\},
\\
B_1&=\{j:j<c,\ r<w^{-1}(j)\},\quand 
\\
B_2
&=
\{j:j<c,\ k<w^{-1}(j)<r\}
\sqcup
\{j:c<j,\ w^{-1}(j)<k\}.
\ea
\]
It is immediate that $A_2\cap B_1=\varnothing$.

For $A_1\cap B_2$, all intersections are clearly empty
except possibly between the second set in $A_1$
and the first set in $B_2$.
An element in this intersection would be some $j$ satisfying
$p<j<c$ and $k<w^{-1}(j)<r$.
Then $(w^{-1}(j),r)$ is an inversion of $w$,
since $w(w^{-1}(j))=j>p=w(r)$.
Its first coordinate is larger than $k$,
contradicting the fact that
$(k,w^{-1}(p))=(k,r)$ is the maximal inversion of $w$
by Lemma~\ref{L: maximal Rothe and inv}.
Thus, $A_1\cap B_2=\varnothing$ as well.
\end{proof}

\subsection{Transition operators}
\label{SS: transition on InvBPD}

Fix a non-identity involution $1\neq z\in \Invol$
with maximal Rothe element $(q,p)$,
and set $p'=z(p)$ and $q'=z(q)$ as in Section~2.4.
The transition formula~\eqref{EQ: iS Tran},
together with the base case $\iS_1=1$,
recursively determines every involution Schubert polynomial.
Thus, to prove Theorem~\ref{main-thm},
it suffices to show that the weight generating functions
of bumpless involution pipe dreams satisfy the same recurrence.
To this end, 
set
\[
\calT(z)
=
\InvBPD(\eta(z))
\bigsqcup
\left(
\bigsqcup_{v\in\hat\Phi^{-}(\eta(z),p)}
\InvBPD(v)
\right).
\]
We now define a transition operator on bumpless involution pipe dreams. This is 
analogous to an operator on square bumpless pipe dreams 
considered in \cite{Weigandt}. 
Recall our Cases (A), (B), and (C) from Lemma~\ref{L: ABC pq}.

\begin{definition}
\label{D: Transition operator}
For $D\in\InvBPD(z)$,
define $\bft(D)$ by the following three steps.
\begin{enumerate}
\item[1.](Cut)
In Cases (A) and (B),
replace the $\PDIAG$ tile in $D$ at $(p',p')$
by  $\EDIAG$.

In \CaseC,
replace the $\PTILE$ tile in $D$ at $(p',q')$
by $\ETILE$.

\item[2.](Undroop)
Perform the corresponding undroop move from the tile modified in the Cut step to $(q,p)$.

That is, undroop
from $(p',p')$ to $(q,p)$ in Cases (A) and (B),
and
from $(p',q')$ to $(q,p)$ in \CaseC.

\item[3.](Fuse)
After the Undroop step,
the tile at $(q,p)$ is 
$\CDIAG$, $\CTILE$, $\EDIAG$, or $\ETILE$.

When it is $\CDIAG$ or $\CTILE$, do nothing.
If the tile is $\EDIAG$ or $\ETILE$, then
replace it by $\PDIAG$ or $\PTILE$, respectively.
\end{enumerate}
We call $\bft$ the \defn{transition operator}.
\end{definition}

Here is this section's main result.
Recall that $1\neq z \in \Invol$ is fixed
with maximal Rothe element $(q,p)$.
Also recall for pipe dreams $D$ that
$\weight(D) = \prod_{(i,j) \in \blank(D)} 
 x_{(i,j)}$ where $x_{(i,j)} = 2^{|\{i,j\}|-2} (x_i +x_j)$.

\begin{theorem}
\label{T: transition}
The  operator $\bft$ is a bijection
$
\InvBPD(z)\longrightarrow\calT(z).
$
Moreover, if $D\in\InvBPD(z)$ then
\be\label{wt-t-eq}
\weight(D)
=
\begin{cases}
 \weight(\bft(D))x_{(p,q)}
&\text{if }
\bft(D)\in\InvBPD(\eta(z)),
\\
\weight(\bft(D))
&\text{if }
\bft(D)\in\InvBPD(v)
\text{ for some }
v\in\hat\Phi^{-}(\eta(z),p).
\end{cases}
\ee
\end{theorem}

We prove this bijectivity by constructing an explicit inverse.

\begin{definition}
\label{D: Reverse Transition operator}
For $D\in\calT(z)$,
define $\bftrev(D)$ by the following three steps.
\begin{enumerate}
\item[1.] (Cut)
If the tile in $D$ at $(q,p)$ is $\PDIAG$
or $\PTILE$, replace it by $\EDIAG$
or $\ETILE$, respectively.

If the tile is $\CDIAG$ or $\CTILE$,
do nothing.

\item[2.] (Droop)
Perform the strong droop move
from $(q,p)$ to $(p',p')$ in Cases (A) and (B).

In \CaseC, perform the strong droop
move from  $(q,p)$ to $(p',q')$.

\item[3.] (Fuse)
In Cases (A) and (B),
replace the $\EDIAG$ tile at $(p',p')$
by $\PDIAG$.

In \CaseC,
replace the $\ETILE$ tile at $(p',q')$
by $\PTILE$.
\end{enumerate}
We call $\bftrev$ the \defn{reverse transition operator}.
\end{definition}

The remainder of this section is devoted to proving
Theorem~\ref{T: transition}.
The key step is to show that
both operators are well-defined. 

\subsection{Consistency of the transition operator}
\label{SS: trans well defined}

We first prove that the transition operator $\bft$ is well-defined
and satisfies the weight identities in Theorem~\ref{T: transition}.
We divide the proof into several statements, explain how they give the
result, and then prove them individually.
Throughout this subsection, fix $D\in\InvBPD(z)$ and write
$\DBox=\double(D)$. By Proposition~\ref{sym-bij-prop},
$\DBox\in\kBPD(z)$, although $\DBox$ need not be reduced.

\begin{lemma}
\label{L: t well-defined}
Let $D_1$, $D_2$, and $D_3$ be the pipe dreams obtained after the Cut,
Undroop, and Fuse steps of Definition~\ref{D: Transition operator}, respectively.
Then the following properties hold:
\begin{enumerate}
\item[(1)] In Cases (A) and (B), position $(p',p')$ of $D$ has a $\PDIAG$ tile
involving pipe $p$ and pipe $p'$.

In \CaseC, position $(p',q')$ of $D$ has a $\PTILE$ tile involving pipe $p$
and pipe $q'$.

\item[(2)] The pipe dream $D_1$ is reduced.

\item[(3)] We have $D_1\in\InvPD(\eta(z))$.

\item[(4)] The undroop in Definition~\ref{D: Transition operator} acts non-trivially
on $D_1$ and the result $D_2\in\InvPD(\eta(z))$.

\item[(5)] The tile at position $(q,p)$ of $D_2$ is  $\CDIAG$, $\CTILE$,
$\EDIAG$, or $\ETILE$, and $D_2$ has no $\ETILE$ or $\EDIAG$ tiles elsewhere.

The pipe of $D_2$ entering $(q,p)$ from the bottom is pipe $p$.

\item[(6)] If a Fuse step is performed, then $D_3\in\InvBPD(v)$ for some
$v\in\hat\Phi^{-}(\eta(z),p)$.
\end{enumerate}
\end{lemma}

We first explain how the lemma gives the desired result.

\begin{proposition}
\label{P: forward transition well defined}
The transition operator $\bft$ is a well-defined map from
$\InvBPD(z)$ to $\calT(z)$ and satisfies the weight identities
in Theorem~\ref{T: transition}.
\end{proposition}

\begin{proof}
Statements~(1), (4), and~(5) of Lemma~\ref{L: t well-defined} show
that the Cut, Undroop, and Fuse steps are valid operations.
By statements~(3) and~(4), both $D_1$ and $D_2$ belong to
$\InvPD(\eta(z))$.
If no fuse is performed, part~(5) shows that $D_2$ is bumpless,
so $D_3=D_2\in\InvBPD(\eta(z))$.
Otherwise, part~(6) gives $D_3\in\InvBPD(v)$ for some
$v\in\hat\Phi^{-}(\eta(z),p)$.
Thus, $\bft(D)\in\calT(z)$.

It remains to compare weights.
Throughout the algorithm, the only
position where a blank tile can appear
or disappear is $(q,p)$
during the Undroop step. 
If that happens, $(q,p)$ changes from a $\BDIAG$ or $\BTILE$ tile
in $D$ to a $\CDIAG$ or $\CTILE$ tile in $D_3$.
Hence $\blank(D)=\blank(D_3)\sqcup\{(q,p)\}$ and
$\weight(D)=x_{(p,q)}\weight(D_3)$.
This is the case when the fuse step does nothing, so  $D_3\in\InvBPD(\eta(z))$ 
by the arguments above. 
Conversely, if $D_3\notin\InvBPD(\eta(z))$,
then we know a fuse step happens,
so position $(q,p)$ in $D_2$ contains $\ETILE$ or
$\EDIAG$. 
Then no blank tile appears or disappears,
so $\blank(D) = \blank(D_3)$.
\end{proof}

We now prove the six statements separately.

\begin{proof}[Proof of Lemma~\ref{L: t well-defined}, part~~(1)]
First suppose that we are in Case (A) or (B).
Since $p'>q$, the position $(p',p')$ lies in the fixed region of
$\DBox$. We have $z(p')=p<p'$ and $z^{-1}(p')=p<p'$.
Thus, Corollary~\ref{C: BPD crossing fixed region} shows that this
position is a $\PTILE$ tile in $\DBox$ involving pipe $p$ and pipe $p'$.
Consequently, it is a $\PDIAG$ tile in $D$.

Now suppose that we are in \CaseC.
The position $(p',q')$ also lies in the fixed region, since $p'>q$.
Moreover, $z(p')=p<q'$ and $z^{-1}(q')=q<p'$.
Corollary~\ref{C: BPD crossing fixed region} shows that this position
is a $\PTILE$ tile involving pipe $p$ and pipe $q'$.
\end{proof}

\begin{proof}[Proof of Lemma~\ref{L: t well-defined}, part~(2)]
Since $D\in\InvBPD(z)$, the pipe dream $D$ is reduced.
In Cases (A) and (B), the cut resolves a crossing on the diagonal.
The two pipes terminate there, so this step cannot introduce any
double crossing.

In \CaseC, consider cutting the $\PTILE$ tile at $(p',q')$ in $\DBox$.
Since $p'=z^{-1}(p)$ and $q'=z(q)$, this is the crossing in
Lemma~\ref{L: No double crossing after cut}.
After the cut, neither affected pipe crosses any other pipe twice.
Thus, these two pipes do not cross any pipe twice in $D_1$,
so $D_1$ is reduced.
\end{proof}

\begin{proof}[Proof of Lemma~\ref{L: t well-defined}, part~(3)]
Let $u=\pi(D)$ and $u_1=\pi(D_1)$.
Since $D\in\InvBPD(z)$, we have $u\in\cW(z)$.
We need to show that $u_1\in\cW(\eta(z))$.
By part~(1) and Lemma~\ref{L: local reconnection},
$u_1=t_{p,p'}u$ in Cases (A) and (B), while $u_1=t_{pq'}u$ in \CaseC.
Since $D$ and $D_1$ are reduced and the cut removes one crossing,
we also have $u_1\lessdot u$.

In Cases (A) and (B), pipes $p$ and $p'$ meet at $(p',p')$.
The vertical pipe is pipe $p'$, and the horizontal pipe is pipe $p$,
so $p'$ appears immediately before $p$ in $u$.
In \CaseB, pipe $q$ passes vertically through row $p'$ in $\DBox$,
since $z(q)=q<p'$.
It therefore reaches the diagonal before $(p',p')$, so $q$ appears
before $p'$ in $u$.
Lemma~\ref{L: AB eta(w)} gives $u_1\in\cW(\eta(z))$.

In \CaseC, the relation $u_1=t_{pq'}u\lessdot u$ implies that
$q'$ appears before $p$ in $u$.
Lemma~\ref{L: C eta(w)} therefore gives $u_1\in\cW(\eta(z))$.
\end{proof}

\begin{proof}[Proof of Lemma~\ref{L: t well-defined}, part~(4)]
Set $c=p'$ in Cases (A) and (B), and $c=q'$ in \CaseC.
Let $R=[q,p']\times[p,c]$.
We show that $D_1$ is the output of a non-trivial strong droop
from $(q,p)$ to $(p',c)$ in Definition~\ref{extended-droop-def}(a).
Only positions in $\LT_n$ are considered in $R$.

Corollary~\ref{C: j or blank}, Lemma~\ref{L: BPD fixed bottom},
and the Cut step show that the tiles in $D_1$ at the four corners of $R$ are described by the following table:
\[
\begingroup
\renewcommand{\arraystretch}{1.5}
\begin{array}{c|c|c|c}
\text{position} & \text{\CaseA} & \text{Case (B)} & \text{Case (C)} \\
\hline
(q,p) & \BDIAG\text{ or }\JDIAG & \BTILE\text{ or }\JTILE
& \BTILE\text{ or }\JTILE \\
(p',p) & \CTILE & \CTILE & \CTILE \\
(p',c) & \EDIAG & \EDIAG & \ETILE \\
(q,c) & \text{outside }\LT_n & \text{outside }\LT_n & \CTILE
\end{array}
\endgroup
\]
In \CaseB, $(q,q)$ is also a $\CDIAG$ tile as $z(q)=q$.
We conclude that the tiles of $D_1$ match the outputs of the transformations in
Definition~\ref{extended-droop-def}(a)
when $(a,b) =(q,p)$ and $(i,j) = (p',c)$.

We next check that $D_1$ has no bend tiles elsewhere in $R$.
For $q<r<p'$, Lemma~\ref{L: on w} gives $z(r)<z(p')=p$.
By Lemma~\ref{L: BPD fixed bottom}, the bend in row $r$ of $\DBox$
is   strictly to the left of $R$.
In row $p'$, the only bends of $D_1$ in $R$ are at the corners
$(p',p)$ and $(p',c)$.
In row $q$, all positions strictly to the right of $(q,p)$ lie in
the fixed region of $\DBox$.
The only additional bend is at $(q,q)$ in \CaseB, or at the corner
$(q,c)$ in \CaseC.
Thus, the only possible interior bend is the $\CDIAG$ tile at $(q,q)$
in \CaseB.

Applying the inverse tile transformations and reconnecting the straight
segments gives a triangular pipe dream $E$.
In \CaseB, these transformations replace the $\CDIAG$ at $(q,q)$
by a $\PDIAG$, and they create no other interior bends.
Thus, $E$ satisfies the corner and interior conditions of
Definition~\ref{extended-droop-def}(a) when $(a,b) =(q,p)$ and $(i,j) = (p',c)$, and applying the strong droop indexed by these positions to $E$ gives $D_1$.
By Definition~\ref{undroop-def}, the corresponding undroop acts non-trivially on $D_1$,
and Proposition~\ref{droops-preserve-prop} gives
$D_2=E\in\InvPD(\eta(z))$.
\end{proof}

\begin{proof}[Proof of Lemma~\ref{L: t well-defined}, part~(5)]
After the undroop described in the previous proof, the tile at $(q,p)$ in $D_2$ is $\CDIAG$, $\CTILE$, $\EDIAG$, or $\ETILE$.
The pipe of $D_2$ entering from the bottom of this 
tile belongs to pipe $p$, which is the pipe we undrooped.
The undroop removes the bump created by the Cut step and creates
no bump tiles other than possibly at $(q,p)$.
\end{proof}

\begin{proof}[Proof of Lemma~\ref{L: t well-defined}, part~(6)]
Suppose that a Fuse step is performed, and set $y=\eta(z)$.
By part~(5), pipe $p$ of $D_2$ enters the bump tile $\ETILE$ at $(q,p)$ from the bottom.
Let $a$ be the label of the other pipe.
This is the pipe in the original $\JTILE$ or $\JDIAG$ tile at position $(q,p)$
of $D$.
Lemma~\ref{L: BPD fixed bottom} says that
in $\DBox$, this pipe enters row $q$ in column $a$ and then travels
to column $p$. Thus, $a<p$.
The pipe turns upward at $(q,p)$, so $z(a)<q$.
Since $a \notin \{p,p',q,q'\}$,
we have $y(a)=z(a)<q=y(p)$.
Definition~\ref{D: tau} therefore gives $\tau_{ap}(y)\ne y$.

Let $u=\pi(D_2)$. By part~(4), $u\in\cW(y)$.
Pipe $a$ of $D_2$ exits the triangle before pipe $p$, so $a$ appears before
$p$ in $u$.
Take $j$ with $a<j<p$.
We claim $j$ cannot appear between $a$ and $p$.
\begin{itemize}
\item If $z(j)>q$, then pipe $j$ exits the triangle below row $q$ in $\DBox$, so $j$ appears after $p$ in $u$.
\item If $z(j)<q$, then pipe $j$ enters row $q$ from below in column $j$
in $\DBox$.
This part of the pipe is unchanged by the cut and undroop
since these steps only affect tiles weakly to the right
of column $p$.
Thus, pipe $j$ enters row $q$ from below
in column $j$ in $D_2$.
Then pipe $a$ crosses pipe $j$ while traveling from column $a$ to column
$p$.
Since $D_2$ is reduced, $j$ appears before $a$ in $u$.
\item The case $z(j)=q$ is impossible, since $z(q)=q'>p$.
\end{itemize}
This proves our claim, so $u\lessdot t_{ap}u$.
Lemma~\ref{L: atom cover transition---missing} gives
$t_{ap}u\in\cW(v)$ and $\ellinv(v)=\ellinv(y)+1$ for
$v=\tau_{ap}(y)$, so $v\in\hat\Phi^{-}(y,p)$.
By Lemma~\ref{L: local reconnection}, the Fuse step turns $D_2$ into a pipe dream
with reading permutation $t_{ap}u$ and one more crossing than $D_2$.
As $u\lessdot t_{ap}u$, the resulting pipe dream $D_3$ is reduced
and thus in $\InvBPD(v)$.
\end{proof}

\subsection{Consistency of the reverse transition operator}
\label{SS: reverse trans well defined}

We now prove that the reverse transition operator is well-defined.
As in the forward direction, we divide the proof into several statements.

\begin{lemma}
\label{L: t' well-defined}
Take $D\in\calT(z)$, and let $D_1$, $D_2$, and $D_3$ be the
pipe dreams obtained after the Cut, Droop, and Fuse steps in Definition~\ref{D: Reverse Transition operator}, respectively.
Then the following statements hold:
\begin{enumerate}
\item[(1)] If $D\in\InvBPD(\eta(z))$, then position $(q,p)$ of $D$ has a $\CDIAG$ or
$\CTILE$ tile, and the Cut step does nothing.

If $D\in\InvBPD(v)$ for $v\in\hat\Phi^{-}(\eta(z),p)$, then
position $(q,p)$ of $D$ has a $\PDIAG$ or $\PTILE$ tile.

In either case, $D_1\in\InvPD(\eta(z))$, and pipe $p$ of $D_1$ enters
$(q,p)$ from the bottom.

\item[(2)] The strong droop in Definition~\ref{D: Reverse Transition operator}
acts non-trivially on $D_1$, and $D_2\in\InvPD(\eta(z))$.

The only bump tile of $D_2$ is at $(p',p')$ in Cases (A) and (B), or at
$(p',q')$ in \CaseC.

\item[(3)] In Cases (A) and (B), position $(p',p')$ of $D_2$ has a $\EDIAG$ tile
involving pipe $p$ and pipe $p'$.

Replacing this tile by $\PDIAG$ gives an element of $\InvBPD(z)$.

\item[(4)] In \CaseC, position $(p',q')$ of $D_2$ has a $\ETILE$ tile
involving pipe $p$ and pipe $q'$.

Replacing this tile by $\PTILE$ gives an element of $\InvBPD(z)$.
\end{enumerate}
\end{lemma}

As in the previous section, we explain the main result which uses this lemma before proving it.

\begin{proposition}
\label{P: reverse transition well defined}
The reverse transition operator $\bftrev$ is a well-defined map
 $\calT(z)\to \InvBPD(z)$.
\end{proposition}

\begin{proof}
Part~(1) of Lemma~\ref{L: t' well-defined} shows that the  Cut
step is valid and gives an element of $\InvPD(\eta(z))$.
Statement~(2) shows that the Droop step is valid and also gives an
element of $\InvPD(\eta(z))$.
Finally, parts~(3) and~(4) show that the Fuse step is valid and
gives an element of $\InvBPD(z)$.
\end{proof}

In the rest of this subsection, set $y=\eta(z)$ and choose
$D\in\InvBPD(v)$, where $v=y$ or $v\in\hat\Phi^{-}(y,p)$.
Again write $\DBox=\double(D)$.
In the case when $v\in\hat\Phi^{-}(y,p)$, set $a=v(q)$.
Lemma~\ref{L: negative transition at qp} implies that 
\[a<p, \quad y(a)<q,
\quad v^{-1}(p)<q, \quand v=\tau_{a,p}(y).\]
In each case, $v(r)=y(r)$ for every $r>q$.
By~\eqref{EQ: y increasing after}, we therefore have
$v(q+1)<\cdots<v(n)$.
Notice that if $v=1$ then
all of the desired claims are immediate. 
Thus, we assume $v$ is not the identity permutation. 

\begin{proof}[Proof of Lemma~\ref{L: t' well-defined}, part~(1)]
First suppose that $v=y$.
We have $y(q)=p$ by Definition~\ref{D: eta(w)}.
Recall that we assume $v \neq 1$. 
Then Lemma~\ref{L: inv transition order} shows that $(q,p)$
lies in the fixed region of $\DBox$.
By Lemma~\ref{L: BPD fixed bottom}, this position is a $\CTILE$
tile containing pipe $p$.
Thus, it is a $\CTILE$ or $\CDIAG$ tile in $D$, according to whether
$p<q$ or $p=q$.
The Cut step does nothing, so $D_1=D\in\InvPD(y)$.

Now suppose that $v\in\hat\Phi^{-}(y,p)$.
Since $v(q+1)<\cdots<v(n)$ and $v(q)=a$, every element of $D(v)$
has first coordinate at most $q$, and those with first coordinate $q$
have second coordinate less than $a$.
Thus, every $(q,j)$ with $j\geq a$ lies in the fixed region of $\DBox$.
Since $a<p$ and $v^{-1}(p)<q$, Corollary~\ref{C: BPD crossing fixed region}
shows that $(q,p)$ is a crossing of pipe $a$ and pipe $p$.
It is a $\PTILE$ tile in $D$ when $p<q$, or a $\PDIAG$ tile when $p=q$.
Let $u=\pi(D)\in\cW(v)$.
Since $D$ is reduced and pipes $a<p$ and $p$ cross at $(q,p)$,
the entry $p$ appears before $a$ in $u$.

We check the additional condition in Lemma~\ref{L: reverse Fuse in Phi}
when $p=q$ and $a<y(a)<p$.
Set $b=y(a)$.
In this case, Definition~\ref{D: tau} replaces the cycle $(a,b)$
by $(a,p)$ and fixes $b$, so $v(b)=b<q$.
The position $(q,b)$ is in the fixed region, since $b>a=v(q)$.
Corollary~\ref{C: BPD crossing fixed region} shows that pipe $b$
crosses pipe $a$ there.
Since $D$ is reduced, $b$ appears before $a$ in $u$.
Lemma~\ref{L: reverse Fuse in Phi} now gives
$t_{a,p}u\in\cW(y)$ and $t_{a,p}u\lessdot u$.

By Lemma~\ref{L: local reconnection}, cutting the crossing gives a
reduced pipe dream with reading permutation $t_{a,p}u$.
Hence $D_1\in\InvPD(y)$.
The segment entering $(q,p)$ from the bottom is unchanged by the cut,
so it still belongs to pipe $p$.
\end{proof}

\begin{proof}[Proof of Lemma~\ref{L: t' well-defined}, statement~(2)]
Set $c=p'$ in Cases (A) and (B), and $c=q'$ in \CaseC.
Let $R=[q,p']\times[p,c]$.
We check the conditions for a strong droop from $(q,p)$ to $(p',c)$
in Definition~\ref{extended-droop-def}(a).
Only positions in $\LT_n$ are considered in $R$.

By statement~(1), Lemma~\ref{L: BPD fixed bottom},
and Corollary~\ref{C: BPD crossing fixed region},
the tiles of $D_1$ at the four corners of $R$ are described by the following table:
\[
\begingroup
\renewcommand{\arraystretch}{1.5}
\begin{array}{c|c|c|c}
\text{position} & \text{\CaseA} & \text{\CaseB} & \text{\CaseC} \\
\hline
(q,p) & \CDIAG\text{ or }\EDIAG
& \CTILE\text{ or }\ETILE
& \CTILE\text{ or }\ETILE \\
(p',p) & \VTILE & \VTILE & \VTILE \\
(p',c) & \CDIAG & \CDIAG & \CTILE \\
(q,c) & \text{outside }\LT_n
& \text{outside }\LT_n & \HTILE
\end{array}
\endgroup
\]
In \CaseB, $(q,q)$ is also a $\PDIAG$ tile, since
$v(q)<q$ and $v^{-1}(q)<q$.
We conclude that the tiles of $D_1$ match the inputs of the transformations in
Definition~\ref{extended-droop-def}(a) when $(a,b) =(q,p)$ and $(i,j) = (p',c)$, including the diagonal
transformation when $p<q<c$.

We next check that $D_1$ has no bend tiles in the interior of $R$.
For $q<r<p'$, we have $v(r)=y(r)=z(r)<p$ by
Lemma~\ref{L: negative transition at qp}, Definition~\ref{D: eta(w)},
and Lemma~\ref{L: on w}.
Thus, Lemma~\ref{L: BPD fixed bottom} places the bend in row $r$
strictly to the left of $R$.
In row $p'$, we have $v(p')=c$ and $v^{-1}(j)<p'$ for every $j<c$,
since $v(q+1)<\cdots<v(n)$.
By Lemma~\ref{L: BPD fixed bottom}, the tile $(p',j)$ of $D_1$
is therefore a $\VTILE$ containing pipe $j$ for every $j<c$.
The only bend in this row inside $R$ is at $(p',c)$.
In row $q$, every position of $R$ lies in the fixed region of $\DBox$,
and $v(q)\leq p$.
By Lemma~\ref{L: BPD fixed bottom}, there is no bend to the right
of column $p$; the cut can create a bend only at $(q,p)$.
Hence $D_1$ has no interior bend tiles in $R$.

Therefore, the strong droop move indexed by $(a,b) =(q,p)$ and $(i,j) = (p',c)$ is non-trivial
when applied to $D_1$.
By Proposition~\ref{droops-preserve-prop},
this results in a pipe dream $D_2\in\InvPD(\eta(z))$.
The transformation removes the possible bump at $(q,p)$
and creates a bump at $(p',c)$ but nowhere else.
\end{proof}

\begin{proof}[Proof of Lemma~\ref{L: t' well-defined}, statement~(3)]
Assume that we are in Case (A) or (B).
By statement~(2), position $(p',p')$ had a $\CDIAG$ tile containing pipe $p'$
before the droop.
The droop turns this position into a $\EDIAG$ tile with pipe $p$ entering from the left
and pipe $p'$ entering from the bottom.
Thus, in $u_2=\pi(D_2)\in\cW(y)$, the entry $p$ appears immediately
before $p'$.
Since $p<p'$, we have $u_2\lessdot t_{p,p'}u_2$.
Definition~\ref{D: tau} gives $\tau_{p,p'}(y)=z\ne y$.
Lemma~\ref{L: atom cover transition---missing} therefore gives
$t_{p,p'}u_2\in\cW(z)$.
By Lemma~\ref{L: local reconnection}, replacing the $\EDIAG$ by a
$\PDIAG$ gives a reduced pipe dream with reading permutation $t_{p,p'}u_2$. This pipe dream is in $\InvBPD(z)$
since the replacement removes the only bump tile.
\end{proof}

\begin{proof}[Proof of Lemma~\ref{L: t' well-defined}, statement~(4)]
Assume that we are in \CaseC.
By statement~(2), position $(p',q')$ had a $\CTILE$ tile containing pipe $q'$
before the droop.
The droop turns it into an $\ETILE$ tile with pipe $p$ entering horizontally
and pipe $q'$ entering vertically.
Let $u_2=\pi(D_2)\in\cW(y)$.
Pipe $p$ of $D_2$ exits the triangle in row $q$, while pipe $q'$ exits in
row $p'$.
Thus, $p$ appears before $q'$ in $u_2$.

Take $j$ with $p<j<q'$.
We claim that $j$ cannot appear between $p$ and $q'$
in $u_2$.
As checked in the proof of statement~(2), position $(p',j)$ had a $\VTILE$ tile containing pipe $j$ before the droop.
The droop makes pipe $p$ cross this pipe while traveling from column
$p$ to column $q'$ in row $p'$.
Since $D_2$ is reduced and $p<j$, the entry $j$ appears before $p$
in $u_2$.
Hence $u_2\lessdot t_{p,q'}u_2$.

Definition~\ref{D: tau} asserts that $\tau_{p,q'}(y)=z\ne y$.
Lemma~\ref{L: atom cover transition---missing} therefore gives
$t_{p,q'}u_2\in\cW(z)$.
By Lemma~\ref{L: local reconnection}, replacing the $\ETILE$ tile at position $(p',q')$ by a
$\PTILE$ tile gives a reduced pipe dream with  reading permutation $t_{p,q'}u_2$.
As the replacement removes the only bump tile, this pipe dream belongs to $\InvBPD(z)$.
\end{proof}

\subsection{Proofs of Theorems~\ref{main-thm} and \ref{T: transition}}

Let us finally indicate how the preceding results let us prove Theorems~\ref{main-thm} and \ref{T: transition}.

\begin{proof}[Proof of Theorem~\ref{T: transition}]
Propositions~\ref{P: forward transition well defined} and \ref{P: reverse transition well defined} show that  $\bft 
$ and $\bftrev$ give well-defined maps $\calT(z) \leftrightarrow \InvBPD(z)$ satisfying the weight relation \eqref{wt-t-eq}.
Each step
in the definition of one operator
reverses the corresponding step
in the definition of the other, so
these maps are inverse bijections.
\end{proof}

\begin{proof}[Proof of Theorem~\ref{main-thm}]
For each $z \in \Invol$ let $f_z = \sum_{D \in \InvBPD(z)}\weight(D)$.
If $z=1$ then $\InvBPD(z)$ has only one element 
and  $f_1=1 = \iS_1$.
If $z\neq 1$ then Theorem~\ref{T: transition}
implies that 
$
f_z
=
x_{(p,q)}f_{\eta(z)}
+
\sum_{v\in\hat\Phi^{-}(\eta(z),p)}
f_v.
$
Since Remark~\ref{rec-rem} shows that $\iS_z$ is  determined by the same recurrence \eqref{EQ: iS Tran}, we have $f_z=\iS_z$ for all $z \in \Invol$.
\end{proof}

\section{Vexillary involutions}\label{vex-sect}

We can give a more refined description of the bumpless involution pipe dreams associated
to vexillary involutions $z \in \Ivex\subseteq \Invol\subseteq S_n$, as we explain in this section.
This leads to some applications outside Schubert calculus in Section~\ref{spp-sect}, and will prepare the way for our proof of Theorem~\ref{fkGO-conj2} in the appendix.
Several of the main results here have analogues in \cite[\S7]{Weigandt}
and our proofs are adapted from methods there.

\subsection{Vexillary simplifications}

Recall that $\Cyc(z)$ denotes the set of 2-cycles of $z \in \Invol$ and 
that
$\cyc(z) = |\Cyc(z)|
= |\{ i \in [n] : z(i) > i\}|
.$
We begin with some results illustrating how the sets $\InvBPD(z)$ and $\kInvBPD(z)$ simplify when $z \in \Ivex$.

\begin{theorem}\label{vex-prop1}
Suppose $z \in \Ivex $. Then
$\InvBPD(z) = \kInvBPD(z)$ and 
\be\label{cyc-bdiag-eq}
\cyc(z) 
= \left|\left\{ (i,j) \in \LT_n : D_{ij} =\PDIAG\right\}\right|
= \left|\left\{ (i,j) \in \LT_n : D_{ij} =\BDIAG\right\}\right|
\quad\text{for all }D\in \InvBPD(z).
\ee
\end{theorem}

\begin{proof}
Let $D \in \kInvBPD(z)$ and $E = \double(D)$. 
We must show that $D\in \InvBPD(z)$.
Weigandt \cite[Lem.~7.2]{Weigandt}
proves that if $w \in S_n$ is vexillary then $\BPD(w)= \kBPD(w)$.
Therefore $E \in \kBPD(z) = \BPD(z)$
by Proposition~\ref{sym-bij-prop} so $D$ and $E$ are both reduced.

It follows that $\ell(z)$ is the number of $\PTILE$ tiles in $E$.
Let $p$ and $q$ be the respective numbers of $\PTILE$ and $\PDIAG$ tiles in $D$,
so that $\ell(\pi(D)) = p+q$ and $\ell(z) = 2p+q$. 
Since $E$ is symmetric and reduced, whenever $a,b\in[n]$ have $a<b=z(a)$, the pipes labeled $a$ and $b$ in $E$ must cross at a diagonal $\PTILE$ tile and nowhere else.
Conversely, the labels of the pipes in any diagonal $\PTILE$ tile in $E$ must form a 2-cycle in $z$. 
As $q$ also counts the diagonal $\PTILE$ tiles in $E$, we conclude that
$q =\cyc(z)$ and so 
\[ \ell(\pi(D)) =p+q = \tfrac{(2p+q) + q}{2} = \tfrac{\ell(z) +\cyc(z)}{2} = \ellinv(z).\]
 As $\pi(D) \circ \pi(D)^{-1} = \kpi(D)\circ \kpi(D)^{-1} =z$
we conclude that $\pi(D) \in \cW(z)$ so $D \in \InvBPD(z)$ as needed.

Our observation that $q = \cyc(z)$ gives the first equality in \eqref{cyc-bdiag-eq}.
The second equality holds by \eqref{Mdiag-eq}.
\end{proof}

\begin{corollary}\label{vex-double-cor}
If $z \in \Ivex$ then the operation
$\double $ from \eqref{double-def} is a bijection \[\InvBPD(z) \to \{ D \in \BPD(z) : D = \transpose(D)\}.\]
\end{corollary}

\begin{proof}
This follows from Proposition~\ref{sym-bij-prop}
as  
 $\InvBPD(z) = \kInvBPD(z)$
 and
  $\BPD(z) = \kBPD(z)$.
\end{proof}

\begin{corollary}\label{ivex-is-cor}
If $z \in \Ivex$ then $2^{\cyc(z)}\iS_z =\sum_{D \in \InvBPD(z)} \prod_{(i,j) \in \blank(D)} (x_i + x_j)$.
\end{corollary}

\begin{proof}
This follows from Theorem~\ref{main-thm} given \eqref{cyc-bdiag-eq}.
\end{proof}

We say   two triangular pipe dreams $D$ and $E$ are connected by a
\defn{basic droop move} if $D\neq \droop^{ab}_{ij}(D) = E$
but $D$ and $E$ both have the same number of diagonal $\BDIAG$ tiles.
Basic droop moves are the simple droop moves corresponding to the pictures in  \eqref{simple-eq} but not in \eqref{simple-eq2}.

\begin{proposition}\label{bbs-prop}
Suppose $z \in \Ivex $.
If two distinct elements of $ \InvBPD(z)$ are connected by a droop move, then this move must be basic.
Hence, the set $ \InvBPD(z)$ consists of all triangular pipe dreams obtained by applying sequences of basic droop moves to $\invRothe(z)$.
\end{proposition}

\begin{proof}
We know from \eqref{cyc-bdiag-eq} that 
two elements of $ \InvBPD(z)$ cannot be connected by a simple droop move that is not basic. 
Suppose instead that $(a,b),(i,j) \in \LT_n$ are such that $\complexDroop^{ab}_{ij}(D) \neq D$ as in the picture \eqref{complex-eq}.
Consider the pipes labeled $a+1,a+2,\dots,n$ in $D$.
  Each of these pipes must reach some diagonal position in $D$.
  As $D_{ij}=\BTILE$, this is only possible if at least two of  these pipes meet 
at a $\PDIAG$ tile in $D$, in which case their labels form a 2-cycle in $z$ by Theorem~\ref{cW-lem}.
But this would create a $2143$-pattern in $z$, contradicting
  the assumption that $z$ is vexillary,
  as the labels of the pipes that meet at position $(a,a)$ in $D$ form another 2-cycle in $z$.

The contradiction proves the first claim in the proposition.
The second claim follows by Theorem~\ref{droop-thm}.
\end{proof}

\subsection{Flagged plane partitions}\label{flagged-sect}

We now explain how the elements of $\InvBPD(z)$ for $z\in \Ivex$ 
may be identified with certain \defn{shifted plane partitions}.
Define a subset $S \subset \PP\times \PP$ to be \defn{northwest aligned} if whenever $(i,j) \in \PP\times \PP$ has $(i,j) \notin S$ then also
$(i+1,j+1) \notin S$. 
Construct the \defn{northwest alignment} of a finite subset $S \subset \PP\times \PP$ 
as the unique northwest aligned set $\alignNW(S) \subset \PP\times \PP$ that has
the same number of positions as $S$ on each diagonal.

Suppose $v \in S_n$ is vexillary.
Weigandt \cite[Lem.~7.3]{Weigandt} proves  that there is 
 a unique partition $\mu^{(v)}$ with
  \be\D_{\mu^{(v)}} = \alignNW( D(v)).\ee
This partition can be computed by sorting the Lehmer code 
of $v$.
Let $\lambda^{(v)}$ be the minimal partition under containment of Young diagrams
that has $D(v) \subseteq \D_{\lambda^{(v)}}$.
Finally, suppose $\mu^{(v)}$ has $r$ nonzero parts
and let \be\label{bff-eq}\bff^{(v)} = (f_1^{(v)}, f_2^{(v)}, f_3^{(v)},\dots,f_r^{(v)})\ee
be the integer sequence such that $f_i^{(v)}$ is the maximal row index
of any position in $ \D_{\lambda^{(v)}}$ on the same diagonal as 
the box $(i,\mu_i^{(v)}) \in \D_{\mu^{(v)}}$.
One always has $i \leq f_i^{(v)} $.

When $\nu $ is any partition, let $\half(\nu)$
be the strict partition with parts
$
\half(\nu)_i = \max\{ 0, \nu_i - (i-1)\}.
$
The shifted diagram of a strict partition $\lambda$ is the set 
$\SD_\lambda = \{ (i,j+i-1) : (i,j) \in \D_\lambda \} 
$
and we have
\[
\SD_{\half(\nu)} = \{(i,j) \in \D_\nu : i\leq j\}.\]
 Suppose $v \in \Ivex$ is a vexillary involution. Then
  $\mu^{(v)}$ is symmetric and the strict partition
  \be
  \hat\mu^{(v)} := \half(\mu^{(v)})
  \ee
   is the unique one with 
$
\left\{ (j,i) : (i,j) \in \SD_{\hat\mu^{(v)}}\right\} = \alignNW( \iD(v)).
$
In this case $\lambda^{(v)}$ is also symmetric and 
\be\hat\lambda^{(v)} := \half(\lambda^{(v)})\ee is the minimal strict partition
such that the transpose of $\iD(v)$ is contained in $\SD_{\hat\lambda^{(v)}}$.

\begin{example}
Suppose $v=526413 = (1,5)(3,6) \in S_6$.
Then $v=v^{-1}$ is vexillary with
\[
D(v) = 
\left\{\begin{smallmatrix}
\square & \square & \square & \square & 1 & \cdot \\
\square & 1 &  \cdot& \cdot & \cdot & \cdot  \\ 
\square &  \cdot& \square& \square& \cdot & 1 \\
\square &  \cdot& \square & 1 & \cdot & \cdot \\
1 & \cdot  & \cdot & \cdot & \cdot & \cdot \\ 
\cdot & \cdot & 1 & \cdot & \cdot & \cdot
\end{smallmatrix}\right\}
\quad\text{so}\quad
\mu^{(v)} = (4,3,2,1) =
\ytabc{0.25cm}{ \ &\  &\  & \ \\ \ & \ & \ \\ \ & \  \\ \ }
\subset
\lambda^{(v)} = (4,4,4,3) =
\ytabc{0.25cm}{ \ &\  &\  & \ \\ \ & \ & \ & \ \\ \ & \  & \ & \  \\ \ & \ & \ }
\]
and $\bff^{(v)} = (1,3,4,4)$.
In this case $v=v^{-1}$ so
$\hat\mu^{(v)} = (4,2)$ and $\hat\lambda^{(v)} = (4,3,2)$. 
\end{example}

Assume $\lambda$ is a partition.
A \defn{plane partition} (respectively, \defn{shifted plane partition}) of shape $\lambda$ is a filling of
 $\D_\lambda$ 
 (respectively, $\SD_\lambda$)
by nonnegative 
integers that weakly increase along rows and columns.
When dealing with shifted plane partitions of shape $\lambda$,
we require $\lambda$ to be a strict partition.

For each vexillary $v \in S_n$ let $\mathsf{PP}(v)$ be the set of
plane partitions of shape $\mu^{(v)}$ 
whose entries in row $i$ are at most $ f^{(v)}_i - i$ for all $i$.
Then for each $z \in \Ivex$
let $\SPP(z)$ be the set of 
shifted plane partitions of shape $\hat\mu^{(z)}$ 
whose entries in row $i$ are all at most $ f^{(z)}_i - i$. For example, we have
\be\label{vex-spp-eq}
\SPP(526413) = 
\left\{
\ytabc{0.4cm}{ 0 & 0  & 0 & 0 \\ \none & 0  & 0 }, \
\ytabc{0.4cm}{ 0 & 0  & 0 & 0 \\ \none & 0  & 1 }, \
\ytabc{0.4cm}{ 0 & 0  & 0 & 0 \\ \none & 1  & 1 }
\right\}.
\ee
Given a (shifted or unshifted) plane partition $P$, let
 $\displace(P) = \{ (i+P_{ij}, j+P_{ij}): (i,j) \in P\}\subset \PP\times\PP$.
 
\begin{lemma}\label{displace-lem}
Fix $z \in \Ivex$ and 
suppose $P$ is either a plane partition of shape $\mu^{(z)}$ or a shifted plane partition of shape $\hat\mu^{(z)}$.
Then $P$ respectively belongs to $\mathsf{PP}(z)$ 
or $\SPP(z)$
if and only if $\displace(P)\subseteq \D_{\lambda^{(z)}}$.
\end{lemma}

\begin{proof}
The entry bounds that define 
$\mathsf{PP}(z)$ 
and $\SPP(z)$ only need to be imposed on the rightmost positions in each (weakly increasing) row.
Likewise, we have $\displace(P)\subseteq \D_{\lambda^{(z)}}$ if and only if $(i+P_{ij},j+P_{ij}) \in \D_{\lambda^{(z)}}$ for the rightmost positions $(i,j)$ in each row.
But for such positions, the definition of $\bff^{(z)}$   implies that we have the bound $P_{ij} \leq f^{(z)}_i-i$ 
precisely when $(i+P_{ij},j+P_{ij}) \in \D_{\lambda^{(z)}}$.
\end{proof}

If $\mu = \mu^\top$ is a symmetric partition 
then there is a doubling operation $\double$ 
that bijectively maps shifted plane partitions of shape $ \half(\mu)$ to symmetric plane partitions of shape $\mu$. This operation acts as 
\[
\double :\ytabc{0.5cm}{ a & b  & c & d \\ \none & e  & f  & g}
\mapsto 
\ytabc{0.5cm}{a & b  & c & d \\ b & e  & f & g \\ c & f \\ d & g }
\]
when $\mu = (4,4,2,2)$ and has a similar definition for arbitrary $\mu=\mu^\top$. 

\begin{lemma}\label{vex-double-lem}
If $z \in \Ivex$ then $\double$ is a bijection $\SPP(z) \to \{ \text{symmetric elements of } \mathsf{PP}(z) \}$.
\end{lemma}

\begin{proof}
The operation $\double$ is a bijection from $\SPP(z)$ to a set of symmetric plane partitions of shape $\mu^{(z)}$.
Since $\displace(\double(P)) = \displace(P) \cup \displace(P)^\top$ and since $\lambda^{(z)}$ is a symmetric partition, Lemma~\ref{displace-lem}
makes it clear that this set consists of all 
symmetric elements of $\mathsf{PP}(z)$.
\end{proof}

\begin{theorem}\label{bij-thm}
Suppose  $z \in\Ivex $. Then there is a unique bijection
\[
\Upsilon : \SPP(z) \to \InvBPD(z) 
\]
such that  if $\Upsilon(T) = D$ then
$ \blank(D) = 
 \left\{(b+T_{ab}, a+T_{ab}) : (a,b) \in \SD_{\hat\mu^{(z)}}\right\}.
$
\end{theorem}

The parameters $a$ and $b$ trade places in the last expression because shifted diagrams are upper triangular while involution pipe dreams are lower triangular.
This bijection sends the elements in \eqref{vex-spp-eq} to the corresponding pipe dreams in 
Example~\ref{fkGO-ex}.

\begin{proof}
Let $\Phi$ act on elements of $\mathsf{PP}(z)$ by adding $i$ to every entry in row $i$.
This operator defines a bijection from $\mathsf{PP}(z)$ to the set of \defn{flagged semistandard tableaux} $\mathsf{FSYT}(z)$ defined in \cite[\S7.1]{Weigandt}. 

Composing $\Phi$ with the map $\gamma$ in \cite[Thm~7.7]{Weigandt} 
gives a bijection $\gamma \circ \Phi: \mathsf{PP}(z) \to \BPD(z)$
such that if $\gamma\circ \Phi(T) = D$
then 
$ \blank(D) = 
 \left\{(a+T_{ab}, b+T_{ab}) : (a,b) \in \D_{\mu^{(z)}}\right\}.
$
Notice that we can recover $T$ from $\blank(D)$ and that $D$ is symmetric if and only if $\blank(D)$ is invariant under transpose.

Thus $\gamma \circ \Phi$
restricts to a bijection from the symmetric elements of $ \mathsf{PP}(z)$ 
to the symmetric pipe dreams in $\BPD(z)$. The desired bijection
$\Upsilon$ is therefore obtained by composing $\gamma \circ \Phi$ with the map in Lemma~\ref{vex-double-lem} on the right and the inverse
of the map in Corollary~\ref{vex-double-cor} on the left.
\end{proof}

\begin{corollary}
\label{spp-ivex-cor}
  If $z\in \Ivex$ 
then
$ |\SPP(z)| = |\InvBPD(z)| = 2^{\cyc(z)} \iS_z(\tfrac{1}{2},\tfrac{1}{2},\tfrac{1}{2},\dots) .$
\end{corollary}

For general $z \in \Invol$ the size of $\InvBPD(z)$ is not a straightforward specialization of $\iS_z$ as \eqref{cyc-bdiag-eq} may fail.
\begin{proof}
The first equality is immediate
from Theorem~\ref{bij-thm} and the second holds by Corollary~\ref{ivex-is-cor}.
\end{proof}

If $z \in \Ivex$ then $\InvBPD(z)$ is also in bijection with the set of \defn{type B excited diagrams} of shape $\hat\lambda^{(z)}/\hat\mu^{(z)}$, which are  defined as follows. 
 Suppose $\mu\subseteq \lambda$ are strict partitions. An element $(i,j)$ of a set $S\subseteq \SD_\lambda$ is \defn{active} if $i\neq j$ and $(i,j+1),(i+1,j),(i+1,j+1) \in \SD_\lambda\setminus S$, or if $i=j$ and $(i,i+1),(i+1,i+1) \in \SD_\lambda\setminus S$. An \defn{excited move} 
 replaces an active cell $(i,j)\in S$ by $(i+1,j+1)$: 
\[
\ytabc{1.5em}{ *(gray) & \\ *(white) &  } \mapsto \ytabc{1.5em}{*(white) & \\ *(white)  & *(gray)  } \qquad \qquad 
\ytabc{1.5em}{   *(gray)  & \\  \none  & } \mapsto  \ytabc{1.5em}{ *(white)    & \\  \none  &  *(gray)}. 
\]
A \defn{type B excited diagram} of shape $\lambda/\mu$ is any subset of $\SD_\lambda$ obtained 
by applying a sequence of (zero or more) excited moves to the initial set $\SD_\mu$.

Let $\mathcal{E}_B(\lambda/\mu)$ be the set of all type B excited diagrams of $\lambda/\mu$. This set appears in formulas of equivariant Schubert classes of the isotropic Grassmannian \cite{GrahamKreimanB,IkedaNaruseB} and in the Naruse hook length formula to count standard shifted tableaux of skew shape \cite{KonvalinkaNaruseBD,NaruseOkada}.

\begin{example}\label{exc-ex} When $z=526413 \in I^{\mathsf{vex}}_6$, the excited diagrams of shape $\hat\lambda^{(z)}/\hat\mu^{(z)}=432/42$ are 
\begin{align*}
\ytabc{1.5em}{ *(gray) & *(gray) & *(gray) & *(gray) \\ \none & *(gray) & *(gray) & \\ \none & \none & & } & \qquad   \ytabc{1.5em}{ *(gray) & *(gray) & *(gray) & *(gray) \\ \none & *(gray) & & \\ \none & \none & & *(gray) }  \qquad \ytabc{1.5em}{ *(gray) & *(gray) & *(gray) & *(gray) \\ \none & & & \\ \none & \none & *(gray) & *(gray) }\\
\intertext{while the bumpless involution pipe dreams of $z$ are}
\begin{tikzpicture}[x=\bpdwidth,y=\bpdwidth,line cap=round,line join=round,baseline=(z.base)]
\node at (0,2.75) (z) {};
\draw[gray!70, line width=0.7pt](0,0)--(6,0);
\draw[gray!70, line width=0.7pt](0,0)--(0,6);
\draw[gray!70, line width=0.7pt](0,1)--(6,1);
\draw[gray!70, line width=0.7pt](1,0)--(1,6);
\draw[gray!70, line width=0.7pt](0,2)--(6,2);
\draw[gray!70, line width=0.7pt](2,0)--(2,6);
\draw[gray!70, line width=0.7pt](0,3)--(6,3);
\draw[gray!70, line width=0.7pt](3,0)--(3,6);
\draw[gray!70, line width=0.7pt](0,4)--(6,4);
\draw[gray!70, line width=0.7pt](4,0)--(4,6);
\draw[gray!70, line width=0.7pt](0,5)--(6,5);
\draw[gray!70, line width=0.7pt](5,0)--(5,6);
\draw[gray!70, line width=0.7pt](0,6)--(6,6);
\draw[gray!70, line width=0.7pt](6,0)--(6,6);
\draw[thick,rounded corners,color=blue] (3.5,2)--(3.5,3);
\draw[thick,rounded corners,color=blue] (0.5,1)--(0.5,1.5)--(1,1.5);
\draw[thick,rounded corners,color=blue] (2.5,0)--(2.5,0.5)--(3,0.5);
\draw[thick,rounded corners,color=blue] (2.5,3)--(2.5,2.5)--(2,2.5);
\draw[thick,rounded corners,color=blue] (2.5,3)--(2.5,4);
\draw[thick,rounded corners,color=blue] (1.5,2)--(1.5,2.5)--(2,2.5);
\draw[thick,rounded corners,color=blue] (0.5,0)--(0.5,1);
\draw[thick,rounded corners,color=blue] (1,1.5)--(2,1.5);
\draw[thick,rounded corners,color=blue] (1.5,1)--(1.5,2);
\draw[thick,rounded corners,color=blue] (2,1.5)--(3,1.5);
\draw[thick,rounded corners,color=blue] (3,1.5)--(4,1.5);
\draw[thick,rounded corners,color=blue] (3.5,1)--(3.5,2);
\draw[thick,rounded corners,color=blue] (4,1.5)--(5,1.5);
\draw[thick,rounded corners,color=blue] (4.5,1)--(4.5,2);
\draw[thick,rounded corners,color=blue] (1.5,0)--(1.5,1);
\draw[thick,rounded corners,color=blue] (3,0.5)--(4,0.5);
\draw[thick,rounded corners,color=blue] (3.5,0)--(3.5,1);
\draw[thick,rounded corners,color=blue] (4,0.5)--(5,0.5);
\draw[thick,rounded corners,color=blue] (4.5,0)--(4.5,1);
\draw[thick,rounded corners,color=blue] (5,0.5)--(6,0.5);
\draw[thick,rounded corners,color=blue] (5.5,0)--(5.5,1);
\draw[white, line width=0.7pt,fill=white] (0,6)--(6,0)--(6,6)--(0,6);
\draw[black, line width=1.6pt] (0,6)--(6,0)--(0,0)--(0,6);
\end{tikzpicture}
& \quad
\begin{tikzpicture}[x=\bpdwidth,y=\bpdwidth,line cap=round,line join=round,baseline=(z.base)]
\node at (0,2.75) (z) {};
\draw[gray!70, line width=0.7pt](0,0)--(6,0);
\draw[gray!70, line width=0.7pt](0,0)--(0,6);
\draw[gray!70, line width=0.7pt](0,1)--(6,1);
\draw[gray!70, line width=0.7pt](1,0)--(1,6);
\draw[gray!70, line width=0.7pt](0,2)--(6,2);
\draw[gray!70, line width=0.7pt](2,0)--(2,6);
\draw[gray!70, line width=0.7pt](0,3)--(6,3);
\draw[gray!70, line width=0.7pt](3,0)--(3,6);
\draw[gray!70, line width=0.7pt](0,4)--(6,4);
\draw[gray!70, line width=0.7pt](4,0)--(4,6);
\draw[gray!70, line width=0.7pt](0,5)--(6,5);
\draw[gray!70, line width=0.7pt](5,0)--(5,6);
\draw[gray!70, line width=0.7pt](0,6)--(6,6);
\draw[gray!70, line width=0.7pt](6,0)--(6,6);
\draw[thick,rounded corners,color=blue] (3.5,2)--(3.5,3);
\draw[thick,rounded corners,color=blue] (0.5,1)--(0.5,1.5)--(1,1.5);
\draw[thick,rounded corners,color=blue] (2.5,0)--(2.5,0.5)--(3,0.5);
\draw[thick,rounded corners,color=blue] (2,3.5)--(3,3.5);
\draw[thick,rounded corners,color=blue] (1.5,3)--(1.5,3.5)--(2,3.5);
\draw[thick,rounded corners,color=blue] (0.5,0)--(0.5,1);
\draw[thick,rounded corners,color=blue] (1,1.5)--(2,1.5);
\draw[thick,rounded corners,color=blue] (1.5,1)--(1.5,2);
\draw[thick,rounded corners,color=blue] (2,1.5)--(3,1.5);
\draw[thick,rounded corners,color=blue] (3,1.5)--(4,1.5);
\draw[thick,rounded corners,color=blue] (3.5,1)--(3.5,2);
\draw[thick,rounded corners,color=blue] (4,1.5)--(5,1.5);
\draw[thick,rounded corners,color=blue] (4.5,1)--(4.5,2);
\draw[thick,rounded corners,color=blue] (1.5,0)--(1.5,1);
\draw[thick,rounded corners,color=blue] (1.5,2)--(1.5,3);
\draw[thick,rounded corners,color=blue] (3,0.5)--(4,0.5);
\draw[thick,rounded corners,color=blue] (3.5,0)--(3.5,1);
\draw[thick,rounded corners,color=blue] (4,0.5)--(5,0.5);
\draw[thick,rounded corners,color=blue] (4.5,0)--(4.5,1);
\draw[thick,rounded corners,color=blue] (5,0.5)--(6,0.5);
\draw[thick,rounded corners,color=blue] (5.5,0)--(5.5,1);
\draw[white, line width=0.7pt,fill=white] (0,6)--(6,0)--(6,6)--(0,6);
\draw[black, line width=1.6pt] (0,6)--(6,0)--(0,0)--(0,6);
\end{tikzpicture} \quad
\begin{tikzpicture}[x=\bpdwidth,y=\bpdwidth,line cap=round,line join=round,baseline=(z.base)]
\node at (0,2.75) (z) {};
\draw[gray!70, line width=0.7pt](0,0)--(6,0);
\draw[gray!70, line width=0.7pt](0,0)--(0,6);
\draw[gray!70, line width=0.7pt](0,1)--(6,1);
\draw[gray!70, line width=0.7pt](1,0)--(1,6);
\draw[gray!70, line width=0.7pt](0,2)--(6,2);
\draw[gray!70, line width=0.7pt](2,0)--(2,6);
\draw[gray!70, line width=0.7pt](0,3)--(6,3);
\draw[gray!70, line width=0.7pt](3,0)--(3,6);
\draw[gray!70, line width=0.7pt](0,4)--(6,4);
\draw[gray!70, line width=0.7pt](4,0)--(4,6);
\draw[gray!70, line width=0.7pt](0,5)--(6,5);
\draw[gray!70, line width=0.7pt](5,0)--(5,6);
\draw[gray!70, line width=0.7pt](0,6)--(6,6);
\draw[gray!70, line width=0.7pt](6,0)--(6,6);
\draw[thick,rounded corners,color=blue] (1.5,4)--(1.5,5);
\draw[thick,rounded corners,color=blue] (3.5,2)--(3.5,3);
\draw[thick,rounded corners,color=blue] (0.5,1)--(0.5,1.5)--(1,1.5);
\draw[thick,rounded corners,color=blue] (2.5,0)--(2.5,0.5)--(3,0.5);
\draw[thick,rounded corners,color=blue] (0.5,0)--(0.5,1);
\draw[thick,rounded corners,color=blue] (1,1.5)--(2,1.5);
\draw[thick,rounded corners,color=blue] (1.5,1)--(1.5,2);
\draw[thick,rounded corners,color=blue] (2,1.5)--(3,1.5);
\draw[thick,rounded corners,color=blue] (3,1.5)--(4,1.5);
\draw[thick,rounded corners,color=blue] (3.5,1)--(3.5,2);
\draw[thick,rounded corners,color=blue] (4,1.5)--(5,1.5);
\draw[thick,rounded corners,color=blue] (4.5,1)--(4.5,2);
\draw[thick,rounded corners,color=blue] (1.5,0)--(1.5,1);
\draw[thick,rounded corners,color=blue] (1.5,2)--(1.5,3);
\draw[thick,rounded corners,color=blue] (1.5,3)--(1.5,4);
\draw[thick,rounded corners,color=blue] (3,0.5)--(4,0.5);
\draw[thick,rounded corners,color=blue] (3.5,0)--(3.5,1);
\draw[thick,rounded corners,color=blue] (4,0.5)--(5,0.5);
\draw[thick,rounded corners,color=blue] (4.5,0)--(4.5,1);
\draw[thick,rounded corners,color=blue] (5,0.5)--(6,0.5);
\draw[thick,rounded corners,color=blue] (5.5,0)--(5.5,1);
\draw[white, line width=0.7pt,fill=white] (0,6)--(6,0)--(6,6)--(0,6);
\draw[black, line width=1.6pt] (0,6)--(6,0)--(0,0)--(0,6);
\end{tikzpicture}
.
\end{align*}

\end{example}

As illustrated in this example, 
there is a simple bijection $\mathcal{E}_B(\hat\lambda^{(z)}/\hat\mu^{(z)})\to \InvBPD(z)$    sending  $S\mapsto D$ where $\blank(D)=\{(j,i) : (i,j) \in S\}$. The existence of such a bijection for general $z \in \Ivex$
is justified by Proposition~\ref{bbs-prop}.
A similar correspondence between \defn{type A excited diagrams} and bumpless pipe dreams of vexillary permutations appears in \cite[\S 7.2]{Weigandt}.

\subsection{Applications to shifted plane partition counts and Ehrhart positivity}
\label{spp-sect}

Given a strict partition $\lambda$, let $\SPP_\lambda(k)$ denote the set of shifted plane partitions of shape $\lambda$
with all entries in $\{0,1,\dots,k\}$. 
In addition, given $k \in \{0,1,\dots,n\}$ and 
some  $w \in S_{n-k}$ let $1^k\times w \in S_n$ be the permutation that fixes each $i \in [k]$
while sending $i +k \mapsto w(i) +k$ for each $i \in [n-k]$.

\begin{theorem} \label{thm: shifted bijection inv bps and plane partitions}
Suppose $0\leq k\leq n$ and $z\in I_{n-k}$ is dominant of shape $\mu=\mu^\top$. Let $\lambda=\half(\mu)$.
Then \[\textstyle
 |\SPP_\lambda(k)|=
|\InvBPD(1^k\times z)| = 2^{\ell(\lambda)}\iS_{1^k\times z}(\tfrac{1}{2},\tfrac{1}{2},\tfrac{1}{2},\dots)=
 2^{\ell(\lambda)-|\lambda|} 
 \sum_{w \in \cW(z)}
\fkS_{1^k\times w^{-1}}(1,1,1,\dots),\]
where $\ell(\lambda)$ is the number of nonzero parts in $\lambda$ and $|\lambda|$ is the sum of these parts.
\end{theorem}

\begin{proof}
The permutation $1^k\times z$ is a vexillary involution by \cite[\S2.2.1]{Manivel},
and its length and number of 2-cycles are $|\mu|$ and $\ell(\lambda)$, respectively \cite[\S4]{HMP6}.
Moreover, we have $\mu^{(1^k\times z)} = \mu$ and $\hat\mu^{(1^k\times z)} = \lambda$ and $f^{(1^k\times z)}_i - i=k$ for all $i$.
The first two equalities therefore hold by Corollary~\ref{spp-ivex-cor}.

The operation $w \mapsto 1^k\times w$ is length-preserving, commutes with inverses, and satisfies
\[(1^k \times w) \circ (1^k \times w)^{-1} = (1^k \times w) \circ (1^k \times w^{-1}) = 1^k \times (w \circ w^{-1}).\]
Therefore $\cW(1^k \times z) = \{ 1^k \times w: w \in \cW(z)\}$
and if $w \in \cW(z)$ then $\fkS_{1^k\times w^{-1}}$ is homogeneous of degree $\ell(1^k\times w^{-1})=\ell(w) = \ellinv(z) = \frac{|\mu| + \ell(\lambda)}{2}$.
Since $\ell(\lambda) - \frac{|\mu| + \ell(\lambda)}{2} = \frac{\ell(\lambda) - |\mu| }{2}=\frac{\ell(\lambda) -(2 |\lambda| - \ell(\lambda)) }{2} =\ell(\lambda) -|\lambda|$, the last equality follows by Corollary~\ref{spp-ivex-cor} and \eqref{atom-eq}.
\end{proof}

Simple product formulas are known \cite{HopkinsLai,Okada,Proctor} for $|\SPP_\lambda(k)|$
when $\lambda$ is a \defn{shifted double staircase partition} of the form $\delta^{(p,q)} = (p,\dots,3,2,1) + (q,\dots,3,2,1)$ for integers $0\leq p\leq q$.
Notice that the reverse permutation $w_0 = n\cdots 321 \in S_n$ is dominant of shape $\mu^{(w_0)}=(n-1,\dots,3,2,1)$
and that the strict partition \[\hat\mu^{(w_0)} = (n-1,n-3,n-5,\dots) = \half((n-1,\dots,3,2,1))\] has the form $\delta^{(p-1,q)}$ for $p = \lceil n/2\rceil $ and $q= \lfloor n/2\rfloor$.
Thus, we can apply a result of Proctor \cite{Proctor} (which is generalized as 
\cite[Thm.~1.1]{HopkinsLai})
to the left side of Theorem~\ref{thm: shifted bijection inv bps and plane partitions}
to obtain the following identity.

\begin{corollary} \label{cor:plane partition product}
If $k \in \NN$ then 
\be\label{ppp-eq}
|\InvBPD(1^k\times n\cdots 321)|
=
\ds\prod_{i=1}^{\lceil n/2\rceil} \prod_{j=1}^{\lfloor n/2\rfloor} \frac{i+j+k-1}{i+j-1}
 = 
|\SPP_{(n-1,n-3,n-5,\dots)}(k)|.
\ee
\end{corollary}

\begin{remark} Proctor \cite{Proctor} shows that the plane partitions of a shifted trapezoid shape 
are counted by MacMahon's box formula,
which is a certain triple product expression.
By Theorem~\ref{thm: shifted bijection inv bps and plane partitions}, this triple product gives a formula for $|\InvBPD(1^k\times z)|$ when $z$ is any 
dominant involution for which $\hat\mu^{(z)}$
 is a shifted trapezoid.
In the special case when $z=n\cdots 321$ this triple product simplifies to
the double product in \eqref{ppp-eq}.
\end{remark}

Our corollary
verifies the conjectural identity \cite[Eq.~(6.3)]{HMP6}.
One can check that this identity is equivalent to the first equality in \eqref{ppp-eq}
by comparing Corollary~\ref{spp-ivex-cor} and \cite[Cor.~6.1]{HMP6}.
Combining Theorem~\ref{bij-thm} with the bijective results in \cite{HopkinsLai}  resolves one half of \cite[Problem~6.4]{HMP6}.

As another application, we can use Theorem~\ref{thm: shifted bijection inv bps and plane partitions} to settle a recent conjecture from \cite{FMP} on \defn{Ehrhart positivity} of order polytopes of posets coming from shifted straight shapes.
In the next two theorems, we write 
$ \delta = (n-1,\dots,3,2,1)$
for the staircase partition of $\binom{n}{2}$. This partition is relevant since every dominant $w \in S_n$ has $D(w) \subseteq \D_\delta$, with equality for the reverse permutation $w=n\cdots 321$.

Given a partition $\lambda$,
let $\mathsf{PP}_\lambda(t)$ denote the set of  plane partitions of shape $\lambda$
with all entries in $\{0,1,\dots,t\}$.
For a permutation $w \in S_n$ of length $m$, let $\cR(w)$ be its set of reduced words $r=(r_1,r_2,\dots,r_m)$.
Our last result in this section is a shifted analogue of the following theorem.

\begin{theorem}[{Fomin--Kirillov \cite[Thm.~2.1]{FominKirillov1997}}] \label{thm:FK}
		Suppose $\lambda$ is a partition contained in $\delta$ and $m=|\lambda|$.
         Then
$$|\mathsf{PP}_{\lambda}(t-1)| = \frac{1}{m!} \sum_{r \in \cR(w)} (t+r_1-1)(t+r_2-1)\cdots (t+r_m-1)\in \QQ_{\geq0}[t]
$$
        where $w\in S_n$ is the dominant permutation with $\D_\lambda= D(w)$.
\end{theorem}

We require one additional lemma.

\begin{lemma}[{Macdonald \cite[Eq.~6.11]{Macdonald}}] \label{lem:macdonald positivity each v}
If $w \in S_{n}$ then $\fkS_{w}(1,1,\ldots,1)=\frac{1}{\ell(w)!} \sum_{r \in \cR(w)} r_1 r_2 \cdots r_{\ell(w)}$.
\end{lemma}

%

When $\mu\subseteq \lambda$ are strict partitions,   let $\SPP_{\lambda/\mu}(t)$ be the set of weakly increasing fillings of the skew shifted shape $\SD_\lambda\setminus \SD_\mu$ with all entries in $\{0,1,\ldots,t\}$. 
Notice that 
\[|\SPP_{\lambda}(t)| = |\SPP_{\lambda/\emptyset}(t)|
\qquand |\mathsf{PP}_{\nu}(t)| = |\SPP_{(\nu+\delta)/\delta}(t)|\] for any partition $\nu$ with at most $n$ parts.
Ferroni--Morales--Panova conjectured in \cite[Conj.~7.10]{FMP} that 
$|\SPP_{\lambda/\mu}(t-1)| \in \QQ_{\geq 0}[t]$. The next result proves the $\mu=\emptyset$ case of this conjecture.

\begin{theorem}\label{fmp-thm}
	Suppose $\lambda$ is a strict partition with $\SD_\lambda \subseteq \D_\delta$. Let $m=|\lambda|$ and $k = \ell(\lambda)$.
    Then
		$$|\SPP_{\lambda}(t-1)| = \frac{2^{k}}{(2m)!!} \sum_{w\in \cW(z)} \sum_{r \in \cR(w)} (t+r_1-1)(t+r_2-1)\cdots (t+r_m-1) \in \QQ_{\geq 0}[t]
		$$
where $z \in \Invol$ is the unique dominant involution with $\SD_\lambda = \{ (i,j) \in D(z) : i\leq j\}$.
\end{theorem}


\begin{proof}
Since $\cR(1^{t-1}\times w^{-1}) = \{ (t+r_m-1,\dots,t+r_1-1) : (r_1,\dots,r_m) \in \cR(w)\}$
when $t \in \PP$, and since every $w \in \cW(z)$ has $\ell(w) = \ellinv(z) = |\lambda|$,
this follows by substituting Lemma~\ref{lem:macdonald positivity each v} into Theorem~\ref{thm: shifted bijection inv bps and plane partitions}.
\end{proof}

Ferroni--Morales--Panova's conjecture that 
$|\SPP_{\lambda/\mu}(t-1)| \in \QQ_{\geq 0}[t]$ arose in the study of the Ehrhart polynomials of polytopes that are (up to integral equivalence) both \defn{order polytopes} of a poset and  \defn{matroid polytopes}. A \defn{matroid polytope} $\mathcal{P}_M$ is the convex hull of the indicator vectors of the bases of the matroid $M$. \defn{Lattice path matroids} are matroids whose bases are given by the up-steps of lattice paths between 
the two lattice paths bordering a skew shape and are examples of \defn{positroids}. 

Among lattice path matroids, \defn{snake matroids} are those lattice path matroids that have at least two elements and are connected and are also \defn{binary} (that is they do not have the uniform matroid $U_{2,4}$ as a minor).  By a result of Benedetti, Knauer, and Valencia-Porras \cite[Thm.~1.1]{BKV-P}, order polytopes are affinely equivalent to  matroid polytopes if and only if the poset is a disjoint union of fence posets and antichains. Fence posets correspond to skew shapes that are \defn{border strips}.

In \cite{FMP}, the authors show that order polytopes of fence posets have Ehrhart polynomials with nonnegative coefficients by showing the stronger statement that $|\mathsf{PP}_{\lambda/\mu}(t-1)| \in \QQ_{\geq 0}[t]$ for all skew shapes $\lambda/\mu$. The proof strategy in \cite{FMP}, which also would apply to the shifted case, reduces the problem to showing that the linear coefficient of the polynomial $|\mathsf{PP}_{\lambda/\mu}(t-1)|$ is nonnegative.  This is then proved by analyzing the Kreweras determinant formula \cite{Kreweras_1965} for $|\mathsf{PP}_{\lambda/\mu}(t)|$. However, for the shifted case it is more challenging to extract the linear coefficient from Pfaffian formulas for $|\SPP_{\lambda/\mu}(t)|$, even when $\mu=\emptyset$ (see \cite[Thm. 2.2]{HopkinsLai}).

\subsection{Set-valued plane partitions}\label{svpp-sect}

We mention an involution analogue of one other result from \cite[\S7]{Weigandt}.
Let $\MInvBPD(z)$ be the set of pairs $(D,S)$ 
where $D \in \InvBPD(z)$ and $S\subseteq \upturn(D)$.
We refer to elements of this set as \defn{(involution) marked bumpless pipe dreams}. 

Next, define a \defn{set-valued shifted plane partition} of strict partition shape $\lambda$
to be a filling of $\SD_\lambda$ by finite, nonempty subsets of nonnegative integers
such that $\max(A) \leq \min(B)$ whenever $A$ is the set assigned to a given cell  $(i,j)$ and $B$ is the set assigned to $(i,j+1)$ or $(i+1,j)$.
When $T $ is such an object we write ``$(i,j) \in T$'' as a shorthand for  ``$(i,j) \in \SD_{\lambda}$'' and let $T_{ij}$ denote the set assigned by $T$ to position $(i,j)$.

Given $z \in \Ivex$, let 
$\SetSPP(z)$ be the set of all set-valued shifted plane partitions of shape $\hat\mu^{(z)}$ whose entries in row $i$ are subsets of $\left\{0,1,2,\dots,f^{(z)}_i-i\right\}$ for all $i$.
Recall the notation $a\oplus b = a + b + \beta ab$.
We define the \defn{$K$-theoretic weights} of elements $(D,S) \in \MInvBPD(z)$ and $T \in \SetSPP(z)$ as
\be \textstyle \kweight(D,S) := \prod_{(i,j) \in \blank(D) \sqcup S} x_i \oplus x_j
\quand
\kweight(T) := \prod_{(i,j) \in T} \prod_{k \in T_{ij}} x_{i+k} \oplus x_{j+k}.
\ee
We also define 
the \defn{degrees} of these objects to be 
\be
\textstyle
|(D,S)| := |\blank(D)| + |S|\quand |T| := \sum_{(i,j) \in T} |T_{ij}|.
\ee
Given $T \in \SetSPP(z)$ let $\flatten(T) \in \SPP(z)$  be the shifted plane partition formed from $T$ by replacing each entry $T_{ij}$ by
its smallest element $\min(T_{ij})$ and define 
\be A_T := \{ (j+k, i+k) : (i,j) \in T \text{ and }\min(T_{ij}) < k \in T_{ij} \}.\ee

\begin{theorem}
\label{upsilon-thm}
Suppose $z \in \Ivex$.
Then the operation defined by 
\[
\overline{\Upsilon}(T) := ( \Upsilon(\flatten(T)) , A_T)
\]
is a weight- and degree-preserving bijection $  \SetSPP(z) \to \MInvBPD(z)$.
\end{theorem}

\begin{proof}
Define
$\SetPP(z)$ to be the set of all set-valued (unshifted) plane partitions of shape $\mu^{(z)}$ whose entries in row $i$ are subsets of $\left\{0,1,2,\dots,f^{(z)}_i-i\right\}$.
Let $\Phi$ be the operator on this set that adds $i$ to every element of every entry in row $i$ for all $i$.

It follows from Lemma~\ref{vex-double-lem} that the obvious doubling operation is a bijection
from $\SetSPP(z)$ to the set of symmetric elements of $\SetPP(z)$.
By Corollary~\ref{vex-double-cor}, the map
\[
\double : (D,S) \mapsto (\double(D), S\cup S^\top)
\]
is likewise a bijection from $\MInvBPD(z)$ to the ``symmetric'' elements of the set $\MPipes(z)$ of
marked bumpless pipe dreams of $z$ defined by Weigandt  in \cite[\S1.2]{Weigandt}.

Composing $\Phi$ with the map $\overline{\gamma}$ in \cite[Thm.~1.6]{Weigandt}
gives a bijection $\overline{\gamma}\circ \Phi :\SetPP(z) \to \MPipes(z)$ that commutes with transpose, and thus sends symmetric elements to symmetric elements.
Inspecting the definition of $\overline{\gamma}$ before \cite[Lem.~7.8]{Weigandt} 
shows that  $\overline{\Upsilon}$, which is weight- and degree-preserving by definition,  is the unique map that makes the diagram
\[
\begin{tikzcd}
\SetSPP(z) \arrow[r, "\overline{\Upsilon}"] \arrow[d, "\double"'] & \MInvBPD(z) \arrow[d, "\double"] \\
\SetPP(z) \arrow[r, "\overline{\gamma}\circ \Phi"']                 & \MPipes(z)
\end{tikzcd}
\]
commute. We conclude that  $\overline{\Upsilon} : \SetSPP(z) \to \MInvBPD(z)$ is a bijection.
\end{proof}

In view of \eqref{kweight-eq} and Corollary~\ref{k-sub-cor}(b),
 when $z \in \Invol$ is any involution we have  
\be\label{fkGO-conj-eq}
\textstyle \sum_{D \in \InvBPD(z)} \kweight(D) = \sum_{(D,S) \in \MInvBPD(z)} \beta^{|(D,S)| - \ellinv(z)}  \kweight(D,S).\ee
Thus, the preceding result shows that Theorem~\ref{fkGO-conj} is equivalent to the following statement.

\begin{theorem}
\label{fkGO-conj2}
For each $z \in \Ivex $ it holds that
$\textstyle
\fkGO_z = \sum_{T \in  \SetSPP(z) } \beta^{|T| - \ellinv(z)}\kweight(T).
$
\end{theorem}

This formula is the natural involution analogue of 
the Knutson--Miller--Yong tableau formula for (double) Grothendieck polynomials of vexillary permutations \cite[Thm.~5.8]{KnutsonMillerYong} 
(see \cite[Thm.~7.1]{Weigandt}),
and can be proved in a similar way.
Specifically, we can derive Theorem~\ref{fkGO-conj2} by studying a Gr\"obner basis for the ideals of \defn{symmetric matrix Schubert varieties} for vexillary involutions. 
As this material is mostly independent of the bumpless pipe dream story at hand, and much more geometric, we present it separately in Appendix~\ref{app-sect2}.

\section{Final remarks and future directions}\label{final-sect}

In this final section we briefly discuss
 the connections between BPDs and alternating sign matrices, and then examine some problems concerning the asymptotics of BPD counts for permutations and involutions.

\subsection{Symmetric frozen ice}\label{ice-sect}

Bumpless pipe dreams (which are not required to be reduced) are in correspondence with alternating sign matrices \cite{L_ice,Weigandt} and this connection has led to key insights into Schubert and Grothendieck polynomials \cite{BuciumasScrimshaw,escobar2025algebrageometryasmweak,MPPY}. Here, we discuss a related correspondence between involution pipe dreams and symmetric alternating sign matrices.

An \defn{alternating sign matrix (ASM)} is a square matrix with entries $0,1,-1$ such that the rows and columns add up to 1 and the nonzero entries alternate in sign. These matrices are of interest in enumerative and algebraic combinatorics, statistical mechanics, and polyhedral combinatorics (e.g., see \cite{BressoudBook}). We denote the set of $n\times n$ alternating sign matrices by $\ASM_n$. By a fundamental result of Zeilberger \cite{ZeilASM} and Kuperberg \cite{Kuperberg} the number of matrices in $\ASM_n$ has the product formula $\prod_{k=0}^{n-1} \frac{(3k+1)!}{(n+k)!}$. Bumpless pipe dreams 
are in bijection with alternating sign matrices via the following correspondence \cite{L_ice,Weigandt}:  
\be
\label{asm-eq}
\JTILE\mapsto -1, \qquad \CTILE \mapsto 1, \qquand \BTILE, \HTILE, \VTILE, \PTILE \mapsto 0.
\ee

If $A$ is an ASM, then 
a zero entry $A_{ij}=0$  is a \defn{balanced zero} when both the sum of the entries above $(i,j)$ in column $j$ and to the left of $(i,j)$ in row $i$ are zero.
Let $\mathrm{BZ}(A)$ be the set of positions $(i,j)$ 
of the balanced zeros in $A$.
Let $M(A)$ be the set of positions $(i,j)$ where $A_{ij}=-1$.

The following identity is similar to the relation between ASMs and tilings of the Aztec diamond observed in \cite{ASMs_original}. 
Let $\kBPD(n) = \bigsqcup_{w \in S_n} \kBPD(w)$ be the set of all  bumpless pipe dreams for permutations in $S_n$.
Then setting $\beta=1$ in Proposition~\ref{prop:partition function inv Grothendieck} 
implies via \cite[Thm.~1.1]{Weigandt} and the correspondence \eqref{asm-eq} that
{\small\be\label{asm-eq2}
\sum_{A \in \ASM_n} \prod_{(i,j) \in M(A)} (1+x_i) \prod_{(i,j) \in \mathrm{BZ}(A)} x_i = \sum_{D \in \kBPD(n)} \prod_{(i,j)\in \upturn(D)} (1+x_i) \prod_{(i,j) \in \blank(D)} x_i =
\prod_{i =1}^{n-1} (1+ x_i)^{n-i}
.
\ee}
Hence (setting $x_i=1$) one also has
$
\sum_{A \in \ASM_n} 2^{|M(A)|} = \sum_{D \in \kBPD(n)} 2^{|\upturn(D)|} = 2^{\binom{n}{2}}.
$

A \defn{diagonally symmetric} alternating sign matrix is an ASM that is symmetric across the main diagonal.
This variant of alternating sign matrices has been studied previously
in \cite{RobbinsSymm,KuperbergSymmetryClasses,Bousquet-Melou-Habsieger}.
We denote the set of diagonally symmetric ASMs of size $n\times n$ by $\DASM_n$. Unlike alternating sign matrices, no simple product formula is known for their number \cite[\href{https://oeis.org/A005163}{A005163}]{OEIS}. However, Behrend--Fischer--Koutschan \cite{BFK} found  a Pfaffian formula which lets one compute the sequence $|\DASM_n|=	1,\: 2,\: 5,\: 16,\: 67,\: 368,\: 2630,\: 24376,\: 293770,\: 4610624,
\ldots
$.

Using the same map \eqref{asm-eq}, diagonally symmetric alternating sign matrices are in bijection with symmetric bumpless pipe dreams. 
There are interesting variants of \eqref{asm-eq2}
that 
hold in this context.
To make these statements a little simpler, define $x_{i\boxplus j}$ by setting $\beta=1$ in $x_{i\oplus j}$ from \eqref{i-oplus-j-eq}, so that
\[
x_{i\boxplus i}=x_i\qquand x_{i\boxplus j} = x_i +  x_j + x_ix_j\text{ when $i\neq j$.}
\]
Finally, let $\kInvBPD(n) = \bigsqcup_{z \in \Invol} \kInvBPD(z)$. 

\begin{proposition}
\label{iG-conj-cor}
If $n \in \PP$ then
$\displaystyle
	   \sum_{D \in \kInvBPD(n)} \prod_{(i,j) \in \upturn(D)} (1+x_{i\boxplus j}) \prod_{(i,j) \in \blank(D)} x_{i\boxplus j}  =
	   \prod_{i =1}^{n-1} (1+ x_i)^{n-i}
	   .
$
	\end{proposition}

This result offers some support for Conjecture~\ref{iG-conj}, since 
Proposition~\ref{prop:partition function inv Grothendieck}
shows that it would be a direct consequence of that conjecture.
The proof of the proposition was found with AI assistance and involves 
a version of of the Yang-Baxter equation \cite{fomin1994grothendieck}. 
Since the argument is technical and independent from the rest of this work,
it is currently omitted. It will be included in a later version of this article.

\begin{remark}
A different weighted sum over $K$-theoretic bumpless involution pipe dreams $D$ results in the same product formula.
Below, define $\cJBdiag(D)$
to be the set of positions $(i,j) \in \upturn(D)\sqcup \blank(D)$ with $i=j$,
then set
$\cJoff(D) = \upturn(D) \setminus \cJBdiag(D)$
and
$\cBoff(D) = \blank(D) \setminus \cJBdiag(D)$.
If $n \in \PP$ then
\be\displaystyle
	   \sum_{D \in \kInvBPD(n)} \prod_{(i,j) \in \cJBdiag(D)} x_i 
       \prod_{(i,j) \in \cJoff(D)} x_ix_j \prod_{(i,j) \in \cBoff(D)} (1+x_i)(1+x_j) = 
        \prod_{i =1}^{n-1} (1+ x_i)^{n-i}
       .
\ee
We do not know if this second identity is related to bumpless pipe dream formulas for involution Grothendieck polynomials. It can be derived as a special case of a more general identity. Fix a vector $\varepsilon=(\varepsilon_1,\varepsilon_2,\dots,\varepsilon_n) \in \{0,1\}^n$ and let $\mathsf{TPD}(n,\varepsilon)$ be the set of all bumpless wirable tilings in $\LT_n$ whose tiles have $\varepsilon_i$ south connections in position $(n,i)$. (In other words, we consider bumpless triangular pipe dreams 
as in Definition~\ref{tpd-def} but with custom requirements for the connections on the south boundary.)
For any $D$ in this set, define 
\[\textstyle
\mathscr{W}(D) = 
\prod_{(i,j) \in \cJBdiag(D)} x_i 
       \prod_{(i,j) \in \cJoff(D)} x_ix_j \prod_{(i,j) \in \cBoff(D)} (1+x_i)(1+x_j)
       .
       \]
Then $ \kInvBPD(n) = \mathsf{TPD}(n,(1,1,\dots,1))$
and one can prove
\be\textstyle
	   \sum_{D \in \mathsf{TPD}(n,\varepsilon)} \mathscr{W}(D) = 
        \prod_{i =1}^{n-1} (1+ x_i)^{n-i} \prod_{i=1}^n \(x_i \prod_{j=i+1}^n (1+x_j)\)^{1-\varepsilon_i}
        \ee
        by induction on $n$ using a transfer-matrix method.
\end{remark}

Specializing $x_i=1$ in the  identities above gives the following enumerative formulas for weighted sums over $K$-theoretic bumpless involution pipe dreams or diagonally symmetric ASMs. Let
\be\ba
 A_n & :=\sum_{D \in \kInvBPD(n)}  2^{|\cJdiag(D)|} 3^{|\cBoff(D)|}4^{|\cJoff(D)|} 
 &&= \sum_{A \in \DASM_n} 2^{\MD(A)}3^{\BZOD(A)}4^{\MOD(A)}  
 \quad\text{and}
 \\
        B_n&:=\sum_{D \in \kInvBPD(n)} 4^{|\cBoff(D)|}  &&=\sum_{A \in \DASM_n} 4^{\BZOD(A)}  
\ea
\ee
    where $\cJdiag(D) 
    = \{(i,j) \in \upturn(D) : i=j\}$, and $\MD(A)$ and $\MOD(A)$ are the numbers of diagonal minus signs and strictly lower triangular minus signs in $A$, respectively, and $\BZOD(A)$ is the number of strictly lower triangular balanced zeros in $A$.
Then we have 
\be A_n=B_n=2^{\binom{n}{2}}.\ee

\subsection{Growth and maximizers of principal specializations}

The growth of the principal specialization of Schubert polynomials was first considered in work of Merzon and Smirnov \cite{Merzon_Smirnov} and a preprint of Stanley \cite{StanleyShenanigans}. In this context, for a permutation $w$, one is interested in the quantities \[
\textstyle
\Upsilon_w:=\fkS_w(1,\ldots,1)=|\BPD(w)| \quand 
 u_n = \max_{w\in S_n}\Upsilon_w\] and the permutations that achieve the maximum $u_n$.

Based on data for $n\leq 13$ (see \cite[\S 1.2]{AndersonPanovaPetrov}), it appeared that the maximum $u_n$ was achieved by a \defn{layered permutation}, meaning a permutation   obtained by taking block diagonal sums \[z(b) = z(b_k,\ldots,b_2,b_1):=w_0(b_k)\oplus \cdots \oplus w_0(b_2)\oplus  w_0(b_1)\] of reverse permutations $w_0(k)=k\cdots 321 \in S_k$. Let $\cL_n$ be the set of layered permutations in $S_n$ and define \[ \textstyle \uPrime_n=\max_{z \in \cL_n} \Upsilon_z.\]
Notice that $\cL_n \subseteq \Invol\subseteq S_n$.
Results in \cite{MPP_schub} identify, among layered permutations, the optimal compositions $b^*$ for each size, their limiting shape, and the value of $\lim_{n \to \infty} (1/n^2)\log_2 \uPrime_n$.

Anderson, Panova, and Petrov discovered in \cite{AndersonPanovaPetrov}  non-layered permutations in $S_{17}$ whose principal specializations exceed the best layered of that size. These values were found by applying (co)transitions on the optimal layered permutation. It remains open to settle a question that Stanley raised in \cite{StanleyShenanigans} about whether the constant $\lim_{n\to \infty} \frac{1}{n^2}\log_2 u_n$ exists and matches the growth constant of the analogous question among layered permutations computed in \cite{MPP_schub}. We summarize the best known bounds for this conjectured constant.

\begin{proposition}[{see \cite[Thm. 1.3]{MPP_schub}
for lower bound and \cite[Rem. 6.8]{MPPY} for upper bound}] \label{prop: bound PS Schub}
We have  
\begin{align*}
0.2932362762\approx \frac{\gamma}{\log 2} \leq \liminf_n \frac{\log_2 u_n}{n^2} &\leq \limsup_n \frac{\log_2 u_n}{n^2} \leq \frac32 \log_2(3) -2  \approx 0.3774437510,
\end{align*}
where $\gamma\approx 0.2032558981$ is a universal constant.
\end{proposition}

Let $\delta_n = (n-1,n-2,\dots,3,2,1)$ and $\delta^\ast_n = \half(\delta_n) = (n-1,n-3,n-5,\dots)$.
\begin{remark} \label{rem: bds max ps Schubs}
The upper bound in the preceding result comes from the fact that we have 
$u_n \leq |\kBPD(n)|= |\ASM_n|$
 and the asymptotics for the number of ASMs: $\lim_{n\to \infty} |\ASM_n|^{1/n^2}=3\sqrt{3}/4$ (see \cite[Eq. (11.5)]{BFK}).
The lower bound  is exactly $\lim_{n \to \infty} \frac{1}{n^2}\log_2 \uPrime_n$ when restricting the problem to layered permutations. One of the reasons this problem is more tractable is the following product formula
    \[
\Upsilon_{w_0(k,n)} = |\mathsf{PP}_{\delta_n}(k)| = \ds\prod_{1\leq i<j\leq n} \frac{2k+i+j-1}{i+j-1}
\qquad\text{for }w_0(k,n):= 1^k\times w_0(n).
\]
The first equality is due to  Fomin and Kirillov in \cite{FominKirillov1997} and the second equality comes from a product formula for plane partitions due to Proctor \cite{Proctor1994}.
We also mention that the principal specialization $\Upsilon_w$ for layered permutations can be computed iteratively \cite[Prop. 2.2]{MPP_schub} as 
\begin{equation} \label{eq: ps layered perm}
{\Upsilon}_{z(b_k,\ldots,b_1)} = {\Upsilon}_{z(b_k,\ldots,b_{2})}\cdot {\Upsilon}_{w_0(b_k+\cdots+b_{2},b_1)}.
\end{equation}
This formula makes it feasible to calculate the optimal value among layered permutations via a dynamic program (see \cite[Appendix]{MPP_schub}). Moreover from \cite[Thm. 1.3]{MPP_schub}, as $n\to \infty$, the maximum value $\uPrime_n$ is achieved at a layered permutation $z(\ldots,b_2,b_1)$ where $b_i \sim \alpha^{i-1}(1-\alpha)n$ for a universal constant $\alpha\approx 0.4331818312$.
\end{remark}

For involutions, consider the 
principal specialization of $\iS_z$. Let 
$\widehat{\Upsilon}_z:
=2^{\cyc(z)}  \iS_{z}(1/2,\ldots,1/2)$. By Corollary~\ref{spp-ivex-cor}, when $z$ is vexillary we have $\widehat{\Upsilon}_z=|\InvBPD(z)|$ where $\cyc(z)$ is the number of $2$-cycles of $z$. For a general involution $z \in \Invol$, we have   
\begin{equation} \label{eq: gen specialization}
\widehat{\Upsilon}_z=\sum_{D \in \InvBPD(z)} 2^{\cyc(z) - \mathsf{db}(D)},
\end{equation}
where $\mathsf{db}(D) $  
is the number of diagonal blank positions,
which is at most $\cyc(z)$ 
 by Corollary~\ref{db-cor}.
Let 
\[\textstyle \hatu_n := \max_{z \in \Invol}\widehat{\Upsilon}_z.\] 
Table~\ref{table:max values PS involutions} lists the values of $\hatu_n$ 
for $2\leq n \leq 12$ and provides involutions that achieve these maxima.
Although the maximizers in these cases are all layered,
there exist non-layered involutions with higher values for $\widehat{\Upsilon}_z$, as in the permutation case.

\begin{example} \label{ex: non layered inv that beats layered}
    For $n=15$, and $z=[1, 3, 2, 8, 6, 5, 15, 4, 14, 13, 12, 11, 10, 9, 7]=s_7\cdot  z(1,2,4,8)\cdot s_7$, we have checked  by a direct calculation using the transition formula \eqref{EQ: iS Tran} for  involution Schubert polynomials that $\widehat{\Upsilon}_{z}=4643513424> 4453792200 = \widehat{\Upsilon}_{z(1,2,4,8)}$, where $z(1,2,4,8)$ is the optimal layered involution. 
\end{example}

\begin{remark}
    By Theorem~\ref{thm: shifted bijection inv bps and plane partitions}, the value $\widehat{\Upsilon}_{w_0(k,n)}$ is exactly the formula in  Corollary~\ref{cor:plane partition product}:
\[
\widehat{\Upsilon}_{w_0(k,n)} = |\SPP_{\delta^*_n}(k)| = \ds\prod_{i=1}^{\lceil n/2\rceil} \prod_{j=1}^{\lfloor n/2\rfloor} \frac{i+j+k-1}{i+j-1}.
\]
Thus, the principal specialization $\widehat{\Upsilon}_z$ for layered involutions can again be computed iteratively as 
\begin{equation} \label{eq:prod layered invol}
\widehat{\Upsilon}_{z(b_k,\ldots,b_1)} = \widehat{\Upsilon}_{z(b_k,\ldots,b_{2})}\cdot \widehat{\Upsilon}_{w_0(b_k+\cdots+b_{2},b_1)},
\end{equation}
and the optimal value among layered involutions can be computed via a dynamic program. Let \[
\textstyle \hatuPrime_n:=\max_{z \in \cL_n} \widehat{\Upsilon}_z.\] 
\end{remark}

It would be interesting to study whether the limit $\lim_{n\to \infty} \frac{1}{n^2}\log_2 \hatu_n$ 
exists and if so, find this constant. We end this section by giving some preliminary bounds that are an analogue and exactly half the bounds of Proposition~\ref{prop: bound PS Schub}.

\begin{proposition}
We have that 
\begin{align*}
0.1466181381 \approx \frac{\gamma}{2\log 2} \leq \lim \inf \frac{\log_2 \hatu_n}{n^2} &\leq \lim \sup \frac{\log_2 \hatu_n}{n^2} \leq \frac34 \log_2(3) -1  \approx 0.1887218755,
\end{align*}
where $\gamma\approx 0.2032558981$ is a universal constant.
\end{proposition}

\begin{proof}
For the upper bound, by \eqref{eq: gen specialization} we have that 
\[
\hatu_n \leq 2^{\lfloor n/2\rfloor} |\kInvBPD(n)|= 2^{\lfloor n/2\rfloor}|\DASM_n|
\]
since $\cyc(z) - \mathsf{db}(D) \leq \lfloor n/2\rfloor$. Behrend, Fischer, and Koutschan found the asymptotics of the  number of diagonally symmetric ASMs \cite[Thm. 11.2]{BFK} to be $\lim_{n\to \infty} |\DASM_n|^{2/n^2} = 3\sqrt{3}/4$. Taking the logarithm and dividing by $n^2$ gives the upper bound.

For the lower bound, let $\SPP_{\delta^*_n}(k;{\bf d})$ be the set of plane partitions in $\SPP_{\delta^*_n}(k)$ with diagonal ${\bf d}$. Analogously, define $\mathsf{PP}_{\delta_n}(k;{\bf d})$.
Since a plane partition in $\mathsf{PP}_{\delta_n}(k;{\bf d})$ can be obtained by gluing two copies of shifted plane partitions in  $\SPP_{\delta^*_n}(k;{\bf d})$ along the diagonal, we have 
\[
|\mathsf{PP}_{\delta_n}(k)| = \sum_{\bf d}  |\mathsf{PP}_{\delta_n}(k;{\bf d})|=\sum_{\bf d} |\SPP_{\delta^*_n}(k;{\bf d})|^2 \leq \left(\sum_{\bf d} |\SPP_{\delta^*_n}(k;{\bf d})|\right)^2= |\SPP_{\delta^*_n}(k)|^2.
\]
 
On the other hand, refining the number of shifted plane partitions by the value on the diagonal and using Cauchy--Schwarz we have that 
\[
|\SPP_{\delta^*_n}(k)|^2 = \left(\sum_{\bf d} |\SPP_{\delta^*_n}(k;{\bf d})|\right)^2 \leq \left(\sum_{\bf d} |\SPP_{\delta^*_n}(k;{\bf d})|^2\right)\cdot |\diag(k,n)|=|\mathsf{PP}_{\delta_n}(k)| \cdot |\diag(k,n)|,
\]
where $\diag(k,n)$ is the set of diagonals ${\bf d}$ appearing in $\SPP_{\delta^*_n}(k)$. 
The main diagonal has $\lfloor n/2\rfloor$ entries, each in $\{0,1,\dots,k\}$. Consequently, we have $| \diag(k,n)|\leq (k+1)^{\lfloor n/2\rfloor}$
along with
\[
|\mathsf{PP}_{\delta_n}(k)| \leq |\SPP_{\delta^*_n}(k)|^2 \leq |\mathsf{PP}_{\delta_n}(k)|\cdot (k+1)^{\lfloor n/2 \rfloor}.
\]
Using these bounds with the product formulas in \eqref{eq: ps layered perm} and \eqref{eq:prod layered invol} gives the layered permutation bound
\[
{\Upsilon}_{z(b_k,\ldots,b_1)} \leq \widehat{\Upsilon}_{z(b_k,\ldots,b_1)}^2 \leq {\Upsilon}_{z(b_k,\ldots,b_1)} \cdot (|b|+1)^{|b|},
\]
where $|b|=b_1+b_2+\cdots + b_k$. 
Since $\lim_{n\to \infty} \tfrac{1}{n^2}\log_2 \uPrime_n$ exists \cite[Thm. 1.3]{MPP_schub},
taking logarithms, dividing by $|b|^2$, and then applying the squeeze theorem gives  $\lim_{n\to \infty} \tfrac{1}{n^2} \log_2 \hatuPrime_n = \tfrac12 \lim_{n\to \infty} \tfrac{1}{n^2}\log_2 \uPrime_n$. Lastly, by Proposition~\ref{prop: bound PS Schub} and Remark~\ref{rem: bds max ps Schubs}, we have $\lim_{n\to \infty} \tfrac{1}{n^2}\log_2 \uPrime_n =\gamma/{\log 2}$, so the result follows.
\end{proof}

\begin{remark}
By an argument similar to the one above, one can also deduce that as $n\to \infty$, the maximum value $\hatuPrime_n$ is achieved at layered permutations with the same shape (not necessarily the same involution) as those  corresponding to $\uPrime_n$: i.e., $z(\ldots,b_2,b_1)$ where $b_i \sim \alpha^{i-1}(1-\alpha)n$ for the same universal constant $\alpha$ as  in Remark~\ref{rem: bds max ps Schubs}.
\end{remark}


\begin{table}
\begin{center}
\begin{tabular}{c|rr}  \hline 
 $n$    
 & $\hatu_n$ &  $z$   \\ \hline
   2  &   1 & [1,2] \\
   3 &  2 & [1,3,2] \\
   4 &  3 & [1,2,4,3] \\
   5 &  8 & [1,3,2,5,4]\\
   6 &  20 & [1,2,6,5,4,3]\\
   7 &  100 & [1,3,2,7,6,5,4]\\
   8 &  350 & [1,3,2,8,7,6,5,4]\\
   9 &  1960 & [1,3,2,9,8,7,6,5,4]\\
   10 & 12348 & [1,2,4,3,10,9,8,7,6,5]\\ 
   11 & 112896& [1,3,2,5,4,11,10,9,8,7,6]\\
   12 & 931392 & [1,3,2,5,4,12,11,10,9,8,7,6] \\\hline 
\end{tabular}
\qquad 
\raisebox{-60pt}{\includegraphics[scale=0.4]{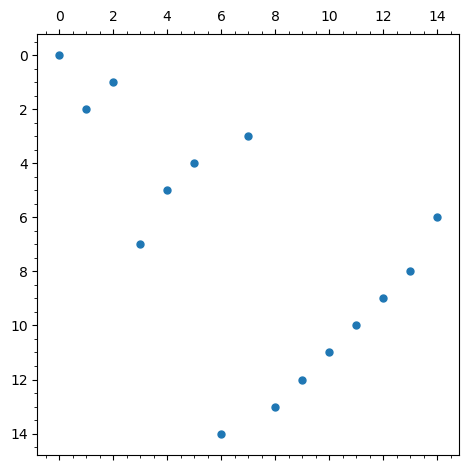}}
\end{center}
\caption{Maximum values of 
$\hatu_n :=\max_{z \in \Invol} 2^{\cyc(z)} \iS_{z}(\frac{1}{2},\ldots,\frac{1}{2})$ and their respective maximizing involutions $z$ up to $n=12$, along with the permutation matrix of $z=[1, 3, 2, 8, 6, 5,15, 4,14, 13, 12, 11, 10, 9, 7]$.}
\label{table:max values PS involutions}
\end{table}

\appendix

\section{Gr\"obner geometry for symmetric matrix Schubert varieties}\label{app-sect2}

This section provides a geometric proof of Theorems~\ref{fkGO-conj} and \ref{fkGO-conj2}.
We have already seen that these statements are equivalent via Theorem~\ref{upsilon-thm} and \eqref{fkGO-conj-eq}. We prove both results here by deriving yet another equivalent reformulation (Theorem~\ref{fkGO-conj3}). 
This material is mostly independent of the rest of this article. Our proof strategy is roughly the same as the one carried out in \cite[\S4.4]{MarPaw2022},
and can be condensed to a few pages 
by using results in \cite{GMN,GrahamKreimanB,KnutsonMillerSubword,KnutsonMillerYong} to derive several key steps.

Let $\Mat_{n\times n}$ be the variety of $n\times n$ matrices over $\CC$ and let $\SymMat_{n\times n} $  be the subvariety of symmetric matrices in $\Mat_{n\times n}$. 
The coordinate rings of $\Mat_{n\times n}$ and $\SymMat_{n\times n}$
are the polynomial algebras 
\[
 \CC[\Mat_{n\times n}] = \CC[x_{ij} : (i,j) \in \square_n]
 \quand
  \CC[\SymMat_{n\times n}] = \CC[u_{ij} : (i,j) \in \UT_n]
 \]
 where $\square_n = [n]\times [n]$ and $\UT_n = \{(i,j) \in \square_n : i\leq j\}$. 
 Let $\pi(x_{ij}) = \pi(x_{ji}) = u_{ij}$ for $i\leq j$ and extend this to a homomorphism between the two coordinate rings.
 
 Fix $w \in S_n$ and $z \in \Invol$.
 Given $v \in S_n$ let $r_v(i,j) = |\{a \in [i] : v(a) \leq j\}|$.
Write $\cX = \mat{{r} x_{ij}}$ and $\cU = \mat{{r} u_{ij}}$
for generic matrices of indeterminates, where $u_{ij}=u_{ji}$.
Then let
\[
\ba
I_w &:= \Bigr\langle \det(\cX_{RC}) : (i,j)\in\square_n,\ R \subseteq [i],\ C\subseteq[j],\ |R|=|C|=1+r_w(i,j) \Bigr\rangle\subseteq \CC[\Mat_{n\times n}],
\\
\SymI_z &:= \Bigr\langle \det(\cU_{RC}) :(i,j)\in\UT_n,\ R \subseteq [i],\ C\subseteq[j],\ |R|=|C|=1+r_z(i,j) \Bigr\rangle\subseteq \CC[\SymMat_{n\times n}].
\ea
\]
Recall that if $D \subseteq \PP\times \PP$ then $\Ess(D)=  \{ (i,j) \in D: (i+1,j)\notin D\text{ and }(i,j+1) \notin D\}.$

\begin{definition}
Let $\cG_w$ and $\cSymG_z$ be the subsets of the minors generating $I_w$ and $\SymI_z$ above 
respectively indexed by the pairs $(i,j) $ in 
$\Ess(D(w))$ and $ \Ess( D(z)\cap \UT_n)= \Ess(D(z)) \cap \UT_n.$
\end{definition}

A   \defn{term order} for a polynomial ring $R$ is a total order on its set of    monomials,
such that $1$ is   minimal and such that if $\mon_1$, $\mon_2$, and $\mon_3$ are monomials with $\mon_1 \leq \mon_2$, 
then $\mon_1  \mon_3 \leq \mon_2 \mon_3$. 
Once a term order has been chosen, the   \defn{initial term} $\init(f)$ of $f \in R$ is either the maximal monomial $\mon$ whose coefficient in $f$ is nonzero, or zero when $f=0$.
The \defn{initial ideal} of an ideal $I\subseteq R$ is
the monomial ideal $\init(I) := \langle\init(f) : f \in I\rangle$.
A \defn{Gr\"obner basis} $G$ for $I$, finally,
is a generating set for $I$ such that $\{ \init(g) :g \in G\}$ generates $\init(I)$.

A term order for $\CC[\Mat_{n\times n}]$ or $\CC[\SymMat_{n\times n}]$ is
 \defn{diagonal} if the leading term of any minor in $\cX$ or $\cU$, respectively, is the product of the terms on the main diagonal of the relevant submatrix.
 From this point on, we fix an arbitrary diagonal term order 
 for each ring and define all initial ideals and Gr\"obner bases relative to this choice. In particular, set $J_w = \init(I_w)$ and $\SymJ_z = \init(\SymI_z)$.
 
 \begin{theorem} 
 \label{gr-thm1}
 Assume $w \in S_n$ and $z \in \Invol$ are vexillary.
Then  for any diagonal term order:
 \ben
 \item[(a)] The sets $\cG_w$ and $\cSymG_z$ are Gr\"obner bases for the respective ideals $I_w$ and $\SymI_z$.

 \item[(b)] 
 The initial ideals $J_w = \init(I_w)$ and $\SymJ_z = \init(\SymI_z)$  are both square-free, with $\SymJ_z = \pi(J_z)$.
 \een
 \end{theorem}
 
 \begin{proof}
The set  $\cG_w$ is a Gr\"obner basis for  $I_w$ by \cite[Thm.~3.8]{KnutsonMillerYong}. The other assertion in part (a) is a direct consequence of \cite[Thm.~4.8]{GMN}, as we now explain.

If $z=1$ then $\cSymG_z=\varnothing$ is trivially a Gr\"obner basis for $\SymI_z =0$.
Assume $1\neq z \in \Ivex$ and recall from Section~\ref{flagged-sect}
that $\lambda = \lambda^{(z)}$ is the smallest partition
whose diagram contains $D(z)$. 
For the rest of this paragraph we follow the notation and terminology of \cite[\S4]{GMN}.
Define 
\[\cL = \D_\lambda, 
\quad\cL^+ = \SD_{\half(\lambda)},
\quand
L = \{u_{ij} : (i,j)\in \cL^+\}.\] Since $\lambda$ is symmetric,  $\cL$ is a \defn{symmetric ladder}
and the \defn{lower border} of $\cL^+$ is just the outermost border strip of $\SD_{\half(\lambda)}$.
Because $z$ is vexillary, the elements of   $\Ess(D(z))$ lie on a piecewise linear curve oriented southwest-to-northeast \cite[Prop.~2.2.8]{Manivel}. This property implies that 
 $E := \Ess(D(z))\cap\UT_n$ is a subset of the lower border of $\cL^+$ that contains all
of its  \defn{lower outside corners}. Moreover, we can list   $E$  
 as 
 \[(v_1,w_1),\dots,(v_s,w_s)\quad\text{where}\quad v_1\leq \dots\leq v_s\leq n\quand n\geq w_1\geq \dots \geq w_s\geq 1.\]
Define $t_i =1+r_z(v_i,w_i)$ and $t=(t_1,\dots,t_s)$.
Then the \defn{symmetric mixed ladder determinantal ideal} $I_t(L)$ from \cite[Def.~4.3]{GMN} is precisely $\SymI_z$ while the set $\cG_t(L)$ from \cite[Not.~4.4]{GMN} is $\cSymG_z$.
As \cite[Thm.~4.8]{GMN} asserts that  $\cG_t(L)$ is a Gr\"obner basis of $I_t(L)$ with respect to any diagonal term order, we get the rest of part (a).

 Under a diagonal term order,
 the initial term of any minor in $\cX$ or $\cU$ is square-free.
 These terms generate
 $J_w$ and $\SymJ_z$ by part (a), so both ideals are square-free.
 We have $\SymJ_z = \pi(J_z)$ since $\pi$ restricts to a surjective map $\{ \init(f) : f \in \cG_z\} \to \{ \init(f) : f\in \cSymG_z\}$
 between generating sets.
 \end{proof}


\begin{corollary}
If $z \in \Ivex$ then $\SymI_z$ is a radical ideal.
\end{corollary}

\begin{proof}
It is a standard fact that any ideal with a square-free initial ideal is radical; see, for example, \cite[Prop.~5.3]{Sturmfels2002} (which works over $\QQ$ but whose proof is valid over $\CC$).
\end{proof}
 
 For a matrix $M$ and subsets $R,C \subseteq [n]$ let $M_{RC}$ be the submatrix in rows $R$ and columns $C$.
The \defn{matrix Schubert variety}  of $w \in S_n$
is then
\[\label{eq:MX-defn}
    \MSV_w = \Bigl\{M \in \Mat_{n\times n} : \text{$\rank M_{[i][j]} \leq   r_w(i,j)$ for all $(i,j) \in \square_n$}\Bigr\}. 
\]
The \defn{symmetric matrix Schubert variety} of $z \in \Invol$ is     $\SymMSV_z 
= \MSV_z  \cap \SymMat_{n\times n}.$
For any $w \in S_n$, the   ideal $I_w$ is known to be prime and equal to the defining radical ideal $I(\MSV_w)$ of $\MSV_w$ \cite[\S3.1]{KnutsonMillerYong}.

\begin{corollary}\label{sym-i-msv-cor}
If $z \in \Ivex$ then $\SymI_z=I(\SymMSV_z)$ is prime and equal to 
the defining ideal of $\SymMSV_z$.
\end{corollary}


\begin{proof}
It follows by standard linear algebra that  $\SymMSV_z$ is equal to the vanishing set of $\SymI_z$, even when $z$ is not vexillary.
When $z$ is vexillary, Theorem~\ref{gr-thm1} implies that $\SymI_z$ is radical so we have $\SymI_z=I(\SymMSV_z)$ by Hilbert's Nullstellensatz. This ideal is prime since the variety $\SymMSV_z$ is irreducible \cite[Lem.~5.2]{BagnoCherniavsky2012}.
\end{proof}

Define
$\displace(P) = \{ (i+k, j+k) : (i,j) \in P\text{ and }k \in P_{ij}\}$
when $P$ is a filling of a subset of $\PP\times \PP$ by
finite sets of nonnegative integers.
When all set-valued entries of $P$ are singletons,
this is the same as the operation used in Lemma~\ref{displace-lem}.
If $P$ is any set-valued (shifted or unshifted) plane partition, then we have 
$   |\displace(P)| = \sum_{(i,j) \in P} |P_{ij}| =:  |P|$.

Fix $z \in \Ivex$ and 
recall 
that $\mu=\mu^{(z)}$ is the partition whose diagram is the northwest alignment of $D(z)$. 
Define $\bff^{(z)}$ as in \eqref{bff-eq}. Then
 $\mathsf{PP}(z)$ and $\SPP(z)$
 respectively
consist of all unshifted and shifted plane partitions $P$ of shapes $\mu$ and $\half(\mu)$ whose entries in row $i$ are at most $f^{(z)}_i - i$.
By Lemma~\ref{displace-lem}, these entry bounds are equivalent to requiring that
$\displace(P)\subseteq \D_\lambda$ for the smallest partition $\lambda =\lambda^{(z)}$ with $\D_\lambda \supseteq D(z)$.

Finally, let $\Delta_z$ and $\SymDelta_z$ be the \defn{Stanley-Reisner complexes} of 
$J_z$ and $\SymJ_z$.
Formally, these are the collections of subsets $S\subseteq \square_n$ and $T \subseteq \UT_n$
such that 
$\textstyle\prod_{(i,j)\in  S} x_{ij} \notin J_z
$ and $ \prod_{(i,j) \in T} u_{ij} \notin \SymJ_z.$

\begin{theorem}
If $z \in \Ivex$ then we can express the initial ideal of $\SymI_z$ as  
\[
\SymJ_z = \bigcap_{P\in \SPP(z)} \bigl\langle u_{ij} : (i,j) \in \displace(P)\bigr\rangle=
\bigcap_{D\in \InvBPD(z)} \bigl\langle u_{ij} : (i,j) \in \blank(D)\bigr\rangle\] and the associated Stanley-Reisner complex is given by 
\be\label{sym-delta-eq}\SymDelta_z = \Bigl\{\UT_n \setminus D   : \displace(P) \subseteq D\subseteq \UT_n\text{ for some }P \in \SPP(z)\Bigr\}.
\ee
\end{theorem}

\begin{proof}
Our two expressions for $\SymJ_z$ are equivalent by Theorem~\ref{bij-thm} since $u_{ij}=u_{ji}$.
By \cite[Thm.~1.7]{MS} one has $J_z = \bigcap_{A \in \Delta_z} \langle x_{ij} :  (i,j) \in \square_n\setminus A\rangle$ and
$\SymJ_z = \bigcap_{A \in \SymDelta_z} \langle u_{ij} :  (i,j) \in \UT_n\setminus A\rangle$.
Thus it is enough to show that the facets of $\SymDelta_z$ are the sets $\UT_n\setminus \displace(P)$ as  $P$ varies over all elements of $\SPP(z)$.
 We can deduce this from the results of Knutson--Miller--Yong \cite[Thms.~4.4 and 5.8]{KnutsonMillerYong}, which imply that the facets of     $\Delta_z$ are the sets $\square_n\setminus \displace(Q)$ for $Q \in \mathsf{PP}(z)$.

Since $\mu$ and $\lambda$ are both symmetric,
it follows from the observations before the theorem statement that  $Q^\top \in \mathsf{PP}(z)$ for all $Q \in \mathsf{PP}(z)$.
 For the same reason, if $Q,Q' \in \mathsf{PP}(z)$ then the entry-wise maximum
 $\max(Q,Q')\in\mathsf{PP}(z)$.
 Finally, notice that if $S \subseteq \UT_n$ then 
 since $\SymJ_z = \pi(J_z)$ 
 we have $S \in\SymDelta_z$ if and only if $S\cup S^\top \in \Delta_z$.

It follows that each set of the form $S = \UT_n\setminus \displace(P) $
for $P \in \SPP(z)$ has $ S \in \SymDelta_z$, since 
if we define $Q = \double(P)\in \mathsf{PP}(z)$  as in Lemma~\ref{vex-double-lem} then
$S\cup S^\top = \square_n\setminus \displace(Q) \in \Delta_z$.
Conversely, if $S \in\SymDelta_z$ so that $S\cup S^\top \in \Delta_z$,
 then there exists $Q \in \mathsf{PP}(z)$ such that $S\cup S^\top  \subseteq \square_n \setminus \displace(Q)$. As $S\cup S^\top$ is symmetric,  
\[
 S\cup S^\top   \subseteq \square_n \setminus (\displace(Q)\cup \displace(Q^\top)) \subseteq \square_n \setminus \displace(\max(Q,Q^\top))
 .\]
But this  means that $S\subseteq \UT_n \setminus \displace(P)$
 for the element $P \in \SPP(z)$ obtained by taking the top half of $\max(Q,Q^\top)$.
 As the sets  $\UT_n\setminus \displace(P) $
for $P \in \SPP(z)$ are all distinct and have the same cardinality,
we conclude that they are precisely
 the facets of $\SymDelta_z$, as desired.
\end{proof}

Continue to fix $z \in \Ivex$ and set 
$ \Sigma := \SD_{\half(\mu)}\subseteq \Gamma := \SD_{\half(\lambda)}\subseteq \UT_n.$
Define $t_0 = (-1,1)$ and $t_i = (i,i+1)(-i,-i-1)$ so that $\{t_0,t_1,\dots,t_{n+1}\}$ is the simple generating set for the Coxeter group of signed permutations $W=W^{\mathsf{B}}_{n+2}$. Write $\circ$ for the Demazure product in $W$ and if $Q=(i_1,i_2,\dots,i_k)$ is a sequence of integers in $\{0,1,\dots,n+1\}$ then let $\delta(Q) = t_{i_1}\circ t_{i_2}\circ \cdots \circ t_{i_k} \in W$.

Label position $(i,j) \in\UT_n$ by either $j - i$
if $(i,j) \in \Gamma$ or else $n+1$ if $(i,j)\notin \Gamma$.
Then form
$Q_D$ for any subset $D\subseteq \UT_n $ by 
reading the labels in the positions of $D$ from right to left in each row, iterating over rows from bottom to top. 
Both $Q_{\Gamma}$ and $Q_{\Sigma}$ are reduced words for elements of $W$ by the discussion preceding \cite[Prop.~5.13]{GrahamKreimanB}.
Let $\sigma_z = \delta(Q_{\Sigma}) \in W $.
Notice that if $D\subseteq \UT_n$ but $D\not\subseteq \Gamma$, then 
\be\label{n2-case} \delta(Q_D) = \delta(Q_{D\cap \Gamma})t_{n+1}.
\ee

For any  $M \subseteq \UT_n$ and  $w \in W$, let  
$\Delta( M, w)$ be the family of subsets $D\subseteq M$ such that  $Q_{M\setminus D}$ contains a subword that is a reduced word for $w$. 
This is an example of Knutson and Miller's \defn{subword complexes}
from \cite{KnutsonMillerSubword}.
Recall that $\SetSPP(z)$ is the family of set-valued plane partitions with the same shape and entry bounds as $\SPP(z)$.
We can derive the following lemma by combining 
some results in \cite{GrahamKreimanB}.

\begin{lemma}
\label{gra-kre-lem}
If $z \in \Ivex$ then the map $P \mapsto \displace(P)$ gives  bijections 
\[ \SetSPP(z) \to \{ D\subseteq \UT_n : \delta(Q_D) = \sigma_z\}
\quand \SPP(z) \to \{ D \subseteq \UT_n : Q_D\text{ is a reduced word for }\sigma_z\}.\]
\end{lemma}

\begin{proof}
Note that since $\ell(\sigma_z) = |\Sigma|= |\half(\mu)|$
is the number of boxes in all elements of $\SPP(z)$,
if the first map is a bijection then the second is also a bijection.
To show that the first map is a bijection, write $\mathcal{E}' $ to denote its codomain.
In view of \eqref{n2-case}, we have $\mathcal{E}'  = \{D\subseteq \Gamma : \delta(Q_D) = \sigma_z\}$.
By
\cite[Prop.~5.13(i)]{GrahamKreimanB}---noting that the generator $s_i$ in \cite{GrahamKreimanB} corresponds to the element $t_{n+2-i}$ used here---the latter is exactly the set of 
\defn{excited shifted Young diagrams} of $\SD_{\half(\mu)}$ in $\SD_{\half(\lambda)}$, as defined in \cite[\S5.3]{GrahamKreimanB} (cf. Example~\ref{exc-ex}).

Now let $\mathcal{T}'$ be the set of set-valued shifted tableaux $T$ (that is, fillings by finite nonempty sets of positive integers, weakly increasing along rows and strictly increasing down columns) of shape $\half(\mu)$ that are \defn{restricted} in the sense that 
for any $(i,j) \in T$ and $k \in T_{ij}$ one has 
$ j - i \leq  \half(\lambda)_k - 1$.
The restriction condition is equivalent to having $j-i \leq \lambda_k-k$, which just means that the map
\be\label{displace-eq2}
T\mapsto \{ (k, k +j- i) : (i,j) \in T\text{ and }k \in T_{ij}\}
\ee
sends every $T \in \mathcal{T}'$ to a subset of $\Gamma$.

As the entry bounds on $P \in \SetSPP(z)$
are likewise equivalent to having $\displace(P)\subseteq \Gamma$, we deduce that  there is a simple bijection $\SetSPP(z)\to \mathcal{T}'$. Namely, just add $i$ to each element of each (set-valued) entry in each row $i$. The operation $\displace$ is exactly the composition of this bijection and \eqref{displace-eq2},
so to prove the lemma we just need to know that \eqref{displace-eq2} is a bijection 
$\mathcal{T}'\to \mathcal{E}'$. This claim is
 \cite[Prop.~5.18]{GrahamKreimanB}.
\end{proof}
We then get the following from  the preceding lemma and \eqref{sym-delta-eq}:

\begin{corollary}\label{sym-delta-cor}
Using the notation introduced above, we have $\SymDelta_z = \Delta(\UT_n,\sigma_z) $ whenever $z \in \Ivex$.
\end{corollary}

Let $T$ be the group of invertible $n\times n$ diagonal matrices over $\CC$. 
Then $\SymMat_{n\times n}$ is a finite-dimensional representation of $T$
under the action $t : M \mapsto tMt$.
Let $a_i : T \to \CC$ for $i \in [n]$ denote the representation $a_i(t) = t_{ii}^{-1}$.
Then each variable $u_{ij}$ is a weight vector in the $T$-module $\CC[\SymMat_{n\times n}]$ with weight $a_ia_j$,
and following the discussion in \cite[\S2.4]{MarPawK} 
we may identify the $T$-equivariant $K$-theory ring 
\[
K_{T}(\SymMat_{n\times n})=\ZZ[a_1^{\pm 1},a_2^{\pm 1}, \ldots, a_n^{\pm 1}].
\]
For any $z \in \Invol$,
we may consider the equivariant $K$-theory class $\iGconj_z := [\SymMSV_z]_T  \in K_{T}(\SymMat_{n\times n}).$
This is \emph{a priori} a Laurent polynomial in the variables $a_1,a_2,\dots,a_n$, but turns out to be an honest polynomial.
The polynomial $\fkGO_z$ is obtained from $\iGconj_z$
by substituting $a_i\mapsto 1 + \beta x_i$  and dividing by $(-\beta)^{\ellinv(z)}$ \cite[Def.~2.18]{MarPawK}.
Thus, to prove Theorem~\ref{fkGO-conj2}, it suffices to derive the following:

\begin{theorem}\label{fkGO-conj3}
If $z \in \Ivex$ then  
$\iGconj_z=\sum_{T \in  \SetSPP(z) } (-1)^{|T|-\ellinv(z)} \prod_{(i,j) \in T} \prod_{k \in T_{ij}} (1 - a_{i+k}a_{j+k}).$
\end{theorem}

 Before starting the proof, we recall from \cite[\S4.4]{MarPaw2022} that
if $R=\CC[x_1,x_2,\dots,x_N]$ is a polynomial ring with a positive grading by a free abelian group $G$,
then an ideal $I \subseteq R$ is \defn{$G$-homogeneous} when it is generated by a set of homogeneous elements. For such an ideal, the quotient ring $R/I$ is a $G$-graded $R$-module with its own \defn{$G$-graded Hilbert series} $\cH(R/I)$,
and the corresponding \defn{$K$-polynomial}  is \[\textstyle\cK(I) = \cH(R/I) \prod_{i=1}^N (1-\deg(x_i)).\]
For the coordinate ring $R = \CC[\SymMat_{n\times n}] = \CC[u_{ij} : (i,j) \in \UT_n]$ we adopt the $G$-grading $\deg(u_{ij} )=a_ia_j$ where $G$ is the multiplicative group generated by $a_1,a_2,\dots,a_n$.
Under this positive grading,
every minor $\det (\cU_{RC})$ is homogeneous so the ideal $\SymI_z$ is $G$-homogeneous.

\begin{proof}
We have $\SymI_z = I(\SymMSV_z)$
by
Corollary~\ref{sym-i-msv-cor},
so, as explained in \cite[\S7.2]{MarPaw2024}, we may also express $\iGconj_z$ as 
 the $K$-polynomial 
 $
 \iGconj_z = \cK(\SymI_z)
 $.  
 By \cite[Lem.~4.17]{MarPaw2022}, we therefore have 
 $ \iGconj_z = \cK(\init(\SymI_z)) = \cK(\SymJ_z).$
Corollary~\ref{sym-delta-cor} shows that $\SymJ_z$ is the \defn{Stanley-Reisner ideal} \cite[Def.~1.6]{MS} of the subword complex $\Delta(\UT_n,\sigma_z)$.
Thus, applying \cite[Thm.~4.1]{KnutsonMillerSubword} (and making the substitution $z_{ij}\mapsto a_ia_j$ as in \cite[Prop.~5.3]{KnutsonMillerSubword}) gives
  \[
 \iGconj_z = \cK(\SymJ_z)=\sum_{\substack{D \subseteq \UT_n \\ \delta(Q_D)=\sigma_z}} (-1)^{|D|-\ell(\sigma_z)} \prod_{(i,j) \in D} (1-a_i a_j).
 \]
 Since $\ell(\sigma_z) = |\Sigma|= |D(z)\cap \UT_n|= \ellinv(z)$, this  transforms to the desired identity
  using Lemma~\ref{gra-kre-lem}.
 \end{proof}

\printbibliography

\end{document}